\documentclass[12pt]{article}
\date{}
\newtheorem{lemma}{Lemma}[section]
\newtheorem{theorem}{Theorem}[section]
\newtheorem{proposition}{Proposition}[section]
\newcommand{\eproof}{~ \hspace*{\fill}\framebox[2.5mm]{}}
\begin{document}
\title{\bf Theta Graph Designs\\}
\author{Anthony D. Forbes, Terry S. Griggs\\
        {\small Department of Mathematics and Statistics, The Open University}\\
        {\small Milton Keynes MK7 6AA, UK}\\
        and\\
        Tamsin J. Forbes\\
        {\small 22 St Albans Road, Kingston upon Thames KT2 5HQ, UK}}
\maketitle
\begin{abstract}
\noindent We solve the design spectrum problem for all theta graphs with 10, 11, 12, 13, 14 and 15 edges.
\end{abstract}


{\noindent {\em E-mail addresses:} anthony.d.forbes@gmail.com, terry.griggs@open.ac.uk, \\tamsin.forbes@gmail.com}

\newcommand{\adfsplit}{\\ \makebox{~~~~~~~~}}
\newcommand{\adfSgap}{\vskip 0.75mm}
\newcommand{\adfLgap}{\vskip 1.5mm}
\newcommand{\adfmod}[1]{~(\mathrm{mod}~#1)}

\section{Introduction}
\label{sec:introduction}
Let $G$ be a simple graph. If the edge set of a simple graph $K$ can be partitioned into edge sets
of graphs each isomorphic to $G$, we say that there exists a {\em decomposition} of $K$ into $G$.
In the case where $K$ is the complete graph $K_n$, we refer to the decomposition as a
$G$ {\em design} of order $n$. The {\em spectrum} of $G$ is the set of non-negative integers $n$ for
which there exists a $G$ design of order $n$.
For completeness we remark that the empty set is a $G$ design of order 0 as well as 1;
these trivial cases will usually be omitted from discussion henceforth.
A complete solution of the spectrum problem often seems to be very difficult.
However it has been achieved in many cases, especially amongst the smaller graphs.
We refer the reader to the survey article of Adams, Bryant and Buchanan, \cite{AdamsBryantBuchanan2008}, and,
for more up to date results, the Web site maintained by Bryant and McCourt, \cite{BryantMcCourt}.

A {\em theta graph} $\Theta(a,b,c)$, $1 \le a \le b \le c$, $b \ge 2$ is the graph with $a+b+c-1$ vertices
and $a+b+c$ edges that consists of three paths of lengths $a$, $b$ and $c$ with common end points
but internally disjoint. Each of the two common end points of the paths has degree 3 and all other vertices have degree 2.
Observe that $\Theta(a,b,c)$ is bipartite if $a \equiv b \equiv c \adfmod{2}$, tripartite otherwise.

The aim of this paper is to establish the design spectrum for
every theta graph with 10, 11, 12, 13, 14 or 15 edges.
The design spectra for all theta graphs with 9 or fewer edges have been determined, \cite{Blinco2003}.
Also there are some partial results for $\Theta(1,b,b)$, odd $b > 1$, \cite{Blinco2001},
and it is known that there exists a $\Theta(a,b,c)$ design of order $2(a+b+c) + 1$, \cite{PunnimPabhapote1987}.
See \cite[Section 5.5]{AdamsBryantBuchanan2008}, specifically Theorems 5.9 and 5.10,
for details and further references.
For convenience, we record the main result of \cite{PunnimPabhapote1987} as follows.
{\theorem \label{thm:2e+1} {\bf (Punnim \& Pabhapote)}
If $1 \le a \le b \le c$ and $b \ge 2$, then there exists a $\Theta(a,b,c)$ design of order $2(a+b+c) + 1$.}\\

It is clear that a $\Theta(a,b,c)$ design of order $n$ can exist only if
(i) $n \le 1$, or $n \ge a+b+c-1$, and (ii) $n(n-1) \equiv 0 \adfmod{2(a+b+c)}$.
These necessary conditions are determined by elementary counting and given explicitly in
Table~\ref{tab:necessary} for some values of $a+b+c$.
In this paper we show that the conditions are sufficient for $n = 10$, 11, 12, 13, 14 and 15.
\begin{table}
\caption{Design existence conditions for theta graphs}
\label{tab:necessary}
\begin{center}
\begin{tabular}{@{}c|l@{}}
 $a+b+c$     & conditions\\
\hline
  10              & $n \equiv 0, 1, 5, 16 \adfmod{20}$, $n \neq 5$ \rule{0mm}{5mm}\\
  12              & $n \equiv 0, 1, 9, 16 \adfmod{24}$, $n \neq 9$\\
  14              & $n \equiv 0, 1, 8, 21 \adfmod{28}$, $n \neq 8$\\
  15              & $n \equiv 0, 1, 6, 10 \adfmod{15}$, $n \neq 6,10$\\
  odd prime power $ \ge 5$  & $n \equiv 0, 1 \adfmod{a+b+c}$
\end{tabular}
\end{center}
\end{table}
\noindent We state our results formally.

{\theorem \label{thm:theta10}
Designs of order $n$ exist for all theta graphs $\Theta(a,b,c)$ with $a+b+c = 10$ if and only if
$n \equiv 0, 1, 5$ or $16 \adfmod{20}$ and $n \neq 5$.}

{\theorem \label{thm:theta11}
Designs of order $n$ exist for all theta graphs $\Theta(a,b,c)$ with $a+b+c = 11$ if and only if
$n \equiv 0$ or $1 \adfmod{11}$.}

{\theorem \label{thm:theta12}
Designs of order $n$ exist for all theta graphs $\Theta(a,b,c)$ with $a+b+c = 12$ if and only if
$n \equiv 0, 1, 9$ or $16 \adfmod{24}$ and $n \neq 9$.}

{\theorem \label{thm:theta13}
Designs of order $n$ exist for all theta graphs $\Theta(a,b,c)$ with $a+b+c = 13$ if and only if
$n \equiv 0$ or $1 \adfmod{13}$.}

{\theorem \label{thm:theta14}
Designs of order $n$ exist for all theta graphs $\Theta(a,b,c)$ with $a+b+c = 14$ if and only if
$n \equiv 0, 1, 8$ or $21 \adfmod{28}$ and $n \neq 8$.}

{\theorem \label{thm:theta15}
Designs of order $n$ exist for all theta graphs $\Theta(a,b,c)$ with $a+b+c = 15$ if and only if
$n \equiv 0$, $1$, $6$ or $10 \adfmod{15}$ and $n \neq 6$, $10$.}

\section{Constructions}
\label{sec:Constructions}

We use Wilson's fundamental construction involving group divisible designs, \cite{WilsonRM1971}.
Recall that a $K$-GDD of type $g_1^{t_1} \dots g_r^{t_r}$ is
an ordered triple ($V,\mathcal{G},\mathcal{B}$)
where $V$ is a set of cardinality $v = t_1 g_1 + \dots + t_r g_r$, $\mathcal{G}$ is a partition of
$V$ into $t_i$ subsets each of cardinality $g_i$, $i = 1, \dots, r$, called \textit{groups} and $\mathcal{B}$
is a collection of subsets of cardinalities $k \in K$, called \textit{blocks}, which collectively have the
property that each pair of elements from different groups occurs in precisely one block but no pair of
elements from the same group occurs at all. A $\{k\}$-GDD is also called a $k$-GDD.
As is well known, whenever $q$ is a prime power there exists a $q$-GDD of type $q^q$ and a $(q+1)$-GDD
of type $q^{q+1}$ (arising from affine and projective planes of order $q$ respectively).
A {\em parallel class} in a group divisible design is a subset of
the block set in which each element of the base set appears exactly once.
A $k$-GDD is called {\em resolvable}, and denoted by $k$-RGDD,
if the entire set of blocks can be partitioned into parallel classes.

Propositions~\ref{prop:bipartite} and \ref{prop:bipartite-alt} are for dealing with bipartite graphs, the latter being used only
for 15-edge bipartite theta graphs, where a suitable decomposition of the form $K_{dr,ds}$ is not available.
Propositions~\ref{prop:prime tripartite} to \ref{prop:theta15} deal with tripartite graphs.

{\proposition \label{prop:bipartite} 
Let $d$, $f$, $r$ and $s$ be positive integers, let $g$ be a non-negative integer, and let $e = 0$ or $1$.
Suppose there exist $G$ designs of orders $fdrs + e$ and $gds + e$, and suppose there exists a decomposition into $G$ of the complete bipartite graph $K_{dr,ds}$.
Then there exist $G$ designs of order $tfdrs + gds + e$ for integers $t \ge 0$.}\\

\noindent\textbf{Proof}~
If $t = 0$, there is nothing to prove; so we may assume $t \ge 1$.
Let $g = 0$ and use induction on $t$.
Assume there exists a $G$ design of order $t fdrs + e$ for some $t \ge 1$.
Take a complete bipartite graph $K_{tfs, fr}$ and
inflate the first part by a factor $dr$ and the second part by a factor of $ds$,
so that the edges become $K_{dr,ds}$ graphs. If $e = 1$, add an extra point, $\infty$.
Overlay the inflated parts, together with $\infty$ when $e = 1$, with $K_{t fdrs + e}$ or $K_{fdrs + e}$ as appropriate.
The result is a complete graph $K_{(t+1) fdrs + e}$ which admits a decomposition into $G$
since there exist decompositions into $G$ of the components, $K_{t fdrs + e}$, $K_{fdrs + e}$ and $K_{dr,ds}$.

For $g \ge 1$, the same construction but starting with a complete bipartite graph $K_{tfs, g}$ yields a design
of order $K_{tfdrs + gds + e}$ since we may now assume that there exists a design of order $tfdrs + e$. \eproof

{\proposition \label{prop:bipartite-alt}
Let $r$ and $s$ be not necessarily distinct positive integers, and let $e = 0$ or $1$.
Suppose there exist $G$ designs of orders $r + e$ and (if $s \neq r$) $s + e$, and suppose there exist decompositions into $G$ of
the complete bipartite graphs $K_{r,r}$ and (if $s \neq r$) $K_{r,s}$.
Then there exist $G$ designs of order $tr + s + e$ for  integers $t \ge 0$.}\\

\noindent\textbf{Proof}~
If $t = 0$, there is nothing to prove; so we may assume $t \ge 1$.
Start with the graph $K_{t+1}$, inflate one point by a factor of $s$ and all others by $r$
so that the original edges become $K_{r,s}$ graphs and, when $t \ge 2$, $K_{r,r}$ graphs.
Add a new point, $\infty$, if $e = 1$.
Overlay the $r$-inflated parts together with $\infty$ if $e = 1$ by $K_{r+e}$.
Overlay the $s$-inflated part together with $\infty$ if $e = 1$ by $K_{s+e}$.
The result is a graph $K_{tr+s+e}$ which admits a decomposition into $G$.\eproof

{\proposition \label{prop:prime tripartite} 
Let $p$ be a positive integer.
Suppose there exist $G$ designs of orders $p$, $p+1$, $2p$ and $2p + 1$.
Suppose also there exist decompositions into $G$ of $K_{p,p,p}$, $K_{p,p,p,p}$ and $K_{p,p,p,p,p}$.
Then there exist $G$ designs of order $n$ for $n \equiv 0,1 \adfmod{p}$.}\\

\noindent\textbf{Proof}~
It is known that there exists a $\{3,4,5\}$-GDD of type $1^t$ for $t \ge 3$, $t \neq 6,8$;
see \cite{AbelBennettGreig2007}.
Inflate each point of the GDD by a factor of $p$, thus expanding the blocks to
complete multipartite graphs $K_{p,p,p}$, $K_{p,p,p,p}$ and $K_{p,p,p,p,p}$.
Let $e = 0$ or $1$. If $e = 1$, add an extra point, $\infty$.
Overlay each inflated group, together with $\infty$ if $e = 1$, with $K_{p + e}$.
Since a design of order $p + e$ and decompositions into $G$ of $K_{p,p,p}$, $K_{p,p,p,p}$ and $K_{p,p,p,p,p}$
exist, this construction creates a design of order $p t +  e$ for $t \ge 3$, $t \ne 6, 8$.

For order $6p + e$, use a $p$-inflated $3$-GDD of type $2^3$, plus an extra point if $e = 1$,
with decompositions into $G$ of $K_{2p+e}$ and $K_{p,p,p}$.
Similarly, for order $8p + e$, use a $3$-GDD of type $2^4$ instead. \eproof

{\proposition \label{prop:2*prime tripartite}
Let $p$ be a positive integer, let $f = 0$ or $1$ and write $f'$ for $1 - f$.
Suppose there exist $G$ designs of orders
$3p+f'$, $4p$, $4p+1$, $5p+f$, $7p+f'$, $8p$, $8p+1$, $9p+f$, $11p+f'$, $13p+f$
and suppose there exist decompositions into $G$ of
$K_{p,p,p,p}$, $K_{p,p,p,p,p}$, $K_{p,p,p,p,3p}$ and $K_{p,p,p,p,4p}$, $K_{2p,2p,2p}$ and $K_{4p,4p,4p,5p}$.
Then there exist $G$ designs of order $n$ for $n \equiv  0, 1, p+f, 3p+f' \adfmod{4p}$, $n \neq p+f$.}\\

\noindent\textbf{Proof}~ Start with a $4$-RGDD of type $4^{3t+1}$, $t \ge 1$, \cite{GeLing2005}; see also \cite{GeMiao2007}.
There are $4t$ parallel classes. Let $x,y,z \ge 0$ and $w = x+y+z \le 4$.
If $w > 0$, add a new group of size $w$ and adjoin each point of this new group to all blocks of a parallel class
to create a \{4,5\}-GDD of type $4^{3t+1} w^1$, $t \ge 1$, which degenerates to a 5-GDD of type $4^5$ if $t = 1$ and $w = 4$.
In the new group inflate $x$ points by a factor of $p$, $y$ points by a factor of $3p$ and $z$ points by a factor of $4p$.
Inflate all points in the other groups by a factor of $p$. So the original blocks become $K_{p,p,p,p}$ graphs,
and new ones $K_{p,p,p,p,p}$ or $K_{p,p,p,p,3p}$ or $K_{p,p,p,p,4p}$ graphs.
Either overlay the groups with $K_{4p}$ and, if $w > 0$, $K_{px + 3py + 4pz}$,
or add a new point and overlay with $K_{4p+1}$ and, if $w > 0$, $K_{px + 3py + 4pz + 1}$.
Using decompositions into $G$ of $K_{4p}$, $K_{4p+1}$, $K_{p,p,p,p}$, $K_{p,p,p,p,p}$, $K_{p,p,p,p,3p}$ and
$K_{p,p,p,p,4p}$ this construction yields a design of order
$12pt + 4p + px + 3py + 4pz + e$ if $t \ge 1$, where $e = 0$ or $1$,
whenever a design of order $px + 3py + 4pz + e$ exists.

As illustrated in Table~\ref{tab:theta10}, all design orders $n \equiv  0, 1, p+f, 3p+f' \adfmod{4p}$,
$n \neq p+f$, are covered except for those values indicated as missing.
\begin{table}
\caption{The construction for Proposition~\ref{prop:2*prime tripartite}}
\label{tab:theta10}
\begin{center}
\begin{tabular}{@{}llll@{}}
   $12p t + 4p$,          & $x = 0$, $y = 0$, $z = 0$, & $e = 0$,  & missing  $4p$\\
   $12p t + 4p + 4p$,     & $x = 0$, $y = 0$, $z = 1$, & $e = 0$,  & missing  $8p$\\
   $12p t + 4p + 8p$,     & $x = 0$, $y = 0$, $z = 2$, & $e = 0$,  & missing  $12p$\\
   $12p t + 4p +  1$,     & $x = 0$, $y = 0$, $z = 0$, & $e = 1$,  & missing  $4p+1$\\
   $12p t + 4p + 4p+1$,   & $x = 0$, $y = 0$, $z = 1$, & $e = 1$,  & missing  $8p+1$\\
   $12p t + 4p + 8p+1$,   & $x = 0$, $y = 0$, $z = 2$, & $e = 1$,  & missing  $12p+1$\\
   $12p t + 4p + 5p+f$,   & $x = 1$, $y = 0$, $z = 1$, & $e = f$,  & missing  $9p+f$\\
   $12p t + 4p + 9p+f$,   & $x = 1$, $y = 0$, $z = 2$, & $e = f$,  & missing  $13p+f$\\
   $12p t + 4p + 13p+f$,  & $x = 1$, $y = 0$, $z = 3$, & $e = f$,  & missing  $5p+f$, $17p+f$\\
   $12p t + 4p + 3p+f'$,  & $x = 0$, $y = 1$, $z = 0$, & $e = f'$, & missing  $7p+f'$\\
   $12p t + 4p + 7p+f'$,  & $x = 0$, $y = 1$, $z = 1$, & $e = f'$, & missing  $11p+f'$\\
   $12p t + 4p + 11p+f'$, & $x = 0$, $y = 1$, $z = 2$, & $e = f'$, & missing  $3p+f'$, $15p+f'$
\end{tabular}
\end{center}
\end{table}
The missing values not assumed as given are handled independently as follows.
See \cite[Tables 4.3 and 4.10]{Ge2007} for lists of small 3- and 4-GDDs.

For ${12p + e}$, use decompositions of $K_{4p + e}$ and $K_{2p,2p,2p}$ with
a $2p$-inflated $3$-GDD of type $2^3$ plus an extra point if $e = 1$.

For ${15p+f'}$, use decompositions of $K_{3p+f'}$ and $K_{p,p,p,p}$  with
a $p$-inflated $4$-GDD of type $3^5$ plus an extra point if $f' = 1$.

For ${17p+f}$, use decompositions of $K_{4p+f}$, $K_{5p+f}$, $K_{4p,4p,4p,5p}$
with the trivial $4$-GDD of type $1^4$ plus an extra point if $f = 1$. \eproof

{\proposition \label{prop:theta10} 
Suppose there exist $G$ designs of orders
$16$, $20$, $21$, $25$, $36$, $40$, $41$, $45$, $56$, $65$
and suppose there exist decompositions into $G$ of $K_{10,10,10}$, $K_{5,5,5,5}$, $K_{20,20,20,25}$,
$K_{5,5,5,5,5}$, $K_{5,5,5,5,15}$ and $K_{5,5,5,5,20}$.
Then there exist $G$ designs of order $n$ for $n \equiv  0, 1, 5, 16 \adfmod{20}$, $n \neq 5$.}\\

\noindent\textbf{Proof}~ Use Proposition \ref{prop:2*prime tripartite} with $p = 5$ and $f = 0$. \eproof

{\proposition \label{prop:theta12} 
Suppose there exist $G$ designs of orders $16$, $24$, $25$, $33$, $40$, $49$, $57$, $81$
and suppose there exist decompositions into $G$ of $K_{8,8,8}$, $K_{8,8,8,8}$ and $K_{8,8,8,24}$.
Then there exist $G$ designs of order $n$ for $n \equiv  0, 1, 9, 16 \adfmod{24}$, $n \neq 9$.}\\

\noindent\textbf{Proof}~ Start with a $3$-RGDD of type $3^{2t+1}$, $t \ge 1$, \cite{Rees1993} (see also \cite{GeMiao2007}),
which has $3t$ parallel classes.
Let $x, y \ge 0$, $x + y \le 3$ and if $x+y > 0$, add a new group of size $x+y$
and adjoin each point of this new group to all blocks of a parallel class
to create a \{3,4\}-GDD of type $3^{2t+1} (x+y)^1$. This degenerates to a 4-GDD of type $3^4$ if $t = 1$ and $x + y = 3$.
In the new group inflate $x$ points by a factor of 8 and $y$ points by a factor of 24. Inflate all other points by a factor of $8$.
The original blocks become $K_{8,8,8}$ graphs and new ones $K_{8,8,8,8}$ or $K_{8,8,8,24}$ graphs.
Either overlay the groups with $K_{24}$ and possibly $K_{8x + 24y}$,
or add a new point and overlay the groups with $K_{25}$ and possibly $K_{8x+24y+1}$.
Using decompositions into $G$ of $K_{24}$, $K_{25}$, $K_{8,8,8}$, $K_{8,8,8,8}$ and $K_{8,8,8,24}$ this
construction yields a design of order $48t + 24 + 8x + 24y + e$, $t \ge 1$, $e = 0$ or $1$,
whenever a design of order $8x + 24y + e$ exists.

As illustrated in Table~\ref{tab:theta12}, all design orders $n \equiv  0, 1, 9, 16 \adfmod{24}$,
$n \neq 9$, are covered except for those values indicated as missing.
\begin{table}
\caption{The construction for Proposition~\ref{prop:theta12}}
\label{tab:theta12}
\begin{center}
\begin{tabular}{@{}llll@{}}
   $48 t + 24     $,& $x = 0$, $y = 0$, & $e = 0$, & missing  24\\
   $48 t + 24 + 24$,& $x = 0$, $y = 1$, & $e = 0$, & missing  48\\
   $48 t + 24 +  1$,& $x = 0$, $y = 0$, & $e = 1$, & missing  25\\
   $48 t + 24 + 25$,& $x = 0$, $y = 1$, & $e = 1$, & missing  49\\
   $48 t + 24 + 33$,& $x = 1$, $y = 1$, & $e = 1$, & missing  57\\
   $48 t + 24 + 57$,& $x = 1$, $y = 2$, & $e = 1$, & missing  33, 81\\
   $48 t + 24 + 16$,& $x = 2$, $y = 0$, & $e = 0$, & missing  40\\
   $48 t + 24 + 40$,& $x = 2$, $y = 1$, & $e = 0$, & missing  16, 64
\end{tabular}
\end{center}
\end{table}
The missing values not assumed as given in the statement of the proposition, namely 48 and 64, are constructed
using decompositions of $K_{16}$ and $K_{8,8,8}$ with an 8-inflated $3$-GDD of type $2^3$ for 48,
or type $2^4$ for 64.\eproof

{\proposition \label{prop:theta14} 
Suppose there exist $G$ designs of orders
$21$, $28$, $29$, $36$, $49$, $56$, $57$, $64$, $77$, $92$
and suppose there exist decompositions into $G$ of $K_{14,14,14}$, $K_{7,7,7,7}$,
$K_{28,28,28,35}$, $K_{7,7,7,7,7}$, $K_{7,7,7,7,21}$ and $K_{7,7,7,7,28}$.
Then there exist $G$ designs of order  $n$ for $n \equiv  0, 1, 8, 21 \adfmod{28}$, $n \neq 8$.}\\

\noindent\textbf{Proof}~ Use Proposition \ref{prop:2*prime tripartite} with $p = 7$ and $f = 1$. \eproof

{\proposition \label{prop:theta15}
Suppose there exist $G$ designs of orders
$15$, $16$, $21$, $25$, $30$, $31$, $36$, $40$, $51$, $55$, $66$ and $70$,
and suppose there exist decompositions into $G$ of
$K_{5,5,5}$ and $K_{5,5,5,5}$.
Then there exist $G$ designs of order $n$ for $n \equiv  0, 1, 6, 10 \adfmod{15}$, $n \neq 6, 10$.}\\

\noindent\textbf{Proof}~ Start with a $3$-RGDD of type $3^{2t+1}$, $t \ge 1$, \cite{Rees1993} (see also \cite{GeMiao2007}),
which has $3t$ parallel classes.
Let $w \ge 0$, $w \le 3t$, and if $w > 0$, add a new group of size $w$ and adjoin each point of this new group to all blocks of a parallel class
to create a \{3,4\}-GDD of type $3^{2t+1} w^1$. This degenerates to a 4-GDD of type $3^4$ if $t = 1$ and $w = 3$.
Inflate all points by a factor of $5$,
so that the original blocks become $K_{5,5,5}$ graphs and new ones $K_{5,5,5,5}$ graphs.
Let $e = 0$ or $1$, and if $e = 1$, add a new point, $\infty$.
Overlay each group, including $\infty$ if $e = 1$, with $K_{15+e}$ or $K_{5w+e}$, as appropriate.
Using decompositions of $K_{15+e}$, $K_{5,5,5}$, and $K_{5,5,5,5}$ this
construction yields a design of order $30t + 15 + 5w + e$, $t \ge w/3$, $e = 0$ or $1$,
whenever a design of order $5w + e$ exists.

As illustrated in Table~\ref{tab:theta15}, all design orders $n \equiv  0, 1, 6, 10 \adfmod{15}$,
$n \neq 6, 10$, are covered except for those values indicated as missing.
\begin{table}
\caption{The construction for Proposition~\ref{prop:theta15}}
\label{tab:theta15}
\begin{center}
\begin{tabular}{@{}llll@{}}
   $30 t + 15$,      & $w = 0$, $e = 0$, & $t \ge 1$, & missing  15\\
   $30 t + 15 + 15$, & $w = 3$, $e = 0$, & $t \ge 1$, & missing  30\\
   $30 t + 15 +  1$, & $w = 0$, $e = 1$, & $t \ge 1$, & missing  16\\
   $30 t + 15 + 16$, & $w = 3$, $e = 1$, & $t \ge 1$, & missing  31\\
   $30 t + 15 + 21$, & $w = 4$, $e = 1$, & $t \ge 2$, & missing  36, 66\\
   $30 t + 15 + 36$, & $w = 7$, $e = 1$, & $t \ge 3$, & missing  21, 51, 81, 111\\
   $30 t + 15 + 25$, & $w = 5$, $e = 0$, & $t \ge 2$, & missing  40, 70\\
   $30 t + 15 + 40$, & $w = 8$, $e = 0$, & $t \ge 3$, & missing  25, 55, 85, 115
\end{tabular}
\end{center}
\end{table}
The missing values not assumed as given are constructed as follows.
See \cite[Table 4.3]{Ge2007} for a list of small 3-GDDs.

For 81, use decompositons of $K_{21}$ and $K_{5,5,5}$ with a 5-inflated 3-GDD of type $4^4$ and an extra point.

For 85, use decompositons of $K_{15}$, $K_{25}$ and $K_{5,5,5}$ with a 5-inflated 3-GDD of type $3^4 5^1$.

For 111, use decompositons of $K_{21}$, $K_{31}$ and $K_{5,5,5}$ with a 5-inflated 3-GDD of type $4^4 6^1$ and an extra point.

For 115, use decompositons of $K_{15}$, $K_{25}$ and $K_{5,5,5}$ with a 5-inflated 3-GDD of type $3^1 5^4$.\eproof

%
%
%
%
%
%
%
%
%
%
%
\newcommand{\ADFvfyParStart}[1]{}

%
%
%

\newcommand{\adfTAABG}{\Theta(1,2,7)}
\newcommand{\adfTAACF}{\Theta(1,3,6)}
\newcommand{\adfTAADE}{\Theta(1,4,5)}
\newcommand{\adfTABBF}{\Theta(2,2,6)}
\newcommand{\adfTABCE}{\Theta(2,3,5)}
\newcommand{\adfTABDD}{\Theta(2,4,4)}
\newcommand{\adfTACCD}{\Theta(3,3,4)}
\newcommand{\adfTBABH}{\Theta(1,2,8)}
\newcommand{\adfTBACG}{\Theta(1,3,7)}
\newcommand{\adfTBADF}{\Theta(1,4,6)}
\newcommand{\adfTBAEE}{\Theta(1,5,5)}
\newcommand{\adfTBBBG}{\Theta(2,2,7)}
\newcommand{\adfTBBCF}{\Theta(2,3,6)}
\newcommand{\adfTBBDE}{\Theta(2,4,5)}
\newcommand{\adfTBCCE}{\Theta(3,3,5)}
\newcommand{\adfTBCDD}{\Theta(3,4,4)}
\newcommand{\adfTCABI}{\Theta(1,2,9)}
\newcommand{\adfTCACH}{\Theta(1,3,8)}
\newcommand{\adfTCADG}{\Theta(1,4,7)}
\newcommand{\adfTCAEF}{\Theta(1,5,6)}
\newcommand{\adfTCBBH}{\Theta(2,2,8)}
\newcommand{\adfTCBCG}{\Theta(2,3,7)}
\newcommand{\adfTCBDF}{\Theta(2,4,6)}
\newcommand{\adfTCBEE}{\Theta(2,5,5)}
\newcommand{\adfTCCCF}{\Theta(3,3,6)}
\newcommand{\adfTCCDE}{\Theta(3,4,5)}
\newcommand{\adfTCDDD}{\Theta(4,4,4)}
\newcommand{\adfTDABJ}{\Theta(1,2,10)}
\newcommand{\adfTDACI}{\Theta(1,3,9)}
\newcommand{\adfTDADH}{\Theta(1,4,8)}
\newcommand{\adfTDAEG}{\Theta(1,5,7)}
\newcommand{\adfTDAFF}{\Theta(1,6,6)}
\newcommand{\adfTDBBI}{\Theta(2,2,9)}
\newcommand{\adfTDBCH}{\Theta(2,3,8)}
\newcommand{\adfTDBDG}{\Theta(2,4,7)}
\newcommand{\adfTDBEF}{\Theta(2,5,6)}
\newcommand{\adfTDCCG}{\Theta(3,3,7)}
\newcommand{\adfTDCDF}{\Theta(3,4,6)}
\newcommand{\adfTDCEE}{\Theta(3,5,5)}
\newcommand{\adfTDDDE}{\Theta(4,4,5)}
\newcommand{\adfTEABK}{\Theta(1,2,11)}
\newcommand{\adfTEACJ}{\Theta(1,3,10)}
\newcommand{\adfTEADI}{\Theta(1,4,9)}
\newcommand{\adfTEAEH}{\Theta(1,5,8)}
\newcommand{\adfTEAFG}{\Theta(1,6,7)}
\newcommand{\adfTEBBJ}{\Theta(2,2,10)}
\newcommand{\adfTEBCI}{\Theta(2,3,9)}
\newcommand{\adfTEBDH}{\Theta(2,4,8)}
\newcommand{\adfTEBEG}{\Theta(2,5,7)}
\newcommand{\adfTEBFF}{\Theta(2,6,6)}
\newcommand{\adfTECCH}{\Theta(3,3,8)}
\newcommand{\adfTECDG}{\Theta(3,4,7)}
\newcommand{\adfTECEF}{\Theta(3,5,6)}
\newcommand{\adfTEDDF}{\Theta(4,4,6)}
\newcommand{\adfTEDEE}{\Theta(4,5,5)}
\newcommand{\adfTFABL}{\Theta(1,2,12)}
\newcommand{\adfTFACK}{\Theta(1,3,11)}
\newcommand{\adfTFADJ}{\Theta(1,4,10)}
\newcommand{\adfTFAEI}{\Theta(1,5,9)}
\newcommand{\adfTFAFH}{\Theta(1,6,8)}
\newcommand{\adfTFAGG}{\Theta(1,7,7)}
\newcommand{\adfTFBBK}{\Theta(2,2,11)}
\newcommand{\adfTFBCJ}{\Theta(2,3,10)}
\newcommand{\adfTFBDI}{\Theta(2,4,9)}
\newcommand{\adfTFBEH}{\Theta(2,5,8)}
\newcommand{\adfTFBFG}{\Theta(2,6,7)}
\newcommand{\adfTFCCI}{\Theta(3,3,9)}
\newcommand{\adfTFCDH}{\Theta(3,4,8)}
\newcommand{\adfTFCEG}{\Theta(3,5,7)}
\newcommand{\adfTFCFF}{\Theta(3,6,6)}
\newcommand{\adfTFDDG}{\Theta(4,4,7)}
\newcommand{\adfTFDEF}{\Theta(4,5,6)}
\newcommand{\adfTFEEE}{\Theta(5,5,5)}


\section{Proofs of the theorems}\label{sec:Proofs}

We deal with each of the theta graphs stated in the theorems. The lemmas in this section assert the
existence of specific graph designs and decompositions of multipartite graphs.
These are used as ingredients for the propositions of Section~\ref{sec:Constructions} to construct the decompositions of
complete graphs required to prove the theorems.
With a few exceptions, the details of the decompositions that constitute the proofs of the lemmas have been deferred to sections in the
rather lengthy Appendix to this paper.
If absent, the Appendix may be obtained from the ArXiv (identifier 1703.01483), or by request from the first author.
In the presentation of our results we represent $\Theta(a,b,c)$ by a
subscripted ordered ($a+b+c-1$)-tuple $(v_1, v_2, \dots, v_{a+b+c-1})_{\Theta(a,b,c)}$,
where $v_1$ and $v_2$ are the vertices of degree 3 and the three paths are
\begin{eqnarray*}
&(\textrm{i})&(v_1,~ v_2) ~~\textrm{if}~a = 1,~~~~ (v_1,~ v_3,~ \dots,~ v_{a + 1},~ v_2)
~~\textrm{if}~ a \ge 2,\\
&(\textrm{ii})&(v_1,~ v_{a + 2},~ \dots,~ v_{a + b},~ v_2),\\
&(\textrm{iii})&(v_1,~ v_{a + b + 1},~ \dots,~ v_{a + b + c - 1},~ v_2).
\end{eqnarray*}
If a graph $G$ has $e$ edges, the numbers of occurrences of $G$ in a decomposition into $G$ of the
complete graph $K_n$, the complete $r$-partite graph $K_{n^r}$ and the complete $(r+1)$-partite
graph $K_{n^r m^1}$ are respectively
$$ \frac{n(n-1)}{2e},~~~ \frac{n^2 r (r-1)}{2e}~~~ \textrm{and}~~~ \frac{nr(n(r - 1) + 2m)}{2e}.$$
The number of theta graphs with $e$ edges is $\lfloor e^2/12 - 1/2\rfloor$ of which
$\lfloor{e^2}/{48} + {(e \mathrm{~mod~} 2)(e-8)}/{8} + 1/2 \rfloor$ are bipartite; see \cite{Forbes2017JulM500} for example.

All decompositions were created by a computer program written in the C language, and
a {\sc Mathematica} program provided final assurance regarding their correctness.
Most of the decompositions were obtained without difficulty. However, there were a few that seemed to present
a challenge, notably $K_{20}$ for theta graphs of 10 edges, $K_{24}$ for theta graphs of 12 edges and
$K_{30}$ for theta graphs of 15 edges.
So, for these three cases, it is appropriate to include the decomposition details in the main body of the paper.

{\lemma \label{lem:theta10 designs}
There exist $\Theta(a,b,c)$ designs of orders $16$, $20$        and $25$ for $a + b + c = 10$.
There exist $\Theta(a,b,c)$ designs of orders $36$, $40$, $41$, $45$, $56$, $65$
for $a + b + c = 10$ with $a$, $b$, $c$ not all even.}\\

\noindent {\bf Proof}~ There are seven theta graphs of which two are bipartite.
We present designs of order 20 here; the remaining decompositions are presented in the Appendix.\\

\noindent {\boldmath $K_{20}$}~
Let the vertex set be $Z_{20}$. The decompositions consist of the graphs

\adfLgap 
$(4,17,0,9,7,1,12,15,13)_{\adfTAABG}$,
$(11,17,18,10,1,13,3,7,16)_{\adfTAABG}$,

$(1,6,18,14,11,12,2,0,15)_{\adfTAABG}$,

\adfSgap
$(12,4,10,18,0,8,15,3,11)_{\adfTAABG}$,
$(6,2,8,11,19,7,0,14,16)_{\adfTAABG}$,

$(2,4,18,7,15,10,14,8,16)_{\adfTAABG}$,
$(6,12,0,10,16,3,18,14,19)_{\adfTAABG}$,

\adfLgap 
$(0,1,9,18,12,15,4,3,5)_{\adfTAACF}$,
$(0,7,17,5,10,14,12,8,2)_{\adfTAACF}$,

$(8,6,14,1,13,3,17,4,19)_{\adfTAACF}$,

\adfSgap
$(3,15,2,1,11,17,18,7,6)_{\adfTAACF}$,
$(3,7,14,11,19,5,6,18,10)_{\adfTAACF}$,

$(7,19,15,18,13,14,2,10,11)_{\adfTAACF}$,
$(15,19,11,2,14,6,3,9,10)_{\adfTAACF}$,

\adfLgap 
$(0,2,3,6,9,11,10,14,8)_{\adfTAADE}$,
$(17,2,7,12,13,9,0,18,11)_{\adfTAADE}$,

$(6,11,5,3,17,18,8,0,7)_{\adfTAADE}$,

\adfSgap
$(8,12,13,15,16,4,3,11,9)_{\adfTAADE}$,
$(1,8,3,15,7,16,13,0,5)_{\adfTAADE}$,

$(1,17,5,7,19,4,9,16,0)_{\adfTAADE}$,
$(12,17,5,9,13,11,19,0,4)_{\adfTAADE}$,

\adfLgap 
$(2,1,8,10,18,7,12,0,6)_{\adfTABBF}$,
$(10,15,7,8,12,11,6,13,5)_{\adfTABBF}$,

$(10,8,11,9,0,16,7,14,17)_{\adfTABBF}$,

\adfSgap
$(5,1,7,3,8,13,11,17,15)_{\adfTABBF}$,
$(1,9,4,5,19,3,7,13,15)_{\adfTABBF}$,

$(5,17,0,19,11,9,13,16,1)_{\adfTABBF}$,
$(9,17,3,12,7,11,15,19,13)_{\adfTABBF}$,

\adfLgap 
$(0,5,18,2,10,4,9,8,17)_{\adfTABCE}$,
$(13,5,2,11,9,3,7,6,12)_{\adfTABCE}$,

$(7,6,0,1,15,2,12,9,10)_{\adfTABCE}$,

\adfSgap
$(11,3,12,16,15,0,8,7,4)_{\adfTABCE}$,
$(3,7,10,0,15,8,16,4,19)_{\adfTABCE}$,

$(7,19,16,14,11,12,4,15,2)_{\adfTABCE}$,
$(11,19,8,3,6,18,15,12,0)_{\adfTABCE}$,

\adfLgap 
$(8,2,7,5,13,9,4,18,17)_{\adfTABDD}$,
$(18,6,10,9,8,0,7,17,4)_{\adfTABDD}$,

$(13,11,15,8,18,0,10,3,5)_{\adfTABDD}$,

\adfSgap
$(7,3,6,13,4,11,15,0,9)_{\adfTABDD}$,
$(3,7,0,2,19,4,8,11,10)_{\adfTABDD}$,

$(15,19,12,3,16,5,8,17,11)_{\adfTABDD}$,
$(15,19,18,1,12,7,14,11,16)_{\adfTABDD}$,

\adfLgap 
$(15,10,18,13,6,12,11,0,8)_{\adfTACCD}$,
$(0,11,16,3,17,8,10,1,5)_{\adfTACCD}$,

$(5,7,0,13,3,8,17,19,14)_{\adfTACCD}$,

\adfSgap
$(6,18,7,17,13,12,0,1,2)_{\adfTACCD}$,
$(6,10,5,4,14,18,2,16,17)_{\adfTACCD}$,

$(10,14,9,8,11,1,2,3,13)_{\adfTACCD}$,
$(14,18,10,6,15,5,2,9,19)_{\adfTACCD}$

\adfLgap
\noindent under the action of the mapping $x \mapsto x + 4$ (mod 20) for only the first three graphs in each design. \eproof

\ADFvfyParStart{20, \{\{20,1,20\}\}, -1, -1, \{3,\{\{20,5,4\}\}\}} 

{\lemma \label{lem:theta10 multipartite}
There exist decompositions into $\Theta(a,b,c)$ of the complete bipartite graph $K_{10,5}$
for $a + b + c = 10$ with $a$, $b$, $c$ all even.
There exist decompositions into $\Theta(a,b,c)$ of the complete multipartite graphs
$K_{10,10,10}$, $K_{5,5,5,5}$, $K_{20,20,20,25}$, $K_{5,5,5,5,5}$, $K_{5,5,5,5,15}$ and $K_{5,5,5,5,20}$
for $a + b + c = 10$ with $a$, $b$, $c$ not all even.}\\

\noindent {\bf Proof}~ The decompositions are presented in the Appendix. \eproof\\

Theorem~\ref{thm:theta10} for $\Theta(a,b,c)$ with $a+b+c=10$ follows from
Theorem~\ref{thm:2e+1}, Lemmas~\ref{lem:theta10 designs} and \ref{lem:theta10 multipartite}, and
Proposition~\ref{prop:bipartite} (with $d = 5$, $r = 2$, $s = 1$, $f = 2$ and $g = 0$, 3 or 5) if $a$, $b$ and $c$ are all even,
Proposition~\ref{prop:theta10} otherwise.\eproof


{\lemma \label{lem:theta11 designs}
There exist $\Theta(a,b,c)$ designs of orders   $11$ and $12$ for $a + b + c = 11$.
There exists a $\Theta(a,b,c)$ design of order  $22$          for $a + b + c = 11$ with $a$, $b$, $c$ not all odd.}\\

\noindent {\bf Proof}~ There are nine theta graphs of which three are bipartite.
The decompositions are presented in the Appendix.\eproof

{\lemma \label{lem:theta11 multipartite}
There exist decompositions into $\Theta(a,b,c)$ of the complete bipartite graph $K_{11,11}$
for $a + b + c = 11$ with $a$, $b$, $c$ all odd.
There exist decompositions into $\Theta(a,b,c)$ of the complete multipartite graphs
$K_{11,11,11}$, $K_{11,11,11,11}$ and $K_{11,11,11,11,11}$ for $a + b + c = 11$ with $a$, $b$, $c$ not all odd.}\\

\noindent {\bf Proof}~ The decompositions are presented in the Appendix.\eproof\\

Theorem~\ref{thm:theta11} for $\Theta(a,b,c)$ with $a+b+c=11$ follows from
Theorem~\ref{thm:2e+1}, Lemmas~\ref{lem:theta11 designs} and \ref{lem:theta11 multipartite}, and
Proposition~\ref{prop:bipartite} (with $d = 11$, $r = s = f = 1$ and $g = 0$) if $a$, $b$ and $c$ are all odd,
Proposition~\ref{prop:prime tripartite} otherwise.\eproof


{\lemma \label{lem:theta12 designs}
There exist $\Theta(a,b,c)$ designs of orders $16$, $24$       and $33$ for $a + b + c = 12$.
There exist $\Theta(a,b,c)$ designs of orders $40$, $49$, $57$ and $81$ for $a + b + c = 12$
with $a$, $b$, $c$ not all even.}\\

\noindent {\bf Proof}~ There are eleven theta graphs of which three are bipartite.
We present designs of order 24 here; the remaining decompositions are presented in the Appendix.\\

\noindent {\boldmath $K_{24}$}~
Let the vertex set be $Z_{24}$. The decompositions consist of

\adfLgap 
$(11,9,16,19,10,12,6,15,7,2,8)_{\adfTCABI}$,
$(8,0,10,7,6,19,18,1,20,22,12)_{\adfTCABI}$,

$(2,6,13,15,1,17,18,14,22,19,21)_{\adfTCABI}$,
$(10,7,1,3,21,8,14,9,19,16,13)_{\adfTCABI}$,

$(16,21,1,7,9,22,15,19,5,17,4)_{\adfTCABI}$,
$(18,4,10,21,17,14,16,6,11,12,19)_{\adfTCABI}$,

\adfSgap
$(21,20,7,2,14,13,18,6,5,4,12)_{\adfTCABI}$,
$(21,12,13,22,10,5,15,20,16,3,9)_{\adfTCABI}$,

$(4,20,0,1,19,23,11,15,3,7,5)_{\adfTCABI}$,
$(4,23,13,8,15,12,7,0,11,17,20)_{\adfTCABI}$,

$(12,23,16,8,19,7,4,15,13,5,21)_{\adfTCABI}$,

\adfLgap 
$(22,18,1,11,8,16,17,13,9,15,2)_{\adfTCACH}$,
$(9,20,10,22,1,3,21,4,19,12,2)_{\adfTCACH}$,

$(7,4,17,5,10,12,21,6,19,16,20)_{\adfTCACH}$,
$(16,2,3,9,21,22,0,18,13,15,11)_{\adfTCACH}$,

$(18,19,14,7,9,12,11,13,0,6,22)_{\adfTCACH}$,
$(12,7,16,11,17,22,4,15,10,13,6)_{\adfTCACH}$,

\adfSgap
$(8,15,17,0,10,21,11,19,5,18,16)_{\adfTCACH}$,
$(5,21,13,7,15,6,11,3,22,9,0)_{\adfTCACH}$,

$(7,23,8,5,22,12,0,2,13,16,1)_{\adfTCACH}$,
$(8,23,1,14,20,6,17,15,7,9,16)_{\adfTCACH}$,

$(13,23,21,15,3,19,14,4,16,7,0)_{\adfTCACH}$,

\adfLgap 
$(22,19,10,9,20,11,16,8,2,21,4)_{\adfTCADG}$,
$(19,9,3,21,18,14,1,7,2,20,4)_{\adfTCADG}$,

$(20,18,15,0,3,6,2,16,21,10,11)_{\adfTCADG}$,
$(21,0,7,20,22,14,18,16,15,3,13)_{\adfTCADG}$,

$(4,11,7,18,15,8,12,0,9,17,13)_{\adfTCADG}$,
$(4,1,22,13,15,18,2,9,11,7,0)_{\adfTCADG}$,

\adfSgap
$(6,22,1,8,14,0,11,17,5,20,21)_{\adfTCADG}$,
$(7,9,15,14,6,22,8,19,1,13,21)_{\adfTCADG}$,

$(1,22,5,21,15,23,7,14,13,9,16)_{\adfTCADG}$,
$(6,23,7,13,5,16,3,9,14,17,22)_{\adfTCADG}$,

$(15,17,23,14,0,6,5,4,13,12,21)_{\adfTCADG}$,

\adfLgap 
$(1,6,22,19,9,8,5,13,12,11,17)_{\adfTCAEF}$,
$(0,7,9,17,19,4,8,3,16,2,21)_{\adfTCAEF}$,

$(5,17,19,16,9,15,6,4,22,12,18)_{\adfTCAEF}$,
$(3,15,22,16,12,10,21,8,20,9,4)_{\adfTCAEF}$,

$(13,4,22,15,7,12,16,20,3,5,18)_{\adfTCAEF}$,
$(15,1,16,6,22,10,0,2,19,18,21)_{\adfTCAEF}$,

\adfSgap
$(18,2,10,19,3,11,1,4,21,14,15)_{\adfTCAEF}$,
$(5,7,0,18,15,11,12,9,2,6,10)_{\adfTCAEF}$,

$(7,13,22,11,19,6,3,14,23,2,8)_{\adfTCAEF}$,
$(10,23,14,18,22,5,17,20,13,15,19)_{\adfTCAEF}$,

$(21,23,16,10,2,22,15,6,7,18,3)_{\adfTCAEF}$,

\adfLgap 
$(11,15,2,13,19,1,0,4,12,3,20)_{\adfTCBBH}$,
$(16,13,18,1,21,17,4,14,15,0,10)_{\adfTCBBH}$,

$(12,11,18,0,14,19,20,5,21,7,15)_{\adfTCBBH}$,
$(13,3,17,14,20,18,15,19,22,0,8)_{\adfTCBBH}$,

$(12,14,2,17,13,19,18,9,22,7,8)_{\adfTCBBH}$,
$(13,20,7,0,3,1,15,8,5,22,17)_{\adfTCBBH}$,

\adfSgap
$(2,17,9,18,22,4,7,19,16,6,10)_{\adfTCBBH}$,
$(14,2,6,10,18,22,15,3,0,17,1)_{\adfTCBBH}$,

$(1,2,7,8,23,17,15,12,6,22,13)_{\adfTCBBH}$,
$(9,10,15,16,7,14,22,8,11,23,18)_{\adfTCBBH}$,

$(14,18,0,5,20,23,6,21,10,9,1)_{\adfTCBBH}$,

\adfLgap 
$(7,11,0,1,14,3,15,10,22,17,16)_{\adfTCBCG}$,
$(1,20,15,17,3,8,18,5,22,14,4)_{\adfTCBCG}$,

$(7,0,13,6,8,16,4,9,2,12,14)_{\adfTCBCG}$,
$(10,14,13,1,5,2,3,21,19,17,18)_{\adfTCBCG}$,

$(9,7,21,20,18,6,19,15,22,4,5)_{\adfTCBCG}$,
$(1,3,19,4,10,16,13,21,2,6,0)_{\adfTCBCG}$,

\adfSgap
$(3,11,12,13,20,18,16,15,4,8,2)_{\adfTCBCG}$,
$(13,15,17,4,7,8,10,19,20,0,18)_{\adfTCBCG}$,

$(5,23,12,0,2,9,7,10,16,21,1)_{\adfTCBCG}$,
$(15,21,12,6,11,23,20,7,22,3,4)_{\adfTCBCG}$,

$(19,23,14,4,0,5,20,16,12,8,7)_{\adfTCBCG}$,

\adfLgap 
$(4,5,18,0,12,3,7,16,6,8,15)_{\adfTCBDF}$,
$(1,16,13,10,6,22,5,11,2,19,8)_{\adfTCBDF}$,

$(22,21,4,12,11,16,5,14,3,17,15)_{\adfTCBDF}$,
$(18,9,17,8,4,5,0,13,21,12,14)_{\adfTCBDF}$,

$(15,17,10,19,12,7,4,1,20,18,6)_{\adfTCBDF}$,
$(3,5,19,15,9,6,1,12,4,10,7)_{\adfTCBDF}$,

\adfSgap
$(10,3,2,8,17,0,11,6,15,22,7)_{\adfTCBDF}$,
$(3,18,22,9,0,2,6,7,15,14,19)_{\adfTCBDF}$,

$(10,18,13,5,2,21,23,0,1,8,7)_{\adfTCBDF}$,
$(16,23,15,1,19,22,9,8,11,14,7)_{\adfTCBDF}$,

$(16,23,19,18,10,14,17,11,15,2,6)_{\adfTCBDF}$,

\adfLgap 
$(15,14,21,3,0,4,11,19,20,1,13)_{\adfTCBEE}$,
$(14,0,18,15,2,4,19,8,16,1,5)_{\adfTCBEE}$,

$(18,1,11,9,8,6,21,19,17,0,15)_{\adfTCBEE}$,
$(10,21,2,13,5,19,8,7,14,0,20)_{\adfTCBEE}$,

$(21,13,4,3,11,17,15,0,7,1,2)_{\adfTCBEE}$,
$(19,22,10,21,7,15,11,14,12,18,1)_{\adfTCBEE}$,

\adfSgap
$(2,7,12,0,10,20,8,6,17,1,4)_{\adfTCBEE}$,
$(10,18,8,15,6,12,4,14,9,20,23)_{\adfTCBEE}$,

$(4,22,17,14,1,12,9,20,15,16,18)_{\adfTCBEE}$,
$(7,12,22,2,16,4,15,20,6,14,23)_{\adfTCBEE}$,

$(20,22,14,12,0,23,4,17,9,1,6)_{\adfTCBEE}$,

\adfLgap 
$(14,13,20,7,16,9,22,17,18,21,0)_{\adfTCCCF}$,
$(21,9,3,5,22,18,2,16,1,17,19)_{\adfTCCCF}$,

$(21,12,4,13,16,18,14,15,20,17,11)_{\adfTCCCF}$,
$(19,2,22,8,15,9,12,1,13,21,11)_{\adfTCCCF}$,

$(7,2,21,10,9,4,14,17,16,22,12)_{\adfTCCCF}$,
$(18,14,7,2,15,1,19,4,16,8,23)_{\adfTCCCF}$,

\adfSgap
$(8,19,7,0,11,6,12,20,16,3,14)_{\adfTCCCF}$,
$(14,3,5,7,12,15,18,11,16,23,0)_{\adfTCCCF}$,

$(0,3,11,19,20,22,4,7,23,15,8)_{\adfTCCCF}$,
$(6,19,4,8,21,23,10,3,11,15,7)_{\adfTCCCF}$,

$(16,22,15,13,19,2,12,4,20,23,11)_{\adfTCCCF}$,

\adfLgap 
$(16,9,17,13,2,11,4,8,22,14,19)_{\adfTCCDE}$,
$(21,19,6,7,22,4,15,2,20,8,11)_{\adfTCCDE}$,

$(3,5,8,7,4,1,19,6,9,14,18)_{\adfTCCDE}$,
$(22,4,5,17,10,11,0,20,6,2,7)_{\adfTCCDE}$,

$(14,11,16,9,8,10,13,7,1,17,15)_{\adfTCCDE}$,
$(22,21,11,5,7,10,20,9,2,1,16)_{\adfTCCDE}$,

\adfSgap
$(10,20,1,15,3,12,4,18,7,21,0)_{\adfTCCDE}$,
$(16,20,10,12,7,0,13,5,23,9,18)_{\adfTCCDE}$,

$(2,16,10,23,15,7,13,4,8,21,12)_{\adfTCCDE}$,
$(7,15,12,5,17,2,8,23,13,4,21)_{\adfTCCDE}$,

$(18,23,0,15,2,19,4,11,20,5,8)_{\adfTCCDE}$,

\adfLgap 
$(20,6,21,4,11,14,18,17,2,10,5)_{\adfTCDDD}$,
$(13,12,10,14,22,17,16,0,2,8,11)_{\adfTCDDD}$,

$(5,19,7,4,2,8,22,0,17,15,18)_{\adfTCDDD}$,
$(8,1,14,15,22,13,21,3,7,2,16)_{\adfTCDDD}$,

$(5,21,14,12,8,19,3,15,20,0,7)_{\adfTCDDD}$,
$(21,7,14,9,3,19,8,1,17,20,15)_{\adfTCDDD}$,

\adfSgap
$(20,11,10,22,7,4,17,1,9,2,14)_{\adfTCDDD}$,
$(11,18,2,15,1,20,7,16,22,19,6)_{\adfTCDDD}$,

$(2,4,0,15,19,12,23,9,17,7,18)_{\adfTCDDD}$,
$(3,10,12,20,1,14,23,8,17,9,19)_{\adfTCDDD}$,

$(3,12,6,15,4,18,9,1,23,10,17)_{\adfTCDDD}$

\adfLgap
\noindent under the action of the mapping $x \mapsto x + 8$ (mod 24) for only the first six graphs in each design. \eproof

\ADFvfyParStart{24, \{\{24,1,24\}\}, -1, -1, \{6,\{\{24,3,8\}\}\}} 

{\lemma \label{lem:theta12 multipartite}
There exist decompositions into $\Theta(a,b,c)$ of the complete bipartite graph $K_{12,8}$
for $a + b + c = 12$ with $a$, $b$, $c$ all even.
There exist decompositions into $\Theta(a,b,c)$ of the complete multipartite graphs
$K_{8,8,8}$, $K_{8,8,8,8}$ and $K_{8,8,8,24}$ for $a + b + c = 12$ with $a$, $b$, $c$ not all even.}\\

\noindent {\bf Proof}~ The decompositions are presented in the Appendix.\eproof\\

Theorem~\ref{thm:theta12} for $\Theta(a,b,c)$ with $a+b+c=12$ follows from
Theorem~\ref{thm:2e+1}, Lemmas~\ref{lem:theta12 designs} and \ref{lem:theta12 multipartite}, and
Proposition~\ref{prop:bipartite} (with $d = 4$, $r = 3$, $s = 2$, $f = 1$ and $g = 0$, 2 or 4) if $a$, $b$ and $c$ are all even,
Proposition~\ref{prop:theta12} otherwise.\eproof


{\lemma \label{lem:theta13 designs}
There exist $\Theta(a,b,c)$ designs of orders   $13$ and $14$ for $a + b + c = 13$.
There exists a $\Theta(a,b,c)$ design of order  $26$          for $a + b + c = 13$ with $a$, $b$, $c$ not all odd.}\\

\noindent {\bf Proof}~ There are thirteen theta graphs of which four are bipartite.
The decompositions are presented in the Appendix.\eproof

{\lemma \label{lem:theta13 multipartite}
There exist decompositions into $\Theta(a,b,c)$ of the complete bipartite graph $K_{13,13}$
for $a + b + c = 13$ with $a$, $b$, $c$ all odd.
There exist decompositions into $\Theta(a,b,c)$ of the complete multipartite graphs
$K_{13,13,13}$, $K_{13,13,13,13}$ and $K_{13,13,13,13,13}$ for $a + b + c = 13$ with $a$, $b$, $c$ not all odd.}\\

\noindent {\bf Proof}~ The decompositions are presented in the Appendix.\eproof\\

Theorem~\ref{thm:theta13} for $\Theta(a,b,c)$ with $a+b+c=13$ follows from
Theorem~\ref{thm:2e+1}, Lemmas~\ref{lem:theta13 designs} and \ref{lem:theta13 multipartite}, and
Proposition~\ref{prop:bipartite} (with $d = 13$, $r = s = f = 1$ and $g = 0$) if $a$, $b$ and $c$ are all odd,
Proposition~\ref{prop:prime tripartite} otherwise.\eproof


{\lemma \label{lem:theta14 designs}
There exist $\Theta(a,b,c)$ designs of order $21$, $28$        and $36$ for $a + b + c = 14$.
There exist $\Theta(a,b,c)$ designs of order $49$, $56$, $57$, $64$, $77$ and $92$
for $a + b + c = 14$ with $a$, $b$, $c$ not all even.}\\

\noindent {\bf Proof}~ There are fifteen theta graphs of which four are bipartite.
The decompositions are presented in the Appendix.\eproof

{\lemma \label{lem:theta14 multipartite}
There exist decompositions into $\Theta(a,b,c)$ of the complete bipartite graph $K_{{14,7}}$
for $a + b + c = 14$ with $a$, $b$, $c$ all even.
There exist decompositions into $\Theta(a,b,c)$ of the complete multipartite graphs
$K_{14,14,14}$, $K_{7,7,7,7}$, $K_{28,28,28,35}$, $K_{7,7,7,7,7}$, $K_{7,7,7,7,21}$ and $K_{7,7,7,7,28}$
for $a + b + c = 14$ with $a$, $b$, $c$ not all even.}\\

\noindent {\bf Proof}~ The decompositions are presented in the Appendix.\eproof\\

Theorem~\ref{thm:theta14} for $\Theta(a,b,c)$ with $a+b+c=14$ follows from
Theorem~\ref{thm:2e+1}, Lemmas~\ref{lem:theta14 designs} and \ref{lem:theta14 multipartite}, and
Proposition~\ref{prop:bipartite} (with $d = 7$, $r = 2$, $s = 1$, $f = 2$ and $g = 0$, 3 or 5) if $a$, $b$ and $c$ are all even,
Proposition~\ref{prop:theta14} otherwise.\eproof


{\lemma \label{lem:theta15 designs}
There exist $\Theta(a,b,c)$ designs of orders $15$, $16$, $21$, and $25$ for $a + b + c = 15$.
There exist $\Theta(a,b,c)$ designs of orders       $30$, $36$, $40$, $51$, $55$, $66$ and $70$
for $a + b + c = 15$ with $a$, $b$, $c$ not all odd.}\\

\noindent {\bf Proof}~
There are eighteen 15-edge theta graphs of which six are bipartite.
We present designs of order 30 here; the remaining decompositions are presented in the Appendix.\\

\noindent {\boldmath $K_{30}$}~
Let the vertex set be $Z_{30}$. The decompositions consist of the graphs

\adfLgap
$(6,2,0,7,22,5,19,15,18,25,23,1,21,12)_{\adfTFABL}$,

$(7,20,15,2,10,9,8,5,6,17,3,18,28,22)_{\adfTFABL}$,

$(20,4,14,19,1,3,5,18,11,16,24,13,2,21)_{\adfTFABL}$,

$(7,13,10,14,6,11,8,4,25,9,18,21,27,23)_{\adfTFABL}$,

$(27,5,16,22,6,18,13,0,14,2,9,21,11,4)_{\adfTFABL}$,

\adfSgap
$(5,17,26,11,23,2,15,22,4,6,10,28,21,8)_{\adfTFABL}$,

$(23,5,14,17,2,11,7,28,24,22,18,16,25,29)_{\adfTFABL}$,

$(17,29,11,13,4,0,28,16,12,10,22,1,5,20)_{\adfTFABL}$,

$(23,29,8,19,10,3,20,11,26,9,16,4,27,14)_{\adfTFABL}$,

\adfLgap
$(19,1,6,11,23,0,27,21,5,18,14,7,15,17)_{\adfTFADJ}$,

$(15,19,0,28,20,1,6,9,18,22,11,16,25,8)_{\adfTFADJ}$,

$(17,26,16,15,8,23,10,28,5,9,20,12,19,21)_{\adfTFADJ}$,

$(0,9,1,10,21,14,12,23,15,25,5,7,2,26)_{\adfTFADJ}$,

$(10,15,4,2,28,25,19,16,13,17,27,26,6,22)_{\adfTFADJ}$,

\adfSgap
$(12,0,5,20,23,24,4,14,11,2,16,6,28,18)_{\adfTFADJ}$,

$(12,18,4,20,10,11,26,23,14,29,6,5,8,17)_{\adfTFADJ}$,

$(0,24,10,26,16,29,2,17,14,28,8,22,12,6)_{\adfTFADJ}$,

$(18,24,11,8,23,6,0,22,2,5,26,29,20,17)_{\adfTFADJ}$,

\adfLgap
$(5,2,26,10,18,8,4,16,25,14,20,17,21,6)_{\adfTFAFH}$,

$(9,21,18,15,4,17,22,11,12,24,13,10,0,14)_{\adfTFAFH}$,

$(10,6,15,2,3,27,13,11,17,0,8,23,12,7)_{\adfTFAFH}$,

$(24,0,22,19,15,20,9,26,14,13,8,28,12,5)_{\adfTFAFH}$,

$(21,13,19,28,4,26,17,18,1,14,7,5,15,23)_{\adfTFAFH}$,

\adfSgap
$(19,7,5,17,3,10,25,11,4,22,15,28,13,1)_{\adfTFAFH}$,

$(7,23,13,3,16,28,5,27,10,22,29,15,25,11)_{\adfTFAFH}$,

$(1,19,17,10,28,21,4,16,9,22,7,29,13,25)_{\adfTFAFH}$,

$(1,23,21,5,13,19,9,25,17,29,11,27,4,16)_{\adfTFAFH}$,

\adfLgap
$(12,18,17,28,24,15,7,19,22,3,23,9,27,2)_{\adfTFBBK}$,

$(18,8,15,7,5,2,24,9,22,1,27,28,11,4)_{\adfTFBBK}$,

$(15,22,17,10,11,19,23,8,6,12,13,26,20,2)_{\adfTFBBK}$,

$(15,13,23,22,25,14,28,6,9,0,4,26,17,20)_{\adfTFBBK}$,

$(5,21,7,19,23,22,20,0,25,12,8,9,15,2)_{\adfTFBBK}$,

\adfSgap
$(13,16,19,10,6,17,24,5,26,1,7,0,28,22)_{\adfTFBBK}$,

$(4,25,10,19,23,17,28,1,24,22,11,18,29,20)_{\adfTFBBK}$,

$(4,13,7,28,6,29,23,14,19,12,5,11,17,8)_{\adfTFBBK}$,

$(16,25,1,18,5,29,10,12,23,0,11,2,7,22)_{\adfTFBBK}$,

\adfLgap
$(12,19,17,15,23,14,4,21,1,25,5,2,6,3)_{\adfTFBCJ}$,

$(22,2,1,6,17,25,28,11,4,0,18,10,16,27)_{\adfTFBCJ}$,

$(6,0,21,20,28,12,2,25,11,15,10,14,26,17)_{\adfTFBCJ}$,

$(11,16,17,2,9,12,5,23,28,26,15,6,13,1)_{\adfTFBCJ}$,

$(22,9,10,12,1,11,21,23,26,20,4,25,0,8)_{\adfTFBCJ}$,

\adfSgap
$(1,13,0,27,9,18,19,2,7,20,25,8,21,15)_{\adfTFBCJ}$,

$(3,1,20,7,26,15,19,8,13,5,21,25,27,14)_{\adfTFBCJ}$,

$(3,9,21,1,23,17,25,14,19,11,27,15,29,7)_{\adfTFBCJ}$,

$(9,13,26,15,2,3,27,21,19,6,7,24,25,12)_{\adfTFBCJ}$,

\adfLgap
$(19,3,9,0,22,7,6,20,26,10,28,17,4,18)_{\adfTFBDI}$,

$(14,6,4,26,19,13,9,20,7,8,16,12,11,24)_{\adfTFBDI}$,

$(8,28,15,18,26,5,21,17,19,14,25,4,7,23)_{\adfTFBDI}$,

$(11,20,23,7,28,21,17,27,6,15,16,22,3,5)_{\adfTFBDI}$,

$(10,4,8,11,3,1,15,12,23,9,21,7,29,24)_{\adfTFBDI}$,

\adfSgap
$(1,25,0,11,20,24,19,7,17,26,23,2,6,3)_{\adfTFBDI}$,

$(25,18,13,5,14,11,7,12,6,29,2,0,24,19)_{\adfTFBDI}$,

$(1,12,9,6,8,5,13,23,0,26,24,17,20,18)_{\adfTFBDI}$,

$(12,18,14,13,21,24,8,29,19,27,0,6,7,15)_{\adfTFBDI}$

\adfLgap
$(5,15,17,7,4,20,28,16,26,25,9,27,18,2)_{\adfTFBEH}$,

$(25,27,24,6,26,5,1,3,2,8,11,18,10,7)_{\adfTFBEH}$,

$(11,24,4,0,12,6,22,19,14,18,17,26,8,1)_{\adfTFBEH}$,

$(26,28,21,7,9,15,24,3,6,14,12,25,13,19)_{\adfTFBEH}$,

$(5,3,18,22,27,28,23,0,25,10,19,2,17,11)_{\adfTFBEH}$,

\adfSgap
$(5,11,25,9,28,22,10,21,2,4,8,27,16,14)_{\adfTFBEH}$,

$(10,17,21,14,3,22,16,4,28,29,15,26,23,7)_{\adfTFBEH}$,

$(5,10,8,4,22,26,28,19,29,13,23,9,20,16)_{\adfTFBEH}$,

$(11,17,1,27,23,22,20,15,4,16,28,2,29,3)_{\adfTFBEH}$,

\adfLgap
$(4,7,25,26,27,15,23,10,9,12,28,14,16,22)_{\adfTFBFG}$,

$(26,12,16,9,19,11,8,1,5,6,24,18,27,13)_{\adfTFBFG}$,

$(3,5,14,22,1,20,26,11,7,17,12,2,9,10)_{\adfTFBFG}$,

$(13,0,11,8,5,23,27,2,18,14,26,22,4,7)_{\adfTFBFG}$,

$(7,18,1,8,24,2,19,5,11,10,17,27,25,3)_{\adfTFBFG}$,

\adfSgap
$(21,22,0,15,12,4,27,24,28,17,3,5,16,9)_{\adfTFBFG}$,

$(22,27,11,12,3,16,18,10,15,28,6,4,23,21)_{\adfTFBFG}$,

$(9,21,18,6,29,15,24,4,23,0,3,10,12,5)_{\adfTFBFG}$,

$(9,27,3,15,17,24,16,6,11,18,28,0,10,29)_{\adfTFBFG}$,

\adfLgap
$(19,17,3,15,13,16,9,23,7,24,12,14,21,2)_{\adfTFCDH}$,

$(11,25,16,17,18,24,2,1,19,28,8,3,6,10)_{\adfTFCDH}$,

$(4,14,9,10,15,11,20,2,5,7,26,1,6,27)_{\adfTFCDH}$,

$(6,15,26,14,17,0,9,2,19,18,7,4,12,28)_{\adfTFCDH}$,

$(11,19,22,12,21,6,9,27,1,3,25,16,29,20)_{\adfTFCDH}$,

\adfSgap
$(22,11,16,4,28,20,17,10,12,26,23,24,8,5)_{\adfTFCDH}$,

$(10,22,16,8,17,23,29,26,4,6,5,0,28,14)_{\adfTFCDH}$,

$(10,29,28,5,4,22,24,2,16,23,18,17,12,11)_{\adfTFCDH}$,

$(11,29,14,0,23,5,17,6,20,4,28,16,18,2)_{\adfTFCDH}$,

\adfLgap
$(21,6,10,11,5,2,16,3,8,24,17,27,23,4)_{\adfTFCFF}$,

$(2,6,18,25,7,27,8,14,26,19,4,28,24,15)_{\adfTFCFF}$,

$(11,16,4,8,0,17,25,22,6,23,15,10,7,26)_{\adfTFCFF}$,

$(4,22,17,13,16,14,27,24,0,3,5,10,19,15)_{\adfTFCFF}$,

$(7,24,8,9,19,18,0,13,20,12,11,26,23,2)_{\adfTFCFF}$,

\adfSgap
$(1,9,17,7,3,15,27,19,25,11,13,23,14,21)_{\adfTFCFF}$,

$(25,1,5,29,11,17,8,15,7,27,13,21,3,9)_{\adfTFCFF}$,

$(9,19,27,21,2,11,5,26,3,15,13,7,23,29)_{\adfTFCFF}$,

$(19,25,17,23,5,7,21,15,1,13,29,20,27,3)_{\adfTFCFF}$,

\adfLgap
$(19,9,5,1,12,26,15,22,16,10,23,21,7,28)_{\adfTFDDG}$,

$(13,3,2,28,18,19,24,7,14,5,10,25,23,15)_{\adfTFDDG}$,

$(10,20,12,3,25,22,11,15,9,4,0,1,21,5)_{\adfTFDDG}$,

$(9,23,7,20,18,26,4,12,0,24,6,22,2,3)_{\adfTFDDG}$,

$(6,1,20,17,23,13,22,19,14,8,15,21,18,11)_{\adfTFDDG}$,

\adfSgap
$(11,14,24,23,16,8,10,26,12,2,6,5,18,17)_{\adfTFDDG}$,

$(14,20,23,22,29,28,26,0,18,8,17,5,2,4)_{\adfTFDDG}$,

$(11,29,10,17,0,23,5,28,2,14,24,20,8,26)_{\adfTFDDG}$,

$(20,29,2,16,17,22,8,12,23,6,26,5,4,11)_{\adfTFDDG}$,

\adfLgap
$(11,7,6,8,15,18,4,9,19,28,27,0,10,26)_{\adfTFDEF}$,

$(24,28,2,17,19,5,14,13,10,25,16,3,8,21)_{\adfTFDEF}$,

$(26,9,28,5,21,8,16,27,13,22,17,18,2,11)_{\adfTFDEF}$,

$(1,22,24,27,18,5,9,7,23,25,12,4,6,2)_{\adfTFDEF}$,

$(28,13,25,20,26,22,7,14,24,17,9,8,11,23)_{\adfTFDEF}$,

\adfSgap
$(23,9,1,6,24,17,11,14,3,26,15,0,12,18)_{\adfTFDEF}$,

$(0,9,25,17,20,6,27,12,29,18,24,11,21,15)_{\adfTFDEF}$,

$(0,5,24,19,11,17,27,3,18,21,6,12,7,29)_{\adfTFDEF}$,

$(5,23,13,18,6,15,24,12,3,8,27,21,2,29)_{\adfTFDEF}$

\adfLgap
\noindent under the action of the mapping $x \mapsto x + 6 \adfmod{30}$
for only the first five graphs in each design. \eproof

\ADFvfyParStart{30, \{\{30,1,30\}\}, -1, -1, \{5,\{\{30,5,6\}\}\}} 

{\lemma \label{lem:theta15 multipartite}
There exist decompositions into $\Theta(a,b,c)$ of the complete bipartite graphs
$K_{15,15}$, $K_{15,20}$ and $K_{15,25}$ for $a + b + c = 15$ with $a$, $b$, $c$ all odd.
There exist decompositions into $\Theta(a,b,c)$ of the complete multipartite graphs
$K_{5,5,5}$ and $K_{5,5,5,5}$ for $a + b + c = 15$ with $a$, $b$, $c$ not all odd.}\\

\noindent {\bf Proof}~ The decompositions are presented in the Appendix. \eproof\\

Theorem~\ref{thm:theta15} follows from Theorem~\ref{thm:2e+1},
Lemmas~\ref{lem:theta15 designs} and \ref{lem:theta15 multipartite}, and
Proposition~\ref{prop:bipartite-alt} (with $(r,s,e)$ = $(15,15,0)$, $(15,15,1)$, $(15,20,1)$ or $(15,25,0)$) if $a$, $b$ and $c$ are all odd,
Proposition~\ref{prop:theta15} otherwise.\eproof


\newcommand{\ADFrefGDD}{\noindent
~[Hanani1975]                 for 3-GDD  $g^u$ \\
~[ColbournHoffmanRees1992]    for 3-GDD  $g^u m^1$ \\
~[BrouwerSchrijverHanani1977] for 4-GDD  $g^u$ \\
~[GeLing2004]                 for 4-GDD  $g^u m^1$ \\
~[Rees1993]                   for 3-RGDD $g^u$ \\
~[GeLing2005]                 for 4-RGDD $g^u$ \\
~[GeLing2005]                 for 5-RGDD $g^u$ \\ }

%
%

\newpage
\appendix

\section{Theta graphs with 10 edges}
\label{sec:theta10}

\noindent {\bf Proof of Lemma \ref{lem:theta10 designs}}~

\noindent {\boldmath $K_{16}$}~
Let the vertex set be $Z_{15} \cup \{\infty\}$. The decompositions consist of the graphs

\adfLgap
$(\infty,11,5,9,2,10,8,7,3)_{\adfTAABG}$,

$(11,12,14,0,2,6,3,8,9)_{\adfTAABG}$,

$(14,1,3,0,4,10,7,2,11)_{\adfTAABG}$,

$(0,3,10,1,9,4,13,7,\infty)_{\adfTAABG}$,

\adfLgap
$(\infty,7,10,1,13,14,4,2,0)_{\adfTAACF}$,

$(1,8,4,12,2,3,9,5,0)_{\adfTAACF}$,

$(8,3,10,6,2,5,4,11,1)_{\adfTAACF}$,

$(0,9,1,\infty,3,14,2,7,11)_{\adfTAACF}$,

\adfLgap
$(\infty,5,8,14,12,1,11,3,6)_{\adfTAADE}$,

$(10,7,4,0,13,1,14,9,6)_{\adfTAADE}$,

$(13,14,9,1,7,3,0,2,\infty)_{\adfTAADE}$,

$(0,11,10,3,7,14,2,12,13)_{\adfTAADE}$,

\adfLgap
$(\infty,6,3,11,4,14,10,1,5)_{\adfTABBF}$,

$(4,5,13,12,2,14,6,9,8)_{\adfTABBF}$,

$(0,8,12,2,5,14,3,11,13)_{\adfTABBF}$,

$(0,12,13,\infty,14,1,2,11,7)_{\adfTABBF}$,

\adfLgap 
$(\infty,12,6,8,4,14,11,3,10)_{\adfTABCE}$,

$(4,0,3,11,10,9,12,1,13)_{\adfTABCE}$,

$(9,14,10,0,12,3,7,2,1)_{\adfTABCE}$,

$(1,8,3,10,2,11,0,\infty,7)_{\adfTABCE}$,

\adfLgap 
$(\infty,14,11,7,8,10,4,5,12)_{\adfTABDD}$,

$(3,7,13,0,11,12,9,14,1)_{\adfTABDD}$,

$(2,0,6,5,13,1,14,3,\infty)_{\adfTABDD}$,

$(1,7,3,6,13,14,9,0,5)_{\adfTABDD}$,

\adfLgap 
$(\infty,5,7,9,8,2,6,1,0)_{\adfTACCD}$,

$(6,13,14,3,10,1,12,5,4)_{\adfTACCD}$,

$(13,9,6,2,0,\infty,12,7,4)_{\adfTACCD}$,

$(2,8,0,9,13,5,1,4,6)_{\adfTACCD}$

\adfLgap
\noindent under the action of the mapping $x \mapsto x + 5$ (mod 15), $\infty \mapsto \infty$.

\ADFvfyParStart{16, \{\{15,3,5\},\{1,1,1\}\}, 15, -1, -1} 



\noindent {\boldmath $K_{25}$}~
Let the vertex set be $Z_{25}$. The decompositions consist of the graphs

\adfLgap 
$(9,10,6,11,19,22,20,7,8)_{\adfTAABG}$,

$(21,20,5,3,12,1,6,23,9)_{\adfTAABG}$,

$(11,4,2,21,8,0,3,1,14)_{\adfTAABG}$,

$(0,4,9,5,2,1,7,3,8)_{\adfTAABG}$,

$(0,7,11,6,3,4,10,2,12)_{\adfTAABG}$,

$(2,14,8,9,17,3,18,5,23)_{\adfTAABG}$,

\adfLgap 
$(0,4,22,8,1,2,3,14,6)_{\adfTAACF}$,

$(14,7,8,2,9,15,3,10,0)_{\adfTAACF}$,

$(13,14,1,12,21,9,6,0,5)_{\adfTAACF}$,

$(8,12,6,22,18,15,17,11,0)_{\adfTAACF}$,

$(13,4,22,19,8,10,6,16,12)_{\adfTAACF}$,

$(0,8,16,11,14,15,7,19,1)_{\adfTAACF}$,

\adfLgap 
$(0,18,21,15,14,8,11,20,7)_{\adfTAADE}$,

$(2,3,22,11,9,19,17,21,14)_{\adfTAADE}$,

$(19,9,14,5,15,7,3,22,21)_{\adfTAADE}$,

$(14,23,10,17,20,0,1,18,11)_{\adfTAADE}$,

$(15,3,7,23,8,13,14,11,1)_{\adfTAADE}$,

$(2,24,0,5,16,12,6,1,17)_{\adfTAADE}$,

\adfLgap 
$(0,8,13,19,7,10,9,20,21)_{\adfTABBF}$,

$(20,6,1,23,12,2,19,7,9)_{\adfTABBF}$,

$(20,3,4,18,6,15,5,0,2)_{\adfTABBF}$,

$(5,23,1,17,13,4,0,18,2)_{\adfTABBF}$,

$(9,18,11,16,13,19,14,4,7)_{\adfTABBF}$,

$(1,4,16,22,7,6,17,12,21)_{\adfTABBF}$,

\adfLgap 
$(0,9,15,8,6,21,5,12,7)_{\adfTABCE}$,

$(2,7,18,21,20,11,0,5,3)_{\adfTABCE}$,

$(16,23,8,9,3,11,19,15,14)_{\adfTABCE}$,

$(4,23,19,3,16,9,11,21,17)_{\adfTABCE}$,

$(12,2,19,9,0,22,11,8,5)_{\adfTABCE}$,

$(0,16,17,18,4,13,12,24,10)_{\adfTABCE}$,

\adfLgap 
$(0,14,11,18,23,13,7,1,10)_{\adfTABDD}$,

$(22,8,0,19,21,20,7,5,17)_{\adfTABDD}$,

$(22,20,9,8,19,23,2,4,10)_{\adfTABDD}$,

$(4,19,9,22,18,1,5,0,6)_{\adfTABDD}$,

$(9,1,3,0,23,22,17,6,13)_{\adfTABDD}$,

$(1,8,11,6,10,2,17,16,24)_{\adfTABDD}$,

\adfLgap 
$(0,9,22,7,10,14,18,6,12)_{\adfTACCD}$,

$(22,13,23,14,8,3,17,16,9)_{\adfTACCD}$,

$(3,5,21,18,9,6,12,10,17)_{\adfTACCD}$,

$(8,16,2,6,19,20,10,5,21)_{\adfTACCD}$,

$(14,18,5,22,16,10,4,15,1)_{\adfTACCD}$,

$(1,7,9,21,17,24,3,0,19)_{\adfTACCD}$

\adfLgap
\noindent under the action of the mapping $x \mapsto x + 5$ (mod 25).

\ADFvfyParStart{25, \{\{25,5,5\}\}, -1, -1, -1} 

\noindent {\boldmath $K_{36}$}~
Let the vertex set be $Z_{36}$. The decompositions consist of the graphs

\adfLgap 
$(2,3,8,17,24,30,25,29,18)_{\adfTAABG}$,

$(5,29,31,20,7,11,34,17,16)_{\adfTAABG}$,

$(12,28,31,10,24,20,34,6,18)_{\adfTAABG}$,

$(2,7,29,27,20,11,19,31,21)_{\adfTAABG}$,

$(20,8,19,12,14,13,4,9,25)_{\adfTAABG}$,

$(0,10,33,15,6,2,5,11,17)_{\adfTAABG}$,

$(1,19,35,29,4,22,3,6,26)_{\adfTAABG}$,

\adfLgap 
$(17,7,20,31,33,26,25,12,15)_{\adfTAACF}$,

$(0,28,14,12,4,33,19,2,26)_{\adfTAACF}$,

$(1,14,15,8,7,12,13,24,29)_{\adfTAACF}$,

$(14,20,4,3,5,1,9,33,2)_{\adfTAACF}$,

$(28,16,13,7,15,33,31,22,2)_{\adfTAACF}$,

$(23,10,26,11,19,17,34,6,35)_{\adfTAACF}$,

$(1,11,26,31,34,30,4,13,32)_{\adfTAACF}$,

\adfLgap 
$(0,15,12,25,31,2,28,18,21)_{\adfTAADE}$,

$(0,8,6,17,7,32,21,2,9)_{\adfTAADE}$,

$(31,21,29,2,26,24,4,23,34)_{\adfTAADE}$,

$(1,19,32,29,15,13,14,30,34)_{\adfTAADE}$,

$(9,28,29,0,3,1,33,11,10)_{\adfTAADE}$,

$(19,32,24,33,18,10,2,7,30)_{\adfTAADE}$,

$(2,19,24,18,11,35,1,16,7)_{\adfTAADE}$,

\adfLgap 
$(0,33,6,2,34,26,19,14,25)_{\adfTABCE}$,

$(19,29,9,2,14,22,23,31,8)_{\adfTABCE}$,

$(14,25,21,1,20,23,29,17,22)_{\adfTABCE}$,

$(34,33,23,12,31,4,18,1,15)_{\adfTABCE}$,

$(1,20,24,23,28,12,5,32,11)_{\adfTABCE}$,

$(24,28,12,25,31,6,14,16,27)_{\adfTABCE}$,

$(0,26,10,7,11,33,35,15,3)_{\adfTABCE}$,

\adfLgap 
$(0,1,22,4,15,30,28,3,33)_{\adfTACCD}$,

$(20,22,21,32,29,11,27,10,1)_{\adfTACCD}$,

$(30,13,16,0,22,20,24,19,33)_{\adfTACCD}$,

$(23,14,33,27,15,19,24,28,26)_{\adfTACCD}$,

$(4,30,14,23,16,10,7,19,17)_{\adfTACCD}$,

$(13,21,2,33,32,19,28,11,31)_{\adfTACCD}$,

$(1,6,2,3,29,10,32,23,9)_{\adfTACCD}$

\adfLgap
\noindent under the action of the mapping $x \mapsto x + 4$ (mod 36).

\ADFvfyParStart{36, \{\{36,9,4\}\}, -1, -1, -1} 

\noindent {\boldmath $K_{40}$}~
Let the vertex set be $Z_{39} \cup \{\infty\}$. The decompositions consist of the graphs

\adfLgap 
$(\infty,25,33,26,38,32,30,21,7)_{\adfTAABG}$,

$(2,33,36,37,22,20,28,24,14)_{\adfTAABG}$,

$(17,22,3,6,38,13,2,20,23)_{\adfTAABG}$,

$(30,31,36,18,35,12,8,28,25)_{\adfTAABG}$,

$(19,35,9,32,2,3,10,27,25)_{\adfTAABG}$,

$(3,31,18,21,37,10,1,23,38)_{\adfTAABG}$,

\adfLgap 
$(\infty,38,28,29,33,3,16,8,9)_{\adfTAACF}$,

$(27,12,0,26,16,33,11,3,13)_{\adfTAACF}$,

$(22,35,34,1,6,29,17,13,33)_{\adfTAACF}$,

$(16,25,31,27,14,8,22,38,18)_{\adfTAACF}$,

$(23,20,34,13,2,9,6,12,37)_{\adfTAACF}$,

$(38,3,14,34,27,6,1,37,8)_{\adfTAACF}$,

\adfLgap 
$(\infty,22,33,0,32,26,15,3,10)_{\adfTAADE}$,

$(31,29,34,16,27,30,8,5,0)_{\adfTAADE}$,

$(8,21,26,1,24,7,29,37,6)_{\adfTAADE}$,

$(19,9,33,29,13,34,28,37,32)_{\adfTAADE}$,

$(23,4,9,30,21,34,2,8,17)_{\adfTAADE}$,

$(1,3,5,24,37,21,29,14,2)_{\adfTAADE}$,

\adfLgap 
$(\infty,26,1,6,31,38,23,25,12)_{\adfTABCE}$,

$(9,5,38,26,0,30,36,19,21)_{\adfTABCE}$,

$(14,20,13,25,10,32,33,34,0)_{\adfTABCE}$,

$(28,27,24,37,16,25,38,2,0)_{\adfTABCE}$,

$(10,13,29,18,33,22,6,38,3)_{\adfTABCE}$,

$(33,5,1,3,14,2,19,25,17)_{\adfTABCE}$,

\adfLgap 
$(\infty,3,22,7,26,5,6,20,37)_{\adfTACCD}$,

$(11,12,13,31,35,32,37,24,22)_{\adfTACCD}$,

$(7,11,6,33,30,12,18,24,21)_{\adfTACCD}$,

$(37,15,23,22,32,6,34,4,8)_{\adfTACCD}$,

$(27,1,13,32,3,8,35,7,34)_{\adfTACCD}$,

$(0,5,22,38,23,32,11,1,9)_{\adfTACCD}$

\adfLgap
\noindent under the action of the mapping $x \mapsto x + 3$ (mod 39), $\infty \mapsto \infty$.

\ADFvfyParStart{40, \{\{39,13,3\},\{1,1,1\}\}, 39, -1, -1} 

\noindent {\boldmath $K_{41}$}~
Let the vertex set be $Z_{41}$. The decompositions consist of the graphs

\adfLgap 
$(0,38,17,18,11,7,1,2,10)_{\adfTAABG}$,

$(0,2,11,5,17,3,18,8,24)_{\adfTAABG}$,

\adfLgap 
$(0,5,7,37,15,32,1,2,8)_{\adfTAACF}$,

$(0,2,4,16,13,5,21,1,20)_{\adfTAACF}$,

\adfLgap 
$(0,21,5,20,33,7,1,2,4)_{\adfTAADE}$,

$(0,3,4,12,21,16,2,24,13)_{\adfTAADE}$,

\adfLgap 
$(0,20,31,12,33,8,1,2,4)_{\adfTABCE}$,

$(0,1,4,5,23,14,8,25,10)_{\adfTABCE}$,

\adfLgap 
$(0,7,19,1,9,24,2,3,10)_{\adfTACCD}$,

$(0,1,4,9,10,21,12,28,14)_{\adfTACCD}$

\adfLgap
\noindent under the action of the mapping $x \mapsto x + 1$ (mod 41).

\ADFvfyParStart{41, \{\{41,41,1\}\}, -1, -1, -1} 

\noindent {\boldmath $K_{45}$}~
Let the vertex set be $Z_{45}$. The decompositions consist of the graphs

\adfLgap 
$(0,12,31,24,39,14,13,17,30)_{\adfTAABG}$,

$(43,41,20,5,14,12,7,19,38)_{\adfTAABG}$,

$(26,3,17,38,21,41,36,32,35)_{\adfTAABG}$,

$(22,32,7,9,14,38,34,1,8)_{\adfTAABG}$,

$(34,17,16,20,1,7,21,3,43)_{\adfTAABG}$,

$(27,19,6,4,39,11,5,15,30)_{\adfTAABG}$,

$(5,0,16,8,19,22,20,11,4)_{\adfTAABG}$,

$(20,19,0,16,8,14,32,38,3)_{\adfTAABG}$,

$(13,12,5,22,6,29,31,30,28)_{\adfTAABG}$,

$(0,18,33,17,6,14,11,1,31)_{\adfTAABG}$,

$(0,29,22,37,3,20,14,23,43)_{\adfTAABG}$,

\adfLgap 
$(0,6,14,5,40,27,42,20,41)_{\adfTAACF}$,

$(21,16,5,41,9,7,10,2,3)_{\adfTAACF}$,

$(23,10,17,21,33,37,8,6,29)_{\adfTAACF}$,

$(4,17,7,28,43,1,8,31,5)_{\adfTAACF}$,

$(3,40,8,34,23,11,41,10,28)_{\adfTAACF}$,

$(8,4,23,31,5,40,18,27,11)_{\adfTAACF}$,

$(4,25,18,9,1,28,35,15,32)_{\adfTAACF}$,

$(2,27,21,34,37,10,12,6,4)_{\adfTAACF}$,

$(40,29,38,9,10,6,17,3,24)_{\adfTAACF}$,

$(27,26,19,9,36,12,24,7,34)_{\adfTAACF}$,

$(0,28,44,29,4,38,21,7,2)_{\adfTAACF}$,

\adfLgap 
$(0,13,10,42,21,24,22,8,31)_{\adfTAADE}$,

$(38,0,35,32,9,33,16,13,30)_{\adfTAADE}$,

$(22,2,12,11,43,29,25,26,20)_{\adfTAADE}$,

$(16,35,12,38,15,26,33,13,34)_{\adfTAADE}$,

$(35,4,8,0,29,40,12,42,30)_{\adfTAADE}$,

$(7,39,41,22,38,34,31,6,21)_{\adfTAADE}$,

$(37,42,21,10,6,29,18,34,3)_{\adfTAADE}$,

$(7,5,8,17,38,18,3,13,39)_{\adfTAADE}$,

$(35,21,29,14,4,6,37,34,26)_{\adfTAADE}$,

$(43,31,4,16,40,0,22,34,29)_{\adfTAADE}$,

$(1,3,14,21,44,7,15,22,39)_{\adfTAADE}$,

\adfLgap 
$(0,21,41,11,38,34,5,43,33)_{\adfTABCE}$,

$(25,40,26,15,34,23,27,37,12)_{\adfTABCE}$,

$(28,12,26,41,34,10,22,19,17)_{\adfTABCE}$,

$(41,42,33,29,3,31,22,40,16)_{\adfTABCE}$,

$(19,21,5,8,4,11,33,34,22)_{\adfTABCE}$,

$(32,17,21,19,34,25,27,5,18)_{\adfTABCE}$,

$(15,9,0,14,38,35,16,10,18)_{\adfTABCE}$,

$(15,39,19,10,33,12,23,43,29)_{\adfTABCE}$,

$(30,28,21,22,43,13,10,34,31)_{\adfTABCE}$,

$(36,37,21,31,29,4,22,35,23)_{\adfTABCE}$,

$(2,4,9,32,26,28,23,7,31)_{\adfTABCE}$,

\adfLgap 
$(0,3,36,24,38,41,11,37,12)_{\adfTACCD}$,

$(13,42,29,11,18,28,24,19,9)_{\adfTACCD}$,

$(39,19,2,37,10,40,5,9,26)_{\adfTACCD}$,

$(33,7,34,36,21,6,14,18,30)_{\adfTACCD}$,

$(25,3,30,39,17,21,1,40,2)_{\adfTACCD}$,

$(30,42,10,0,3,23,38,27,4)_{\adfTACCD}$,

$(8,1,16,41,40,23,6,19,4)_{\adfTACCD}$,

$(23,20,6,16,29,21,7,1,37)_{\adfTACCD}$,

$(27,9,15,34,42,40,0,3,7)_{\adfTACCD}$,

$(15,41,13,28,31,7,1,20,19)_{\adfTACCD}$,

$(0,24,32,27,39,7,23,17,38)_{\adfTACCD}$

\adfLgap
\noindent under the action of the mapping $x \mapsto x + 5$ (mod 45).

\ADFvfyParStart{45, \{\{45,9,5\}\}, -1, -1, -1} 

\noindent {\boldmath $K_{56}$}~
Let the vertex set be $Z_{55} \cup \{\infty\}$. The decompositions consist of the graphs

\adfLgap 
$(\infty,12,18,29,24,51,31,1,7)_{\adfTAABG}$,

$(33,32,39,54,18,10,41,8,51)_{\adfTAABG}$,

$(39,35,41,42,28,21,11,50,27)_{\adfTAABG}$,

$(37,51,30,19,43,29,28,52,36)_{\adfTAABG}$,

$(4,48,2,37,0,51,9,29,35)_{\adfTAABG}$,

$(9,0,24,30,41,43,16,34,15)_{\adfTAABG}$,

$(16,4,30,3,48,8,43,46,14)_{\adfTAABG}$,

$(27,14,0,1,2,5,3,7,11)_{\adfTAABG}$,

$(0,1,10,2,11,3,8,4,12)_{\adfTAABG}$,

$(0,3,20,5,4,11,6,14,25)_{\adfTAABG}$,

$(0,12,22,17,2,13,4,20,38)_{\adfTAABG}$,

$(0,23,30,26,2,14,31,8,35)_{\adfTAABG}$,

$(0,36,\infty,42,3,22,47,19,2)_{\adfTAABG}$,

$(1,18,34,38,9,39,7,27,48)_{\adfTAABG}$,

\adfLgap 
$(\infty,37,23,24,21,50,9,1,43)_{\adfTAACF}$,

$(51,8,16,15,40,48,39,11,17)_{\adfTAACF}$,

$(47,1,45,12,31,44,46,38,27)_{\adfTAACF}$,

$(47,4,49,40,8,34,45,10,48)_{\adfTAACF}$,

$(43,33,20,32,10,53,19,3,9)_{\adfTAACF}$,

$(6,47,13,11,43,25,40,9,5)_{\adfTAACF}$,

$(38,25,7,4,2,40,21,15,52)_{\adfTAACF}$,

$(4,39,37,2,1,5,0,3,6)_{\adfTAACF}$,

$(0,7,10,2,16,1,6,15,14)_{\adfTAACF}$,

$(0,21,25,9,28,1,2,5,46)_{\adfTAACF}$,

$(0,29,31,14,32,1,11,4,12)_{\adfTAACF}$,

$(0,49,53,13,\infty,4,8,3,17)_{\adfTAACF}$,

$(1,18,19,47,22,32,2,17,38)_{\adfTAACF}$,

$(3,28,36,32,26,49,19,9,14)_{\adfTAACF}$,

\adfLgap 
$(\infty,7,54,36,25,41,14,26,17)_{\adfTAADE}$,

$(25,39,45,44,3,34,9,22,27)_{\adfTAADE}$,

$(42,22,48,28,47,6,9,31,43)_{\adfTAADE}$,

$(52,11,23,10,3,44,12,27,21)_{\adfTAADE}$,

$(4,6,37,54,25,3,47,5,8)_{\adfTAADE}$,

$(32,33,5,48,4,1,10,22,51)_{\adfTAADE}$,

$(22,18,4,38,11,13,53,30,9)_{\adfTAADE}$,

$(46,13,43,18,8,11,12,15,36)_{\adfTAADE}$,

$(0,1,49,4,40,22,8,21,16)_{\adfTAADE}$,

$(42,31,14,6,44,50,2,41,45)_{\adfTAADE}$,

$(38,21,0,32,15,54,14,17,51)_{\adfTAADE}$,

$(35,13,20,38,44,43,27,29,9)_{\adfTAADE}$,

$(4,35,20,9,5,11,7,38,40)_{\adfTAADE}$,

$(0,26,28,\infty,50,45,47,54,49)_{\adfTAADE}$,

\adfLgap 
$(\infty,28,1,44,41,17,14,32,36)_{\adfTABCE}$,

$(38,24,49,19,53,23,52,16,47)_{\adfTABCE}$,

$(53,15,47,12,22,3,36,4,39)_{\adfTABCE}$,

$(53,7,29,14,18,45,54,1,2)_{\adfTABCE}$,

$(46,4,15,52,41,30,40,42,5)_{\adfTABCE}$,

$(26,48,31,41,15,23,42,2,10)_{\adfTABCE}$,

$(49,24,34,10,30,44,25,31,41)_{\adfTABCE}$,

$(13,12,40,50,54,4,30,21,41)_{\adfTABCE}$,

$(34,30,0,7,27,33,23,35,47)_{\adfTABCE}$,

$(49,52,43,16,28,2,32,34,6)_{\adfTABCE}$,

$(11,49,32,18,36,25,20,21,37)_{\adfTABCE}$,

$(32,38,45,48,13,53,40,36,6)_{\adfTABCE}$,

$(46,45,48,34,26,32,33,\infty,5)_{\adfTABCE}$,

$(1,12,45,49,8,35,3,5,19)_{\adfTABCE}$,

\adfLgap 
$(\infty,25,28,2,12,52,5,24,10)_{\adfTACCD}$,

$(41,5,9,34,48,1,17,23,47)_{\adfTACCD}$,

$(27,38,10,1,23,33,2,7,52)_{\adfTACCD}$,

$(28,8,17,5,50,19,53,1,23)_{\adfTACCD}$,

$(45,28,54,14,39,26,47,34,7)_{\adfTACCD}$,

$(20,34,30,18,38,14,25,45,29)_{\adfTACCD}$,

$(35,50,6,25,14,16,48,44,1)_{\adfTACCD}$,

$(31,44,12,26,45,43,11,27,36)_{\adfTACCD}$,

$(25,38,36,9,33,10,1,16,37)_{\adfTACCD}$,

$(19,39,41,46,28,47,9,52,17)_{\adfTACCD}$,

$(37,45,11,12,34,27,40,33,46)_{\adfTACCD}$,

$(46,6,3,12,23,51,16,33,17)_{\adfTACCD}$,

$(16,25,12,14,54,17,19,\infty,41)_{\adfTACCD}$,

$(2,14,19,15,50,33,34,28,48)_{\adfTACCD}$

\adfLgap
\noindent under the action of the mapping $x \mapsto x + 5$ (mod 55), $\infty \mapsto \infty$.

\ADFvfyParStart{56, \{\{55,11,5\},\{1,1,1\}\}, 55, -1, -1} 

\noindent {\boldmath $K_{65}$}~
Let the vertex set be $Z_{65}$. The decompositions consist of the graphs

\adfLgap 
$(0,4,41,14,12,22,27,7,35)_{\adfTAABG}$,

$(55,56,52,25,27,54,34,39,26)_{\adfTAABG}$,

$(35,59,8,13,23,55,15,14,20)_{\adfTAABG}$,

$(42,49,60,38,12,11,20,27,62)_{\adfTAABG}$,

$(29,63,52,54,16,18,46,31,7)_{\adfTAABG}$,

$(17,26,6,42,14,31,38,50,29)_{\adfTAABG}$,

$(39,6,22,47,13,49,41,46,24)_{\adfTAABG}$,

$(28,27,3,5,21,2,8,9,13)_{\adfTAABG}$,

$(56,31,10,8,3,52,33,11,17)_{\adfTAABG}$,

$(22,43,56,5,63,4,57,30,40)_{\adfTAABG}$,

$(56,17,14,48,24,33,3,60,40)_{\adfTAABG}$,

$(33,59,44,48,29,20,39,4,6)_{\adfTAABG}$,

$(3,55,21,23,11,4,14,50,61)_{\adfTAABG}$,

$(50,48,21,62,26,0,28,45,30)_{\adfTAABG}$,

$(31,21,35,63,19,3,26,47,18)_{\adfTAABG}$,

$(0,52,60,22,7,39,17,50,34)_{\adfTAABG}$,

\adfLgap 
$(0,17,13,46,61,52,57,6,25)_{\adfTAACF}$,

$(29,52,24,20,15,54,16,21,50)_{\adfTAACF}$,

$(57,61,38,19,9,11,25,48,46)_{\adfTAACF}$,

$(49,55,48,61,17,39,54,41,2)_{\adfTAACF}$,

$(19,12,9,3,7,6,55,23,38)_{\adfTAACF}$,

$(33,27,60,55,8,0,26,44,62)_{\adfTAACF}$,

$(43,51,2,29,55,40,58,31,62)_{\adfTAACF}$,

$(4,61,18,49,15,60,30,8,44)_{\adfTAACF}$,

$(7,9,10,54,47,60,17,62,13)_{\adfTAACF}$,

$(15,24,42,8,25,7,17,9,37)_{\adfTAACF}$,

$(26,33,5,12,61,13,3,44,47)_{\adfTAACF}$,

$(42,13,57,15,21,54,29,60,20)_{\adfTAACF}$,

$(30,21,14,3,49,28,63,26,46)_{\adfTAACF}$,

$(15,43,26,42,44,20,3,46,23)_{\adfTAACF}$,

$(47,58,1,63,6,9,16,28,19)_{\adfTAACF}$,

$(0,1,3,7,31,55,54,19,56)_{\adfTAACF}$,

\adfLgap 
$(0,15,23,13,31,62,42,25,61)_{\adfTAADE}$,

$(22,54,17,14,16,36,52,58,6)_{\adfTAADE}$,

$(57,47,32,63,17,0,51,7,29)_{\adfTAADE}$,

$(2,28,53,40,30,44,60,20,24)_{\adfTAADE}$,

$(11,6,17,10,51,54,20,50,62)_{\adfTAADE}$,

$(2,40,56,48,13,23,27,26,20)_{\adfTAADE}$,

$(20,8,15,2,0,14,46,23,51)_{\adfTAADE}$,

$(6,10,40,58,57,37,53,14,43)_{\adfTAADE}$,

$(29,22,42,10,49,9,58,43,0)_{\adfTAADE}$,

$(23,14,26,33,9,30,51,50,4)_{\adfTAADE}$,

$(24,45,36,46,4,35,44,7,48)_{\adfTAADE}$,

$(62,14,8,28,11,53,5,16,44)_{\adfTAADE}$,

$(15,57,29,17,21,14,58,53,59)_{\adfTAADE}$,

$(31,40,16,18,4,1,9,34,3)_{\adfTAADE}$,

$(56,49,32,3,36,9,28,53,64)_{\adfTAADE}$,

$(1,13,40,12,51,26,7,22,14)_{\adfTAADE}$,

\adfLgap 
$(0,18,23,30,52,28,43,6,1)_{\adfTABCE}$,

$(15,60,57,42,2,34,16,9,58)_{\adfTABCE}$,

$(15,60,36,4,29,32,50,56,24)_{\adfTABCE}$,

$(35,54,27,36,42,34,48,12,32)_{\adfTABCE}$,

$(0,7,46,9,54,51,37,18,11)_{\adfTABCE}$,

$(27,26,36,15,58,60,10,43,9)_{\adfTABCE}$,

$(9,13,57,11,34,37,22,38,27)_{\adfTABCE}$,

$(21,29,41,52,62,2,56,34,55)_{\adfTABCE}$,

$(2,54,46,51,16,0,48,35,43)_{\adfTABCE}$,

$(4,42,39,20,38,33,7,31,6)_{\adfTABCE}$,

$(24,33,9,17,23,11,58,2,30)_{\adfTABCE}$,

$(2,27,62,44,40,59,0,5,29)_{\adfTABCE}$,

$(30,0,56,41,25,26,54,53,14)_{\adfTABCE}$,

$(26,48,13,28,8,29,24,37,36)_{\adfTABCE}$,

$(56,38,41,48,29,13,12,53,11)_{\adfTABCE}$,

$(0,63,10,36,5,45,18,14,4)_{\adfTABCE}$,

\adfLgap 
$(0,61,29,34,36,43,27,4,62)_{\adfTACCD}$,

$(15,56,19,1,60,39,63,53,61)_{\adfTACCD}$,

$(7,13,52,48,33,10,41,4,56)_{\adfTACCD}$,

$(50,7,0,35,52,53,43,1,42)_{\adfTACCD}$,

$(13,35,42,16,45,56,14,60,61)_{\adfTACCD}$,

$(1,34,4,7,12,35,41,57,43)_{\adfTACCD}$,

$(24,56,1,62,52,28,42,37,49)_{\adfTACCD}$,

$(48,61,3,16,17,4,59,29,45)_{\adfTACCD}$,

$(4,33,15,55,37,48,24,63,58)_{\adfTACCD}$,

$(24,1,14,5,26,60,55,12,28)_{\adfTACCD}$,

$(21,23,53,29,35,11,2,31,19)_{\adfTACCD}$,

$(16,26,13,61,25,17,31,52,12)_{\adfTACCD}$,

$(36,19,38,25,4,48,14,63,45)_{\adfTACCD}$,

$(28,43,7,29,22,35,55,4,45)_{\adfTACCD}$,

$(2,23,57,60,20,32,50,40,35)_{\adfTACCD}$,

$(2,29,19,4,60,27,52,44,63)_{\adfTACCD}$

\adfLgap
\noindent under the action of the mapping $x \mapsto x + 5$ (mod 65).\eproof

\ADFvfyParStart{65, \{\{65,13,5\}\}, -1, -1, -1} 

~\\
\noindent {\bf Proof of Lemma \ref{lem:theta10 multipartite}}~

\noindent {\boldmath $K_{5,10}$}~
Let the vertex set be $\{0, 1, \dots, 14\}$ partitioned into $\{0, 3, 6, 9, 12\}$
and \{1, 2, 4, 5, 7, 8, 10, 11, 13, 14\}. The decompositions consist of the graphs

\adfLgap 
$(2,13,9,0,3,8,12,7,6)_{\adfTABBF}$

\adfLgap
\noindent under the action of the mapping $x \mapsto x + 3$ (mod 15), and

\ADFvfyParStart{15, \{\{15,5,3\}\}, -1, \{\{5,\{0,1,1\}\}\}, -1} 


\adfLgap 
$(5,7,0,3,4,12,9,2,6)_{\adfTABDD}$,

$(1,13,6,3,11,12,9,8,0)_{\adfTABDD}$,

$(11,8,6,9,7,3,0,1,12)_{\adfTABDD}$,

$(4,2,0,6,10,3,9,14,12)_{\adfTABDD}$,

$(10,14,0,9,13,3,12,5,6)_{\adfTABDD}$.


\ADFvfyParStart{15, \{\{15,1,15\}\}, -1, \{\{5,\{0,1,1\}\}\}, -1} 

\noindent {\boldmath $K_{10,10,10}$}~
Let the vertex set be $Z_{30}$ partitioned according to residue class modulo 3.
The decompositions consist of the graphs

\adfLgap 
$(15,28,17,1,0,4,9,16,6)_{\adfTAABG}$,

\adfLgap 
$(2,3,9,1,6,11,21,5,16)_{\adfTAACF}$,

\adfLgap 
$(0,17,23,1,3,4,9,8,27)_{\adfTAADE}$,

\adfLgap 
$(0,1,26,17,3,7,8,19,9)_{\adfTABCE}$,

\adfLgap 
$(0,19,28,15,1,8,5,13,3)_{\adfTACCD}$

\adfLgap
\noindent under the action of the mapping $x \mapsto x + 1$ (mod 30).

\ADFvfyParStart{30, \{\{30,30,1\}\}, -1, \{\{10,\{0,1,2\}\}\}, -1} 

\noindent {\boldmath $K_{5,5,5,5}$}~
Let the vertex set be $\{0, 1, \dots, 19\}$ partitioned into
$\{3j + i: j = 0, 1, 2, 3, 4\}$, $i = 0, 1, 2$, and $\{15, 16, 17, 18, 19\}$.
The decompositions consist of the graphs

\adfLgap 
$(18,0,2,1,5,4,11,6,15)_{\adfTAABG}$,

\adfLgap 
$(2,18,1,3,6,11,15,8,0)_{\adfTAACF}$,

\adfLgap 
$(0,17,5,13,2,14,1,15,4)_{\adfTAADE}$,

\adfLgap 
$(0,18,10,14,1,4,11,16,2)_{\adfTABCE}$,

\adfLgap 
$(0,16,7,8,4,2,5,15,9)_{\adfTACCD}$

\adfLgap
\noindent under the action of the mapping $x \mapsto x + 1 \adfmod{15}$ for $x < 15$,
$x \mapsto (x + 1 \adfmod{5}) + 15$ for $x \ge 15$.

\ADFvfyParStart{20, \{\{15,15,1\},\{5,5,1\}\}, -1, \{\{5,\{0,1,2\}\},\{5,\{3\}\}\}, -1} 

\noindent {\boldmath $K_{20,20,20,25}$}~
Let the vertex set be $\{0, 1, \dots, 84\}$ partitioned into
$\{3j + i: j = 0, 1, \dots, 19\}$, $i = 0, 1, 2$, and $\{60, 61, \dots, 84\}$.
The decompositions consist of the graphs
\adfLgap 
$(0,63,44,69,37,12,2,33,16)_{\adfTAABG}$,

$(26,72,1,64,20,80,25,44,10)_{\adfTAABG}$,

$(49,29,0,3,64,27,23,36,31)_{\adfTAABG}$,

$(30,61,52,60,13,47,76,34,53)_{\adfTAABG}$,

$(3,26,40,62,1,45,75,7,24)_{\adfTAABG}$,

$(24,61,16,31,23,83,20,25,47)_{\adfTAABG}$,

$(67,56,57,0,68,43,81,7,65)_{\adfTAABG}$,

$(46,59,6,80,4,8,40,62,49)_{\adfTAABG}$,

$(53,78,25,74,38,37,79,8,19)_{\adfTAABG}$,

\adfLgap 
$(0,10,22,66,46,74,29,78,44)_{\adfTAACF}$,

$(72,41,25,73,47,6,78,42,69)_{\adfTAACF}$,

$(5,43,63,54,4,15,41,57,66)_{\adfTAACF}$,

$(14,61,10,39,60,35,62,19,47)_{\adfTAACF}$,

$(7,54,53,10,70,5,58,80,47)_{\adfTAACF}$,

$(34,53,33,58,54,19,15,78,48)_{\adfTAACF}$,

$(76,39,38,69,22,59,79,57,49)_{\adfTAACF}$,

$(25,8,38,10,45,80,24,32,74)_{\adfTAACF}$,

$(32,72,77,16,3,5,42,60,8)_{\adfTAACF}$,

\adfLgap 
$(0,13,60,51,75,58,83,45,64)_{\adfTAADE}$,

$(58,71,44,78,2,17,6,28,45)_{\adfTAADE}$,

$(73,16,15,61,50,28,32,57,64)_{\adfTAADE}$,

$(21,34,66,48,60,65,24,67,42)_{\adfTAADE}$,

$(61,33,52,53,37,27,17,24,34)_{\adfTAADE}$,

$(72,45,47,79,5,33,52,35,13)_{\adfTAADE}$,

$(46,67,62,45,40,15,49,47,39)_{\adfTAADE}$,

$(4,24,27,70,40,41,46,18,74)_{\adfTAADE}$,

$(15,78,29,40,11,8,79,4,39)_{\adfTAADE}$,

\adfLgap 
$(0,23,19,52,64,60,58,72,42)_{\adfTABCE}$,

$(65,67,1,27,22,53,24,55,15)_{\adfTABCE}$,

$(39,9,31,68,22,65,25,50,71)_{\adfTABCE}$,

$(56,36,19,64,50,76,53,39,29)_{\adfTABCE}$,

$(74,32,39,41,16,43,59,57,81)_{\adfTABCE}$,

$(40,67,39,50,27,81,1,12,13)_{\adfTABCE}$,

$(80,69,30,44,6,2,78,21,53)_{\adfTABCE}$,

$(45,44,68,71,12,28,41,63,24)_{\adfTABCE}$,

$(57,50,46,83,16,23,28,26,72)_{\adfTABCE}$,

\adfLgap 
$(0,48,5,61,14,80,75,55,26)_{\adfTACCD}$,

$(39,51,81,50,56,68,44,79,37)_{\adfTACCD}$,

$(3,51,47,62,23,74,14,24,83)_{\adfTACCD}$,

$(16,11,33,60,69,45,57,32,48)_{\adfTACCD}$,

$(45,63,62,18,49,30,14,1,2)_{\adfTACCD}$,

$(34,23,54,71,32,76,82,2,15)_{\adfTACCD}$,

$(78,29,21,19,23,51,17,65,54)_{\adfTACCD}$,

$(52,38,81,31,20,82,59,77,15)_{\adfTACCD}$,

$(26,41,52,70,64,30,34,38,79)_{\adfTACCD}$

\adfLgap
\noindent under the action of the mapping $x \mapsto x + 2 \adfmod{60}$ for $x < 60$,
$x \mapsto (x - 60 + 5 \adfmod{25}) + 60$ for $x \ge 60$.

\ADFvfyParStart{85, \{\{60,30,2\},\{25,5,5\}\}, -1, \{\{20,\{0,1,2\}\},\{25,\{3\}\}\}, -1} 

\noindent {\boldmath $K_{5,5,5,5,5}$}~
Let the vertex set be $Z_{25}$ partitioned according to residue class modulo 5.
The decompositions consist of the graphs

\adfLgap 
$(4,8,21,1,0,2,9,3,17)_{\adfTAABG}$,

\adfLgap 
$(3,19,0,1,5,9,20,8,2)_{\adfTAACF}$,

\adfLgap 
$(0,19,13,5,1,9,6,4,18)_{\adfTAADE}$,

\adfLgap 
$(0,22,21,2,5,9,3,15,4)_{\adfTABCE}$,

\adfLgap 
$(0,22,2,21,9,1,3,15,4)_{\adfTACCD}$

\adfLgap
\noindent under the action of the mapping $x \mapsto x + 1$ (mod 25).

\ADFvfyParStart{25, \{\{25,25,1\}\}, -1, \{\{5,\{0,1,2,3,4\}\}\}, -1} 

\noindent {\boldmath $K_{5,5,5,5,15}$}~
Let the vertex set be $\{0, 1, \dots, 34\}$ partitioned into
$\{3j + i: j = 0, 1, 2, 3, 4\}$, $i = 0, 1, 2$, $\{15, 16, 17, 18, 19\}$ and $\{20, 21, \dots, 34\}$.
The decompositions consist of the graphs

\adfLgap 
$(5,26,6,15,22,16,4,8,1)_{\adfTAABG}$,

$(21,17,3,13,28,1,32,8,30)_{\adfTAABG}$,

$(15,14,12,20,1,6,32,0,28)_{\adfTAABG}$,

\adfLgap 
$(4,20,23,12,6,24,3,11,10)_{\adfTAACF}$,

$(7,21,16,1,24,12,25,18,6)_{\adfTAACF}$,

$(14,10,15,34,4,26,16,24,18)_{\adfTAACF}$,

\adfLgap 
$(6,24,13,23,9,16,33,14,18)_{\adfTAADE}$,

$(7,5,33,19,6,11,17,32,0)_{\adfTAADE}$,

$(6,27,18,21,5,30,7,20,0)_{\adfTAADE}$,

\adfLgap 
$(23,8,13,18,20,3,29,1,9)_{\adfTABCE}$,

$(10,9,15,17,25,8,29,7,26)_{\adfTABCE}$,

$(33,1,19,0,24,3,7,21,15)_{\adfTABCE}$,

\adfLgap 
$(8,32,31,15,25,3,29,16,13)_{\adfTACCD}$,

$(27,28,14,4,0,13,1,26,10)_{\adfTACCD}$,

$(27,6,7,16,17,5,18,2,10)_{\adfTACCD}$

\adfLgap
\noindent under the action of the mapping $x \mapsto x + 1 \adfmod{15}$ for $x < 15$,
$x \mapsto (x + 1 \adfmod{5}) + 15$ for $15 \le x < 20$,
$x \mapsto (x - 20 + 1 \adfmod{15}) + 20$ for $x \ge 20$.

\ADFvfyParStart{35, \{\{15,15,1\},\{5,5,1\},\{15,15,1\}\}, -1, \{\{5,\{0,1,2\}\},\{5,\{3\}\},\{15,\{4\}\}\}, -1} 

\noindent {\boldmath $K_{5,5,5,5,20}$}~
Let the vertex set be $\{0, 1, \dots, 39\}$ partitioned into
$\{3j + i: j = 0, 1, 2, 3, 4\}$, $i = 0, 1, 2, 3$ and $\{20, 21, \dots, 39\}$.
The decompositions consist of the graphs

\adfLgap 
$(16,37,17,34,4,20,15,26,1)_{\adfTAABG}$,

$(27,10,4,6,21,16,7,13,28)_{\adfTAABG}$,

$(2,34,12,15,28,1,6,9,0)_{\adfTAABG}$,

$(30,1,3,14,20,19,31,13,34)_{\adfTAABG}$,

$(12,27,11,38,19,2,32,18,16)_{\adfTAABG}$,

$(35,8,1,10,30,6,11,4,32)_{\adfTAABG}$,

$(0,33,14,29,6,35,7,27,3)_{\adfTAABG}$,

$(14,16,1,22,5,28,9,21,19)_{\adfTAABG}$,

$(3,37,13,38,15,32,12,29,11)_{\adfTAABG}$,

$(31,9,18,1,19,4,36,14,15)_{\adfTAABG}$,

$(1,29,2,27,8,5,36,7,18)_{\adfTAABG}$,

\adfLgap 
$(14,25,24,0,35,2,26,9,7)_{\adfTAACF}$,

$(17,38,0,3,24,15,20,5,10)_{\adfTAACF}$,

$(37,9,4,32,8,29,10,17,11)_{\adfTAACF}$,

$(8,17,28,6,14,19,36,0,23)_{\adfTAACF}$,

$(12,23,24,11,13,33,6,21,5)_{\adfTAACF}$,

$(33,11,1,31,19,12,10,30,5)_{\adfTAACF}$,

$(4,15,38,0,30,11,37,14,31)_{\adfTAACF}$,

$(21,4,17,31,11,35,13,26,3)_{\adfTAACF}$,

$(12,6,14,3,27,8,13,0,22)_{\adfTAACF}$,

$(27,18,2,36,19,6,7,38,8)_{\adfTAACF}$,

$(3,30,24,18,13,10,9,20,1)_{\adfTAACF}$,

\adfLgap 
$(8,11,19,33,9,15,22,0,32)_{\adfTAADE}$,

$(16,29,36,12,7,2,15,10,0)_{\adfTAADE}$,

$(37,6,0,23,8,18,1,19,38)_{\adfTAADE}$,

$(13,20,22,8,17,19,28,2,3)_{\adfTAADE}$,

$(27,7,15,38,13,5,35,2,8)_{\adfTAADE}$,

$(28,17,15,21,18,3,31,5,25)_{\adfTAADE}$,

$(18,20,33,15,6,23,14,21,5)_{\adfTAADE}$,

$(18,35,36,0,1,38,5,14,11)_{\adfTAADE}$,

$(36,17,8,33,12,6,1,14,35)_{\adfTAADE}$,

$(12,38,23,7,17,31,4,22,14)_{\adfTAADE}$,

$(3,38,33,1,8,34,17,0,2)_{\adfTAADE}$,

\adfLgap 
$(8,13,35,18,12,23,9,33,19)_{\adfTABCE}$,

$(11,16,29,33,3,1,30,7,35)_{\adfTABCE}$,

$(13,31,11,29,8,2,5,24,14)_{\adfTABCE}$,

$(23,2,19,11,22,10,36,15,30)_{\adfTABCE}$,

$(35,37,10,19,0,14,28,16,18)_{\adfTABCE}$,

$(8,33,2,19,5,1,3,36,18)_{\adfTABCE}$,

$(16,21,1,38,12,24,9,2,15)_{\adfTABCE}$,

$(30,31,13,12,1,17,38,8,5)_{\adfTABCE}$,

$(29,28,17,4,19,6,38,11,8)_{\adfTABCE}$,

$(36,34,9,0,18,11,10,32,15)_{\adfTABCE}$,

$(0,1,24,31,2,34,10,15,6)_{\adfTABCE}$,

\adfLgap 
$(29,35,1,0,19,12,18,13,16)_{\adfTACCD}$,

$(13,4,3,38,35,10,33,16,7)_{\adfTACCD}$,

$(36,25,2,0,13,6,8,30,9)_{\adfTACCD}$,

$(10,5,19,32,3,34,20,9,0)_{\adfTACCD}$,

$(20,6,8,7,3,28,4,17,11)_{\adfTACCD}$,

$(18,2,1,25,4,37,31,10,12)_{\adfTACCD}$,

$(4,31,33,7,22,14,28,9,15)_{\adfTACCD}$,

$(11,27,25,13,29,7,16,36,1)_{\adfTACCD}$,

$(4,34,30,11,25,18,35,7,9)_{\adfTACCD}$,

$(13,22,26,2,14,5,15,36,10)_{\adfTACCD}$,

$(2,23,20,11,31,13,19,26,16)_{\adfTACCD}$

\adfLgap
\noindent under the action of the mapping $x \mapsto x + 4 \adfmod{20}$ for $x < 20$,
$x \mapsto (x + 4 \adfmod{20}) + 20$ for $x \ge 20$. \eproof

\ADFvfyParStart{40, \{\{20,5,4\},\{20,5,4\}\}, -1, \{\{5,\{0,1,2,3\}\},\{20,\{4\}\}\}, -1} 

\section{Theta graphs with 11 edges}
\label{sec:theta11}


\noindent {\bf Proof of Lemma \ref{lem:theta11 designs}}~

\noindent {\boldmath $K_{11}$}~
Let the vertex set be $Z_{10} \cup \{\infty\}$. The decompositions consist of the graphs

\adfLgap 
$(0,2,1,3,7,5,8,4,9,\infty)_{\adfTBABH}$,

\adfLgap 
$(0,2,1,5,4,3,6,\infty,9,7)_{\adfTBACG}$,

\adfLgap 
$(0,2,1,3,6,9,5,\infty,4,7)_{\adfTBADF}$,

\adfLgap 
$(0,2,1,3,4,7,6,\infty,5,9)_{\adfTBAEE}$,

\adfLgap 
$(0,2,1,4,3,6,\infty,5,9,7)_{\adfTBBBG}$,

\adfLgap 
$(0,2,1,3,6,5,9,7,\infty,4)_{\adfTBBCF}$,

\adfLgap 
$(0,2,1,3,5,9,4,6,\infty,7)_{\adfTBBDE}$,

\adfLgap 
$(0,2,1,4,3,8,9,5,7,\infty)_{\adfTBCCE}$,

\adfLgap 
$(0,2,1,4,3,9,7,6,5,\infty)_{\adfTBCDD}$

\adfLgap
\noindent under the action of the mapping $x \mapsto x + 2$ (mod 10), $\infty \mapsto \infty$.

\ADFvfyParStart{11, \{\{10,5,2\},\{1,1,1\}\}, 10, -1, -1} 

\noindent {\boldmath $K_{12}$}~
Let the vertex set be $Z_{12}$. The decompositions consist of

\adfLgap 
$(0,1,2,3,4,8,5,7,6,9)_{\adfTBABH}$,

$(2,7,9,4,10,6,3,11,5,0)_{\adfTBABH}$,

\adfLgap 
$(0,1,2,3,4,7,8,5,6,9)_{\adfTBACG}$,

$(2,7,4,9,8,3,6,10,5,11)_{\adfTBACG}$,

\adfLgap 
$(0,1,2,3,4,5,6,8,11,7)_{\adfTBADF}$,

$(2,7,6,3,5,9,1,10,4,0)_{\adfTBADF}$,

\adfLgap 
$(0,1,2,3,4,7,5,6,9,11)_{\adfTBAEE}$,

$(2,7,4,0,6,10,9,5,8,3)_{\adfTBAEE}$,

\adfLgap 
$(0,1,2,3,4,5,8,6,7,9)_{\adfTBBBG}$,

$(2,3,6,7,8,1,10,5,11,4)_{\adfTBBBG}$,

\adfLgap 
$(0,1,2,3,4,5,7,6,8,9)_{\adfTBBCF}$,

$(2,3,6,5,10,8,0,7,1,11)_{\adfTBBCF}$,

\adfLgap 
$(0,1,2,3,4,5,6,9,11,7)_{\adfTBBDE}$,

$(2,3,10,9,7,6,4,0,5,8)_{\adfTBBDE}$,

\adfLgap 
$(0,1,2,3,4,5,6,8,7,10)_{\adfTBCCE}$,

$(2,7,1,4,6,11,9,3,5,0)_{\adfTBCCE}$,

\adfLgap 
$(0,1,2,3,4,5,6,7,11,8)_{\adfTBCDD}$,

$(2,3,5,9,6,0,10,11,1,4)_{\adfTBCDD}$

\adfLgap
\noindent under the action of the mapping $x \mapsto x + 4$ (mod 12).

\ADFvfyParStart{12, \{\{12,3,4\}\}, -1, -1, -1} 

\noindent {\boldmath $K_{22}$}~
Let the vertex set of be $Z_{21} \cup \{\infty\}$. The decompositions consist of

\adfLgap 
$(\infty,0,1,2,3,5,4,6,9,13)_{\adfTBABH}$,

$(0,5,7,6,1,4,8,2,9,17)_{\adfTBABH}$,

$(0,10,17,9,20,2,7,1,11,19)_{\adfTBABH}$,

\adfLgap 
$(\infty,0,1,2,3,5,7,4,6,8)_{\adfTBADF}$,

$(0,1,4,8,11,5,9,2,7,12)_{\adfTBADF}$,

$(0,11,6,13,4,9,1,16,8,2)_{\adfTBADF}$,

\adfLgap 
$(\infty,0,1,2,3,6,4,5,7,10)_{\adfTBBBG}$,

$(0,1,5,6,4,10,2,3,11,8)_{\adfTBBBG}$,

$(0,2,7,9,11,17,8,12,4,13)_{\adfTBBBG}$,

\adfLgap 
$(\infty,0,1,2,3,6,4,5,7,10)_{\adfTBBCF}$,

$(0,1,5,2,6,4,10,3,9,17)_{\adfTBBCF}$,

$(0,2,9,11,5,13,1,8,20,10)_{\adfTBBCF}$,

\adfLgap 
$(\infty,0,1,2,3,5,6,4,7,8)_{\adfTBBDE}$,

$(0,1,6,3,7,5,9,2,8,11)_{\adfTBBDE}$,

$(0,2,7,10,1,14,11,4,19,6)_{\adfTBBDE}$,

\adfLgap 
$(\infty,0,1,2,3,4,6,5,7,10)_{\adfTBCDD}$,

$(0,1,3,7,5,2,6,8,4,11)_{\adfTBCDD}$,

$(0,2,9,8,7,16,11,13,5,12)_{\adfTBCDD}$

\adfLgap
\noindent under the action of the mapping $x \mapsto x + 3$ (mod 21), $\infty \mapsto \infty$. \eproof

\ADFvfyParStart{22, \{\{21,7,3\},\{1,1,1\}\}, 21, -1, -1} 


~\\
\noindent {\bf Proof of Lemma \ref{lem:theta11 multipartite}}~

\noindent {\boldmath $K_{11,11}$}~
Let the vertex set be $Z_{22}$ partitioned according to residue classes modulo 2.
The decompositions consist of

\adfLgap 
$(0,1,3,4,5,16,7,20,13,6)_{\adfTBACG}$,

\adfLgap 
$(0,1,3,4,9,12,15,2,11,16)_{\adfTBAEE}$,

\adfLgap 
$(0,1,3,2,5,8,7,12,21,10)_{\adfTBCCE}$

\adfLgap
\noindent under the action of the mapping $x \mapsto x + 2$ (mod 22).

\ADFvfyParStart{22, \{\{22,11,2\}\}, -1, \{\{11,\{0,1\}\}\}, -1} 

\noindent {\boldmath $K_{11,11,11}$}~
Let the vertex set be $Z_{33}$ partitioned according to residue classes modulo 3.
The decompositions consist of

\adfLgap 
$(0,1,5,2,9,17,3,14,4,21)_{\adfTBABH}$,

\adfLgap 
$(0,1,2,6,11,8,19,5,21,14)_{\adfTBADF}$,

\adfLgap 
$(0,1,2,5,7,15,26,10,24,11)_{\adfTBBBG}$,

\adfLgap 
$(0,1,2,4,9,7,20,3,22,11)_{\adfTBBCF}$,

\adfLgap 
$(0,1,2,4,9,17,7,18,5,15)_{\adfTBBDE}$,

\adfLgap 
$(0,1,2,6,7,15,14,10,26,12)_{\adfTBCDD}$

\adfLgap
\noindent under the action of the mapping $x \mapsto x + 1$ (mod 33).

\ADFvfyParStart{33, \{\{33,33,1\}\}, -1, \{\{11,\{0,1,2\}\}\}, -1} 

\noindent {\boldmath $K_{11,11,11,11}$}~
Let the vertex set be $\{0, 1, \dots, 43\}$ partitioned into
$\{3j + i: j = 0, 1, \dots, 10\}$, $i = 0, 1, 2$, and $\{33, 34, \dots, 43\}$.
The decompositions consist of

\adfLgap 
$(0,10,32,37,8,38,14,31,24,35)_{\adfTBABH}$,

$(5,37,3,0,4,12,26,13,40,9)_{\adfTBABH}$,

\adfLgap 
$(0,33,13,30,29,28,35,19,27,16)_{\adfTBADF}$,

$(27,17,1,3,37,8,4,34,9,41)_{\adfTBADF}$,

\adfLgap 
$(0,27,20,17,2,36,4,32,34,5)_{\adfTBBBG}$,

$(17,35,9,13,3,4,34,6,37,10)_{\adfTBBBG}$,

\adfLgap 
$(0,31,17,26,27,22,24,33,4,39)_{\adfTBBCF}$,

$(33,6,1,5,42,0,20,10,2,40)_{\adfTBBCF}$,

\adfLgap 
$(0,12,42,26,28,32,33,9,34,2)_{\adfTBBDE}$,

$(40,3,11,1,6,17,2,18,41,4)_{\adfTBBDE}$,

\adfLgap 
$(0,21,38,29,33,32,1,40,5,4)_{\adfTBCDD}$,

$(38,8,1,12,2,21,40,10,3,13)_{\adfTBCDD}$

\adfLgap
\noindent under the action of the mapping $x \mapsto x + 1$ (mod 33) for $x < 33$,
$x \mapsto (x + 1 \adfmod{11}) + 33$ for $x \ge 33$.

\ADFvfyParStart{44, \{\{33,33,1\},\{11,11,1\}\}, -1, \{\{11,\{0,1,2\}\},\{11,\{3\}\}\}, -1} 

\noindent {\boldmath $K_{11,11,11,11,11}$}~
Let the vertex set be $Z_{55}$ partitioned according to residue class modulo 5.
The decompositions consist of

\adfLgap 
$(0,46,2,26,38,1,20,21,37,15)_{\adfTBABH}$,

$(7,4,0,1,9,37,14,48,31,17)_{\adfTBABH}$,

\adfLgap 
$(0,38,41,47,20,31,15,26,14,12)_{\adfTBADF}$,

$(34,27,0,1,4,12,16,24,5,14)_{\adfTBADF}$,

\adfLgap 
$(0,21,32,18,4,51,44,28,30,43)_{\adfTBBBG}$,

$(34,3,15,17,5,6,12,21,0,27)_{\adfTBBBG}$,

\adfLgap 
$(0,42,38,41,19,3,1,10,21,48)_{\adfTBBCF}$,

$(6,21,5,13,0,24,12,4,23,47)_{\adfTBBCF}$,

\adfLgap 
$(0,2,16,22,45,46,27,25,28,40)_{\adfTBBDE}$,

$(20,26,2,1,9,5,7,0,6,35)_{\adfTBBDE}$,

\adfLgap 
$(0,26,46,34,39,37,50,27,5,12)_{\adfTBCDD}$,

$(47,34,3,2,10,4,0,18,1,37)_{\adfTBCDD}$

\adfLgap
\noindent under the action of the mapping $x \mapsto x + 1 \adfmod{55}$. \eproof

\ADFvfyParStart{55, \{\{55,55,1\}\}, -1, \{\{11,\{0,1,2,3,4\}\}\}, -1} 

\section{Theta graphs with 12 edges}
\label{sec:theta12}


\noindent {\bf Proof of Lemma \ref{lem:theta12 designs}}~

\noindent {\boldmath $K_{16}$}~
Let the vertex set be $Z_{15} \cup \{\infty\}$. The decompositions consist of the graphs

\adfLgap 
$(\infty,0,1,2,3,5,4,6,9,13,7)_{\adfTCABI}$,

$(1,5,8,4,2,7,12,3,14,6,11)_{\adfTCABI}$,

\adfLgap 
$(\infty,0,1,2,5,6,3,4,7,11,8)_{\adfTCACH}$,

$(1,8,3,10,6,0,5,9,13,4,14)_{\adfTCACH}$,

\adfLgap 
$(\infty,0,1,2,3,5,7,4,6,8,11)_{\adfTCADG}$,

$(1,5,0,4,10,9,3,8,14,7,12)_{\adfTCADG}$,

\adfLgap 
$(\infty,0,1,2,3,4,5,7,9,6,8)_{\adfTCAEF}$,

$(1,5,4,9,0,11,7,2,10,3,8)_{\adfTCAEF}$,

\adfLgap 
$(\infty,0,1,2,3,6,4,5,7,10,14)_{\adfTCBBH}$,

$(1,2,6,7,8,5,0,4,12,3,11)_{\adfTCBBH}$,

\adfLgap 
$(\infty,0,1,2,3,6,4,5,7,11,8)_{\adfTCBCG}$,

$(1,2,6,4,0,7,14,5,10,3,12)_{\adfTCBCG}$,

\adfLgap 
$(\infty,0,1,2,3,5,6,4,7,8,10)_{\adfTCBDF}$,

$(1,2,5,8,0,11,12,3,6,13,7)_{\adfTCBDF}$,

\adfLgap 
$(\infty,0,1,2,3,5,4,6,9,7,10)_{\adfTCBEE}$,

$(1,2,5,8,0,6,11,10,3,14,4)_{\adfTCBEE}$,

\adfLgap 
$(\infty,0,1,2,3,4,5,6,9,7,10)_{\adfTCCCF}$,

$(1,2,5,7,8,11,10,3,9,14,6)_{\adfTCCCF}$,

\adfLgap 
$(\infty,0,1,2,3,4,6,5,7,10,14)_{\adfTCCDE}$,

$(1,2,6,9,7,3,10,11,0,5,8)_{\adfTCCDE}$,

\adfLgap 
$(\infty,0,1,2,3,5,7,4,6,8,11)_{\adfTCDDD}$,

$(1,2,0,5,10,6,4,13,9,3,11)_{\adfTCDDD}$

\adfLgap
\noindent under the action of the mapping $x \mapsto x + 3$ (mod 15), $\infty \mapsto \infty$.

\ADFvfyParStart{16, \{\{15,5,3\},\{1,1,1\}\}, 15, -1, -1} 



\noindent {\boldmath $K_{33}$}~
Let the vertex set be $Z_{33}$. The decompositions consist of

\adfLgap 
$(0,4,31,10,18,27,9,11,19,29,20)_{\adfTCABI}$,

$(5,24,17,0,13,14,1,27,28,26,11)_{\adfTCABI}$,

$(28,25,4,2,13,30,18,14,8,5,6)_{\adfTCABI}$,

$(7,2,12,11,0,6,3,20,1,16,27)_{\adfTCABI}$,

\adfLgap 
$(0,3,19,25,6,1,2,16,29,21,4)_{\adfTCACH}$,

$(11,18,5,14,19,6,16,7,10,28,3)_{\adfTCACH}$,

$(1,30,23,20,5,10,26,6,17,2,21)_{\adfTCACH}$,

$(10,8,3,1,22,32,11,6,23,21,9)_{\adfTCACH}$,

\adfLgap 
$(0,18,30,26,10,24,1,3,17,19,13)_{\adfTCADG}$,

$(20,10,2,27,14,3,19,12,23,29,1)_{\adfTCADG}$,

$(9,11,8,20,19,22,10,13,27,16,2)_{\adfTCADG}$,

$(10,6,9,2,12,17,28,13,26,23,18)_{\adfTCADG}$,

\adfLgap 
$(0,11,24,1,6,5,4,3,28,15,18)_{\adfTCAEF}$,

$(19,22,13,15,20,2,4,5,28,6,29)_{\adfTCAEF}$,

$(20,23,16,5,29,4,3,30,13,8,27)_{\adfTCAEF}$,

$(3,5,10,19,0,18,15,23,7,28,26)_{\adfTCAEF}$,

\adfLgap 
$(0,1,23,21,29,27,3,10,28,18,15)_{\adfTCBBH}$,

$(25,14,28,27,31,29,30,8,26,0,4)_{\adfTCBBH}$,

$(17,12,1,4,16,8,3,20,25,30,13)_{\adfTCBBH}$,

$(8,2,11,14,0,18,12,26,19,7,31)_{\adfTCBBH}$,

\adfLgap 
$(0,15,2,30,23,28,5,26,10,19,21)_{\adfTCBCG}$,

$(1,0,5,7,29,8,17,27,10,22,15)_{\adfTCBCG}$,

$(0,13,9,19,5,12,23,25,20,3,16)_{\adfTCBCG}$,

$(2,5,8,1,24,17,4,12,13,28,6)_{\adfTCBCG}$,

\adfLgap 
$(0,5,27,4,2,18,14,9,16,6,29)_{\adfTCBDF}$,

$(20,5,17,10,22,16,25,11,7,15,31)_{\adfTCBDF}$,

$(18,1,21,19,0,2,11,3,15,30,26)_{\adfTCBDF}$,

$(28,2,8,4,1,17,0,9,7,25,3)_{\adfTCBDF}$,

\adfLgap 
$(0,16,17,19,29,13,6,15,10,1,23)_{\adfTCBEE}$,

$(28,11,30,15,19,31,25,8,4,1,9)_{\adfTCBEE}$,

$(3,17,6,24,15,5,12,9,29,4,2)_{\adfTCBEE}$,

$(14,2,5,6,7,24,23,19,1,12,8)_{\adfTCBEE}$,

\adfLgap 
$(0,18,10,31,23,3,27,15,17,26,1)_{\adfTCCCF}$,

$(19,8,15,13,0,7,5,17,2,25,24)_{\adfTCCCF}$,

$(4,17,22,19,11,6,26,29,30,27,23)_{\adfTCCCF}$,

$(21,0,2,9,13,22,29,25,31,11,28)_{\adfTCCCF}$,

\adfLgap 
$(0,18,29,11,15,26,24,19,8,20,10)_{\adfTCCDE}$,

$(8,14,4,6,5,18,30,2,3,19,1)_{\adfTCCDE}$,

$(6,31,28,24,11,21,18,1,4,5,14)_{\adfTCCDE}$,

$(24,1,0,10,5,7,13,25,11,4,29)_{\adfTCCDE}$,

\adfLgap 
$(0,20,26,27,13,31,7,4,18,12,2)_{\adfTCDDD}$,

$(17,5,3,31,8,9,25,13,28,1,21)_{\adfTCDDD}$,

$(28,30,24,15,17,3,14,18,29,31,2)_{\adfTCDDD}$,

$(20,1,7,30,0,11,16,27,26,5,19)_{\adfTCDDD}$

\adfLgap
\noindent under the action of the mapping $x \mapsto x + 3$ (mod 33).

\ADFvfyParStart{33, \{\{33,11,3\}\}, -1, -1, -1} 

\noindent {\boldmath $K_{40}$}~
Let the vertex set be $Z_{39} \cup \{\infty\}$. The decompositions consist of

\adfLgap 
$(\infty,31,0,2,27,17,16,19,37,11,26)_{\adfTCABI}$,

$(2,4,19,29,20,38,36,25,30,8,16)_{\adfTCABI}$,

$(33,20,0,12,24,29,23,19,15,31,9)_{\adfTCABI}$,

$(4,30,29,20,10,38,19,28,3,18,21)_{\adfTCABI}$,

$(21,14,17,5,36,4,24,22,12,13,7)_{\adfTCABI}$,

\adfLgap 
$(\infty,26,31,12,24,2,22,33,9,29,28)_{\adfTCACH}$,

$(11,23,14,21,0,36,1,28,4,35,17)_{\adfTCACH}$,

$(13,27,16,25,8,0,6,7,30,22,1)_{\adfTCACH}$,

$(17,21,28,5,12,22,15,27,14,38,0)_{\adfTCACH}$,

$(14,24,1,7,4,8,22,5,35,28,33)_{\adfTCACH}$,

\adfLgap 
$(\infty,20,25,37,31,6,28,9,5,26,17)_{\adfTCADG}$,

$(0,23,16,11,1,29,37,27,28,36,22)_{\adfTCADG}$,

$(35,15,16,23,10,27,12,14,28,30,33)_{\adfTCADG}$,

$(32,16,26,11,7,34,10,13,24,23,12)_{\adfTCADG}$,

$(5,0,12,3,17,32,6,13,34,21,27)_{\adfTCADG}$,

\adfLgap 
$(\infty,2,18,31,20,19,13,25,6,32,33)_{\adfTCAEF}$,

$(26,5,21,4,11,38,24,1,34,8,10)_{\adfTCAEF}$,

$(9,23,18,12,8,4,16,19,37,13,38)_{\adfTCAEF}$,

$(23,33,20,28,0,25,3,6,1,30,34)_{\adfTCAEF}$,

$(9,7,2,11,0,17,24,3,30,14,37)_{\adfTCAEF}$,

\adfLgap 
$(\infty,29,19,11,7,0,21,18,13,8,12)_{\adfTCBCG}$,

$(37,23,5,28,10,4,30,13,11,9,18)_{\adfTCBCG}$,

$(15,29,23,22,37,30,3,31,19,35,36)_{\adfTCBCG}$,

$(27,10,33,2,38,25,0,29,14,34,35)_{\adfTCBCG}$,

$(0,1,4,10,9,23,3,29,2,32,21)_{\adfTCBCG}$,

\adfLgap 
$(\infty,1,10,24,19,20,14,32,11,36,15)_{\adfTCBEE}$,

$(13,32,33,37,16,28,30,17,25,11,15)_{\adfTCBEE}$,

$(29,15,7,22,11,18,24,31,8,3,19)_{\adfTCBEE}$,

$(29,12,19,14,26,23,15,3,25,24,0)_{\adfTCBEE}$,

$(0,1,20,10,5,15,4,23,32,21,34)_{\adfTCBEE}$,

\adfLgap 
$(\infty,14,17,3,4,12,6,23,31,7,26)_{\adfTCCCF}$,

$(37,38,35,13,7,10,4,36,32,25,0)_{\adfTCCCF}$,

$(28,27,1,14,9,0,38,5,8,4,21)_{\adfTCCCF}$,

$(13,18,9,2,8,38,36,10,21,31,33)_{\adfTCCCF}$,

$(15,2,5,20,12,33,10,11,6,7,25)_{\adfTCCCF}$,

\adfLgap 
$(\infty,30,5,1,33,38,19,28,9,17,14)_{\adfTCCDE}$,

$(36,38,15,22,13,1,29,24,0,3,25)_{\adfTCCDE}$,

$(31,35,32,6,22,36,8,38,17,0,33)_{\adfTCCDE}$,

$(38,12,14,8,28,27,26,16,18,22,7)_{\adfTCCDE}$,

$(29,2,4,7,31,18,9,37,19,13,21)_{\adfTCCDE}$

\adfLgap
\noindent under the action of the mapping $x \mapsto x + 3$ (mod 39), $\infty \mapsto \infty$.

\ADFvfyParStart{40, \{\{39,13,3\},\{1,1,1\}\}, 39, -1, -1} 

\noindent {\boldmath $K_{49}$}~
Let the vertex set be $Z_{49}$. The decompositions consist of

\adfLgap 
$(0,27,37,26,8,6,19,47,23,30,41)_{\adfTCABI}$,

$(14,5,6,8,3,0,15,11,27,7,24)_{\adfTCABI}$,

\adfLgap 
$(0,26,44,41,45,28,12,10,47,22,15)_{\adfTCACH}$,

$(24,2,3,11,4,5,19,0,31,18,8)_{\adfTCACH}$,

\adfLgap 
$(0,26,16,46,14,11,47,32,7,42,35)_{\adfTCADG}$,

$(23,1,2,0,4,3,8,7,15,9,19)_{\adfTCADG}$,

\adfLgap 
$(0,40,43,21,24,17,19,37,13,45,6)_{\adfTCAEF}$,

$(28,0,8,3,1,12,14,6,2,15,16)_{\adfTCAEF}$,

\adfLgap 
$(0,33,13,10,34,23,2,18,22,5,27)_{\adfTCBCG}$,

$(25,0,7,6,8,10,1,4,9,23,11)_{\adfTCBCG}$,

\adfLgap 
$(0,22,36,16,10,1,3,5,20,2,39)_{\adfTCBEE}$,

$(36,0,7,8,4,1,11,9,10,2,25)_{\adfTCBEE}$,

\adfLgap 
$(0,20,28,26,39,3,35,13,29,14,21)_{\adfTCCCF}$,

$(47,0,1,5,6,18,9,28,3,29,20)_{\adfTCCCF}$,

\adfLgap 
$(0,42,13,31,16,38,35,45,17,47,7)_{\adfTCCDE}$,

$(13,0,1,2,3,8,23,19,27,7,24)_{\adfTCCDE}$

\adfLgap
\noindent under the action of the mapping $x \mapsto x + 1$ (mod 49).

\ADFvfyParStart{49, \{\{49,49,1\}\}, -1, -1, -1} 

\noindent {\boldmath $K_{57}$}~
Let the vertex set be $Z_{57}$. The decompositions consist of

\adfLgap 
$(0,49,36,51,16,40,53,35,46,15,17)_{\adfTCABI}$,

$(11,15,19,35,54,52,17,10,44,47,18)_{\adfTCABI}$,

$(47,53,10,0,40,46,16,11,49,20,9)_{\adfTCABI}$,

$(16,15,0,31,35,30,12,26,38,45,23)_{\adfTCABI}$,

$(13,44,53,49,52,7,18,2,3,33,21)_{\adfTCABI}$,

$(20,35,19,41,11,13,27,51,22,32,15)_{\adfTCABI}$,

$(0,7,25,9,4,13,33,8,39,16,54)_{\adfTCABI}$,

\adfLgap 
$(0,16,45,24,6,35,33,15,11,37,7)_{\adfTCACH}$,

$(17,24,55,27,31,34,15,37,50,47,22)_{\adfTCACH}$,

$(28,23,7,1,21,53,42,0,34,49,54)_{\adfTCACH}$,

$(12,20,52,53,50,43,32,42,15,35,41)_{\adfTCACH}$,

$(51,50,18,23,17,2,15,40,39,49,36)_{\adfTCACH}$,

$(20,2,10,22,42,46,3,34,17,19,11)_{\adfTCACH}$,

$(46,7,0,41,5,1,21,30,47,2,31)_{\adfTCACH}$,

\adfLgap 
$(0,33,54,4,3,42,26,53,15,34,41)_{\adfTCADG}$,

$(2,43,6,16,46,5,17,11,33,44,54)_{\adfTCADG}$,

$(13,0,48,28,51,7,19,45,38,21,4)_{\adfTCADG}$,

$(33,45,25,16,31,47,46,8,32,55,11)_{\adfTCADG}$,

$(12,10,40,8,50,33,49,54,36,41,43)_{\adfTCADG}$,

$(2,11,3,23,36,22,53,32,37,55,19)_{\adfTCADG}$,

$(34,5,20,2,33,38,16,26,37,12,3)_{\adfTCADG}$,

\adfLgap 
$(0,11,54,21,30,4,17,48,35,50,18)_{\adfTCAEF}$,

$(24,20,45,28,42,1,16,32,2,19,17)_{\adfTCAEF}$,

$(31,8,41,29,24,52,21,1,3,23,14)_{\adfTCAEF}$,

$(51,53,29,21,6,12,13,40,18,25,28)_{\adfTCAEF}$,

$(16,4,48,14,42,41,21,22,23,47,52)_{\adfTCAEF}$,

$(55,51,2,37,3,21,41,5,15,53,7)_{\adfTCAEF}$,

$(28,4,2,20,49,10,7,56,42,30,19)_{\adfTCAEF}$,

\adfLgap 
$(0,39,52,32,8,15,34,31,26,29,21)_{\adfTCBCG}$,

$(51,8,55,44,45,42,22,20,1,5,21)_{\adfTCBCG}$,

$(12,1,42,9,38,34,10,23,5,3,49)_{\adfTCBCG}$,

$(48,47,5,29,20,25,7,0,17,1,41)_{\adfTCBCG}$,

$(49,26,38,9,3,24,36,0,31,5,34)_{\adfTCBCG}$,

$(48,47,26,24,22,34,6,16,50,49,4)_{\adfTCBCG}$,

$(45,9,56,37,10,41,19,13,28,7,14)_{\adfTCBCG}$,

\adfLgap 
$(0,24,30,9,53,17,14,3,1,32,13)_{\adfTCBEE}$,

$(52,4,44,39,55,16,29,9,48,49,34)_{\adfTCBEE}$,

$(37,46,9,41,25,6,44,43,40,19,24)_{\adfTCBEE}$,

$(40,36,5,26,9,50,51,35,11,41,2)_{\adfTCBEE}$,

$(38,26,24,32,52,43,51,4,0,20,19)_{\adfTCBEE}$,

$(40,21,9,7,24,49,42,28,38,47,54)_{\adfTCBEE}$,

$(0,3,8,10,33,44,56,35,46,17,32)_{\adfTCBEE}$,

\adfLgap 
$(0,41,24,25,28,9,54,38,19,2,4)_{\adfTCCCF}$,

$(7,53,33,35,14,36,38,45,40,34,44)_{\adfTCCCF}$,

$(32,36,2,54,24,49,53,19,12,55,52)_{\adfTCCCF}$,

$(24,26,53,39,47,2,3,12,6,36,1)_{\adfTCCCF}$,

$(27,55,7,25,17,6,15,0,10,23,22)_{\adfTCCCF}$,

$(36,8,34,43,40,11,19,42,23,28,14)_{\adfTCCCF}$,

$(47,0,1,46,21,20,5,13,28,49,53)_{\adfTCCCF}$,

\adfLgap 
$(0,54,29,32,30,19,20,3,25,28,49)_{\adfTCCDE}$,

$(18,45,11,6,43,42,5,54,3,10,33)_{\adfTCCDE}$,

$(39,43,6,2,49,51,41,24,55,29,10)_{\adfTCCDE}$,

$(35,22,51,8,43,28,10,25,5,14,3)_{\adfTCCDE}$,

$(47,46,35,33,41,52,9,13,21,29,2)_{\adfTCCDE}$,

$(43,40,15,47,16,12,11,52,36,19,44)_{\adfTCCDE}$,

$(43,3,8,29,26,50,12,49,44,5,47)_{\adfTCCDE}$

\adfLgap
\noindent under the action of the mapping $x \mapsto x + 3$ (mod 57).

\ADFvfyParStart{57, \{\{57,19,3\}\}, -1, -1, -1} 

\noindent {\boldmath $K_{81}$}~
Let the vertex set be $Z_{81}$. The decompositions consist of

\adfLgap 
$(0,23,32,38,13,9,54,70,27,67,60)_{\adfTCABI}$,

$(28,23,56,43,77,57,63,73,53,27,1)_{\adfTCABI}$,

$(55,18,53,14,65,28,51,24,7,31,77)_{\adfTCABI}$,

$(40,60,31,27,22,36,17,20,44,2,8)_{\adfTCABI}$,

$(26,36,47,24,33,74,27,52,23,40,21)_{\adfTCABI}$,

$(19,29,37,40,13,17,30,14,28,78,74)_{\adfTCABI}$,

$(65,12,10,77,11,78,79,46,0,18,42)_{\adfTCABI}$,

$(40,59,41,34,50,33,72,64,13,16,5)_{\adfTCABI}$,

$(61,68,19,12,60,13,25,3,2,27,55)_{\adfTCABI}$,

$(30,35,27,18,11,42,21,10,46,8,58)_{\adfTCABI}$,

\adfLgap 
$(0,11,4,77,65,43,42,50,8,41,10)_{\adfTCACH}$,

$(4,47,36,58,27,63,7,41,21,18,37)_{\adfTCACH}$,

$(53,59,29,52,24,7,60,45,63,54,10)_{\adfTCACH}$,

$(4,37,59,73,1,31,50,70,64,52,57)_{\adfTCACH}$,

$(24,70,23,28,72,47,11,78,49,65,30)_{\adfTCACH}$,

$(71,31,74,57,43,36,15,53,41,62,58)_{\adfTCACH}$,

$(78,7,41,28,27,59,32,50,54,14,21)_{\adfTCACH}$,

$(58,2,29,57,56,26,36,41,54,4,19)_{\adfTCACH}$,

$(18,65,57,19,10,28,12,35,44,63,6)_{\adfTCACH}$,

$(17,61,15,4,48,1,3,57,45,39,52)_{\adfTCACH}$,

\adfLgap 
$(0,65,70,36,53,26,51,42,19,46,45)_{\adfTCADG}$,

$(54,16,22,15,11,49,25,9,8,61,31)_{\adfTCADG}$,

$(25,4,69,24,52,45,78,49,10,27,30)_{\adfTCADG}$,

$(30,70,77,40,56,48,38,14,13,9,78)_{\adfTCADG}$,

$(5,64,59,21,23,67,54,62,79,77,38)_{\adfTCADG}$,

$(0,29,14,27,77,31,49,58,33,26,48)_{\adfTCADG}$,

$(8,53,58,38,35,37,25,3,27,13,48)_{\adfTCADG}$,

$(77,73,70,67,35,31,76,65,33,68,2)_{\adfTCADG}$,

$(55,49,47,41,62,32,60,54,15,36,51)_{\adfTCADG}$,

$(0,51,41,32,62,19,9,63,5,61,14)_{\adfTCADG}$,

\adfLgap 
$(0,24,62,52,38,20,76,71,60,14,2)_{\adfTCAEF}$,

$(30,18,9,16,19,33,67,1,54,47,63)_{\adfTCAEF}$,

$(12,43,76,72,62,5,45,74,71,50,10)_{\adfTCAEF}$,

$(51,73,64,9,23,4,28,49,20,53,8)_{\adfTCAEF}$,

$(7,32,29,52,12,1,13,4,66,8,25)_{\adfTCAEF}$,

$(19,55,51,48,49,29,76,33,50,12,56)_{\adfTCAEF}$,

$(71,58,43,63,65,16,29,37,55,57,56)_{\adfTCAEF}$,

$(1,30,72,42,29,2,36,44,40,77,21)_{\adfTCAEF}$,

$(76,41,49,15,62,21,46,30,3,66,72)_{\adfTCAEF}$,

$(29,24,23,70,59,50,14,65,33,58,66)_{\adfTCAEF}$,

\adfLgap 
$(0,40,68,78,37,24,10,6,50,17,4)_{\adfTCBCG}$,

$(70,60,69,5,24,50,64,51,33,16,59)_{\adfTCBCG}$,

$(73,1,21,75,16,49,44,25,30,19,74)_{\adfTCBCG}$,

$(33,5,41,60,77,7,44,51,56,59,57)_{\adfTCBCG}$,

$(63,14,20,31,29,16,34,76,5,66,41)_{\adfTCBCG}$,

$(64,79,5,53,4,55,48,59,58,10,44)_{\adfTCBCG}$,

$(23,47,64,27,17,57,71,30,42,77,24)_{\adfTCBCG}$,

$(17,27,46,67,57,13,34,4,69,47,70)_{\adfTCBCG}$,

$(1,22,45,70,78,36,15,43,16,41,72)_{\adfTCBCG}$,

$(51,41,57,36,62,3,42,74,5,23,65)_{\adfTCBCG}$,

\adfLgap 
$(0,72,9,29,59,6,31,50,7,26,5)_{\adfTCBEE}$,

$(55,78,70,72,67,14,71,79,46,58,2)_{\adfTCBEE}$,

$(79,11,0,8,72,68,27,69,40,10,7)_{\adfTCBEE}$,

$(53,63,17,9,16,32,51,62,75,54,31)_{\adfTCBEE}$,

$(7,54,3,62,56,33,76,16,10,73,18)_{\adfTCBEE}$,

$(37,9,48,35,20,38,11,8,22,43,24)_{\adfTCBEE}$,

$(0,66,56,46,5,53,16,61,29,30,33)_{\adfTCBEE}$,

$(26,0,13,65,6,30,38,37,76,77,1)_{\adfTCBEE}$,

$(64,20,17,71,39,76,40,78,72,25,8)_{\adfTCBEE}$,

$(27,16,74,47,1,28,0,54,80,7,38)_{\adfTCBEE}$,

\adfLgap 
$(0,63,14,55,54,73,22,15,58,71,33)_{\adfTCCCF}$,

$(74,56,53,9,22,59,18,31,43,46,71)_{\adfTCCCF}$,

$(32,66,43,72,31,64,45,49,35,13,78)_{\adfTCCCF}$,

$(71,45,73,47,32,42,17,62,53,43,74)_{\adfTCCCF}$,

$(51,73,1,22,60,49,44,12,45,53,58)_{\adfTCCCF}$,

$(16,5,77,36,44,52,23,18,60,24,1)_{\adfTCCCF}$,

$(15,42,49,26,68,25,41,65,19,63,59)_{\adfTCCCF}$,

$(45,63,4,58,73,28,31,25,44,62,74)_{\adfTCCCF}$,

$(59,46,39,7,53,72,29,75,55,30,45)_{\adfTCCCF}$,

$(18,43,77,61,17,34,0,57,36,59,11)_{\adfTCCCF}$,

\adfLgap 
$(0,19,28,14,26,63,34,8,25,38,50)_{\adfTCCDE}$,

$(27,6,11,39,58,56,38,24,7,41,49)_{\adfTCCDE}$,

$(75,10,71,51,55,65,28,5,58,52,61)_{\adfTCCDE}$,

$(73,34,70,48,75,36,61,52,76,17,23)_{\adfTCCDE}$,

$(78,24,60,36,69,2,40,26,46,13,39)_{\adfTCCDE}$,

$(47,11,77,56,38,73,24,23,0,76,34)_{\adfTCCDE}$,

$(48,77,12,62,29,75,74,78,44,71,38)_{\adfTCCDE}$,

$(17,6,16,52,49,65,27,58,48,61,60)_{\adfTCCDE}$,

$(46,75,42,35,71,52,64,20,16,39,31)_{\adfTCCDE}$,

$(55,72,21,77,36,14,78,48,38,67,74)_{\adfTCCDE}$

\adfLgap
\noindent under the action of the mapping $x \mapsto x + 3$ (mod 81). \eproof

\ADFvfyParStart{81, \{\{81,27,3\}\}, -1, -1, -1} 

~\\
\noindent {\bf Proof of Lemma \ref{lem:theta12 multipartite}}~

\noindent {\boldmath $K_{8,12}$}~
Let the vertex set be $\{0, 1, \dots, 19\}$ partitioned into $\{0, 1, \dots, 7\}$
and $\{8, 9, \dots, 19\}$. The decompositions consist of

\adfLgap 
$(0,2,8,9,10,1,14,3,13,4,12)_{\adfTCBBH}$,

\adfLgap 
$(0,2,8,9,1,10,11,6,16,5,18)_{\adfTCBDF}$,

\adfLgap 
$(0,2,8,1,9,10,3,13,11,7,18)_{\adfTCDDD}$

\adfLgap
\noindent under the action of the mapping $x \mapsto x + 1 \adfmod{8}$ for $x < 8$,
$x \mapsto (x - 8 + 3 \adfmod{12}) + 8$ for $x \ge 8$.

\ADFvfyParStart{20, \{\{8,8,1\},\{12,4,3\}\}, -1, \{\{8,\{0\}\},\{12,\{1\}\}\}, -1} 

\noindent {\boldmath $K_{8,8,8}$}~
Let the vertex set be $Z_{24}$ partitioned according to residue classes modulo 3.
The decompositions consist of

\adfLgap 
$(0,1,2,4,6,5,7,11,3,8,12)_{\adfTCABI}$,

$(0,7,14,10,2,13,8,22,6,23,12)_{\adfTCABI}$,

\adfLgap 
$(0,1,2,3,4,5,7,11,6,13,8)_{\adfTCACH}$,

$(0,13,8,18,10,2,6,17,4,14,21)_{\adfTCACH}$,

\adfLgap 
$(0,1,2,3,8,4,5,7,9,16,6)_{\adfTCADG}$,

$(0,13,8,4,23,11,19,3,20,6,2)_{\adfTCADG}$,

\adfLgap 
$(0,1,2,3,7,5,8,12,4,6,11)_{\adfTCAEF}$,

$(0,7,10,2,1,12,11,4,23,6,20)_{\adfTCAEF}$,

\adfLgap 
$(0,1,2,4,3,5,6,13,8,10,14)_{\adfTCBCG}$,

$(0,1,8,11,15,13,5,19,3,20,6)_{\adfTCBCG}$,

\adfLgap 
$(0,1,2,4,3,8,9,7,5,10,12)_{\adfTCBEE}$,

$(0,1,8,10,2,9,5,14,3,22,11)_{\adfTCBEE}$,

\adfLgap 
$(0,1,2,3,4,5,7,6,11,13,8)_{\adfTCCCF}$,

$(0,1,8,12,16,6,11,21,14,4,17)_{\adfTCCCF}$,

\adfLgap 
$(0,1,2,3,4,5,9,7,6,16,8)_{\adfTCCDE}$,

$(0,1,5,12,8,4,14,19,17,6,20)_{\adfTCCDE}$

\adfLgap
\noindent under the action of the mapping $x \mapsto x + 3$ (mod 24).

\ADFvfyParStart{24, \{\{24,8,3\}\}, -1, \{\{8,\{0,1,2\}\}\}, -1} 

\noindent {\boldmath $K_{8,8,8,8}$}~
Let the vertex set be $Z_{32}$ partitioned according to residue classes modulo 4.
The decompositions consist of
\adfLgap 
$(0,2,3,5,11,1,8,21,4,22,13)_{\adfTCABI}$,

\adfLgap 
$(0,2,1,7,3,10,19,4,14,27,13)_{\adfTCACH}$,

\adfLgap 
$(0,2,1,4,9,6,15,5,16,3,17)_{\adfTCADG}$,

\adfLgap 
$(0,2,1,4,9,15,7,16,26,5,19)_{\adfTCAEF}$,

\adfLgap 
$(0,2,3,5,11,7,9,22,8,23,12)_{\adfTCBCG}$,

\adfLgap 
$(0,2,3,5,7,1,11,13,20,6,17)_{\adfTCBEE}$,

\adfLgap 
$(0,2,1,4,5,11,7,20,6,23,12)_{\adfTCCCF}$,

\adfLgap 
$(0,2,1,4,5,11,20,7,17,6,19)_{\adfTCCDE}$.

\adfLgap
\noindent under the action of the mapping $x \mapsto x + 1$ (mod 32).

\ADFvfyParStart{32, \{\{32,32,1\}\}, -1, \{\{8,\{0,1,2,3\}\}\}, -1} 

\noindent {\boldmath $K_{8,8,8,24}$}~
Let the vertex set be $Z_{48}$ partitioned according to residue classes 0, 1, 2 and $\{3,4,5\}$ modulo 6.
The decompositions consist of

\adfLgap 
$(0,35,32,31,29,7,2,10,6,4,36)_{\adfTCABI}$,

$(27,20,43,42,19,9,38,39,2,1,6)_{\adfTCABI}$,

$(24,27,2,5,25,10,7,15,18,8,1)_{\adfTCABI}$,

$(11,20,18,2,35,12,17,44,40,8,46)_{\adfTCABI}$,

$(46,12,1,7,3,36,28,18,19,5,8)_{\adfTCABI}$,

$(31,40,18,45,14,43,24,32,46,44,13)_{\adfTCABI}$,

$(0,17,7,27,36,45,25,46,18,29,13)_{\adfTCABI}$,

$(19,11,32,34,14,25,27,6,26,39,44)_{\adfTCABI}$,

\adfLgap 
$(0,22,10,6,33,19,32,41,42,38,31)_{\adfTCACH}$,

$(15,38,26,0,44,12,25,24,3,30,19)_{\adfTCACH}$,

$(22,43,36,17,25,15,37,41,7,38,28)_{\adfTCACH}$,

$(23,12,30,35,6,31,11,8,22,44,3)_{\adfTCACH}$,

$(20,12,5,14,15,13,23,7,27,26,9)_{\adfTCACH}$,

$(17,19,20,40,44,30,39,24,43,35,8)_{\adfTCACH}$,

$(36,43,8,4,16,14,46,6,1,26,39)_{\adfTCACH}$,

$(21,13,18,16,37,4,44,1,2,17,30)_{\adfTCACH}$,

\adfLgap 
$(0,22,29,20,43,13,46,7,23,44,30)_{\adfTCADG}$,

$(32,5,4,31,36,39,38,6,22,12,43)_{\adfTCADG}$,

$(15,24,44,29,31,7,4,37,39,18,23)_{\adfTCADG}$,

$(14,29,3,31,38,10,8,36,23,0,37)_{\adfTCADG}$,

$(6,10,14,9,12,39,1,35,42,28,2)_{\adfTCADG}$,

$(32,42,16,26,43,29,25,44,1,33,20)_{\adfTCADG}$,

$(26,30,29,7,10,3,20,18,45,19,33)_{\adfTCADG}$,

$(27,18,12,23,43,31,20,19,32,46,37)_{\adfTCADG}$,

\adfLgap 
$(0,14,16,1,44,12,40,37,36,5,25)_{\adfTCAEF}$,

$(29,31,20,3,13,42,18,40,36,23,2)_{\adfTCAEF}$,

$(35,2,31,23,18,39,42,41,7,29,19)_{\adfTCAEF}$,

$(22,14,32,11,30,39,24,27,25,46,36)_{\adfTCAEF}$,

$(27,42,14,15,24,21,20,35,38,0,7)_{\adfTCAEF}$,

$(34,32,43,15,19,35,38,13,45,30,10)_{\adfTCAEF}$,

$(27,8,19,36,13,22,32,23,0,20,40)_{\adfTCAEF}$,

$(20,13,19,4,25,39,12,46,1,14,18)_{\adfTCAEF}$,

\adfLgap 
$(0,28,43,1,32,37,17,36,2,18,14)_{\adfTCBCG}$,

$(35,12,2,37,23,30,32,16,38,18,22)_{\adfTCBCG}$,

$(5,8,31,12,29,43,22,0,35,44,33)_{\adfTCBCG}$,

$(10,7,24,13,22,18,33,42,45,12,9)_{\adfTCBCG}$,

$(42,40,19,3,32,15,2,24,22,31,12)_{\adfTCBCG}$,

$(28,41,26,8,1,12,11,2,3,20,25)_{\adfTCBCG}$,

$(6,23,19,29,26,33,37,38,28,25,44)_{\adfTCBCG}$,

$(7,2,45,33,13,15,31,44,36,43,9)_{\adfTCBCG}$,

\adfLgap 
$(0,8,21,16,24,19,4,41,25,9,7)_{\adfTCBEE}$,

$(30,42,29,2,13,8,33,14,21,32,3)_{\adfTCBEE}$,

$(2,5,42,6,21,7,36,29,43,12,38)_{\adfTCBEE}$,

$(44,15,30,37,34,36,43,5,20,17,7)_{\adfTCBEE}$,

$(12,36,37,1,23,32,15,34,44,42,22)_{\adfTCBEE}$,

$(4,32,1,42,17,25,34,31,9,12,16)_{\adfTCBEE}$,

$(30,22,43,35,7,11,8,20,28,2,31)_{\adfTCBEE}$,

$(45,5,7,1,26,27,24,2,37,4,32)_{\adfTCBEE}$,

\adfLgap 
$(0,45,26,43,1,24,34,18,3,20,30)_{\adfTCCCF}$,

$(8,45,16,7,43,44,24,19,30,4,26)_{\adfTCCCF}$,

$(0,22,45,18,41,1,19,34,14,46,8)_{\adfTCCCF}$,

$(10,20,7,23,31,9,12,40,38,11,13)_{\adfTCCCF}$,

$(34,31,42,8,38,3,24,5,19,23,0)_{\adfTCCCF}$,

$(10,38,19,33,13,45,25,35,2,31,18)_{\adfTCCCF}$,

$(2,18,17,26,7,29,27,36,39,31,20)_{\adfTCCCF}$,

$(24,1,20,29,11,42,23,44,41,36,45)_{\adfTCCCF}$,

\adfLgap 
$(0,44,45,6,14,25,3,15,20,21,13)_{\adfTCCDE}$,

$(21,7,25,24,1,34,14,36,4,8,23)_{\adfTCCDE}$,

$(39,36,8,31,1,22,20,18,27,43,46)_{\adfTCCDE}$,

$(25,22,18,7,27,38,30,12,23,44,0)_{\adfTCCDE}$,

$(29,18,1,46,24,22,43,12,15,36,11)_{\adfTCCDE}$,

$(12,20,38,23,31,35,25,46,8,6,41)_{\adfTCCDE}$,

$(4,23,44,24,26,10,1,0,29,31,32)_{\adfTCCDE}$,

$(14,32,33,19,23,37,29,39,26,7,46)_{\adfTCCDE}$

\adfLgap
\noindent under the action of the mapping $x \mapsto x + 6 \adfmod{48}$. \eproof

\ADFvfyParStart{48, \{\{48,8,6\}\}, -1, \{\{8,\{0,1,2,4,4,4\}\}\}, -1} 

\section{Theta graphs with 13 edges}
\label{sec:theta13}


\noindent {\bf Proof of Lemma \ref{lem:theta13 designs}}~

\noindent {\boldmath $K_{13}$}~
Let the vertex set be $Z_{12} \cup \{\infty\}$. The decompositions consist of the graphs

\adfLgap 
$(0,1,2,3,4,8,5,7,6,10,\infty,11)_{\adfTDABJ}$,

$(0,5,10,6,3,9,1,\infty,4,11,7,2)_{\adfTDABJ}$,

\adfLgap 
$(1,0,2,3,4,6,8,7,10,5,11,\infty)_{\adfTDACI}$,

$(0,5,6,2,4,11,1,9,\infty,10,3,7)_{\adfTDACI}$,

\adfLgap 
$(0,1,2,3,4,6,5,7,10,\infty,8,11)_{\adfTDADH}$,

$(0,5,4,2,9,7,1,10,6,11,3,\infty)_{\adfTDADH}$,

\adfLgap 
$(1,0,2,3,4,6,5,7,10,\infty,11,8)_{\adfTDAEG}$,

$(0,5,7,1,4,\infty,10,2,9,6,11,3)_{\adfTDAEG}$,

\adfLgap 
$(0,1,2,3,4,7,5,6,9,10,\infty,11)_{\adfTDAFF}$,

$(0,5,4,1,6,8,\infty,7,3,10,2,11)_{\adfTDAFF}$,

\adfLgap 
$(0,1,2,3,4,5,8,6,7,10,\infty,11)_{\adfTDBBI}$,

$(0,1,5,6,11,4,\infty,9,3,7,2,10)_{\adfTDBBI}$,

\adfLgap 
$(0,1,2,3,4,5,7,9,6,10,11,\infty)_{\adfTDBCH}$,

$(0,1,6,4,5,10,3,8,\infty,2,11,7)_{\adfTDBCH}$,

\adfLgap 
$(0,1,2,3,4,5,6,7,9,11,\infty,10)_{\adfTDBDG}$,

$(0,1,7,5,8,\infty,4,2,11,3,10,6)_{\adfTDBDG}$,

\adfLgap 
$(0,1,2,3,4,5,7,6,8,\infty,10,11)_{\adfTDBEF}$,

$(0,1,4,5,9,2,10,7,3,6,11,\infty)_{\adfTDBEF}$,

\adfLgap 
$(1,0,2,3,4,5,6,8,10,7,11,\infty)_{\adfTDCCG}$,

$(0,1,4,3,6,9,7,2,10,\infty,5,11)_{\adfTDCCG}$,

\adfLgap 
$(0,1,2,3,4,5,6,7,8,11,\infty,10)_{\adfTDCDF}$,

$(0,5,10,3,9,4,\infty,6,2,11,7,1)_{\adfTDCDF}$,

\adfLgap 
$(1,0,2,3,4,5,7,8,6,11,\infty,10)_{\adfTDCEE}$,

$(0,1,2,5,6,10,7,11,\infty,9,3,8)_{\adfTDCEE}$,

\adfLgap 
$(0,1,2,3,4,5,6,9,7,10,\infty,11)_{\adfTDDDE}$,

$(0,1,3,11,6,4,5,7,10,2,8,\infty)_{\adfTDDDE}$

\adfLgap
\noindent under the action of the mapping $x \mapsto x + 4$ (mod 12), $\infty \mapsto \infty$.

\ADFvfyParStart{13, \{\{12,3,4\},\{1,1,1\}\}, 12, -1, -1} 

\noindent {\boldmath $K_{14}$}~
Let the vertex set be $Z_{14}$. The decompositions consist of

\adfLgap 
$(0,1,2,3,5,8,4,9,13,6,12,7)_{\adfTDABJ}$,

\adfLgap 
$(0,1,2,5,4,3,8,13,6,12,9,7)_{\adfTDACI}$,

\adfLgap 
$(0,1,2,5,6,4,10,3,12,9,11,7)_{\adfTDADH}$,

\adfLgap 
$(0,1,2,5,3,4,8,12,7,13,6,11)_{\adfTDAEG}$,

\adfLgap 
$(0,1,2,5,3,4,8,6,11,7,12,9)_{\adfTDAFF}$,

\adfLgap 
$(0,1,2,3,4,9,12,6,13,5,10,11)_{\adfTDBBI}$,

\adfLgap 
$(0,1,2,3,5,6,7,13,4,9,12,8)_{\adfTDBCH}$,

\adfLgap 
$(0,1,2,3,5,4,8,12,7,13,6,11)_{\adfTDBDG}$,

\adfLgap 
$(0,1,2,3,5,4,8,6,11,7,12,9)_{\adfTDBEF}$,

\adfLgap 
$(0,1,2,3,4,7,8,13,9,12,5,6)_{\adfTDCCG}$,

\adfLgap 
$(0,1,2,3,4,7,8,5,13,9,12,6)_{\adfTDCDF}$,

\adfLgap 
$(0,1,2,3,4,7,12,9,8,13,6,5)_{\adfTDCEE}$,

\adfLgap 
$(0,1,2,3,4,5,12,6,10,13,9,7)_{\adfTDDDE}$

\adfLgap
\noindent under the action of the mapping $x \mapsto x + 2$ (mod 14).

\ADFvfyParStart{14, \{\{14,7,2\}\}, -1, -1, -1} 

\noindent {\boldmath $K_{26}$}~
Let the vertex set be $Z_{25} \cup \{\infty\}$. The decompositions consist of

\adfLgap 
$(\infty,5,2,6,7,23,17,21,3,14,20,1)_{\adfTDABJ}$,

$(6,18,16,15,19,10,12,23,4,7,21,1)_{\adfTDABJ}$,

$(22,15,16,9,4,12,7,14,21,19,17,5)_{\adfTDABJ}$,

$(18,14,15,17,0,5,19,16,4,13,3,6)_{\adfTDABJ}$,

$(20,3,8,13,14,4,\infty,18,22,7,16,5)_{\adfTDABJ}$,

\adfLgap 
$(\infty,14,11,1,6,23,22,2,15,21,5,8)_{\adfTDADH}$,

$(22,12,16,9,19,3,15,7,5,24,23,21)_{\adfTDADH}$,

$(10,8,21,7,11,9,4,20,16,3,14,22)_{\adfTDADH}$,

$(11,14,24,3,10,12,9,22,18,2,5,0)_{\adfTDADH}$,

$(17,19,10,18,3,\infty,15,0,1,8,13,21)_{\adfTDADH}$,

\adfLgap 
$(\infty,1,3,11,5,15,16,22,7,6,20,12)_{\adfTDAFF}$,

$(19,13,6,3,20,23,9,21,14,22,4,15)_{\adfTDAFF}$,

$(18,13,11,14,6,15,12,2,21,0,20,4)_{\adfTDAFF}$,

$(18,22,0,12,21,17,15,8,9,24,20,19)_{\adfTDAFF}$,

$(8,20,2,7,19,24,\infty,6,11,23,12,14)_{\adfTDAFF}$,

\adfLgap 
$(\infty,23,9,15,8,4,22,19,14,10,0,1)_{\adfTDBBI}$,

$(17,6,1,0,15,18,2,13,11,9,10,7)_{\adfTDBBI}$,

$(19,3,9,5,17,18,1,13,7,22,14,23)_{\adfTDBBI}$,

$(5,8,12,1,14,21,10,22,2,16,6,15)_{\adfTDBBI}$,

$(21,6,2,9,4,3,13,0,5,24,12,\infty)_{\adfTDBBI}$,

\adfLgap 
$(\infty,21,19,0,14,11,12,24,15,17,5,6)_{\adfTDBCH}$,

$(22,8,17,14,4,12,6,1,10,21,0,19)_{\adfTDBCH}$,

$(9,14,8,12,19,5,23,13,17,3,0,15)_{\adfTDBCH}$,

$(6,20,2,22,3,17,18,5,10,8,13,1)_{\adfTDBCH}$,

$(9,1,18,6,3,21,4,2,5,22,\infty,8)_{\adfTDBCH}$,

\adfLgap 
$(\infty,18,11,20,8,9,22,2,21,19,15,12)_{\adfTDBDG}$,

$(3,22,20,8,4,12,14,17,19,5,11,21)_{\adfTDBDG}$,

$(0,23,11,20,12,13,1,10,19,9,16,21)_{\adfTDBDG}$,

$(14,3,9,6,10,0,2,18,20,19,16,7)_{\adfTDBDG}$,

$(8,1,22,\infty,19,12,16,4,10,17,5,23)_{\adfTDBDG}$,

\adfLgap 
$(\infty,21,23,11,6,2,18,19,14,1,16,15)_{\adfTDBEF}$,

$(15,4,11,5,13,9,10,6,20,22,2,14)_{\adfTDBEF}$,

$(12,13,20,15,3,21,7,9,8,23,14,10)_{\adfTDBEF}$,

$(13,12,21,17,3,16,19,8,14,5,10,22)_{\adfTDBEF}$,

$(17,8,10,16,22,24,7,\infty,15,4,21,19)_{\adfTDBEF}$,

\adfLgap 
$(\infty,7,5,16,11,13,6,23,20,8,12,14)_{\adfTDCDF}$,

$(15,0,11,8,2,22,7,21,17,6,19,1)_{\adfTDCDF}$,

$(14,1,16,6,2,19,13,4,7,15,0,18)_{\adfTDCDF}$,

$(3,24,17,8,4,18,20,13,7,5,10,19)_{\adfTDCDF}$,

$(0,4,22,\infty,19,5,21,24,3,8,7,1)_{\adfTDCDF}$,

\adfLgap 
$(\infty,17,8,13,5,10,20,19,9,4,18,6)_{\adfTDDDE}$,

$(21,19,16,14,7,17,15,13,12,5,9,6)_{\adfTDDDE}$,

$(4,5,10,13,21,8,2,1,12,15,16,18)_{\adfTDDDE}$,

$(14,11,17,12,8,21,6,0,5,19,9,18)_{\adfTDDDE}$,

$(16,3,22,5,10,24,23,12,\infty,2,17,18)_{\adfTDDDE}$

\adfLgap
\noindent under the action of the mapping $x \mapsto x + 5$ (mod 25), $\infty \mapsto \infty$. \eproof

\ADFvfyParStart{26, \{\{25,5,5\},\{1,1,1\}\}, 25, -1, -1} 


~\\
\noindent {\bf Proof of Lemma \ref{lem:theta13 multipartite}}~

\noindent {\boldmath $K_{13,13}$}~
Let the vertex set be $Z_{26}$ partitioned according to residue classes modulo 2.
The decompositions consist of

\adfLgap 
$(0,17,11,4,1,2,5,8,13,18,7,24)_{\adfTDACI}$,

\adfLgap 
$(0,11,21,8,17,10,3,4,9,16,5,14)_{\adfTDAEG}$,

\adfLgap 
$(0,19,11,14,21,2,1,12,15,16,23,10)_{\adfTDCCG}$,

\adfLgap 
$(0,1,3,8,15,2,7,10,11,12,19,18)_{\adfTDCEE}$

\adfLgap
\noindent under the action of the mapping $x \mapsto x + 2$ (mod 26).

\ADFvfyParStart{26, \{\{26,13,2\}\}, -1, \{\{13,\{0,1\}\}\}, -1} 

\noindent {\boldmath $K_{13,13,13}$}~
Let the vertex set be $Z_{39}$ partitioned according to residue classes modulo 3.
The decompositions consist of

\adfLgap 
$(0,1,5,2,9,17,3,13,26,6,28,12)_{\adfTDABJ}$,

\adfLgap 
$(0,1,2,6,11,7,15,26,3,17,4,21)_{\adfTDADH}$,

\adfLgap 
$(0,1,2,6,11,3,14,7,23,9,28,18)_{\adfTDAFF}$,

\adfLgap 
$(0,1,2,5,7,15,4,14,28,8,30,17)_{\adfTDBBI}$,

\adfLgap 
$(0,1,2,4,9,7,17,3,14,27,5,21)_{\adfTDBCH}$,

\adfLgap 
$(0,1,2,4,9,17,7,27,10,35,22,11)_{\adfTDBDG}$,

\adfLgap 
$(0,1,2,4,9,16,24,10,23,6,26,12)_{\adfTDBEF}$,

\adfLgap 
$(0,1,2,6,7,8,18,11,3,17,4,20)_{\adfTDCDF}$,

\adfLgap 
$(0,1,2,3,8,4,12,23,10,24,37,17)_{\adfTDDDE}$

\adfLgap
\noindent under the action of the mapping $x \mapsto x + 1$ (mod 39).

\ADFvfyParStart{39, \{\{39,39,1\}\}, -1, \{\{13,\{0,1,2\}\}\}, -1} 

\noindent {\boldmath $K_{13,13,13,13}$}~
Let the vertex set be $\{0, 1, \dots, 51\}$ partitioned into
$\{3j + i: j = 0, 1, \dots, 12\}$, $i = 0, 1, 2$, and $\{39, 40, \dots, 51\}$.
The decompositions consist of

\adfLgap 
$(0,11,49,34,41,2,43,27,44,20,16,3)_{\adfTDABJ}$,

$(49,1,2,4,6,13,3,17,0,16,36,41)_{\adfTDABJ}$,

\adfLgap 
$(0,38,46,12,1,13,23,9,14,6,29,47)_{\adfTDADH}$,

$(7,43,14,44,16,49,11,15,32,12,51,6)_{\adfTDADH}$,

\adfLgap 
$(0,49,14,6,28,23,36,41,37,8,46,38)_{\adfTDAFF}$,

$(1,17,33,37,35,41,36,2,13,42,9,44)_{\adfTDAFF}$,

\adfLgap 
$(0,12,5,22,11,30,14,27,42,38,45,13)_{\adfTDBBI}$,

$(21,42,13,7,39,15,19,44,8,10,49,2)_{\adfTDBBI}$,

\adfLgap 
$(0,46,20,41,1,5,37,6,10,24,26,13)_{\adfTDBCH}$,

$(34,28,6,45,27,11,1,39,4,40,10,44)_{\adfTDBCH}$,

\adfLgap 
$(0,44,29,39,7,2,25,14,6,26,19,20)_{\adfTDBDG}$,

$(21,11,34,19,42,15,39,1,18,40,7,41)_{\adfTDBDG}$,

\adfLgap 
$(0,25,47,23,12,40,15,17,19,18,41,38)_{\adfTDBEF}$,

$(10,45,14,18,25,20,0,24,41,1,47,8)_{\adfTDBEF}$,

\adfLgap 
$(0,31,26,48,46,17,41,2,30,38,40,14)_{\adfTDCDF}$,

$(38,44,6,10,19,35,30,9,8,33,45,0)_{\adfTDCDF}$,

\adfLgap 
$(0,45,8,48,10,5,47,21,28,35,40,38)_{\adfTDDDE}$,

$(9,7,19,32,30,28,24,41,8,25,11,43)_{\adfTDDDE}$

\adfLgap
\noindent under the action of the mapping $x \mapsto x + 1$ (mod 39) for $x < 39$,
$x \mapsto (x + 1 \adfmod{13}) + 39$ for $x \ge 39$.

\ADFvfyParStart{52, \{\{39,39,1\},\{13,13,1\}\}, -1, \{\{13,\{0,1,2\}\},\{13,\{3\}\}\}, -1} 

\noindent {\boldmath $K_{13,13,13,13,13}$}~
Let the vertex set be $Z_{65}$ partitioned according to residue class modulo 5.
The decompositions consist of

\adfLgap 
$(0,17,29,22,48,11,34,23,61,43,19,25)_{\adfTDABJ}$,

$(27,29,60,40,47,3,19,0,4,7,6,15)_{\adfTDABJ}$,

\adfLgap 
$(0,9,2,20,13,32,44,61,19,11,63,12)_{\adfTDADH}$,

$(17,1,6,0,22,16,40,3,37,8,46,20)_{\adfTDADH}$,

\adfLgap 
$(0,43,7,53,2,19,10,39,16,44,55,37)_{\adfTDAFF}$,

$(17,51,21,18,2,46,22,5,6,8,0,13)_{\adfTDAFF}$,

\adfLgap 
$(0,5,21,38,34,63,45,6,48,36,23,27)_{\adfTDBBI}$,

$(5,31,42,12,13,7,4,2,1,10,24,48)_{\adfTDBBI}$,

\adfLgap 
$(0,4,13,37,23,59,30,61,54,62,15,16)_{\adfTDBCH}$,

$(57,11,55,9,35,30,19,22,0,4,20,43)_{\adfTDBCH}$,

\adfLgap 
$(0,45,54,46,12,51,17,9,58,61,10,57)_{\adfTDBDG}$,

$(13,53,40,57,29,31,6,2,3,26,50,17)_{\adfTDBDG}$,

\adfLgap 
$(0,43,37,34,58,46,25,19,41,18,35,32)_{\adfTDBEF}$,

$(17,9,25,24,28,57,48,3,1,2,15,42)_{\adfTDBEF}$,

\adfLgap 
$(0,13,49,21,52,14,31,6,35,2,9,32)_{\adfTDCDF}$,

$(20,8,19,5,29,60,62,16,28,4,25,47)_{\adfTDCDF}$,

\adfLgap 
$(0,9,46,45,1,39,3,27,38,6,62,25)_{\adfTDDDE}$,

$(16,23,33,31,17,19,7,0,5,1,14,45)_{\adfTDDDE}$

\adfLgap
\noindent under the action of the mapping $x \mapsto x + 1 \adfmod{65}$. \eproof

\ADFvfyParStart{65, \{\{65,65,1\}\}, -1, \{\{13,\{0,1,2,3,4\}\}\}, -1} 

\section{Theta graphs with 14 edges}
\label{sec:theta14}


\noindent {\bf Proof of Lemma \ref{lem:theta14 designs}}~

\noindent {\boldmath $K_{21}$}~
Let the vertex set be $Z_{20} \cup \{\infty\}$. The decompositions consist of the graphs

\adfLgap 
$(0,1,2,3,4,8,5,7,6,9,13,11,14)_{\adfTEABK}$,

$(0,5,10,6,2,4,11,1,7,3,8,16,\infty)_{\adfTEABK}$,

$(2,10,15,8,1,9,0,11,3,14,5,19,\infty)_{\adfTEABK}$,

\adfLgap 
$(0,1,2,3,4,7,8,5,6,9,13,11,14)_{\adfTEACJ}$,

$(0,5,6,10,7,1,8,2,4,12,3,11,\infty)_{\adfTEACJ}$,

$(1,9,10,0,11,2,7,12,\infty,14,6,19,3)_{\adfTEACJ}$,

\adfLgap 
$(0,1,2,3,4,5,6,8,11,7,9,13,10)_{\adfTEADI}$,

$(0,4,6,1,11,8,2,7,5,12,3,9,\infty)_{\adfTEADI}$,

$(0,9,10,3,17,15,6,14,11,19,\infty,18,2)_{\adfTEADI}$,

\adfLgap 
$(0,1,2,3,4,7,5,6,8,12,9,11,13)_{\adfTEAEH}$,

$(0,6,7,2,5,1,8,3,9,16,10,17,\infty)_{\adfTEAEH}$,

$(1,10,11,0,\infty,3,12,2,18,6,15,19,7)_{\adfTEAEH}$,

\adfLgap 
$(0,1,2,3,4,7,5,6,8,12,9,10,13)_{\adfTEAFG}$,

$(0,5,7,1,6,2,12,8,3,9,18,4,\infty)_{\adfTEAFG}$,

$(3,7,5,15,2,10,19,12,1,14,11,6,\infty)_{\adfTEAFG}$,

\adfLgap 
$(0,1,2,3,4,5,8,6,7,9,13,10,14)_{\adfTEBBJ}$,

$(0,1,6,7,5,12,2,8,3,4,13,\infty,10)_{\adfTEBBJ}$,

$(3,7,10,11,12,0,\infty,15,6,14,19,5,17)_{\adfTEBBJ}$,

\adfLgap 
$(0,1,2,3,4,5,7,6,8,9,13,10,14)_{\adfTEBCI}$,

$(0,1,6,4,11,8,2,7,3,5,12,\infty,9)_{\adfTEBCI}$,

$(2,6,15,12,3,10,1,7,13,4,19,11,\infty)_{\adfTEBCI}$,

\adfLgap 
$(0,1,2,3,4,5,6,7,9,11,14,10,12)_{\adfTEBDH}$,

$(0,1,4,5,2,7,8,3,9,14,6,12,\infty)_{\adfTEBDH}$,

$(2,6,13,11,3,\infty,15,19,12,5,17,7,16)_{\adfTEBDH}$,

\adfLgap 
$(0,1,2,3,4,5,7,6,8,12,9,13,10)_{\adfTEBEG}$,

$(0,1,8,5,3,2,6,7,10,15,4,14,\infty)_{\adfTEBEG}$,

$(3,7,11,8,2,9,17,14,6,19,5,16,\infty)_{\adfTEBEG}$,

\adfLgap 
$(0,1,2,3,4,5,7,6,8,12,9,13,10)_{\adfTEBFF}$,

$(0,1,7,5,3,6,2,8,9,17,10,4,\infty)_{\adfTEBFF}$,

$(3,7,11,8,6,15,4,14,9,19,\infty,10,2)_{\adfTEBFF}$,

\adfLgap 
$(0,1,2,3,4,5,6,8,7,9,10,13,16)_{\adfTECCH}$,

$(0,1,3,6,7,11,8,2,9,15,4,14,\infty)_{\adfTECCH}$,

$(3,7,8,\infty,10,2,11,17,4,13,5,14,18)_{\adfTECCH}$,

\adfLgap 
$(0,1,2,3,4,5,6,7,8,11,13,9,12)_{\adfTECDG}$,

$(0,1,5,11,6,2,8,10,3,7,12,4,\infty)_{\adfTECDG}$,

$(2,6,4,15,5,17,11,9,18,10,7,19,\infty)_{\adfTECDG}$,

\adfLgap 
$(0,1,2,3,4,5,6,8,7,9,12,11,14)_{\adfTECEF}$,

$(0,1,3,7,5,2,6,11,8,14,4,15,\infty)_{\adfTECEF}$,

$(2,6,8,17,10,19,4,\infty,15,3,9,13,1)_{\adfTECEF}$,

\adfLgap 
$(0,1,2,3,4,5,6,8,7,9,11,14,10)_{\adfTEDDF}$,

$(0,1,3,7,12,4,5,9,6,11,2,8,\infty)_{\adfTEDDF}$,

$(1,5,6,\infty,11,14,2,15,18,8,16,7,19)_{\adfTEDDF}$,

\adfLgap 
$(0,1,2,3,4,5,6,8,9,7,10,13,11)_{\adfTEDEE}$,

$(0,1,3,7,5,4,10,2,6,8,17,11,\infty)_{\adfTEDEE}$,

$(0,4,10,3,17,11,2,7,19,14,5,18,\infty)_{\adfTEDEE}$

\adfLgap
\noindent under the action of the mapping $x \mapsto x + 4$ (mod 20), $\infty \mapsto \infty$.

\ADFvfyParStart{21, \{\{20,5,4\},\{1,1,1\}\}, 20, -1, -1} 

\noindent {\boldmath $K_{28}$}~
Let the vertex set be $Z_{27} \cup \{\infty\}$. The decompositions consist of the graphs

\adfLgap 
$(\infty,0,1,2,3,5,4,6,9,13,7,10,8)_{\adfTEABK}$,

$(0,5,9,6,1,8,2,7,11,14,3,10,18)_{\adfTEABK}$,

$(0,10,20,12,1,13,2,11,21,7,16,8,23)_{\adfTEABK}$,

\adfLgap 
$(\infty,0,1,2,5,6,3,4,7,9,13,8,10)_{\adfTEACJ}$,

$(0,5,6,1,7,13,2,8,11,3,12,4,14)_{\adfTEACJ}$,

$(0,13,11,4,12,1,14,2,6,23,3,17,25)_{\adfTEACJ}$,

\adfLgap 
$(\infty,0,1,2,3,5,7,4,6,8,11,15,9)_{\adfTEADI}$,

$(0,1,4,8,13,5,11,2,9,16,3,14,6)_{\adfTEADI}$,

$(0,14,10,1,7,12,2,13,23,11,25,6,22)_{\adfTEADI}$,

\adfLgap 
$(\infty,0,1,2,3,4,5,7,9,6,8,11,15)_{\adfTEAEH}$,

$(0,5,6,1,4,10,7,11,2,8,15,23,13)_{\adfTEAEH}$,

$(0,10,9,1,8,19,13,26,12,2,17,6,22)_{\adfTEAEH}$,

\adfLgap 
$(\infty,0,1,2,3,4,6,5,7,10,14,8,11)_{\adfTEAFG}$,

$(0,2,3,7,1,6,11,8,12,4,13,5,15)_{\adfTEAFG}$,

$(0,7,9,2,13,1,18,12,25,11,16,26,14)_{\adfTEAFG}$,

\adfLgap 
$(\infty,0,1,2,3,6,4,5,7,10,14,8,9)_{\adfTEBBJ}$,

$(0,1,6,7,4,11,2,5,9,14,3,13,23)_{\adfTEBBJ}$,

$(0,2,12,13,8,15,4,16,7,20,5,18,10)_{\adfTEBBJ}$,

\adfLgap 
$(\infty,0,1,2,3,6,4,5,7,10,14,8,11)_{\adfTEBCI}$,

$(0,1,6,2,7,4,11,3,8,12,5,13,21)_{\adfTEBCI}$,

$(0,2,12,9,19,13,1,10,21,8,17,5,16)_{\adfTEBCI}$,

\adfLgap 
$(\infty,0,1,2,3,5,6,4,7,8,10,14,11)_{\adfTEBDH}$,

$(0,1,6,3,7,13,8,2,9,5,10,19,11)_{\adfTEBDH}$,

$(0,2,12,7,20,13,9,1,15,4,21,8,17)_{\adfTEBDH}$,

\adfLgap 
$(\infty,0,1,2,3,5,4,6,9,7,10,8,11)_{\adfTEBEG}$,

$(0,1,5,6,2,7,12,8,14,4,10,3,13)_{\adfTEBEG}$,

$(0,2,9,12,4,17,10,13,22,11,20,5,15)_{\adfTEBEG}$,

\adfLgap 
$(\infty,0,1,2,3,5,4,6,9,12,7,10,8)_{\adfTEBFF}$,

$(0,1,5,4,10,2,6,13,9,19,3,14,8)_{\adfTEBFF}$,

$(0,2,12,13,4,14,11,16,20,5,23,7,15)_{\adfTEBFF}$,

\adfLgap 
$(\infty,0,1,2,3,4,5,6,9,7,10,8,11)_{\adfTECCH}$,

$(0,1,5,9,6,13,8,2,7,11,4,10,15)_{\adfTECCH}$,

$(0,2,10,19,18,11,12,1,14,25,17,5,15)_{\adfTECCH}$,

\adfLgap 
$(\infty,0,1,2,3,4,6,5,7,10,14,8,9)_{\adfTECDG}$,

$(0,1,3,7,5,2,6,8,13,4,11,19,9)_{\adfTECDG}$,

$(0,2,7,17,11,21,9,13,26,8,19,4,15)_{\adfTECDG}$,

\adfLgap 
$(\infty,0,1,2,3,4,6,5,8,10,7,11,14)_{\adfTECEF}$,

$(0,1,3,7,6,2,8,13,9,4,11,18,10)_{\adfTECEF}$,

$(0,2,7,17,8,16,3,13,12,23,10,21,11)_{\adfTECEF}$,

\adfLgap 
$(\infty,0,1,2,3,5,7,4,6,8,11,15,9)_{\adfTEDDF}$,

$(0,1,5,10,3,8,2,9,11,4,13,7,6)_{\adfTEDDF}$,

$(0,2,10,20,16,12,1,13,14,22,9,26,11)_{\adfTEDDF}$,

\adfLgap 
$(\infty,0,1,2,3,5,7,4,6,9,10,14,8)_{\adfTEDEE}$,

$(0,1,2,5,9,4,10,3,8,11,16,6,15)_{\adfTEDEE}$,

$(0,2,12,1,10,17,3,23,11,22,8,25,13)_{\adfTEDEE}$

\adfLgap
\noindent under the action of the mapping $x \mapsto x + 3$ (mod 27), $\infty \mapsto \infty$.

\ADFvfyParStart{28, \{\{27,9,3\},\{1,1,1\}\}, 27, -1, -1} 
\noindent {\boldmath $K_{36}$}~
Let the vertex set be $Z_{36}$. The decompositions consist of

\adfLgap 
$(21,2,19,3,9,1,23,7,29,16,24,15,27)_{\adfTEABK}$,

$(24,26,6,21,14,28,9,18,15,10,2,23,0)_{\adfTEABK}$,

$(3,7,33,4,8,17,2,32,6,1,26,19,0)_{\adfTEABK}$,

$(31,29,30,22,18,24,27,16,1,5,12,28,6)_{\adfTEABK}$,

$(1,34,0,12,7,17,29,9,4,16,31,23,10)_{\adfTEABK}$,

\adfLgap 
$(0,23,30,7,21,6,10,26,9,8,25,15,13)_{\adfTEACJ}$,

$(30,27,2,12,1,34,3,5,32,4,20,31,6)_{\adfTEACJ}$,

$(2,19,13,33,4,8,14,26,27,31,22,20,25)_{\adfTEACJ}$,

$(33,0,2,29,3,31,19,1,5,18,4,26,25)_{\adfTEACJ}$,

$(2,20,28,11,31,0,3,4,16,29,1,13,27)_{\adfTEACJ}$,

\adfLgap 
$(0,1,30,5,2,21,14,19,22,13,33,4,24)_{\adfTEADI}$,

$(3,20,9,13,1,27,14,26,30,32,34,24,10)_{\adfTEADI}$,

$(23,0,25,3,24,17,27,9,12,21,29,6,7)_{\adfTEADI}$,

$(3,12,13,24,23,14,34,29,10,27,25,11,16)_{\adfTEADI}$,

$(12,18,4,7,3,17,2,10,28,6,35,19,27)_{\adfTEADI}$,

\adfLgap 
$(0,19,9,6,18,14,16,22,33,29,10,5,34)_{\adfTEAEH}$,

$(26,28,27,21,33,19,17,0,11,34,31,22,2)_{\adfTEAEH}$,

$(34,12,23,28,17,24,16,15,27,31,21,0,8)_{\adfTEAEH}$,

$(24,21,18,17,31,1,25,12,14,7,27,4,19)_{\adfTEAEH}$,

$(6,16,14,0,3,21,29,31,5,13,34,15,23)_{\adfTEAEH}$,

\adfLgap 
$(0,9,7,26,1,16,6,15,3,5,27,34,30)_{\adfTEAFG}$,

$(15,6,21,26,8,19,11,12,29,13,1,27,14)_{\adfTEAFG}$,

$(0,30,34,18,3,4,16,23,13,27,9,26,20)_{\adfTEAFG}$,

$(22,9,19,23,7,24,20,34,32,25,17,10,11)_{\adfTEAFG}$,

$(16,13,0,1,5,11,22,7,12,20,33,28,14)_{\adfTEAFG}$,

\adfLgap 
$(0,13,17,33,8,28,32,25,19,4,9,26,2)_{\adfTEBBJ}$,

$(21,13,6,3,30,15,28,34,1,11,27,16,19)_{\adfTEBBJ}$,

$(30,26,32,27,12,14,18,11,3,5,20,29,34)_{\adfTEBBJ}$,

$(22,2,12,27,31,32,20,6,26,3,7,24,19)_{\adfTEBBJ}$,

$(15,2,1,24,3,17,4,5,16,23,21,13,25)_{\adfTEBBJ}$,

\adfLgap 
$(0,24,17,19,22,28,13,25,14,29,33,6,30)_{\adfTEBCI}$,

$(17,15,10,19,3,7,18,1,29,32,6,2,12)_{\adfTEBCI}$,

$(29,3,16,7,21,34,4,26,24,33,8,13,11)_{\adfTEBCI}$,

$(30,4,11,22,15,3,8,24,28,6,9,25,19)_{\adfTEBCI}$,

$(25,3,2,24,12,26,6,27,14,32,31,35,29)_{\adfTEBCI}$,

\adfLgap 
$(0,30,14,5,20,3,22,10,25,7,6,29,19)_{\adfTEBDH}$,

$(34,3,33,2,20,26,5,10,15,17,9,19,31)_{\adfTEBDH}$,

$(31,10,0,15,24,12,14,5,27,23,30,2,13)_{\adfTEBDH}$,

$(19,3,16,12,5,1,8,25,6,27,33,13,4)_{\adfTEBDH}$,

$(11,9,32,14,12,8,33,21,24,18,28,0,20)_{\adfTEBDH}$,

\adfLgap 
$(0,19,16,5,32,29,13,22,24,17,18,7,33)_{\adfTEBEG}$,

$(34,28,2,27,12,20,3,9,33,14,6,21,32)_{\adfTEBEG}$,

$(11,22,2,26,32,14,27,10,8,31,17,19,29)_{\adfTEBEG}$,

$(24,17,14,9,18,6,0,15,20,21,13,11,4)_{\adfTEBEG}$,

$(15,2,19,3,6,29,33,7,13,31,11,12,24)_{\adfTEBEG}$,

\adfLgap 
$(0,2,25,28,1,21,20,18,10,17,26,15,6)_{\adfTEBFF}$,

$(9,28,7,3,29,15,23,5,34,1,13,16,30)_{\adfTEBFF}$,

$(28,20,3,15,16,11,7,1,10,18,33,14,8)_{\adfTEBFF}$,

$(7,28,22,31,14,13,5,19,4,30,18,23,33)_{\adfTEBFF}$,

$(20,0,4,5,1,3,6,7,13,18,31,15,22)_{\adfTEBFF}$,

\adfLgap 
$(0,7,24,8,21,5,28,11,17,9,26,13,25)_{\adfTECCH}$,

$(32,25,26,30,21,14,22,16,7,23,1,24,10)_{\adfTECCH}$,

$(34,31,26,2,14,19,1,10,12,8,11,3,33)_{\adfTECCH}$,

$(34,6,23,27,15,5,7,20,21,12,17,31,24)_{\adfTECCH}$,

$(11,0,14,15,16,33,21,17,24,26,12,23,10)_{\adfTECCH}$,

\adfLgap 
$(0,19,24,9,8,2,13,10,34,7,25,22,23)_{\adfTECDG}$,

$(9,23,25,15,3,34,2,7,28,8,13,12,30)_{\adfTECDG}$,

$(28,20,5,22,27,14,3,18,34,21,23,32,13)_{\adfTECDG}$,

$(29,9,33,2,14,22,0,21,32,28,15,1,10)_{\adfTECDG}$,

$(34,3,0,6,15,27,32,20,23,12,9,33,19)_{\adfTECDG}$,

\adfLgap 
$(0,31,27,17,15,25,2,12,29,1,6,4,7)_{\adfTECEF}$,

$(13,5,7,6,1,18,31,27,9,28,23,24,2)_{\adfTECEF}$,

$(1,9,32,18,26,11,19,30,24,8,14,22,2)_{\adfTECEF}$,

$(7,1,27,21,5,16,20,28,0,30,18,22,4)_{\adfTECEF}$,

$(14,3,4,5,7,10,15,28,31,22,24,0,21)_{\adfTECEF}$,

\adfLgap 
$(0,2,8,30,17,7,32,22,1,15,5,11,27)_{\adfTEDDF}$,

$(15,7,29,19,6,17,34,16,3,0,33,24,12)_{\adfTEDDF}$,

$(2,12,1,20,14,26,30,33,23,21,8,4,18)_{\adfTEDDF}$,

$(8,7,13,31,3,33,21,14,6,1,29,23,24)_{\adfTEDDF}$,

$(17,2,1,5,30,24,4,19,26,16,3,6,11)_{\adfTEDDF}$,

\adfLgap 
$(0,22,19,33,7,31,32,17,11,27,4,28,1)_{\adfTEDEE}$,

$(6,25,7,19,12,13,15,33,32,17,28,10,9)_{\adfTEDEE}$,

$(10,17,22,28,31,33,2,34,7,12,23,26,20)_{\adfTEDEE}$,

$(12,15,2,30,11,26,29,10,0,16,24,4,9)_{\adfTEDEE}$,

$(0,5,17,21,29,2,7,30,14,22,15,23,3)_{\adfTEDEE}$

\adfLgap
\noindent under the action of the mapping $x \mapsto x + 4$ (mod 36).

\ADFvfyParStart{36, \{\{36,9,4\}\}, -1, -1, -1} 

\noindent {\boldmath $K_{49}$}~
Let the vertex set be $Z_{49}$. The decompositions consist of

\adfLgap 
$(42,34,3,38,41,5,27,14,19,4,2,35,20)_{\adfTEABK}$,

$(47,13,37,4,8,0,32,9,28,34,25,26,35)_{\adfTEABK}$,

$(24,36,23,47,43,20,30,5,15,27,33,4,0)_{\adfTEABK}$,

$(13,14,34,25,24,17,33,45,47,30,28,7,36)_{\adfTEABK}$,

$(38,35,23,30,33,47,11,36,6,25,28,37,17)_{\adfTEABK}$,

$(9,6,8,39,24,11,16,3,27,17,25,47,22)_{\adfTEABK}$,

$(40,22,29,43,24,8,46,39,18,32,6,5,7)_{\adfTEABK}$,

$(21,10,39,35,5,47,15,17,3,29,32,27,38)_{\adfTEABK}$,

$(38,6,44,34,30,19,10,37,9,1,23,28,4)_{\adfTEABK}$,

$(34,29,43,18,9,5,46,24,7,19,40,45,39)_{\adfTEABK}$,

$(16,32,42,23,41,34,5,28,35,10,40,7,1)_{\adfTEABK}$,

$(1,18,30,10,15,9,40,48,16,2,36,8,7)_{\adfTEABK}$,

\adfLgap 
$(0,25,33,29,26,5,35,22,17,39,19,20,1)_{\adfTEACJ}$,

$(35,31,40,37,7,9,34,43,20,30,25,23,5)_{\adfTEACJ}$,

$(30,19,27,3,21,31,43,14,44,11,24,15,23)_{\adfTEACJ}$,

$(7,45,46,16,18,35,9,20,2,11,19,28,44)_{\adfTEACJ}$,

$(33,46,31,38,22,21,35,47,15,4,39,12,34)_{\adfTEACJ}$,

$(29,8,11,43,32,6,2,12,17,15,16,21,18)_{\adfTEACJ}$,

$(22,45,3,41,14,13,34,40,10,7,5,4,24)_{\adfTEACJ}$,

$(29,19,41,43,36,42,24,33,13,21,39,45,31)_{\adfTEACJ}$,

$(47,29,32,13,1,28,45,6,24,8,30,16,44)_{\adfTEACJ}$,

$(30,38,3,6,13,26,32,4,20,18,44,8,2)_{\adfTEACJ}$,

$(36,21,41,34,16,9,24,25,18,30,5,40,6)_{\adfTEACJ}$,

$(7,0,11,6,27,41,32,17,42,20,12,5,37)_{\adfTEACJ}$,

\adfLgap 
$(0,45,38,35,23,2,32,4,22,16,6,47,3)_{\adfTEADI}$,

$(47,1,17,32,19,39,3,6,25,5,40,44,22)_{\adfTEADI}$,

$(21,28,1,46,3,27,13,17,29,6,42,10,44)_{\adfTEADI}$,

$(10,27,24,3,2,47,16,21,39,35,14,22,29)_{\adfTEADI}$,

$(19,43,9,1,17,15,30,2,31,47,45,6,34)_{\adfTEADI}$,

$(13,40,1,10,33,4,2,25,41,26,27,0,20)_{\adfTEADI}$,

$(35,21,20,16,36,8,7,17,25,27,30,43,26)_{\adfTEADI}$,

$(40,19,18,45,25,21,20,15,10,41,46,14,6)_{\adfTEADI}$,

$(20,13,26,37,44,45,5,21,19,16,3,1,39)_{\adfTEADI}$,

$(40,25,1,18,9,14,23,0,11,12,37,2,13)_{\adfTEADI}$,

$(9,46,4,29,32,26,15,28,18,17,23,6,36)_{\adfTEADI}$,

$(23,42,31,1,36,22,6,17,35,11,18,21,33)_{\adfTEADI}$,

\adfLgap 
$(0,1,43,6,39,16,47,17,14,35,32,18,31)_{\adfTEAEH}$,

$(43,21,27,25,23,47,12,36,29,19,10,28,39)_{\adfTEAEH}$,

$(22,11,41,26,39,2,14,34,0,32,17,3,15)_{\adfTEAEH}$,

$(41,9,1,27,5,3,42,29,24,15,43,10,20)_{\adfTEAEH}$,

$(13,38,31,26,18,10,19,2,14,29,43,45,30)_{\adfTEAEH}$,

$(12,18,2,0,41,37,42,33,32,4,31,34,17)_{\adfTEAEH}$,

$(2,38,21,5,1,32,15,44,26,11,8,19,45)_{\adfTEAEH}$,

$(11,20,35,30,27,22,36,46,17,33,1,2,44)_{\adfTEAEH}$,

$(43,37,14,18,28,40,17,16,45,35,44,34,4)_{\adfTEAEH}$,

$(35,0,28,45,7,13,10,14,30,8,25,18,23)_{\adfTEAEH}$,

$(9,47,23,2,43,40,5,34,11,6,41,45,33)_{\adfTEAEH}$,

$(32,20,5,6,28,33,12,40,48,30,3,41,13)_{\adfTEAEH}$,

\adfLgap 
$(0,17,28,34,8,43,14,19,30,35,24,5,36)_{\adfTEAFG}$,

$(21,20,39,10,9,28,2,45,11,38,42,26,41)_{\adfTEAFG}$,

$(30,23,27,38,17,15,18,39,19,25,37,10,36)_{\adfTEAFG}$,

$(28,22,5,27,11,24,44,26,38,3,42,2,30)_{\adfTEAFG}$,

$(4,14,36,0,2,16,46,18,42,20,38,47,11)_{\adfTEAFG}$,

$(45,13,38,39,46,22,10,4,37,19,29,24,15)_{\adfTEAFG}$,

$(14,26,34,4,19,22,6,18,8,0,42,27,15)_{\adfTEAFG}$,

$(6,11,5,37,47,45,9,15,4,25,31,26,20)_{\adfTEAFG}$,

$(19,6,3,13,43,15,14,28,10,33,12,39,41)_{\adfTEAFG}$,

$(35,23,21,36,3,44,27,46,34,26,33,25,29)_{\adfTEAFG}$,

$(33,8,32,9,24,30,15,19,16,41,34,29,12)_{\adfTEAFG}$,

$(40,16,35,2,22,4,27,44,34,10,6,42,15)_{\adfTEAFG}$,

\adfLgap 
$(0,30,28,24,42,7,13,29,6,45,39,44,10)_{\adfTEBCI}$,

$(25,18,9,40,30,3,26,31,8,4,34,21,1)_{\adfTEBCI}$,

$(46,18,34,44,43,17,20,29,47,12,41,42,37)_{\adfTEBCI}$,

$(35,18,12,20,11,24,43,32,0,4,3,9,15)_{\adfTEBCI}$,

$(13,25,11,42,43,37,15,44,12,14,31,19,47)_{\adfTEBCI}$,

$(46,42,26,20,23,28,25,38,10,18,3,27,9)_{\adfTEBCI}$,

$(43,28,15,9,22,29,21,44,45,35,5,4,42)_{\adfTEBCI}$,

$(5,19,27,12,29,23,14,6,8,18,13,26,30)_{\adfTEBCI}$,

$(24,18,26,21,46,33,36,38,19,15,8,27,44)_{\adfTEBCI}$,

$(37,45,23,41,3,9,20,13,19,28,33,34,29)_{\adfTEBCI}$,

$(33,5,30,8,42,22,34,20,41,45,36,24,38)_{\adfTEBCI}$,

$(6,1,28,16,9,38,2,10,41,26,39,0,45)_{\adfTEBCI}$,

\adfLgap 
$(0,43,44,45,2,29,42,37,3,14,13,36,12)_{\adfTEBEG}$,

$(31,43,22,39,28,47,8,9,30,23,34,41,2)_{\adfTEBEG}$,

$(26,27,11,6,8,20,45,13,5,31,23,39,0)_{\adfTEBEG}$,

$(44,1,45,29,40,43,11,32,8,46,34,9,35)_{\adfTEBEG}$,

$(3,31,15,35,8,39,7,41,23,12,2,37,1)_{\adfTEBEG}$,

$(27,0,14,13,16,36,43,21,18,32,31,41,8)_{\adfTEBEG}$,

$(35,2,19,16,47,18,25,27,42,40,36,30,34)_{\adfTEBEG}$,

$(42,2,35,46,25,26,38,3,1,4,31,33,47)_{\adfTEBEG}$,

$(33,44,42,30,25,20,46,7,19,40,8,31,24)_{\adfTEBEG}$,

$(18,32,38,42,45,14,23,37,12,13,47,25,34)_{\adfTEBEG}$,

$(1,17,20,21,26,18,31,41,24,3,19,32,13)_{\adfTEBEG}$,

$(40,5,10,13,8,4,47,31,46,28,0,20,48)_{\adfTEBEG}$,

\adfLgap 
$(0,20,28,34,13,5,18,37,33,6,25,32,4)_{\adfTECCH}$,

$(43,33,2,15,32,38,6,46,29,36,28,24,10)_{\adfTECCH}$,

$(42,10,5,43,19,6,35,36,1,45,23,13,34)_{\adfTECCH}$,

$(29,5,26,1,20,44,4,39,23,38,28,10,46)_{\adfTECCH}$,

$(3,18,6,21,5,30,44,39,10,36,46,19,28)_{\adfTECCH}$,

$(28,13,9,30,6,8,43,27,32,7,12,41,31)_{\adfTECCH}$,

$(12,46,31,20,32,1,13,2,0,23,7,42,34)_{\adfTECCH}$,

$(0,23,11,37,32,17,36,38,39,41,21,24,43)_{\adfTECCH}$,

$(28,38,26,5,16,47,27,1,29,12,18,33,14)_{\adfTECCH}$,

$(37,7,6,2,22,16,24,41,34,15,42,26,36)_{\adfTECCH}$,

$(37,38,5,17,30,29,34,40,26,33,46,43,0)_{\adfTECCH}$,

$(26,44,23,24,37,46,25,17,39,35,3,10,22)_{\adfTECCH}$,

\adfLgap 
$(0,45,7,12,8,2,41,22,47,35,21,39,16)_{\adfTECDG}$,

$(32,2,7,35,12,6,24,46,30,16,40,47,19)_{\adfTECDG}$,

$(44,47,4,7,42,3,1,32,2,27,15,0,20)_{\adfTECDG}$,

$(13,30,32,43,4,1,40,14,20,45,10,16,34)_{\adfTECDG}$,

$(38,1,12,25,4,9,21,43,29,10,13,33,32)_{\adfTECDG}$,

$(14,45,6,17,37,19,25,40,29,3,10,44,33)_{\adfTECDG}$,

$(3,14,21,17,28,45,46,15,37,44,36,19,27)_{\adfTECDG}$,

$(36,26,34,24,45,32,11,0,4,8,37,7,41)_{\adfTECDG}$,

$(1,17,7,34,8,25,47,5,10,32,24,23,26)_{\adfTECDG}$,

$(38,36,2,13,32,6,5,30,29,44,23,20,27)_{\adfTECDG}$,

$(44,20,40,32,27,13,18,46,39,29,1,6,36)_{\adfTECDG}$,

$(14,11,9,0,12,48,27,15,31,42,26,40,21)_{\adfTECDG}$,

\adfLgap 
$(0,5,31,14,3,41,15,44,43,13,8,40,29)_{\adfTECEF}$,

$(35,32,34,20,19,29,21,8,14,7,4,11,45)_{\adfTECEF}$,

$(42,33,27,17,32,18,22,29,37,23,26,41,35)_{\adfTECEF}$,

$(34,26,44,17,14,28,5,4,25,0,22,1,13)_{\adfTECEF}$,

$(0,35,10,16,29,3,20,9,38,42,8,11,23)_{\adfTECEF}$,

$(7,32,25,30,20,27,6,9,12,44,1,17,31)_{\adfTECEF}$,

$(30,2,14,38,22,23,21,33,37,9,5,35,18)_{\adfTECEF}$,

$(0,34,1,36,9,33,47,26,41,16,5,24,43)_{\adfTECEF}$,

$(45,0,33,4,20,39,6,11,17,16,43,31,27)_{\adfTECEF}$,

$(4,10,40,33,44,8,46,18,22,6,24,17,39)_{\adfTECEF}$,

$(41,4,23,6,38,32,47,27,40,22,37,45,36)_{\adfTECEF}$,

$(34,5,30,10,12,18,37,3,40,43,38,36,46)_{\adfTECEF}$,

\adfLgap 
$(0,3,15,2,34,26,27,16,31,33,39,13,28)_{\adfTEDEE}$,

$(19,23,32,27,17,10,15,4,25,9,30,47,45)_{\adfTEDEE}$,

$(9,38,0,21,22,34,18,8,32,13,10,11,7)_{\adfTEDEE}$,

$(38,23,18,44,37,21,11,30,7,2,0,13,8)_{\adfTEDEE}$,

$(24,45,23,19,31,22,28,17,41,39,36,30,42)_{\adfTEDEE}$,

$(30,2,35,0,5,20,27,19,46,1,40,39,21)_{\adfTEDEE}$,

$(3,45,41,13,40,23,12,29,15,7,39,31,8)_{\adfTEDEE}$,

$(36,4,33,5,19,14,43,6,8,37,45,32,39)_{\adfTEDEE}$,

$(35,19,5,23,11,45,36,12,21,32,15,37,34)_{\adfTEDEE}$,

$(22,14,13,29,6,33,21,28,39,41,47,18,20)_{\adfTEDEE}$,

$(3,44,10,27,13,25,6,26,19,33,17,40,36)_{\adfTEDEE}$,

$(4,0,44,11,20,22,15,28,27,41,5,36,8)_{\adfTEDEE}$

\adfLgap
\noindent under the action of the mapping $x \mapsto x + 7$ (mod 49).

\ADFvfyParStart{49, \{\{49,7,7\}\}, -1, -1, -1} 

\noindent {\boldmath $K_{56}$}~
Let the vertex set be $Z_{55} \cup \{\infty\}$. The decompositions consist of the graphs

\adfLgap 
$(\infty,38,40,19,10,30,15,46,42,35,18,44,32)_{\adfTEABK}$,

$(26,27,7,41,46,1,48,13,5,50,29,14,10)_{\adfTEABK}$,

$(15,38,20,40,32,30,43,0,27,28,4,41,23)_{\adfTEABK}$,

$(44,3,1,52,11,50,24,47,25,53,14,46,7)_{\adfTEABK}$,

$(12,23,29,43,48,15,26,38,35,31,28,51,24)_{\adfTEABK}$,

$(53,26,17,31,51,25,14,16,49,3,24,52,37)_{\adfTEABK}$,

$(36,29,19,15,12,35,47,14,1,10,4,52,18)_{\adfTEABK}$,

$(26,1,43,52,46,22,47,29,15,39,42,21,20)_{\adfTEABK}$,

$(12,14,30,17,27,18,25,11,5,4,23,13,9)_{\adfTEABK}$,

$(43,13,27,17,4,1,9,29,54,35,22,\infty,6)_{\adfTEABK}$,

\adfLgap 
$(\infty,26,5,25,48,13,8,33,15,45,35,51,4)_{\adfTEACJ}$,

$(22,42,53,10,14,40,49,51,32,26,16,13,12)_{\adfTEACJ}$,

$(45,40,34,43,37,20,2,42,1,8,24,22,39)_{\adfTEACJ}$,

$(50,33,28,14,6,35,4,3,30,43,7,34,29)_{\adfTEACJ}$,

$(27,38,48,40,51,33,18,2,39,13,7,12,44)_{\adfTEACJ}$,

$(4,37,10,44,21,42,13,27,31,46,34,14,30)_{\adfTEACJ}$,

$(49,1,26,43,19,5,28,16,35,21,53,7,17)_{\adfTEACJ}$,

$(51,53,0,2,20,26,52,39,8,36,54,50,43)_{\adfTEACJ}$,

$(4,14,1,29,46,21,43,54,2,5,32,10,50)_{\adfTEACJ}$,

$(2,14,1,48,11,6,15,27,16,51,30,17,\infty)_{\adfTEACJ}$,

\adfLgap 
$(\infty,17,35,48,52,13,33,2,45,40,18,1,24)_{\adfTEADI}$,

$(51,23,46,25,9,34,31,18,45,47,7,54,0)_{\adfTEADI}$,

$(31,41,4,0,30,44,33,7,23,29,54,11,47)_{\adfTEADI}$,

$(44,47,54,18,33,53,21,13,34,38,17,42,48)_{\adfTEADI}$,

$(13,15,38,31,33,1,27,20,7,40,2,12,29)_{\adfTEADI}$,

$(22,36,44,0,10,17,25,31,53,50,3,15,11)_{\adfTEADI}$,

$(19,34,39,37,10,1,9,30,15,51,17,6,22)_{\adfTEADI}$,

$(17,8,54,6,26,53,15,46,45,27,31,29,24)_{\adfTEADI}$,

$(35,15,26,40,34,12,44,2,5,21,54,23,24)_{\adfTEADI}$,

$(11,2,8,3,13,12,36,21,\infty,19,48,0,29)_{\adfTEADI}$,

\adfLgap 
$(\infty,10,9,7,40,11,27,45,4,14,23,6,44)_{\adfTEAEH}$,

$(48,44,54,6,8,13,23,30,43,53,41,2,19)_{\adfTEAEH}$,

$(37,38,51,29,9,17,16,19,50,42,48,11,35)_{\adfTEAEH}$,

$(10,12,52,8,21,9,14,32,25,53,38,7,3)_{\adfTEAEH}$,

$(49,31,54,14,53,21,36,43,44,46,26,23,25)_{\adfTEAEH}$,

$(33,45,0,29,20,10,53,26,30,13,32,31,39)_{\adfTEAEH}$,

$(0,50,32,47,17,39,40,2,34,22,29,30,14)_{\adfTEAEH}$,

$(28,9,39,52,42,37,6,21,2,30,53,27,36)_{\adfTEAEH}$,

$(13,50,27,15,6,11,5,30,51,7,1,32,36)_{\adfTEAEH}$,

$(12,15,28,14,43,51,32,6,29,8,\infty,1,26)_{\adfTEAEH}$,

\adfLgap 
$(\infty,48,2,40,33,10,42,54,21,16,47,7,23)_{\adfTEAFG}$,

$(7,12,16,1,40,28,31,11,13,6,4,49,0)_{\adfTEAFG}$,

$(51,10,34,41,15,23,11,38,25,12,26,45,1)_{\adfTEAFG}$,

$(28,18,6,9,20,26,47,11,36,1,35,10,44)_{\adfTEAFG}$,

$(26,14,3,47,40,20,34,36,9,33,29,6,53)_{\adfTEAFG}$,

$(42,31,23,41,4,46,47,24,32,26,34,13,35)_{\adfTEAFG}$,

$(22,40,25,47,0,31,3,34,33,42,28,43,23)_{\adfTEAFG}$,

$(22,42,39,0,19,32,9,29,5,15,14,44,17)_{\adfTEAFG}$,

$(2,4,3,5,14,33,19,47,20,15,44,38,7)_{\adfTEAFG}$,

$(24,20,13,8,12,33,5,19,28,25,27,1,\infty)_{\adfTEAFG}$,

\adfLgap 
$(\infty,42,15,51,45,12,26,24,32,53,8,23,9)_{\adfTEBCI}$,

$(21,43,44,13,36,9,19,17,3,5,31,51,6)_{\adfTEBCI}$,

$(24,4,19,8,40,0,17,29,30,9,28,2,21)_{\adfTEBCI}$,

$(40,15,25,34,23,42,43,50,30,31,33,37,48)_{\adfTEBCI}$,

$(41,11,45,14,6,47,24,49,28,22,15,4,32)_{\adfTEBCI}$,

$(47,35,22,34,5,36,54,23,20,25,11,24,44)_{\adfTEBCI}$,

$(49,28,12,43,0,3,45,6,22,17,27,7,53)_{\adfTEBCI}$,

$(40,7,16,3,31,32,0,39,43,12,52,14,11)_{\adfTEBCI}$,

$(17,27,5,35,26,53,15,19,48,31,16,49,1)_{\adfTEBCI}$,

$(13,5,16,8,41,1,31,18,38,\infty,9,12,19)_{\adfTEBCI}$,

\adfLgap 
$(\infty,39,27,44,13,0,46,38,50,26,16,7,6)_{\adfTEBEG}$,

$(28,41,11,51,37,45,43,50,16,20,9,30,40)_{\adfTEBEG}$,

$(48,1,27,12,6,49,41,20,44,39,35,34,43)_{\adfTEBEG}$,

$(9,4,3,27,16,45,29,24,8,33,13,12,37)_{\adfTEBEG}$,

$(33,11,47,28,54,34,32,18,10,30,17,15,9)_{\adfTEBEG}$,

$(48,44,17,25,22,42,2,21,24,27,19,38,34)_{\adfTEBEG}$,

$(24,16,5,50,45,4,21,38,41,37,14,52,23)_{\adfTEBEG}$,

$(50,13,7,22,6,29,1,11,30,36,3,24,17)_{\adfTEBEG}$,

$(18,12,36,35,5,22,45,25,40,47,2,48,30)_{\adfTEBEG}$,

$(28,1,36,12,7,18,15,38,49,31,44,35,\infty)_{\adfTEBEG}$,

\adfLgap 
$(\infty,32,52,18,34,49,43,15,19,6,5,45,37)_{\adfTECCH}$,

$(30,42,1,31,6,51,36,25,35,15,13,20,11)_{\adfTECCH}$,

$(45,50,20,33,31,19,29,6,46,26,42,53,9)_{\adfTECCH}$,

$(23,11,16,37,53,13,18,17,31,15,9,36,44)_{\adfTECCH}$,

$(32,20,31,29,3,38,35,12,39,8,41,36,54)_{\adfTECCH}$,

$(1,34,7,4,49,42,43,9,6,2,48,12,14)_{\adfTECCH}$,

$(8,15,4,18,51,48,53,47,27,2,20,7,38)_{\adfTECCH}$,

$(14,4,33,37,15,27,13,25,44,54,16,20,22)_{\adfTECCH}$,

$(41,49,23,7,20,37,22,15,51,28,54,4,43)_{\adfTECCH}$,

$(44,18,15,10,53,26,0,22,7,17,45,\infty,1)_{\adfTECCH}$,

\adfLgap 
$(\infty,48,43,54,32,0,8,14,31,20,23,42,50)_{\adfTECDG}$,

$(8,22,7,35,45,20,2,12,6,10,36,3,48)_{\adfTECDG}$,

$(5,27,39,43,29,1,14,7,34,30,52,41,17)_{\adfTECDG}$,

$(53,16,34,31,30,13,25,19,41,36,29,7,10)_{\adfTECDG}$,

$(47,33,24,19,35,29,40,9,39,21,3,28,0)_{\adfTECDG}$,

$(18,7,21,37,5,10,26,31,33,38,11,19,4)_{\adfTECDG}$,

$(19,14,21,53,12,24,34,45,30,31,10,9,0)_{\adfTECDG}$,

$(48,23,41,3,49,6,47,27,20,10,30,14,32)_{\adfTECDG}$,

$(34,40,43,2,3,11,21,8,52,44,42,37,16)_{\adfTECDG}$,

$(2,0,1,\infty,8,51,19,11,31,27,42,16,41)_{\adfTECDG}$,

\adfLgap 
$(\infty,19,54,21,22,25,53,14,13,36,18,17,26)_{\adfTECEF}$,

$(37,7,17,30,8,14,41,49,43,45,46,51,26)_{\adfTECEF}$,

$(1,2,10,49,25,33,24,27,12,48,7,39,13)_{\adfTECEF}$,

$(21,54,44,10,11,24,42,33,0,9,50,7,14)_{\adfTECEF}$,

$(19,14,1,4,49,48,31,2,25,45,16,38,3)_{\adfTECEF}$,

$(37,9,33,45,21,27,5,7,10,25,32,47,13)_{\adfTECEF}$,

$(2,54,16,13,40,1,45,27,52,19,31,25,50)_{\adfTECEF}$,

$(43,18,48,25,10,0,47,46,40,45,14,49,11)_{\adfTECEF}$,

$(13,21,53,6,50,49,20,1,37,41,45,3,33)_{\adfTECEF}$,

$(4,0,28,38,23,7,17,41,42,21,8,6,\infty)_{\adfTECEF}$,

\adfLgap 
$(\infty,7,51,17,34,54,2,8,9,5,25,44,53)_{\adfTEDEE}$,

$(38,8,25,0,36,18,14,48,3,17,5,21,51)_{\adfTEDEE}$,

$(28,3,21,30,32,9,36,41,43,12,27,40,2)_{\adfTEDEE}$,

$(25,39,11,24,2,29,18,20,51,54,1,4,19)_{\adfTEDEE}$,

$(50,38,17,48,9,28,5,11,21,6,39,1,30)_{\adfTEDEE}$,

$(27,2,34,4,53,38,1,16,50,14,24,35,30)_{\adfTEDEE}$,

$(6,29,14,22,11,52,20,2,43,38,51,18,15)_{\adfTEDEE}$,

$(44,3,35,29,30,20,16,51,22,32,7,42,11)_{\adfTEDEE}$,

$(45,22,29,34,41,33,8,11,17,30,23,40,32)_{\adfTEDEE}$,

$(37,3,51,50,40,14,46,42,\infty,36,52,0,34)_{\adfTEDEE}$

\adfLgap
\noindent under the action of the mapping $x \mapsto x + 5$ (mod 55), $\infty \mapsto \infty$.

\ADFvfyParStart{56, \{\{55,11,5\},\{1,1,1\}\}, 55, -1, -1} 

\noindent {\boldmath $K_{57}$}~
Let the vertex set be $Z_{57}$. The decompositions consist of

\adfLgap 
$(0,1,3,4,9,2,8,16,5,14,24,6,18)_{\adfTEABK}$,

$(0,13,27,15,31,2,21,1,23,49,24,3,36)_{\adfTEABK}$,

\adfLgap 
$(0,1,2,5,6,11,3,10,19,4,14,25,37)_{\adfTEACJ}$,

$(0,13,14,30,18,37,2,22,45,12,41,9,39)_{\adfTEACJ}$,

\adfLgap 
$(0,1,2,5,9,6,11,4,13,3,14,26,39)_{\adfTEADI}$,

$(0,14,15,31,48,18,38,2,27,1,30,3,36)_{\adfTEADI}$,

\adfLgap 
$(0,1,2,5,9,14,6,13,3,11,20,4,15)_{\adfTEAEH}$,

$(0,12,15,32,1,30,19,41,11,43,20,53,33)_{\adfTEAEH}$,

\adfLgap 
$(0,1,2,5,9,3,8,10,18,4,13,24,36)_{\adfTEAFG}$,

$(0,13,15,31,1,18,36,19,39,3,34,5,37)_{\adfTEAFG}$,

\adfLgap 
$(0,1,2,3,7,5,12,4,13,23,6,17,29)_{\adfTEBCI}$,

$(0,1,14,15,31,18,37,2,22,45,20,46,25)_{\adfTEBCI}$,

\adfLgap 
$(0,1,2,3,7,12,18,8,15,4,13,23,35)_{\adfTEBEG}$,

$(0,1,14,15,31,2,20,21,41,6,30,5,32)_{\adfTEBEG}$,

\adfLgap 
$(0,1,2,5,6,11,7,8,16,3,12,23,35)_{\adfTECCH}$,

$(0,1,14,29,16,33,18,38,5,32,6,41,20)_{\adfTECCH}$,

\adfLgap 
$(0,1,2,5,6,7,12,8,15,3,13,4,17)_{\adfTECDG}$,

$(0,1,14,29,17,38,19,20,42,8,32,2,27)_{\adfTECDG}$,

\adfLgap 
$(0,1,2,5,6,7,12,19,8,17,3,13,24)_{\adfTECEF}$,

$(0,1,12,25,15,31,2,22,19,44,13,40,18)_{\adfTECEF}$,

\adfLgap 
$(0,1,2,3,6,4,10,17,25,9,19,5,16)_{\adfTEDEE}$,

$(0,1,12,25,41,18,37,2,29,20,43,7,32)_{\adfTEDEE}$

\adfLgap
\noindent under the action of the mapping $x \mapsto x + 1$ (mod 57).

\ADFvfyParStart{57, \{\{57,57,1\}\}, -1, -1, -1} 

\noindent {\boldmath $K_{64}$}~
Let the vertex set be $Z_{63} \cup \{\infty\}$. The decompositions consist of the graphs

\adfLgap 
$(\infty,21,23,59,31,15,27,1,48,52,16,26,61)_{\adfTEABK}$,

$(4,30,0,23,58,29,44,34,31,49,2,9,32)_{\adfTEABK}$,

$(48,2,10,19,62,20,49,54,11,12,37,51,8)_{\adfTEABK}$,

$(44,23,10,59,16,5,2,32,49,11,45,53,47)_{\adfTEABK}$,

$(53,35,55,61,0,12,49,33,51,13,4,56,25)_{\adfTEABK}$,

$(46,61,11,37,6,59,17,0,55,48,4,50,31)_{\adfTEABK}$,

$(28,27,57,2,33,37,32,35,54,13,26,17,21)_{\adfTEABK}$,

$(46,58,13,50,8,18,32,48,47,62,31,56,0)_{\adfTEABK}$,

$(6,8,12,9,22,60,45,52,38,43,15,42,28)_{\adfTEABK}$,

$(46,20,43,7,8,57,58,12,50,38,40,45,56)_{\adfTEABK}$,

$(47,40,24,33,4,27,38,8,39,45,21,9,3)_{\adfTEABK}$,

$(46,10,1,41,56,24,21,2,13,8,49,14,27)_{\adfTEABK}$,

$(37,59,41,14,29,22,12,42,0,50,33,54,36)_{\adfTEABK}$,

$(2,57,18,41,13,62,26,29,20,32,53,31,30)_{\adfTEABK}$,

$(55,36,47,31,57,0,54,3,1,33,52,53,\infty)_{\adfTEABK}$,

$(3,46,53,56,1,23,14,55,\infty,19,43,32,5)_{\adfTEABK}$,

\adfLgap 
$(\infty,11,61,57,48,38,35,31,0,2,16,60,46)_{\adfTEACJ}$,

$(2,34,20,41,11,10,40,27,37,22,44,33,49)_{\adfTEACJ}$,

$(55,39,47,7,61,5,29,34,19,9,51,25,45)_{\adfTEACJ}$,

$(9,24,49,47,15,8,19,37,21,29,18,10,46)_{\adfTEACJ}$,

$(40,9,20,34,22,38,15,18,62,10,31,43,39)_{\adfTEACJ}$,

$(50,18,48,20,33,39,23,58,29,10,35,13,36)_{\adfTEACJ}$,

$(36,7,55,0,14,56,28,42,24,39,46,26,4)_{\adfTEACJ}$,

$(17,61,41,42,13,38,24,58,45,47,10,19,33)_{\adfTEACJ}$,

$(10,26,16,18,13,45,8,20,9,21,41,40,25)_{\adfTEACJ}$,

$(28,23,43,15,61,40,52,57,8,10,0,39,26)_{\adfTEACJ}$,

$(31,60,13,18,22,35,40,50,12,46,23,49,20)_{\adfTEACJ}$,

$(57,32,37,41,54,18,55,33,16,15,42,6,22)_{\adfTEACJ}$,

$(51,12,34,37,42,4,49,29,62,32,45,9,48)_{\adfTEACJ}$,

$(3,10,49,38,58,2,60,48,42,46,22,43,30)_{\adfTEACJ}$,

$(33,62,42,36,1,0,53,2,3,20,29,\infty,49)_{\adfTEACJ}$,

$(0,24,30,57,11,28,9,5,42,54,59,37,\infty)_{\adfTEACJ}$,

\adfLgap 
$(\infty,10,62,15,1,46,22,4,3,27,50,16,37)_{\adfTEADI}$,

$(17,36,60,21,39,48,26,23,52,3,13,9,59)_{\adfTEADI}$,

$(51,8,0,27,46,56,42,61,45,4,23,3,29)_{\adfTEADI}$,

$(51,29,45,3,62,47,52,54,60,55,1,26,12)_{\adfTEADI}$,

$(40,28,39,9,54,55,48,22,61,32,34,26,15)_{\adfTEADI}$,

$(8,49,27,11,22,40,6,57,9,48,34,2,60)_{\adfTEADI}$,

$(11,46,51,7,49,4,14,8,6,32,15,10,33)_{\adfTEADI}$,

$(11,19,3,55,30,47,48,61,7,31,0,46,13)_{\adfTEADI}$,

$(49,24,31,42,16,52,17,47,30,44,57,45,23)_{\adfTEADI}$,

$(15,9,16,32,7,22,26,19,60,48,6,52,27)_{\adfTEADI}$,

$(18,61,50,22,14,24,39,35,0,16,40,9,19)_{\adfTEADI}$,

$(6,24,14,29,28,23,61,41,21,19,55,42,11)_{\adfTEADI}$,

$(7,37,36,44,34,41,18,8,26,56,62,21,30)_{\adfTEADI}$,

$(12,38,31,4,8,3,10,58,46,44,16,56,36)_{\adfTEADI}$,

$(56,61,33,5,8,10,39,53,44,55,52,36,\infty)_{\adfTEADI}$,

$(2,39,\infty,0,48,29,34,6,7,14,24,12,60)_{\adfTEADI}$,

\adfLgap 
$(\infty,9,49,16,5,35,47,46,57,20,10,19,23)_{\adfTEAEH}$,

$(2,11,36,0,17,18,50,49,55,9,40,10,26)_{\adfTEAEH}$,

$(41,29,53,10,54,16,42,14,22,0,3,9,39)_{\adfTEAEH}$,

$(33,53,30,4,39,34,10,22,5,41,14,18,59)_{\adfTEAEH}$,

$(21,14,7,23,46,1,45,20,6,53,47,34,19)_{\adfTEAEH}$,

$(36,44,40,26,57,1,41,49,27,6,62,2,9)_{\adfTEAEH}$,

$(6,58,56,40,29,48,51,32,21,23,24,14,12)_{\adfTEAEH}$,

$(62,45,30,12,10,38,54,18,47,6,26,49,11)_{\adfTEAEH}$,

$(11,42,35,55,24,2,43,45,7,52,9,17,38)_{\adfTEAEH}$,

$(7,59,1,62,27,4,60,57,34,33,40,49,29)_{\adfTEAEH}$,

$(40,18,34,9,4,43,22,6,30,51,55,46,20)_{\adfTEAEH}$,

$(50,36,35,1,60,9,29,19,58,0,44,59,10)_{\adfTEAEH}$,

$(53,51,3,61,24,29,39,60,44,50,59,36,17)_{\adfTEAEH}$,

$(6,24,35,54,33,8,10,25,7,53,36,27,60)_{\adfTEAEH}$,

$(62,29,59,47,60,5,33,61,51,0,42,10,\infty)_{\adfTEAEH}$,

$(0,12,9,45,34,49,60,\infty,55,31,44,43,15)_{\adfTEAEH}$,

\adfLgap 
$(\infty,51,43,62,52,33,30,42,16,53,4,5,60)_{\adfTEAFG}$,

$(33,36,60,2,34,32,44,58,42,26,27,41,37)_{\adfTEAFG}$,

$(46,20,8,18,50,26,2,56,60,14,40,58,35)_{\adfTEAFG}$,

$(46,25,40,62,7,41,13,18,33,10,43,39,17)_{\adfTEAFG}$,

$(38,35,62,53,46,28,8,55,24,25,47,7,3)_{\adfTEAFG}$,

$(23,33,59,53,21,12,20,45,41,28,49,25,54)_{\adfTEAFG}$,

$(13,55,18,15,26,46,44,22,57,42,29,11,14)_{\adfTEAFG}$,

$(22,17,45,20,26,9,5,0,58,2,36,48,42)_{\adfTEAFG}$,

$(60,16,30,17,36,57,45,49,47,29,33,8,32)_{\adfTEAFG}$,

$(32,17,45,31,52,34,1,9,58,47,16,41,12)_{\adfTEAFG}$,

$(43,42,23,21,28,57,9,6,25,12,7,56,50)_{\adfTEAFG}$,

$(47,15,1,6,33,26,17,12,45,2,29,19,62)_{\adfTEAFG}$,

$(3,57,10,16,31,49,11,32,48,25,55,2,43)_{\adfTEAFG}$,

$(40,24,3,46,57,31,7,54,42,30,27,0,35)_{\adfTEAFG}$,

$(0,10,44,57,50,48,45,9,55,62,47,\infty,46)_{\adfTEAFG}$,

$(3,5,2,10,\infty,6,49,55,16,51,57,34,35)_{\adfTEAFG}$,

\adfLgap 
$(\infty,45,22,44,10,19,25,49,46,60,52,9,11)_{\adfTEBCI}$,

$(58,44,53,6,40,34,49,41,43,16,21,7,36)_{\adfTEBCI}$,

$(52,14,3,33,51,32,13,55,9,44,62,36,42)_{\adfTEBCI}$,

$(34,7,10,11,50,38,39,51,45,57,24,6,52)_{\adfTEBCI}$,

$(56,60,37,4,33,24,14,55,36,27,59,46,15)_{\adfTEBCI}$,

$(50,53,46,30,27,28,52,11,21,47,35,23,16)_{\adfTEBCI}$,

$(3,39,41,12,11,19,9,42,2,52,46,16,54)_{\adfTEBCI}$,

$(22,23,40,16,38,36,57,62,15,7,5,12,13)_{\adfTEBCI}$,

$(58,24,27,59,47,16,7,3,13,61,30,8,17)_{\adfTEBCI}$,

$(10,45,12,35,18,62,5,27,11,1,61,47,60)_{\adfTEBCI}$,

$(41,26,13,28,15,53,36,21,25,46,8,11,6)_{\adfTEBCI}$,

$(26,23,37,5,36,4,13,51,50,62,0,47,39)_{\adfTEBCI}$,

$(27,6,20,7,39,50,52,8,3,40,21,57,33)_{\adfTEBCI}$,

$(1,47,57,29,44,30,28,49,11,54,26,31,14)_{\adfTEBCI}$,

$(31,14,15,58,19,23,7,41,37,18,0,6,\infty)_{\adfTEBCI}$,

$(0,20,56,54,50,40,15,4,\infty,24,3,29,12)_{\adfTEBCI}$,

\adfLgap 
$(\infty,52,57,32,10,31,38,19,14,25,8,7,24)_{\adfTEBEG}$,

$(37,5,62,56,23,41,22,13,60,58,48,0,57)_{\adfTEBEG}$,

$(28,12,50,11,25,18,56,7,43,30,35,24,21)_{\adfTEBEG}$,

$(12,46,43,47,39,50,16,15,35,11,36,17,13)_{\adfTEBEG}$,

$(12,30,2,41,32,36,22,27,45,33,13,24,28)_{\adfTEBEG}$,

$(32,23,1,24,50,2,40,47,5,56,11,13,39)_{\adfTEBEG}$,

$(18,33,38,21,29,6,9,19,10,34,2,11,53)_{\adfTEBEG}$,

$(61,47,34,35,50,59,57,22,45,48,1,6,25)_{\adfTEBEG}$,

$(5,19,16,38,25,53,59,52,28,35,3,37,49)_{\adfTEBEG}$,

$(24,11,51,25,35,22,40,49,15,57,18,41,30)_{\adfTEBEG}$,

$(58,54,52,32,17,53,47,40,56,9,55,51,50)_{\adfTEBEG}$,

$(0,43,23,34,8,38,55,41,62,42,7,52,37)_{\adfTEBEG}$,

$(4,38,22,55,54,5,13,35,49,26,52,62,6)_{\adfTEBEG}$,

$(40,52,44,62,4,38,51,53,43,15,48,35,9)_{\adfTEBEG}$,

$(28,13,27,37,25,21,48,34,36,43,16,2,\infty)_{\adfTEBEG}$,

$(2,3,44,60,33,35,\infty,58,26,45,57,48,56)_{\adfTEBEG}$,

\adfLgap 
$(\infty,15,57,7,3,19,41,40,8,9,23,33,13)_{\adfTECCH}$,

$(10,57,37,4,23,44,39,6,62,2,5,41,27)_{\adfTECCH}$,

$(33,42,40,18,28,37,4,48,0,25,30,49,24)_{\adfTECCH}$,

$(55,28,46,32,35,42,30,18,21,31,13,11,57)_{\adfTECCH}$,

$(0,46,53,47,15,43,55,49,12,24,58,11,22)_{\adfTECCH}$,

$(1,34,16,6,15,58,32,40,61,38,55,14,60)_{\adfTECCH}$,

$(25,41,4,54,16,33,53,26,1,57,40,6,15)_{\adfTECCH}$,

$(59,36,3,58,35,38,45,23,15,55,50,14,10)_{\adfTECCH}$,

$(57,42,53,59,17,16,35,5,7,23,55,4,60)_{\adfTECCH}$,

$(42,44,36,12,13,9,22,40,62,24,4,45,15)_{\adfTECCH}$,

$(57,2,3,35,51,45,6,58,56,40,21,42,14)_{\adfTECCH}$,

$(27,36,6,17,43,33,32,59,4,44,40,2,9)_{\adfTECCH}$,

$(19,26,58,41,52,61,8,50,3,45,30,4,20)_{\adfTECCH}$,

$(9,41,11,42,49,38,53,5,56,6,51,33,31)_{\adfTECCH}$,

$(39,24,19,5,57,52,28,55,31,26,37,\infty,11)_{\adfTECCH}$,

$(3,31,4,27,18,43,34,15,14,54,0,\infty,5)_{\adfTECCH}$,

\adfLgap 
$(\infty,62,37,40,22,56,42,21,46,14,45,10,36)_{\adfTECDG}$,

$(53,28,58,12,37,19,7,54,18,41,62,38,55)_{\adfTECDG}$,

$(9,36,21,28,17,39,60,29,50,40,13,42,25)_{\adfTECDG}$,

$(50,5,17,46,57,13,11,18,33,19,42,47,29)_{\adfTECDG}$,

$(20,9,5,8,39,31,62,46,18,28,50,15,11)_{\adfTECDG}$,

$(12,59,62,34,52,61,41,56,32,19,53,36,44)_{\adfTECDG}$,

$(22,1,16,59,60,51,30,20,8,12,5,31,10)_{\adfTECDG}$,

$(58,39,44,28,8,59,45,25,32,51,24,22,35)_{\adfTECDG}$,

$(61,13,6,10,28,25,22,35,8,14,44,7,51)_{\adfTECDG}$,

$(29,60,45,42,13,47,59,54,23,16,5,36,46)_{\adfTECDG}$,

$(10,17,25,51,26,31,56,0,61,2,13,60,6)_{\adfTECDG}$,

$(23,16,55,54,58,22,7,1,41,31,18,13,20)_{\adfTECDG}$,

$(19,13,17,28,54,11,62,38,42,24,60,27,14)_{\adfTECDG}$,

$(14,19,37,13,5,60,9,57,40,61,50,55,58)_{\adfTECDG}$,

$(28,45,34,1,23,8,22,0,17,24,2,20,\infty)_{\adfTECDG}$,

$(3,55,16,14,46,1,0,2,59,26,\infty,32,9)_{\adfTECDG}$,

\adfLgap 
$(\infty,5,45,19,32,51,3,30,16,9,11,26,15)_{\adfTECEF}$,

$(13,6,39,53,40,33,52,46,30,44,36,51,7)_{\adfTECEF}$,

$(1,44,5,56,22,45,8,13,25,14,34,50,38)_{\adfTECEF}$,

$(38,57,30,10,40,55,36,34,17,50,39,47,60)_{\adfTECEF}$,

$(58,0,15,46,21,61,30,19,48,52,50,4,57)_{\adfTECEF}$,

$(24,14,21,22,6,30,31,59,54,19,18,44,43)_{\adfTECEF}$,

$(24,19,31,6,49,11,47,22,32,53,18,36,56)_{\adfTECEF}$,

$(1,49,37,19,36,0,30,41,34,39,7,20,10)_{\adfTECEF}$,

$(40,23,17,20,56,10,59,7,6,5,51,4,33)_{\adfTECEF}$,

$(36,3,17,8,27,48,39,44,5,26,50,28,55)_{\adfTECEF}$,

$(10,25,5,7,42,21,28,32,27,41,60,54,49)_{\adfTECEF}$,

$(31,44,4,47,35,25,13,48,18,43,26,38,53)_{\adfTECEF}$,

$(23,39,29,43,58,13,46,17,44,21,16,41,61)_{\adfTECEF}$,

$(10,14,44,42,34,32,31,62,53,2,15,27,1)_{\adfTECEF}$,

$(11,62,44,40,51,42,13,5,14,0,15,8,\infty)_{\adfTECEF}$,

$(0,1,10,23,41,34,28,19,\infty,47,38,54,30)_{\adfTECEF}$,

\adfLgap 
$(\infty,44,25,22,20,58,50,17,43,12,26,60,53)_{\adfTEDEE}$,

$(10,23,\infty,48,19,56,3,35,30,41,34,28,1)_{\adfTEDEE}$,

$(58,19,33,23,0,43,62,11,21,39,14,7,41)_{\adfTEDEE}$,

$(9,45,38,54,1,11,36,28,39,52,49,34,20)_{\adfTEDEE}$,

$(13,5,17,53,24,46,49,43,3,14,58,15,29)_{\adfTEDEE}$,

$(30,38,28,44,26,1,10,32,22,24,42,20,5)_{\adfTEDEE}$,

$(25,26,24,0,54,4,47,1,22,51,52,32,13)_{\adfTEDEE}$,

$(44,29,5,4,22,6,32,9,62,23,41,18,46)_{\adfTEDEE}$,

$(11,29,50,62,60,22,7,54,49,44,33,27,17)_{\adfTEDEE}$,

$(42,20,18,35,57,12,49,58,23,17,3,53,37)_{\adfTEDEE}$,

$(19,36,45,10,23,12,4,40,13,20,48,35,49)_{\adfTEDEE}$,

$(21,59,48,39,51,41,50,49,7,61,19,4,62)_{\adfTEDEE}$,

$(11,22,33,41,57,5,36,7,19,60,2,16,27)_{\adfTEDEE}$,

$(50,2,45,60,28,5,43,54,34,52,6,48,24)_{\adfTEDEE}$,

$(3,12,10,31,58,26,17,21,25,48,32,36,9)_{\adfTEDEE}$,

$(0,34,35,53,45,51,3,61,30,\infty,15,9,42)_{\adfTEDEE}$

\adfLgap
\noindent under the action of the mapping $x \mapsto x + 7$ (mod 63), $\infty \mapsto \infty$.

\ADFvfyParStart{64, \{\{63,9,7\},\{1,1,1\}\}, 63, -1, -1} 

\noindent {\boldmath $K_{77}$}~
Let the vertex set be $Z_{76} \cup \{\infty\}$. The decompositions consist of the graphs

\adfLgap 
$(0,1,\infty,68,54,56,10,7,42,40,71,27,73)_{\adfTEABK}$,

$(32,61,0,45,6,25,49,38,11,36,48,31,40)_{\adfTEABK}$,

$(9,35,36,1,41,42,3,25,23,47,13,75,61)_{\adfTEABK}$,

$(7,2,61,67,40,1,70,23,0,69,36,6,56)_{\adfTEABK}$,

$(63,20,15,70,61,4,38,41,9,68,75,29,62)_{\adfTEABK}$,

$(12,26,62,23,36,18,28,25,9,37,17,59,50)_{\adfTEABK}$,

$(51,4,13,43,18,65,2,23,19,37,47,41,52)_{\adfTEABK}$,

$(3,67,34,1,24,2,29,14,37,49,18,69,11)_{\adfTEABK}$,

$(51,12,57,61,26,18,71,16,21,42,43,7,22)_{\adfTEABK}$,

$(5,40,56,10,38,58,70,66,50,44,63,46,2)_{\adfTEABK}$,

$(67,52,32,54,43,40,0,4,28,22,3,\infty,70)_{\adfTEABK}$,

\adfLgap 
$(0,1,\infty,33,67,36,53,49,26,72,45,19,38)_{\adfTEACJ}$,

$(43,63,47,57,51,48,44,18,10,9,4,20,40)_{\adfTEACJ}$,

$(38,40,2,74,47,42,26,3,37,72,8,66,69)_{\adfTEACJ}$,

$(9,1,24,66,26,14,43,31,18,49,4,22,13)_{\adfTEACJ}$,

$(41,8,64,63,44,4,48,42,29,13,33,20,59)_{\adfTEACJ}$,

$(19,22,5,39,30,4,52,61,9,50,74,35,65)_{\adfTEACJ}$,

$(23,39,75,12,64,7,36,50,21,70,13,18,32)_{\adfTEACJ}$,

$(45,75,38,32,8,67,73,25,4,6,31,29,27)_{\adfTEACJ}$,

$(21,28,43,50,57,67,27,26,54,29,48,40,18)_{\adfTEACJ}$,

$(10,0,30,15,51,69,54,8,33,71,39,13,24)_{\adfTEACJ}$,

$(12,50,7,29,23,2,34,3,17,35,62,19,\infty)_{\adfTEACJ}$,

\adfLgap 
$(0,1,\infty,6,45,30,65,49,64,24,41,58,15)_{\adfTEADI}$,

$(5,75,35,52,38,10,26,16,69,67,12,14,15)_{\adfTEADI}$,

$(14,54,36,26,34,10,60,25,0,31,75,24,5)_{\adfTEADI}$,

$(9,36,33,35,40,67,5,53,16,49,69,18,63)_{\adfTEADI}$,

$(57,34,3,62,43,19,66,9,27,36,74,10,64)_{\adfTEADI}$,

$(64,67,62,14,13,9,40,11,38,20,19,71,52)_{\adfTEADI}$,

$(17,57,13,60,65,20,55,49,62,51,28,37,47)_{\adfTEADI}$,

$(11,62,24,35,7,61,6,59,67,21,54,2,20)_{\adfTEADI}$,

$(33,30,2,43,1,42,53,11,67,60,28,4,74)_{\adfTEADI}$,

$(54,59,61,49,23,67,63,66,52,72,24,36,20)_{\adfTEADI}$,

$(37,24,44,52,58,22,7,14,40,29,\infty,3,67)_{\adfTEADI}$,

\adfLgap 
$(0,1,\infty,6,5,26,47,67,27,14,54,60,35)_{\adfTEAEH}$,

$(9,57,69,50,70,51,42,64,72,21,11,34,65)_{\adfTEAEH}$,

$(16,22,74,7,46,8,25,49,72,53,39,14,37)_{\adfTEAEH}$,

$(49,37,15,72,0,64,5,27,44,42,16,53,51)_{\adfTEAEH}$,

$(68,20,14,1,36,70,52,5,63,9,39,66,59)_{\adfTEAEH}$,

$(26,10,54,43,12,55,29,49,23,27,20,53,14)_{\adfTEAEH}$,

$(12,36,47,15,5,54,48,33,16,9,74,30,27)_{\adfTEAEH}$,

$(47,20,14,60,30,31,41,28,72,58,49,19,17)_{\adfTEAEH}$,

$(41,3,36,31,59,67,37,42,47,24,27,40,61)_{\adfTEAEH}$,

$(31,55,46,58,6,14,47,32,52,51,72,41,34)_{\adfTEAEH}$,

$(32,42,22,56,45,71,34,5,41,6,23,\infty,25)_{\adfTEAEH}$,

\adfLgap 
$(0,1,\infty,43,15,23,37,60,62,57,26,70,24)_{\adfTEAFG}$,

$(38,15,53,55,75,56,3,0,71,45,4,70,61)_{\adfTEAFG}$,

$(23,14,34,52,35,8,63,66,38,39,7,32,0)_{\adfTEAFG}$,

$(49,25,28,9,69,50,16,45,20,12,27,6,75)_{\adfTEAFG}$,

$(50,24,19,34,16,66,10,13,2,74,21,18,35)_{\adfTEAFG}$,

$(17,60,46,16,56,61,48,42,2,54,15,34,1)_{\adfTEAFG}$,

$(71,11,74,0,72,6,33,42,50,75,8,56,63)_{\adfTEAFG}$,

$(20,14,9,41,34,39,26,19,32,5,15,11,13)_{\adfTEAFG}$,

$(8,5,39,45,51,9,17,32,3,41,69,15,12)_{\adfTEAFG}$,

$(4,39,47,1,40,69,49,46,24,30,52,21,3)_{\adfTEAFG}$,

$(71,37,57,70,35,17,58,32,12,73,14,30,\infty)_{\adfTEAFG}$,

\adfLgap 
$(0,1,\infty,65,42,74,72,50,71,58,44,67,27)_{\adfTEBCI}$,

$(6,25,5,20,60,53,0,58,39,17,22,63,67)_{\adfTEBCI}$,

$(65,66,23,69,0,28,55,72,64,68,49,31,7)_{\adfTEBCI}$,

$(39,56,11,60,63,62,68,10,55,20,8,42,2)_{\adfTEBCI}$,

$(13,20,63,62,53,49,21,65,64,32,2,31,70)_{\adfTEBCI}$,

$(59,53,50,72,34,8,37,47,49,62,18,11,10)_{\adfTEBCI}$,

$(47,46,61,9,54,41,29,7,72,12,37,45,40)_{\adfTEBCI}$,

$(19,54,51,24,33,5,30,10,38,43,70,44,65)_{\adfTEBCI}$,

$(64,3,63,30,18,49,4,23,34,38,62,32,71)_{\adfTEBCI}$,

$(52,3,32,49,33,25,55,67,42,58,65,41,59)_{\adfTEBCI}$,

$(68,17,7,20,11,44,54,37,24,41,2,\infty,19)_{\adfTEBCI}$,

\adfLgap 
$(0,1,\infty,22,11,26,8,36,71,3,54,40,2)_{\adfTEBEG}$,

$(43,42,32,14,63,22,18,0,30,72,12,57,69)_{\adfTEBEG}$,

$(50,49,59,34,11,39,51,19,56,53,3,72,21)_{\adfTEBEG}$,

$(40,1,30,13,75,4,58,23,52,29,7,68,44)_{\adfTEBEG}$,

$(48,51,60,27,43,1,47,0,42,34,36,23,24)_{\adfTEBEG}$,

$(18,16,53,74,31,38,2,9,56,12,6,34,1)_{\adfTEBEG}$,

$(26,40,63,24,75,69,10,7,53,74,19,5,43)_{\adfTEBEG}$,

$(54,1,57,66,71,27,53,61,48,56,52,32,37)_{\adfTEBEG}$,

$(5,45,13,65,59,40,49,10,60,25,8,14,27)_{\adfTEBEG}$,

$(64,61,46,19,29,63,7,9,28,2,65,42,59)_{\adfTEBEG}$,

$(49,62,18,60,61,3,\infty,10,7,47,23,22,73)_{\adfTEBEG}$,

\adfLgap 
$(0,1,\infty,18,57,54,21,44,38,65,2,71,30)_{\adfTECCH}$,

$(27,32,4,7,14,12,46,26,22,43,35,69,20)_{\adfTECCH}$,

$(67,68,12,21,43,66,24,74,44,55,27,39,9)_{\adfTECCH}$,

$(29,26,51,20,67,71,3,50,32,73,22,47,56)_{\adfTECCH}$,

$(25,57,9,53,64,56,27,54,69,33,0,32,39)_{\adfTECCH}$,

$(45,48,23,72,43,53,71,35,1,64,61,5,6)_{\adfTECCH}$,

$(58,1,20,10,48,8,6,38,66,71,27,64,13)_{\adfTECCH}$,

$(51,4,66,58,57,62,34,20,36,31,74,53,46)_{\adfTECCH}$,

$(75,63,60,45,74,7,65,51,64,53,42,5,57)_{\adfTECCH}$,

$(12,21,38,35,74,62,8,56,27,68,14,33,66)_{\adfTECCH}$,

$(1,67,29,21,32,51,\infty,3,4,63,26,42,2)_{\adfTECCH}$,

\adfLgap 
$(0,1,\infty,27,22,12,69,64,46,71,58,59,19)_{\adfTECDG}$,

$(69,68,47,74,50,53,18,49,1,48,11,13,45)_{\adfTECDG}$,

$(44,3,34,4,15,49,22,68,11,41,27,72,74)_{\adfTECDG}$,

$(30,58,42,5,7,51,45,67,71,6,46,14,73)_{\adfTECDG}$,

$(25,41,12,46,49,63,68,59,71,1,61,50,29)_{\adfTECDG}$,

$(32,30,0,37,53,13,16,49,40,55,73,68,51)_{\adfTECDG}$,

$(52,26,50,59,37,68,10,31,57,6,35,0,43)_{\adfTECDG}$,

$(64,62,11,14,44,52,53,29,67,1,23,13,58)_{\adfTECDG}$,

$(7,40,35,15,14,55,31,22,2,36,8,62,10)_{\adfTECDG}$,

$(6,31,12,39,49,2,47,14,52,1,8,70,20)_{\adfTECDG}$,

$(11,28,24,31,57,68,32,9,5,6,\infty,1,44)_{\adfTECDG}$,

\adfLgap 
$(0,1,\infty,35,75,46,20,51,48,72,37,28,30)_{\adfTECEF}$,

$(33,21,3,32,54,16,75,23,42,41,14,2,34)_{\adfTECEF}$,

$(63,68,22,37,19,65,59,3,10,40,21,17,62)_{\adfTECEF}$,

$(34,33,4,25,27,31,26,11,71,16,30,72,28)_{\adfTECEF}$,

$(70,24,5,44,19,56,64,60,27,15,14,38,47)_{\adfTECEF}$,

$(74,8,66,62,64,30,19,68,52,67,9,56,5)_{\adfTECEF}$,

$(16,54,4,70,22,39,49,29,23,56,38,73,21)_{\adfTECEF}$,

$(50,13,70,30,64,7,10,49,23,51,41,59,61)_{\adfTECEF}$,

$(69,37,5,43,7,38,19,44,25,41,27,11,3)_{\adfTECEF}$,

$(33,42,55,68,56,73,23,63,38,45,44,62,65)_{\adfTECEF}$,

$(11,4,20,53,74,26,\infty,17,8,6,25,40,75)_{\adfTECEF}$,

\adfLgap 
$(0,1,\infty,15,7,2,73,26,27,55,3,75,19)_{\adfTEDEE}$,

$(4,8,71,36,24,13,66,50,57,34,74,16,44)_{\adfTEDEE}$,

$(59,54,69,12,6,15,68,39,65,0,53,38,67)_{\adfTEDEE}$,

$(71,73,44,12,36,40,48,15,57,46,33,35,74)_{\adfTEDEE}$,

$(45,52,13,46,9,62,53,67,58,10,27,17,32)_{\adfTEDEE}$,

$(25,50,12,33,71,17,28,69,29,55,56,52,42)_{\adfTEDEE}$,

$(18,26,11,54,40,62,52,47,50,75,64,3,0)_{\adfTEDEE}$,

$(70,56,14,36,54,20,65,48,49,32,51,10,61)_{\adfTEDEE}$,

$(61,54,30,35,23,64,10,63,12,41,38,50,39)_{\adfTEDEE}$,

$(29,21,25,6,33,15,31,9,67,71,1,48,73)_{\adfTEDEE}$,

$(21,36,49,27,75,19,55,48,6,\infty,46,50,23)_{\adfTEDEE}$

\adfLgap
\noindent under the action of the mapping $x \mapsto x + 4$ (mod 76), $\infty \mapsto \infty$.

\ADFvfyParStart{77, \{\{76,19,4\},\{1,1,1\}\}, 76, -1, -1} 

\noindent {\boldmath $K_{92}$}~
Let the vertex set be $Z_{92}$. The decompositions consist of the graphs

\adfLgap 
$(63,36,7,78,58,73,49,39,28,1,80,0,67)_{\adfTEABK}$,

$(80,83,37,43,2,56,8,33,84,69,14,67,38)_{\adfTEABK}$,

$(9,71,75,28,38,7,79,64,19,53,6,13,63)_{\adfTEABK}$,

$(30,12,6,43,25,16,66,9,29,64,46,68,84)_{\adfTEABK}$,

$(19,81,69,21,74,16,11,82,38,78,89,10,73)_{\adfTEABK}$,

$(3,43,71,10,88,4,27,2,86,69,50,73,75)_{\adfTEABK}$,

$(22,1,55,50,45,17,78,82,81,25,58,42,15)_{\adfTEABK}$,

$(69,64,11,48,84,81,21,73,44,13,89,35,34)_{\adfTEABK}$,

$(39,80,23,83,24,58,72,33,29,88,22,49,12)_{\adfTEABK}$,

$(65,72,76,14,26,31,40,8,48,70,59,81,78)_{\adfTEABK}$,

$(39,22,13,40,17,23,14,52,24,43,86,84,58)_{\adfTEABK}$,

$(16,14,46,6,49,74,19,42,23,41,51,64,71)_{\adfTEABK}$,

$(2,91,48,52,5,4,21,7,13,61,83,3,20)_{\adfTEABK}$,

\adfLgap 
$(0,60,9,62,10,17,88,49,84,37,80,76,55)_{\adfTEACJ}$,

$(21,37,57,24,4,67,74,63,48,32,52,51,34)_{\adfTEACJ}$,

$(83,28,63,78,47,66,88,42,2,80,50,75,15)_{\adfTEACJ}$,

$(49,51,79,76,82,57,7,47,35,28,65,11,20)_{\adfTEACJ}$,

$(69,39,54,26,34,40,59,31,53,51,55,9,90)_{\adfTEACJ}$,

$(83,73,56,8,35,21,55,71,20,31,74,1,88)_{\adfTEACJ}$,

$(11,46,50,54,3,62,17,81,56,64,57,43,66)_{\adfTEACJ}$,

$(20,1,63,45,9,5,74,56,44,8,60,90,69)_{\adfTEACJ}$,

$(80,12,52,34,78,44,67,22,27,37,66,0,59)_{\adfTEACJ}$,

$(9,69,60,26,46,34,2,33,15,81,84,27,61)_{\adfTEACJ}$,

$(45,33,65,39,18,66,40,30,84,53,76,23,50)_{\adfTEACJ}$,

$(36,37,19,89,42,51,14,64,50,79,25,66,87)_{\adfTEACJ}$,

$(51,27,57,58,50,14,1,6,22,46,84,21,30)_{\adfTEACJ}$,

\adfLgap 
$(0,21,44,84,60,81,71,31,50,43,65,51,14)_{\adfTEADI}$,

$(50,70,29,63,8,34,80,86,42,38,3,49,2)_{\adfTEADI}$,

$(15,54,64,46,60,65,22,57,42,14,80,71,3)_{\adfTEADI}$,

$(9,4,3,15,35,25,72,81,85,77,21,80,2)_{\adfTEADI}$,

$(64,71,86,6,75,28,56,53,78,20,67,23,8)_{\adfTEADI}$,

$(27,0,8,42,75,69,28,44,36,35,13,76,61)_{\adfTEADI}$,

$(23,64,51,45,14,56,46,84,88,61,68,79,76)_{\adfTEADI}$,

$(45,17,6,15,26,33,3,28,78,77,75,30,55)_{\adfTEADI}$,

$(32,57,67,37,34,71,5,19,55,87,38,65,33)_{\adfTEADI}$,

$(79,5,63,62,60,14,35,38,57,0,73,72,49)_{\adfTEADI}$,

$(37,26,70,62,52,42,12,72,85,19,11,16,87)_{\adfTEADI}$,

$(10,32,73,86,49,62,26,31,29,9,67,57,19)_{\adfTEADI}$,

$(44,24,62,91,73,82,41,81,6,66,51,34,47)_{\adfTEADI}$,

\adfLgap 
$(0,33,7,57,44,52,19,42,12,83,59,53,10)_{\adfTEAEH}$,

$(31,11,4,10,39,85,88,78,21,2,86,8,80)_{\adfTEAEH}$,

$(56,2,73,29,43,80,59,18,31,60,22,50,74)_{\adfTEAEH}$,

$(19,7,58,32,9,63,15,67,78,90,53,22,66)_{\adfTEAEH}$,

$(9,77,31,78,61,55,25,57,52,50,28,40,13)_{\adfTEAEH}$,

$(89,59,7,24,73,2,1,37,44,60,19,68,42)_{\adfTEAEH}$,

$(53,0,42,11,14,70,65,55,29,63,16,10,89)_{\adfTEAEH}$,

$(48,79,0,40,61,21,23,78,83,39,53,27,6)_{\adfTEAEH}$,

$(68,55,5,43,45,53,33,88,11,2,51,59,87)_{\adfTEAEH}$,

$(1,78,16,44,86,53,0,56,5,55,70,10,49)_{\adfTEAEH}$,

$(37,34,17,64,30,25,28,10,62,55,39,14,61)_{\adfTEAEH}$,

$(68,57,70,77,78,32,35,44,39,0,18,76,87)_{\adfTEAEH}$,

$(77,16,26,10,6,48,59,60,36,40,50,51,78)_{\adfTEAEH}$,

\adfLgap 
$(0,16,46,75,50,45,42,1,22,8,36,14,89)_{\adfTEAFG}$,

$(78,24,50,46,88,15,83,1,66,32,34,12,41)_{\adfTEAFG}$,

$(62,2,71,22,32,26,76,55,33,43,46,36,19)_{\adfTEAFG}$,

$(84,47,76,51,57,50,20,83,63,28,40,45,2)_{\adfTEAFG}$,

$(33,7,75,15,27,62,56,70,86,29,63,30,83)_{\adfTEAFG}$,

$(20,33,64,87,22,48,45,89,30,43,84,71,79)_{\adfTEAFG}$,

$(84,66,41,16,52,45,4,87,83,14,75,35,89)_{\adfTEAFG}$,

$(49,31,30,18,5,14,45,28,8,41,52,43,36)_{\adfTEAFG}$,

$(35,14,50,90,3,84,55,73,87,57,46,38,77)_{\adfTEAFG}$,

$(71,21,27,25,45,15,60,34,90,22,24,63,73)_{\adfTEAFG}$,

$(9,85,33,64,29,76,8,35,17,77,41,32,18)_{\adfTEAFG}$,

$(56,29,4,36,32,70,25,71,65,1,2,75,77)_{\adfTEAFG}$,

$(67,3,9,1,42,43,64,11,22,2,46,80,25)_{\adfTEAFG}$,

\adfLgap 
$(0,49,66,12,72,81,73,3,43,15,75,29,65)_{\adfTEBCI}$,

$(67,22,50,57,63,8,55,64,40,37,79,81,69)_{\adfTEBCI}$,

$(72,28,58,62,24,81,22,16,71,51,87,5,30)_{\adfTEBCI}$,

$(59,52,30,72,1,74,6,18,10,89,66,29,47)_{\adfTEBCI}$,

$(45,78,19,38,42,64,84,29,26,0,35,41,36)_{\adfTEBCI}$,

$(19,49,81,35,58,67,4,52,13,33,14,30,76)_{\adfTEBCI}$,

$(54,70,85,88,67,56,49,24,51,28,47,86,21)_{\adfTEBCI}$,

$(71,65,17,66,14,88,36,72,54,1,29,0,7)_{\adfTEBCI}$,

$(65,9,61,70,83,79,21,22,23,35,57,15,56)_{\adfTEBCI}$,

$(90,63,71,58,67,55,79,14,49,84,10,52,24)_{\adfTEBCI}$,

$(40,86,43,25,24,54,7,73,60,75,82,48,64)_{\adfTEBCI}$,

$(63,22,42,25,1,61,47,10,73,24,85,55,66)_{\adfTEBCI}$,

$(65,7,32,48,56,54,23,24,34,40,2,27,88)_{\adfTEBCI}$,

\adfLgap 
$(0,78,73,15,90,66,22,46,69,87,39,82,65)_{\adfTEBEG}$,

$(68,40,81,55,21,73,46,4,79,48,51,31,30)_{\adfTEBEG}$,

$(22,55,60,30,61,6,83,54,24,63,38,17,85)_{\adfTEBEG}$,

$(80,86,48,25,19,85,1,74,65,15,84,14,83)_{\adfTEBEG}$,

$(45,32,59,41,88,31,7,33,44,77,14,19,12)_{\adfTEBEG}$,

$(83,61,84,22,56,52,13,43,51,40,25,41,44)_{\adfTEBEG}$,

$(0,46,57,10,74,2,60,83,3,73,52,88,58)_{\adfTEBEG}$,

$(50,51,87,83,7,34,38,59,5,30,81,47,1)_{\adfTEBEG}$,

$(86,5,79,29,14,12,36,68,23,60,18,21,11)_{\adfTEBEG}$,

$(89,68,46,42,41,71,20,75,43,69,13,11,84)_{\adfTEBEG}$,

$(79,49,69,8,67,29,56,90,24,53,34,36,44)_{\adfTEBEG}$,

$(22,27,74,61,12,52,66,72,60,86,51,30,46)_{\adfTEBEG}$,

$(11,23,74,7,56,81,21,40,41,13,4,78,45)_{\adfTEBEG}$,

\adfLgap 
$(0,44,48,14,47,57,10,30,24,39,21,6,64)_{\adfTECCH}$,

$(3,82,64,9,19,73,1,72,31,26,21,27,0)_{\adfTECCH}$,

$(4,6,72,57,27,63,67,42,43,15,24,90,75)_{\adfTECCH}$,

$(38,6,44,86,76,77,25,21,67,47,55,10,65)_{\adfTECCH}$,

$(62,28,59,56,32,71,79,25,20,77,33,76,17)_{\adfTECCH}$,

$(82,42,6,14,81,79,23,11,41,32,7,55,45)_{\adfTECCH}$,

$(73,11,39,7,31,37,6,84,70,88,57,1,48)_{\adfTECCH}$,

$(9,8,28,46,51,58,36,19,38,86,21,53,1)_{\adfTECCH}$,

$(85,6,71,58,16,76,9,70,68,66,62,36,53)_{\adfTECCH}$,

$(23,30,75,66,12,19,56,40,18,71,42,77,54)_{\adfTECCH}$,

$(58,54,9,1,90,61,12,16,28,23,41,30,3)_{\adfTECCH}$,

$(40,52,32,51,65,87,27,58,15,77,80,29,0)_{\adfTECCH}$,

$(1,5,13,33,15,39,69,47,21,38,59,80,44)_{\adfTECCH}$,

\adfLgap 
$(0,21,23,8,42,9,62,68,51,58,66,69,38)_{\adfTECDG}$,

$(73,57,24,34,55,31,59,33,40,6,62,64,85)_{\adfTECDG}$,

$(85,65,60,64,38,3,68,83,50,67,22,1,21)_{\adfTECDG}$,

$(5,56,81,38,30,19,11,62,44,63,33,9,43)_{\adfTECDG}$,

$(5,82,87,34,31,45,59,6,27,33,62,63,22)_{\adfTECDG}$,

$(42,52,82,86,40,24,90,67,55,48,28,14,8)_{\adfTECDG}$,

$(67,22,19,34,17,75,3,21,10,25,37,76,90)_{\adfTECDG}$,

$(27,8,76,79,52,51,81,82,56,3,59,85,77)_{\adfTECDG}$,

$(27,88,36,48,14,11,61,56,41,46,26,31,13)_{\adfTECDG}$,

$(58,48,27,70,86,3,79,31,68,26,64,75,16)_{\adfTECDG}$,

$(43,52,8,44,47,85,80,33,65,46,68,13,17)_{\adfTECDG}$,

$(21,9,84,74,57,16,46,14,60,30,36,87,47)_{\adfTECDG}$,

$(90,7,14,1,27,32,77,81,48,17,8,89,67)_{\adfTECDG}$,

\adfLgap 
$(0,66,12,68,89,24,33,22,81,90,43,67,36)_{\adfTECEF}$,

$(51,6,45,35,71,24,73,85,53,75,43,16,1)_{\adfTECEF}$,

$(41,35,33,1,64,40,21,88,48,68,87,4,74)_{\adfTECEF}$,

$(88,63,34,5,6,25,76,48,39,69,52,22,40)_{\adfTECEF}$,

$(15,64,44,48,4,70,39,82,66,29,28,72,77)_{\adfTECEF}$,

$(3,49,90,82,20,80,23,42,30,34,10,77,50)_{\adfTECEF}$,

$(72,45,38,76,79,39,75,31,37,85,42,40,7)_{\adfTECEF}$,

$(51,47,63,38,37,22,86,44,18,3,78,55,29)_{\adfTECEF}$,

$(21,54,90,4,5,43,80,34,25,45,47,69,51)_{\adfTECEF}$,

$(89,61,38,26,7,3,2,33,27,6,46,20,19)_{\adfTECEF}$,

$(59,21,51,77,66,82,85,61,22,5,76,55,80)_{\adfTECEF}$,

$(0,14,84,74,79,15,2,72,40,85,39,44,50)_{\adfTECEF}$,

$(59,42,70,56,2,73,44,81,8,61,55,39,89)_{\adfTECEF}$,

\adfLgap 
$(0,2,4,88,6,82,38,8,83,26,65,77,30)_{\adfTEDEE}$,

$(70,66,45,21,16,38,71,46,26,72,42,28,74)_{\adfTEDEE}$,

$(49,15,11,73,42,36,83,24,1,82,4,77,13)_{\adfTEDEE}$,

$(82,48,8,20,19,41,9,57,37,33,14,68,17)_{\adfTEDEE}$,

$(89,20,56,40,75,18,9,0,54,39,8,65,64)_{\adfTEDEE}$,

$(6,43,3,37,22,69,79,23,15,5,68,12,19)_{\adfTEDEE}$,

$(86,10,7,48,87,35,32,53,17,69,4,26,8)_{\adfTEDEE}$,

$(84,4,41,48,75,35,66,23,32,79,62,78,10)_{\adfTEDEE}$,

$(20,60,72,17,83,5,31,78,84,45,0,53,26)_{\adfTEDEE}$,

$(70,1,82,83,9,0,79,63,6,75,85,55,23)_{\adfTEDEE}$,

$(58,76,35,49,87,19,41,1,51,16,88,15,73)_{\adfTEDEE}$,

$(56,15,88,13,89,71,31,79,6,30,33,29,35)_{\adfTEDEE}$,

$(19,3,7,9,66,23,78,65,10,25,2,38,49)_{\adfTEDEE}$

\adfLgap
\noindent under the action of the mapping $x \mapsto x + 4$ (mod 92). \eproof

\ADFvfyParStart{92, \{\{92,23,4\}\}, -1, -1, -1} 

~\\
\noindent {\bf Proof of Lemma \ref{lem:theta14 multipartite}}~

\noindent {\boldmath $K_{14,7}$}~
Let the vertex set be $\{0, 1, \dots, 20\}$ partitioned into
$\{0, 1, \dots, 13\}$ and $\{14, 15, \dots, 20\}$. The decompositions consist of

\adfLgap 
$(0,1,14,16,15,4,19,2,20,11,18,13,17)_{\adfTEBBJ}$,

\adfLgap 
$(0,9,15,16,8,17,14,6,20,13,18,5,19)_{\adfTEBDH}$,

\adfLgap 
$(0,9,18,15,4,16,13,19,17,10,14,5,20)_{\adfTEBFF}$,

\adfLgap 
$(0,11,17,4,18,19,10,16,20,7,15,1,14)_{\adfTEDDF}$

\adfLgap
\noindent under the action of the mapping $x \mapsto x + 2 \adfmod{14}$ for $x < 14$,
$x \mapsto (x + 1 \adfmod{7}) + 14$ for $x \ge 14$.

\ADFvfyParStart{21, \{\{14,7,2\},\{7,7,1\}\}, -1, \{\{14,\{0\}\},\{7,\{1\}\}\}, -1} 

\noindent {\boldmath $K_{14,14,14}$}~
Let the vertex set be $Z_{42}$ partitioned according to residue classes modulo 3.
The decompositions consist of

\adfLgap 
$(0,17,16,10,18,11,30,34,3,1,2,4,9)_{\adfTEABK}$,

$(0,7,20,13,3,8,22,5,1,12,31,9,23)_{\adfTEABK}$,

\adfLgap 
$(0,20,17,39,37,27,13,8,1,2,3,5,9)_{\adfTEACJ}$,

$(0,2,4,12,7,15,31,6,19,8,22,9,28)_{\adfTEACJ}$,

\adfLgap 
$(0,14,10,27,25,34,18,7,2,1,5,4,33)_{\adfTEADI}$,

$(0,2,4,11,25,13,3,8,1,9,31,5,22)_{\adfTEADI}$,

\adfLgap 
$(0,29,10,9,23,21,1,11,27,5,40,26,4)_{\adfTEAEH}$,

$(21,32,10,2,19,6,16,11,15,22,3,26,28)_{\adfTEAEH}$,

\adfLgap 
$(0,19,16,29,34,32,39,35,12,37,6,17,18)_{\adfTEAFG}$,

$(13,8,9,17,7,20,0,15,29,3,28,14,4)_{\adfTEAFG}$,

\adfLgap 
$(0,23,28,17,31,5,34,27,26,6,38,30,13)_{\adfTEBCI}$,

$(28,13,17,30,14,15,4,0,7,5,25,6,29)_{\adfTEBCI}$,

\adfLgap 
$(0,12,40,31,14,27,8,34,18,19,32,15,35)_{\adfTEBEG}$,

$(29,33,25,13,18,38,1,15,8,7,5,12,22)_{\adfTEBEG}$,

\adfLgap 
$(0,23,25,27,40,33,1,12,31,15,34,29,9)_{\adfTECCH}$,

$(28,39,20,31,32,22,0,7,2,15,16,36,10)_{\adfTECCH}$,

\adfLgap 
$(0,3,1,20,28,12,23,29,30,5,10,32,37)_{\adfTECDG}$,

$(39,34,29,36,23,25,11,1,12,2,6,13,0)_{\adfTECDG}$,

\adfLgap 
$(0,2,37,36,23,30,8,40,13,11,12,28,33)_{\adfTECEF}$,

$(31,13,5,27,14,0,7,3,2,4,15,38,21)_{\adfTECEF}$,

\adfLgap 
$(0,23,8,25,24,11,22,29,39,37,32,12,10)_{\adfTEDEE}$,

$(40,39,27,35,16,8,4,18,1,14,7,5,19)_{\adfTEDEE}$

\adfLgap
\noindent under the action of the mapping $x \mapsto x + 2$ (mod 42).

\ADFvfyParStart{42, \{\{42,21,2\}\}, -1, \{\{14,\{0,1,2\}\}\}, -1} 

\noindent {\boldmath $K_{7,7,7,7}$}~
Let the vertex set be $\{0, 1, \dots, 27\}$ partitioned into
$\{3j + i: j = 0, 1, \dots, 6\}$, $i = 0, 1, 2$, and $\{21, 22, \dots, 27\}$.
The decompositions consist of

\adfLgap 
$(0,23,13,7,17,1,2,21,3,5,9,22,8)_{\adfTEABK}$,

\adfLgap 
$(0,25,1,11,2,21,12,8,13,6,14,22,5)_{\adfTEACJ}$,

\adfLgap 
$(0,24,14,18,19,13,11,26,1,6,16,22,8)_{\adfTEADI}$,

\adfLgap 
$(0,23,1,21,6,16,17,3,8,25,7,5,18)_{\adfTEAEH}$,

\adfLgap 
$(0,22,26,17,1,18,4,21,8,6,7,15,5)_{\adfTEAFG}$,

\adfLgap 
$(0,21,13,5,9,19,20,26,2,12,23,7,14)_{\adfTEBCI}$,

\adfLgap 
$(0,22,17,11,10,12,4,5,19,21,6,23,2)_{\adfTEBEG}$,

\adfLgap 
$(0,21,26,20,11,18,19,3,2,6,24,1,14)_{\adfTECCH}$,

\adfLgap 
$(0,21,10,12,24,3,17,16,8,4,5,25,6)_{\adfTECDG}$,

\adfLgap 
$(0,23,24,2,10,11,15,7,14,9,22,3,5)_{\adfTECEF}$,

\adfLgap 
$(0,25,4,15,1,23,8,21,7,2,3,11,6)_{\adfTEDEE}$

\adfLgap
\noindent under the action of the mapping $x \mapsto x + 1$ (mod 21) for $x < 21$,
$x \mapsto (x + 1 \adfmod{7}) + 21$ for $x \ge 21$.

\ADFvfyParStart{28, \{\{21,21,1\},\{7,7,1\}\}, -1, \{\{7,\{0,1,2\}\},\{7,\{3\}\}\}, -1} 

\noindent {\boldmath $K_{28,28,28,35}$}~
Let the vertex set be $\{0, 1, \dots, 118\}$ partitioned into
$\{3j + i: j = 0, 1, \dots, 27\}$, $i = 0, 1, 2$, and $\{84, 85, \dots, 118\}$.
The decompositions consist of

\adfLgap 
$(0,87,65,2,113,20,36,84,74,101,12,52,11)_{\adfTEABK}$,

$(45,26,100,20,25,83,61,54,87,44,81,116,19)_{\adfTEABK}$,

$(35,6,93,18,106,30,90,66,40,17,84,19,62)_{\adfTEABK}$,

$(79,29,116,35,40,41,39,31,3,99,49,69,98)_{\adfTEABK}$,

$(75,111,71,117,25,15,50,1,108,72,99,29,67)_{\adfTEABK}$,

$(24,2,94,32,0,10,81,65,95,11,90,20,6)_{\adfTEABK}$,

$(69,80,16,117,44,93,4,15,90,33,103,21,92)_{\adfTEABK}$,

$(112,55,72,20,67,108,76,75,43,85,68,34,99)_{\adfTEABK}$,

$(116,68,30,28,41,10,69,46,17,3,94,6,110)_{\adfTEABK}$,

\adfLgap 
$(0,25,86,81,73,65,60,111,12,38,109,27,50)_{\adfTEACJ}$,

$(73,90,75,5,88,49,59,55,6,29,82,80,9)_{\adfTEACJ}$,

$(86,76,43,32,73,78,11,88,10,63,115,52,94)_{\adfTEACJ}$,

$(0,74,89,45,92,6,71,72,79,95,28,69,88)_{\adfTEACJ}$,

$(13,71,53,113,68,40,24,84,55,5,70,23,102)_{\adfTEACJ}$,

$(77,110,102,44,0,4,12,13,85,18,117,53,66)_{\adfTEACJ}$,

$(111,73,33,104,37,117,76,95,82,32,67,94,51)_{\adfTEACJ}$,

$(27,7,79,92,11,101,24,10,108,78,102,62,45)_{\adfTEACJ}$,

$(15,52,103,74,56,4,106,6,26,64,88,20,89)_{\adfTEACJ}$,

\adfLgap 
$(0,38,19,78,95,77,109,80,63,102,25,36,70)_{\adfTEADI}$,

$(10,66,97,48,62,104,39,11,9,113,8,30,85)_{\adfTEADI}$,

$(17,96,27,79,80,9,35,70,112,60,73,56,30)_{\adfTEADI}$,

$(68,76,63,111,15,25,109,6,105,31,95,57,101)_{\adfTEADI}$,

$(0,89,88,76,77,65,36,83,97,16,23,46,57)_{\adfTEADI}$,

$(116,69,42,104,1,8,10,50,70,23,91,61,102)_{\adfTEADI}$,

$(60,84,110,11,36,107,22,57,35,81,76,93,79)_{\adfTEADI}$,

$(19,39,113,60,70,15,105,8,86,32,3,73,88)_{\adfTEADI}$,

$(30,73,46,77,93,88,66,53,13,63,95,9,102)_{\adfTEADI}$,

\adfLgap 
$(0,19,91,23,93,62,70,96,7,51,40,78,112)_{\adfTEAEH}$,

$(68,66,96,50,75,14,40,53,51,55,85,76,90)_{\adfTEAEH}$,

$(114,72,36,1,15,5,75,70,98,54,84,17,117)_{\adfTEAEH}$,

$(16,38,90,59,89,64,11,87,6,103,68,55,27)_{\adfTEAEH}$,

$(75,16,110,82,113,23,102,8,58,3,87,1,111)_{\adfTEAEH}$,

$(1,8,105,27,11,61,23,86,53,79,15,91,45)_{\adfTEAEH}$,

$(31,39,103,46,45,93,100,44,91,32,75,28,8)_{\adfTEAEH}$,

$(47,1,88,31,80,115,104,69,14,30,40,0,99)_{\adfTEAEH}$,

$(38,104,61,9,26,7,34,35,87,0,76,117,8)_{\adfTEAEH}$,

\adfLgap 
$(0,80,62,42,100,79,12,40,71,96,41,85,66)_{\adfTEAFG}$,

$(69,68,90,30,107,15,43,11,76,89,73,77,106)_{\adfTEAFG}$,

$(81,32,56,103,62,6,87,4,109,76,90,78,89)_{\adfTEAFG}$,

$(63,2,94,74,116,61,86,43,59,13,14,67,87)_{\adfTEAFG}$,

$(6,64,77,12,91,24,103,113,15,4,59,92,30)_{\adfTEAFG}$,

$(58,66,96,5,113,17,115,9,70,59,15,85,79)_{\adfTEAFG}$,

$(14,4,9,76,107,17,103,66,68,92,5,109,50)_{\adfTEAFG}$,

$(64,69,27,20,49,41,7,93,43,112,79,38,111)_{\adfTEAFG}$,

$(18,55,61,71,73,59,99,90,9,108,11,63,104)_{\adfTEAFG}$,

\adfLgap 
$(0,94,56,55,27,5,105,14,1,50,28,111,35)_{\adfTEBCI}$,

$(111,114,66,51,37,72,1,48,73,6,79,32,70)_{\adfTEBCI}$,

$(106,111,20,52,11,69,92,50,58,62,102,4,27)_{\adfTEBCI}$,

$(74,76,110,31,103,39,55,15,53,19,88,32,42)_{\adfTEBCI}$,

$(83,116,19,61,62,1,24,64,9,13,78,76,77)_{\adfTEBCI}$,

$(36,102,56,97,47,10,27,19,100,78,46,39,44)_{\adfTEBCI}$,

$(36,49,89,105,17,94,32,113,11,95,75,112,68)_{\adfTEBCI}$,

$(80,81,108,10,98,49,23,109,74,6,13,99,7)_{\adfTEBCI}$,

$(117,54,29,33,85,59,108,4,35,24,118,11,105)_{\adfTEBCI}$,

\adfLgap 
$(0,39,7,74,49,93,23,25,95,29,34,108,28)_{\adfTEBEG}$,

$(9,105,55,92,83,109,76,91,77,19,96,8,7)_{\adfTEBEG}$,

$(34,33,98,107,4,65,52,32,25,91,55,86,31)_{\adfTEBEG}$,

$(82,76,27,117,48,97,33,14,6,73,63,77,54)_{\adfTEBEG}$,

$(100,36,7,68,54,107,70,37,15,55,51,85,10)_{\adfTEBEG}$,

$(47,22,108,30,94,82,42,55,21,89,69,68,111)_{\adfTEBEG}$,

$(55,70,35,94,60,84,57,83,46,0,47,42,38)_{\adfTEBEG}$,

$(64,60,116,23,87,5,114,11,117,57,103,24,111)_{\adfTEBEG}$,

$(6,32,90,79,109,68,95,34,3,98,0,103,51)_{\adfTEBEG}$,

\adfLgap 
$(0,108,16,21,62,58,91,12,11,84,10,88,6)_{\adfTECCH}$,

$(42,25,102,69,109,74,43,41,104,70,17,39,71)_{\adfTECCH}$,

$(35,21,103,13,86,79,94,30,20,100,14,6,52)_{\adfTECCH}$,

$(62,56,60,107,49,33,73,18,116,53,43,80,63)_{\adfTECCH}$,

$(60,96,10,75,5,12,32,69,50,76,105,54,77)_{\adfTECCH}$,

$(93,32,7,91,39,102,60,17,28,8,67,72,87)_{\adfTECCH}$,

$(86,60,42,105,39,89,46,98,1,95,45,25,8)_{\adfTECCH}$,

$(44,13,9,100,89,17,58,83,49,84,41,87,0)_{\adfTECCH}$,

$(55,56,69,115,27,113,92,5,85,9,107,59,16)_{\adfTECCH}$,

\adfLgap 
$(0,72,84,44,95,19,65,22,63,61,27,58,32)_{\adfTECDG}$,

$(24,70,19,109,117,57,32,84,49,69,50,30,86)_{\adfTECDG}$,

$(89,79,28,80,73,83,109,35,4,75,71,88,56)_{\adfTECDG}$,

$(30,94,59,43,5,115,40,31,53,27,116,72,38)_{\adfTECDG}$,

$(64,17,68,45,110,2,112,111,9,113,43,51,97)_{\adfTECDG}$,

$(50,15,31,103,113,42,26,33,80,100,37,96,10)_{\adfTECDG}$,

$(72,30,79,65,70,33,20,58,50,117,62,27,85)_{\adfTECDG}$,

$(52,29,23,12,102,76,97,88,16,113,3,106,81)_{\adfTECDG}$,

$(48,43,7,110,115,49,102,111,4,113,13,53,96)_{\adfTECDG}$,

\adfLgap 
$(0,34,26,33,87,43,110,24,95,83,1,114,62)_{\adfTECEF}$,

$(78,71,64,88,5,93,44,48,58,3,22,98,34)_{\adfTECEF}$,

$(45,105,53,25,1,68,104,47,76,81,8,49,66)_{\adfTECEF}$,

$(60,91,111,62,55,85,78,13,68,89,67,102,19)_{\adfTECEF}$,

$(78,36,44,112,65,117,20,106,80,98,51,109,1)_{\adfTECEF}$,

$(90,111,37,57,8,63,5,75,70,114,61,48,25)_{\adfTECEF}$,

$(54,69,76,17,13,112,14,93,29,39,46,92,65)_{\adfTECEF}$,

$(13,70,88,38,51,50,95,33,75,101,36,20,104)_{\adfTECEF}$,

$(38,59,78,25,96,67,109,24,0,97,13,29,108)_{\adfTECEF}$,

\adfLgap 
$(0,38,61,86,58,76,47,89,42,84,36,79,100)_{\adfTEDEE}$,

$(62,18,100,54,101,79,115,13,107,76,42,29,43)_{\adfTEDEE}$,

$(31,30,51,32,84,78,103,66,65,106,46,69,91)_{\adfTEDEE}$,

$(109,58,76,5,103,31,102,33,97,61,23,16,27)_{\adfTEDEE}$,

$(99,32,5,22,97,77,88,57,76,31,27,104,70)_{\adfTEDEE}$,

$(25,19,96,12,44,88,40,103,81,9,17,108,53)_{\adfTEDEE}$,

$(12,54,68,13,112,70,111,27,95,53,48,105,32)_{\adfTEDEE}$,

$(46,38,77,67,92,98,82,102,27,47,0,49,54)_{\adfTEDEE}$,

$(39,71,116,83,27,115,75,17,19,46,44,54,110)_{\adfTEDEE}$

\adfLgap
\noindent under the action of the mapping $x \mapsto x + 2$ (mod 84) for $x < 84$,
$x \mapsto (x - 84 + 5 \adfmod{35}) + 84$ for $x \ge 84$.

\ADFvfyParStart{119, \{\{84,42,2\},\{35,7,5\}\}, -1, \{\{28,\{0,1,2\}\},\{35,\{3\}\}\}, -1} 

\noindent {\boldmath $K_{7,7,7,7,7}$}~
Let the vertex set be $Z_{35}$ partitioned according to residue class modulo 5.
The decompositions consist of

\adfLgap 
$(0,9,23,3,16,24,18,20,1,8,4,21,10)_{\adfTEABK}$,

\adfLgap 
$(0,32,1,3,22,29,6,33,2,16,5,24,15)_{\adfTEACJ}$,

\adfLgap 
$(0,9,18,12,20,13,14,28,5,1,3,19,16)_{\adfTEADI}$,

\adfLgap 
$(0,22,14,18,27,11,33,6,3,4,10,17,5)_{\adfTEAEH}$,

\adfLgap 
$(0,23,13,15,31,3,14,1,7,4,18,10,6)_{\adfTEAFG}$,

\adfLgap 
$(0,2,29,4,21,28,15,17,3,6,5,14,25)_{\adfTEBCI}$,

\adfLgap 
$(0,26,18,7,4,33,14,13,2,1,5,3,12)_{\adfTEBEG}$,

\adfLgap 
$(0,14,9,12,16,3,27,10,4,5,17,21,7)_{\adfTECCH}$,

\adfLgap 
$(0,12,27,9,23,7,16,14,1,2,13,11,5)_{\adfTECDG}$,

\adfLgap 
$(0,5,6,22,26,3,11,18,1,4,28,7,9)_{\adfTECEF}$,

\adfLgap 
$(0,16,12,14,17,7,29,2,20,6,22,1,25)_{\adfTEDEE}$

\adfLgap
\noindent under the action of the mapping $x \mapsto x + 1 \adfmod{35}$.

\ADFvfyParStart{35, \{\{35,35,1\}\}, -1, \{\{7,\{0,1,2,3,4\}\}\}, -1} 

\noindent {\boldmath $K_{7,7,7,7,21}$}~
Let the vertex set be $\{0, 1, \dots, 48\}$ partitioned into
$\{3j + i: j = 0, 1, \dots, 6\}$, $i = 0, 1, 2$, $\{21, 22, \dots, 27\}$ and $\{28, 29, \dots, 48\}$.
The decompositions consist of the graphs

\adfLgap 
$(0,4,5,35,14,30,15,21,39,22,45,1,25)_{\adfTEABK}$,

$(42,17,10,16,8,19,32,3,44,22,36,9,25)_{\adfTEABK}$,

$(18,30,20,40,21,16,27,33,9,28,10,38,0)_{\adfTEABK}$,

\adfLgap 
$(0,32,25,20,34,7,36,23,14,22,40,17,12)_{\adfTEACJ}$,

$(18,11,21,15,5,45,2,30,14,44,23,32,1)_{\adfTEACJ}$,

$(32,15,17,41,11,25,42,3,1,0,46,24,29)_{\adfTEACJ}$,

\adfLgap 
$(0,4,10,45,3,21,18,43,13,5,23,37,6)_{\adfTEADI}$,

$(21,43,32,20,15,44,4,30,19,34,0,45,2)_{\adfTEADI}$,

$(45,26,21,2,39,1,24,4,31,13,6,35,18)_{\adfTEADI}$,

\adfLgap 
$(0,19,26,5,13,34,31,17,25,46,8,3,35)_{\adfTEAEH}$,

$(27,33,10,36,24,8,31,1,26,35,7,30,11)_{\adfTEAEH}$,

$(30,22,18,7,0,39,6,10,11,31,2,29,16)_{\adfTEAEH}$,

\adfLgap 
$(0,34,20,1,26,28,3,43,6,19,30,9,24)_{\adfTEAFG}$,

$(13,32,22,35,27,45,6,9,47,26,19,42,18)_{\adfTEAFG}$,

$(8,13,38,2,31,19,25,1,12,22,34,6,33)_{\adfTEAFG}$,

\adfLgap 
$(0,30,7,47,16,21,46,1,40,26,13,25,8)_{\adfTEBCI}$,

$(6,17,34,46,9,47,24,32,5,3,7,37,16)_{\adfTEBCI}$,

$(31,17,26,19,32,21,10,0,16,41,9,38,25)_{\adfTEBCI}$,

\adfLgap 
$(0,18,39,41,8,23,43,5,40,26,45,22,20)_{\adfTEBEG}$,

$(3,38,14,21,32,6,10,27,35,1,45,15,16)_{\adfTEBEG}$,

$(41,0,14,10,46,8,40,4,17,22,1,30,27)_{\adfTEBEG}$,

\adfLgap 
$(0,6,24,4,36,14,42,19,45,26,1,18,30)_{\adfTECCH}$,

$(0,20,26,29,22,37,5,40,15,31,13,12,39)_{\adfTECCH}$,

$(1,3,35,24,11,31,15,44,12,43,22,28,26)_{\adfTECCH}$,

\adfLgap 
$(0,41,20,6,39,5,1,32,4,47,11,42,17)_{\adfTECDG}$,

$(40,11,13,23,22,44,25,27,36,3,33,12,28)_{\adfTECDG}$,

$(43,31,2,24,26,6,8,14,1,12,7,23,5)_{\adfTECDG}$,

\adfLgap 
$(0,1,5,26,40,22,29,2,24,37,7,43,15)_{\adfTECEF}$,

$(1,42,5,25,14,34,11,27,44,9,38,20,18)_{\adfTECEF}$,

$(10,34,22,0,36,27,46,9,31,14,4,35,2)_{\adfTECEF}$,

\adfLgap 
$(0,43,28,16,2,42,4,38,19,39,3,20,22)_{\adfTEDEE}$,

$(20,12,9,44,26,7,24,19,17,45,21,29,23)_{\adfTEDEE}$,

$(41,17,10,47,18,15,45,1,23,9,36,27,39)_{\adfTEDEE}$

\adfLgap
\noindent under the action of the mapping $x \mapsto x + 1 \adfmod{21}$ for $x < 21$,
$x \mapsto (x + 1 \adfmod{7}) + 21$ for $21 \le x < 28$,
$x \mapsto (x - 28 + 1 \adfmod{21}) + 28$ for $x \ge 28$.

\ADFvfyParStart{49, \{\{21,21,1\},\{7,7,1\},\{21,21,1\}\}, -1, \{\{7,\{0,1,2\}\},\{7,\{3\}\},\{21,\{4\}\}\}, -1} 

\noindent {\boldmath $K_{7,7,7,7,28}$}~
Let the vertex set be $\{0, 1, \dots, 55\}$ partitioned into
$\{3j + i: j = 0, 1, \dots, 6\}$, $i = 0, 1, 2, 3$ and $\{28, 29, \dots, 55\}$.
The decompositions consist of the graphs

\adfLgap 
$(0,18,35,6,43,4,28,5,40,27,34,16,42)_{\adfTEABK}$,

$(21,47,23,16,54,25,12,5,27,29,0,42,6)_{\adfTEABK}$,

$(16,44,18,27,49,22,51,19,9,14,1,7,10)_{\adfTEABK}$,

$(38,17,6,21,30,3,51,23,42,27,14,16,49)_{\adfTEABK}$,

$(4,11,47,35,8,3,42,9,36,25,18,17,51)_{\adfTEABK}$,

$(37,2,11,16,32,27,33,18,36,17,41,24,46)_{\adfTEABK}$,

$(21,34,11,29,17,37,27,48,26,38,8,28,3)_{\adfTEABK}$,

$(1,47,0,12,26,29,13,43,27,44,8,45,22)_{\adfTEABK}$,

$(28,26,16,18,7,4,45,20,54,1,51,0,27)_{\adfTEABK}$,

$(13,10,16,44,15,52,9,7,49,2,50,22,43)_{\adfTEABK}$,

$(49,10,3,21,7,16,48,6,11,24,5,43,1)_{\adfTEABK}$,

\adfLgap 
$(0,44,42,25,2,8,45,22,16,23,28,10,3)_{\adfTEACJ}$,

$(13,49,10,19,39,21,28,0,32,12,47,18,20)_{\adfTEACJ}$,

$(52,23,3,26,1,8,39,22,34,4,37,2,9)_{\adfTEACJ}$,

$(15,40,46,14,24,6,44,1,10,49,8,23,0)_{\adfTEACJ}$,

$(35,2,25,39,22,28,1,42,8,34,19,20,30)_{\adfTEACJ}$,

$(1,41,18,19,51,27,17,12,34,26,45,13,3)_{\adfTEACJ}$,

$(47,19,13,42,11,26,29,3,1,12,9,0,53)_{\adfTEACJ}$,

$(30,5,9,29,19,33,12,23,32,10,37,22,11)_{\adfTEACJ}$,

$(8,26,53,25,9,24,43,27,1,32,18,5,42)_{\adfTEACJ}$,

$(34,10,15,35,1,31,19,37,21,32,24,27,5)_{\adfTEACJ}$,

$(28,4,26,47,11,38,3,35,8,54,6,20,43)_{\adfTEACJ}$,

\adfLgap 
$(0,32,28,26,25,40,22,39,13,46,23,29,16)_{\adfTEADI}$,

$(50,6,22,30,3,1,43,12,37,5,48,21,8)_{\adfTEADI}$,

$(0,23,10,32,22,42,18,33,16,36,3,51,17)_{\adfTEADI}$,

$(19,0,43,3,52,34,2,49,26,40,11,12,6)_{\adfTEADI}$,

$(27,24,44,8,26,16,42,4,11,36,10,23,51)_{\adfTEADI}$,

$(0,47,21,38,8,33,2,16,25,36,27,49,3)_{\adfTEADI}$,

$(31,10,22,49,13,21,39,6,1,52,18,30,17)_{\adfTEADI}$,

$(28,9,15,54,4,25,0,43,8,42,7,38,18)_{\adfTEADI}$,

$(50,21,25,27,10,3,33,13,43,18,23,14,7)_{\adfTEADI}$,

$(37,25,11,47,24,13,41,20,38,1,12,49,7)_{\adfTEADI}$,

$(26,33,37,8,23,13,3,5,55,14,43,11,17)_{\adfTEADI}$,

\adfLgap 
$(0,48,53,6,30,22,41,3,39,4,14,24,18)_{\adfTEAEH}$,

$(53,20,3,54,16,44,22,43,10,25,23,14,39)_{\adfTEAEH}$,

$(17,18,33,4,21,32,44,10,1,49,14,27,28)_{\adfTEAEH}$,

$(11,25,48,16,7,0,39,15,41,20,18,23,2)_{\adfTEAEH}$,

$(50,7,6,46,11,5,23,12,42,0,43,9,37)_{\adfTEAEH}$,

$(50,5,9,49,12,43,18,45,22,51,24,10,28)_{\adfTEAEH}$,

$(23,0,28,7,41,2,20,43,27,5,48,3,50)_{\adfTEAEH}$,

$(32,1,10,17,47,8,7,6,3,42,22,35,26)_{\adfTEAEH}$,

$(51,25,19,31,9,44,5,22,11,53,10,4,30)_{\adfTEAEH}$,

$(51,9,3,16,25,24,6,34,0,36,20,32,19)_{\adfTEAEH}$,

$(25,49,50,4,3,13,54,5,42,11,29,12,17)_{\adfTEAEH}$,

\adfLgap 
$(0,31,46,17,42,5,15,35,9,27,41,21,4)_{\adfTEAFG}$,

$(48,0,23,38,22,3,37,24,2,42,1,50,5)_{\adfTEAFG}$,

$(2,39,5,3,40,10,0,29,17,4,30,19,18)_{\adfTEAFG}$,

$(1,22,36,25,0,32,9,18,17,39,20,48,19)_{\adfTEAFG}$,

$(0,26,41,17,23,54,19,40,2,47,9,51,21)_{\adfTEAFG}$,

$(19,39,43,0,15,37,11,32,18,52,2,49,3)_{\adfTEAFG}$,

$(11,28,12,46,14,37,1,30,20,42,18,9,23)_{\adfTEAFG}$,

$(54,21,6,37,16,9,52,24,18,33,17,49,20)_{\adfTEAFG}$,

$(1,3,35,22,29,27,20,47,15,12,2,41,8)_{\adfTEAFG}$,

$(51,18,22,48,12,28,7,11,26,54,3,41,16)_{\adfTEAFG}$,

$(38,2,11,6,20,1,48,24,15,21,49,4,55)_{\adfTEAFG}$,

\adfLgap 
$(0,30,10,48,17,50,23,12,52,22,37,18,1)_{\adfTEBCI}$,

$(11,37,5,42,17,49,13,52,24,18,45,12,9)_{\adfTEBCI}$,

$(11,0,46,37,13,29,12,1,6,45,8,49,19)_{\adfTEBCI}$,

$(20,16,22,23,47,10,35,7,2,9,31,8,32)_{\adfTEBCI}$,

$(36,28,21,9,10,14,45,16,43,25,27,47,2)_{\adfTEBCI}$,

$(16,25,23,9,51,41,20,54,3,53,19,44,15)_{\adfTEBCI}$,

$(26,18,39,15,28,0,9,54,11,2,44,8,50)_{\adfTEBCI}$,

$(13,43,7,26,11,50,10,54,0,32,15,36,3)_{\adfTEBCI}$,

$(37,23,2,25,24,21,51,22,46,27,41,18,9)_{\adfTEBCI}$,

$(1,31,26,48,25,24,39,4,19,18,52,15,20)_{\adfTEBCI}$,

$(42,26,12,21,54,4,51,14,11,1,39,23,34)_{\adfTEBCI}$,

\adfLgap 
$(0,53,17,37,24,41,6,44,25,39,5,31,27)_{\adfTEBEG}$,

$(0,19,9,23,32,5,48,38,21,26,52,18,29)_{\adfTEBEG}$,

$(47,18,9,8,41,27,40,16,26,46,0,32,1)_{\adfTEBEG}$,

$(32,20,21,25,34,4,39,22,16,52,12,9,6)_{\adfTEBEG}$,

$(27,7,14,45,1,53,10,47,18,43,13,35,17)_{\adfTEBEG}$,

$(31,39,16,15,36,19,2,23,35,18,17,46,12)_{\adfTEBEG}$,

$(46,14,23,18,25,3,24,21,19,50,7,54,15)_{\adfTEBEG}$,

$(29,26,4,17,24,7,42,27,54,12,34,1,50)_{\adfTEBEG}$,

$(54,43,22,13,36,12,10,0,5,33,1,14,19)_{\adfTEBEG}$,

$(24,25,45,52,19,20,11,26,53,2,43,15,12)_{\adfTEBEG}$,

$(0,40,15,48,6,46,10,29,23,1,3,53,22)_{\adfTEBEG}$,

\adfLgap 
$(0,36,15,17,40,14,26,7,21,6,41,1,8)_{\adfTECCH}$,

$(31,35,12,2,9,15,19,32,18,33,5,46,1)_{\adfTECCH}$,

$(29,11,6,40,26,21,16,42,7,14,33,27,4)_{\adfTECCH}$,

$(48,1,21,30,10,47,4,50,3,2,20,31,18)_{\adfTECCH}$,

$(34,5,3,35,26,38,23,49,10,39,13,29,15)_{\adfTECCH}$,

$(37,25,17,8,20,49,19,6,0,5,28,4,18)_{\adfTECCH}$,

$(39,18,15,38,14,46,24,27,16,1,33,0,34)_{\adfTECCH}$,

$(29,23,27,32,8,44,2,51,6,25,22,52,4)_{\adfTECCH}$,

$(33,13,25,28,23,48,11,47,12,44,27,10,7)_{\adfTECCH}$,

$(38,21,16,46,24,23,17,20,29,4,42,15,31)_{\adfTECCH}$,

$(27,47,22,20,55,24,42,18,36,25,26,4,5)_{\adfTECCH}$,

\adfLgap 
$(0,23,50,2,3,18,47,28,21,34,7,32,17)_{\adfTECDG}$,

$(20,37,30,14,39,13,11,29,16,40,10,35,5)_{\adfTECDG}$,

$(42,10,25,48,27,32,19,13,24,2,21,52,15)_{\adfTECDG}$,

$(7,9,28,22,12,3,47,52,5,45,6,30,24)_{\adfTECDG}$,

$(48,37,12,23,9,14,16,21,24,10,0,51,3)_{\adfTECDG}$,

$(6,0,47,7,38,1,53,39,4,37,21,49,2)_{\adfTECDG}$,

$(0,10,21,7,30,25,45,40,24,44,26,48,16)_{\adfTECDG}$,

$(50,34,3,16,15,33,11,19,47,26,1,24,23)_{\adfTECDG}$,

$(12,8,41,26,39,25,50,51,19,2,47,13,34)_{\adfTECDG}$,

$(42,12,2,13,6,43,21,14,45,18,25,27,29)_{\adfTECDG}$,

$(32,27,18,17,3,41,5,6,7,43,12,55,9)_{\adfTECDG}$,

\adfLgap 
$(0,52,18,27,37,14,30,12,25,4,10,1,3)_{\adfTECEF}$,

$(18,4,3,36,13,30,16,33,37,1,47,6,5)_{\adfTECEF}$,

$(10,36,39,5,7,47,19,2,3,33,17,29,9)_{\adfTECEF}$,

$(9,19,39,14,41,8,21,30,34,10,36,17,45)_{\adfTECEF}$,

$(7,0,34,5,41,19,42,2,45,27,21,3,48)_{\adfTECEF}$,

$(29,41,26,27,5,31,4,14,8,43,12,38,16)_{\adfTECEF}$,

$(51,32,3,8,27,46,10,17,14,42,12,6,9)_{\adfTECEF}$,

$(5,54,46,11,18,39,24,23,43,21,15,17,6)_{\adfTECEF}$,

$(38,3,0,47,5,16,23,12,6,49,8,52,13)_{\adfTECEF}$,

$(13,43,48,20,4,2,44,7,27,24,11,40,10)_{\adfTECEF}$,

$(2,28,41,6,47,5,54,20,16,55,23,36,18)_{\adfTECEF}$,

\adfLgap 
$(0,54,32,9,6,36,21,10,5,19,34,25,16)_{\adfTEDEE}$,

$(9,12,51,26,49,41,24,36,19,2,23,1,31)_{\adfTEDEE}$,

$(25,33,41,26,20,47,10,43,8,44,9,40,6)_{\adfTEDEE}$,

$(14,0,46,21,29,13,23,31,27,53,6,12,39)_{\adfTEDEE}$,

$(17,28,16,30,12,26,36,23,14,45,19,44,7)_{\adfTEDEE}$,

$(53,33,22,27,21,7,35,9,11,5,44,26,0)_{\adfTEDEE}$,

$(26,18,8,25,54,47,7,30,3,16,38,5,19)_{\adfTEDEE}$,

$(29,36,8,42,7,15,17,30,14,19,35,1,16)_{\adfTEDEE}$,

$(7,0,24,47,23,17,20,48,15,41,11,31,2)_{\adfTEDEE}$,

$(6,38,46,15,12,41,17,44,14,32,8,22,19)_{\adfTEDEE}$,

$(0,10,5,12,55,30,2,43,25,31,21,38,27)_{\adfTEDEE}$

\adfLgap
\noindent under the action of the mapping $x \mapsto x + 4 \adfmod{28}$ for $x < 28$,
$x \mapsto (x + 4 \adfmod{28}) + 28$ for $x \ge 28$. \eproof

\ADFvfyParStart{56, \{\{28,7,4\},\{28,7,4\}\}, -1, \{\{7,\{0,1,2,3\}\},\{28,\{4\}\}\}, -1} 

\section{Theta graphs with 15 edges}
\label{sec:theta15}


\noindent {\bf Proof of Lemma \ref{lem:theta15 designs}}

\noindent {\boldmath $K_{15}$}~
Let the vertex set be $Z_{14} \cup \{\infty\}$. The decompositions consist of the graphs

\adfLgap
$(0,2,1,3,5,8,4,9,\infty,12,6,13,7,11)_{\adfTFABL}$,

\adfLgap
$(0,2,1,5,4,3,6,12,7,13,8,\infty,11,9)_{\adfTFACK}$,

\adfLgap
$(0,2,1,3,6,5,9,10,13,7,12,4,11,\infty)_{\adfTFADJ}$,

\adfLgap
$(0,2,1,3,4,7,6,10,\infty,13,9,12,5,11)_{\adfTFAEI}$,

\adfLgap
$(0,2,1,3,4,7,10,9,5,11,6,13,\infty,12)_{\adfTFAFH}$,

\adfLgap
$(0,2,1,3,4,7,11,6,9,12,5,13,\infty,8)_{\adfTFAGG}$,

\adfLgap
$(0,2,1,4,3,5,8,\infty,11,6,12,7,13,9)_{\adfTFBBK}$,

\adfLgap
$(0,2,1,3,6,5,7,13,9,\infty,8,10,4,11)_{\adfTFBCJ}$,

\adfLgap
$(0,2,1,3,5,8,4,6,11,\infty,12,7,13,9)_{\adfTFBDI}$,

\adfLgap
$(0,2,1,3,5,8,4,6,11,\infty,12,7,13,9)_{\adfTFBEH}$,

\adfLgap
$(0,2,1,3,5,8,4,9,12,6,\infty,13,7,11)_{\adfTFBFG}$,

\adfLgap
$(0,2,1,4,3,7,6,10,5,13,11,12,\infty,9)_{\adfTFCCI}$,

\adfLgap
$(0,2,1,4,3,5,9,6,\infty,7,13,8,12,11)_{\adfTFCDH}$,

\adfLgap
$(0,2,1,4,3,5,6,11,7,13,9,\infty,12,8)_{\adfTFCEG}$,

\adfLgap
$(0,2,1,4,3,5,6,11,\infty,8,12,7,13,9)_{\adfTFCFF}$,

\adfLgap
$(0,2,1,3,4,5,8,11,6,10,\infty,7,13,9)_{\adfTFDDG}$,

\adfLgap
$(0,2,1,3,4,5,8,\infty,9,6,10,13,7,11)_{\adfTFDEF}$

\adfLgap
$(0,2,1,3,4,6,5,10,\infty,9,8,11,7,13)_{\adfTFEEE}$

\adfLgap
\noindent under the action of the mapping $x \mapsto x + 2 \adfmod{14}$, $\infty \mapsto \infty$.

\ADFvfyParStart{15, \{\{14,7,2\},\{1,1,1\}\}, 14, -1, -1} 

\noindent {\boldmath $K_{16}$}~
Let the vertex set be $Z_{16}$. The decompositions consist of the graphs

\adfLgap
$(0,1,15,8,2,9,12,7,3,4,10,5,14,11)_{\adfTFABL}$,

$(2,3,14,6,9,11,4,7,13,8,5,15,12,1)_{\adfTFABL}$,

$(4,5,0,13,11,8,6,1,14,9,3,15,10,12)_{\adfTFABL}$,

$(6,7,5,15,8,4,12,0,13,3,10,9,1,2)_{\adfTFABL}$,

$(8,9,7,10,6,14,4,1,5,11,0,2,13,15)_{\adfTFABL}$,

$(10,11,7,14,0,9,5,3,8,1,13,6,12,2)_{\adfTFABL}$,

$(12,13,14,3,0,7,1,10,2,15,11,6,4,9)_{\adfTFABL}$,

$(14,15,7,8,12,11,3,6,0,10,13,5,2,4)_{\adfTFABL}$,

\adfLgap
$(0,1,15,11,13,5,3,7,14,8,10,2,12,6)_{\adfTFACK}$,

$(2,3,8,12,7,1,4,14,10,0,6,13,9,15)_{\adfTFACK}$,

$(4,5,7,9,11,2,1,3,0,12,10,6,8,15)_{\adfTFACK}$,

$(6,7,11,0,5,1,9,4,3,10,15,12,14,13)_{\adfTFACK}$,

$(8,9,4,6,0,5,11,14,1,13,10,7,15,2)_{\adfTFACK}$,

$(10,11,4,13,5,2,0,14,6,3,9,12,1,8)_{\adfTFACK}$,

$(12,13,4,2,7,11,3,8,5,14,9,10,1,15)_{\adfTFACK}$,

$(14,15,2,6,3,13,8,7,5,12,11,9,0,4)_{\adfTFACK}$,

\adfLgap
$(0,1,9,12,11,5,7,10,14,8,6,4,3,15)_{\adfTFADJ}$,

$(2,3,11,6,9,15,13,0,7,14,12,1,4,8)_{\adfTFADJ}$,

$(4,5,2,10,8,12,6,13,9,14,1,3,0,11)_{\adfTFADJ}$,

$(6,7,5,12,8,3,13,11,4,15,0,10,9,2)_{\adfTFADJ}$,

$(8,9,11,14,4,1,7,13,2,6,10,12,15,5)_{\adfTFADJ}$,

$(10,11,13,1,9,4,0,2,5,3,14,6,15,7)_{\adfTFADJ}$,

$(12,13,3,10,5,7,9,15,8,0,6,1,2,14)_{\adfTFADJ}$,

$(14,15,5,1,10,0,12,2,8,13,4,7,3,11)_{\adfTFADJ}$,

\adfLgap
$(0,1,15,13,4,7,12,2,11,6,5,8,14,10)_{\adfTFAEI}$,

$(2,3,14,5,11,9,15,12,1,6,4,10,0,7)_{\adfTFAEI}$,

$(4,5,11,8,2,13,15,10,7,9,6,14,1,3)_{\adfTFAEI}$,

$(6,7,8,10,9,13,15,3,11,12,14,0,2,5)_{\adfTFAEI}$,

$(8,9,13,10,3,4,0,11,1,15,5,12,6,2)_{\adfTFAEI}$,

$(10,11,5,1,9,14,6,13,0,3,8,4,12,7)_{\adfTFAEI}$,

$(12,13,9,0,6,3,10,2,4,1,8,15,7,14)_{\adfTFAEI}$,

$(14,15,4,0,5,9,3,12,8,7,2,1,13,11)_{\adfTFAEI}$,

\adfLgap
$(0,1,7,13,11,4,10,14,9,2,15,5,3,6)_{\adfTFAFH}$,

$(2,3,12,6,10,13,14,8,7,4,0,5,11,9)_{\adfTFAFH}$,

$(4,5,1,13,6,15,12,14,10,8,0,11,3,7)_{\adfTFAFH}$,

$(6,7,5,1,2,10,12,0,13,9,4,15,8,11)_{\adfTFAFH}$,

$(8,9,1,14,11,2,7,13,4,12,0,3,15,10)_{\adfTFAFH}$,

$(10,11,0,9,1,3,12,5,8,4,6,2,14,15)_{\adfTFAFH}$,

$(12,14,9,15,0,2,5,13,3,8,6,11,1,7)_{\adfTFAFH}$,

$(13,15,2,4,3,10,7,5,9,6,14,8,12,1)_{\adfTFAFH}$,

\adfLgap
$(0,1,3,6,8,4,12,5,14,7,10,9,2,15)_{\adfTFAGG}$,

$(2,3,1,8,13,0,12,11,4,15,10,5,6,9)_{\adfTFAGG}$,

$(4,5,1,9,15,12,13,11,14,2,6,10,3,7)_{\adfTFAGG}$,

$(6,7,11,0,10,12,8,15,13,3,1,14,9,4)_{\adfTFAGG}$,

$(8,9,14,13,10,1,7,0,11,2,12,3,15,5)_{\adfTFAGG}$,

$(10,11,2,13,7,8,3,4,14,5,0,6,12,9)_{\adfTFAGG}$,

$(12,14,1,13,4,6,15,11,7,2,0,8,5,3)_{\adfTFAGG}$,

$(13,15,5,2,8,10,4,0,9,7,11,1,6,14)_{\adfTFAGG}$,

\adfLgap
$(0,1,7,4,12,10,5,6,14,13,8,3,9,11)_{\adfTFBBK}$,

$(2,3,7,10,0,8,12,4,5,13,1,14,9,15)_{\adfTFBBK}$,

$(4,5,14,8,15,2,12,13,11,7,6,0,3,1)_{\adfTFBBK}$,

$(6,7,8,9,10,15,1,2,3,11,4,13,0,5)_{\adfTFBBK}$,

$(8,9,1,2,11,6,4,10,7,15,12,14,3,5)_{\adfTFBBK}$,

$(10,11,14,0,1,6,2,5,12,7,4,9,13,15)_{\adfTFBBK}$,

$(12,13,3,6,1,0,9,8,10,11,5,15,14,2)_{\adfTFBBK}$,

$(14,15,0,8,7,13,10,9,12,11,2,4,3,6)_{\adfTFBBK}$,

\adfLgap
$(0,1,6,10,5,7,4,13,8,12,15,3,11,14)_{\adfTFBCJ}$,

$(2,3,10,14,6,15,4,8,9,11,7,1,0,13)_{\adfTFBCJ}$,

$(4,5,9,10,13,3,1,8,15,11,6,7,12,2)_{\adfTFBCJ}$,

$(6,7,2,13,15,5,14,12,4,11,10,9,0,3)_{\adfTFBCJ}$,

$(8,9,2,5,7,10,6,4,1,13,14,15,0,12)_{\adfTFBCJ}$,

$(10,11,12,14,0,7,8,6,9,3,5,4,2,1)_{\adfTFBCJ}$,

$(12,15,6,1,10,13,2,11,5,0,8,3,14,9)_{\adfTFBCJ}$,

$(13,14,7,11,8,9,1,15,5,12,3,2,0,4)_{\adfTFBCJ}$,

\adfLgap
$(0,1,3,14,4,9,6,15,13,7,5,2,10,8)_{\adfTFBDI}$,

$(2,3,11,12,5,8,14,9,0,15,4,13,1,7)_{\adfTFBDI}$,

$(4,5,11,0,2,1,3,14,12,9,6,8,13,10)_{\adfTFBDI}$,

$(6,7,10,4,5,9,3,2,13,11,8,12,1,15)_{\adfTFBDI}$,

$(8,9,15,4,10,11,14,1,6,2,7,0,12,3)_{\adfTFBDI}$,

$(10,11,1,0,13,12,3,5,14,6,7,8,2,15)_{\adfTFBDI}$,

$(12,14,10,4,7,11,6,5,0,8,9,13,3,15)_{\adfTFBDI}$,

$(13,15,5,14,7,12,6,11,0,1,4,2,9,10)_{\adfTFBDI}$,

\adfLgap
$(0,1,6,5,8,12,13,9,15,2,11,4,10,3)_{\adfTFBEH}$,

$(2,3,8,0,1,7,6,5,4,9,11,14,12,15)_{\adfTFBEH}$,

$(4,5,6,0,8,1,10,7,14,3,11,13,9,12)_{\adfTFBEH}$,

$(6,7,10,12,4,1,9,8,11,15,0,13,14,2)_{\adfTFBEH}$,

$(8,9,10,15,5,1,2,7,11,6,14,4,13,3)_{\adfTFBEH}$,

$(10,11,12,0,3,7,5,13,2,4,8,14,15,1)_{\adfTFBEH}$,

$(12,13,7,3,4,15,6,2,10,11,0,14,9,5)_{\adfTFBEH}$,

$(14,15,10,1,12,0,7,5,3,2,6,9,8,13)_{\adfTFBEH}$,

\adfLgap
$(0,1,11,6,9,5,7,10,14,4,12,8,13,2)_{\adfTFBFG}$,

$(2,3,6,15,7,8,1,13,5,12,9,11,10,14)_{\adfTFBFG}$,

$(4,5,3,6,11,14,8,15,7,9,10,0,12,13)_{\adfTFBFG}$,

$(6,7,1,5,14,15,13,0,8,9,4,2,10,3)_{\adfTFBFG}$,

$(8,9,0,10,13,7,12,1,11,5,4,15,3,2)_{\adfTFBFG}$,

$(10,11,4,5,8,2,0,3,6,7,14,1,15,12)_{\adfTFBFG}$,

$(12,13,6,3,1,5,0,4,14,2,7,11,15,9)_{\adfTFBFG}$,

$(14,15,6,13,11,2,12,10,9,3,8,4,1,0)_{\adfTFBFG}$,

\adfLgap
$(0,1,10,2,8,9,6,15,7,11,12,3,4,5)_{\adfTFCCI}$,

$(2,3,12,15,4,1,14,13,10,7,8,11,5,0)_{\adfTFCCI}$,

$(4,5,8,3,0,13,9,10,14,6,7,1,15,2)_{\adfTFCCI}$,

$(6,7,11,4,5,14,10,8,15,0,2,13,12,9)_{\adfTFCCI}$,

$(8,9,13,3,6,2,5,12,14,0,1,10,15,11)_{\adfTFCCI}$,

$(10,11,3,14,4,13,12,8,2,7,5,9,6,1)_{\adfTFCCI}$,

$(12,15,1,14,6,4,0,9,13,7,3,11,10,5)_{\adfTFCCI}$,

$(13,14,1,8,15,9,6,3,2,11,0,7,12,4)_{\adfTFCCI}$,

\adfLgap
$(0,1,2,15,7,12,11,10,14,13,5,8,6,4)_{\adfTFCDH}$,

$(2,3,9,10,11,6,1,13,7,4,15,12,8,0)_{\adfTFCDH}$,

$(4,5,13,9,2,14,12,8,3,15,7,10,11,0)_{\adfTFCDH}$,

$(6,7,0,1,5,14,8,15,9,3,11,4,10,2)_{\adfTFCDH}$,

$(8,9,10,6,1,12,0,13,15,11,5,2,3,14)_{\adfTFCDH}$,

$(10,11,5,7,12,9,8,1,13,6,3,4,0,14)_{\adfTFCDH}$,

$(12,15,4,5,3,13,0,6,14,7,9,1,2,8)_{\adfTFCDH}$,

$(13,14,10,15,11,9,4,12,2,6,7,3,5,1)_{\adfTFCDH}$,

\adfLgap
$(0,1,14,7,6,8,15,9,11,2,4,12,13,5)_{\adfTFCEG}$,

$(2,3,10,12,1,0,7,13,9,8,4,6,11,15)_{\adfTFCEG}$,

$(4,5,7,6,1,8,10,14,11,3,0,2,13,15)_{\adfTFCEG}$,

$(6,7,10,9,1,12,0,15,2,14,8,11,5,3)_{\adfTFCEG}$,

$(8,9,13,14,5,0,10,11,3,4,15,2,12,6)_{\adfTFCEG}$,

$(10,11,5,12,3,9,13,1,7,2,8,0,4,14)_{\adfTFCEG}$,

$(12,15,14,6,9,4,13,10,8,7,5,2,3,1)_{\adfTFCEG}$,

$(13,14,6,3,11,7,12,15,0,9,5,4,10,1)_{\adfTFCEG}$,

\adfLgap
$(0,1,2,3,15,6,12,11,7,8,14,4,13,5)_{\adfTFCFF}$,

$(2,3,4,7,11,8,9,12,14,6,10,0,1,13)_{\adfTFCFF}$,

$(4,5,3,15,1,12,10,9,2,0,6,7,13,8)_{\adfTFCFF}$,

$(6,7,13,15,9,1,2,14,10,8,12,4,11,5)_{\adfTFCFF}$,

$(8,9,15,4,7,2,10,11,0,3,5,6,1,14)_{\adfTFCFF}$,

$(10,11,15,14,1,8,4,5,9,13,12,0,3,6)_{\adfTFCFF}$,

$(12,13,7,14,2,8,10,5,0,3,9,15,1,11)_{\adfTFCFF}$,

$(14,15,5,12,0,7,9,13,2,6,4,10,3,11)_{\adfTFCFF}$,

\adfLgap
$(0,1,12,10,14,2,7,15,4,8,9,5,11,6)_{\adfTFDDG}$,

$(2,3,8,5,15,12,1,11,13,4,10,7,14,0)_{\adfTFDDG}$,

$(4,5,3,1,7,2,14,13,12,6,15,0,9,10)_{\adfTFDDG}$,

$(6,7,9,2,11,3,12,8,13,0,1,10,15,4)_{\adfTFDDG}$,

$(8,9,11,4,14,10,13,15,1,5,6,2,3,7)_{\adfTFDDG}$,

$(10,11,2,1,9,6,7,13,0,5,3,8,14,12)_{\adfTFDDG}$,

$(12,14,7,0,6,9,4,5,15,8,13,3,10,11)_{\adfTFDDG}$,

$(13,15,9,3,14,12,5,2,1,4,6,8,0,11)_{\adfTFDDG}$,

\adfLgap
$(0,1,2,11,6,13,10,7,3,5,15,4,14,12)_{\adfTFDEF}$,

$(2,3,15,14,13,9,12,8,11,5,10,1,0,6)_{\adfTFDEF}$,

$(4,5,7,1,8,13,11,10,14,6,12,0,3,9)_{\adfTFDEF}$,

$(6,7,5,3,15,14,0,9,8,10,4,12,13,2)_{\adfTFDEF}$,

$(8,9,13,5,1,14,7,0,4,10,2,12,11,15)_{\adfTFDEF}$,

$(10,11,3,14,1,9,7,5,4,15,6,2,8,0)_{\adfTFDEF}$,

$(12,13,7,11,9,3,4,2,1,10,0,15,8,6)_{\adfTFDEF}$,

$(14,15,11,5,12,9,6,7,13,2,3,8,4,1)_{\adfTFDEF}$

\adfLgap
$(0,1,10,9,6,2,7,12,14,11,15,4,5,8)_{\adfTFEEE}$,

$(2,3,15,13,8,7,11,6,1,12,10,4,9,5)_{\adfTFEEE}$,

$(4,5,1,7,9,14,8,10,15,12,3,0,2,13)_{\adfTFEEE}$,

$(6,7,15,8,3,10,0,14,2,4,12,9,11,13)_{\adfTFEEE}$,

$(8,9,14,15,11,3,6,13,4,0,2,7,5,1)_{\adfTFEEE}$,

$(10,11,13,0,8,12,6,7,14,4,1,3,2,5)_{\adfTFEEE}$,

$(12,13,0,5,15,1,4,6,14,3,10,11,8,9)_{\adfTFEEE}$,

$(14,15,1,0,11,7,10,5,6,3,13,12,2,9)_{\adfTFEEE}$.

\adfLgap
\ADFvfyParStart{16, \{\{16,1,16\}\}, -1, -1, -1} 

\noindent {\boldmath $K_{21}$}~
Let the vertex set be $Z_{21}$. The decompositions consist of the graphs

\adfLgap
$(0,9,7,13,18,19,5,10,2,1,11,12,15,17)_{\adfTFABL}$,

$(5,8,15,9,3,7,16,6,11,17,13,19,1,20)_{\adfTFABL}$,

\adfLgap
$(0,14,4,18,9,10,16,6,19,11,13,1,8,5)_{\adfTFACK}$,

$(12,11,14,6,15,9,7,2,1,19,8,18,13,17)_{\adfTFACK}$,

\adfLgap
$(0,11,8,6,16,5,2,13,12,15,14,1,18,10)_{\adfTFADJ}$,

$(0,14,9,5,7,15,1,3,19,4,13,16,20,8)_{\adfTFADJ}$,

\adfLgap
$(0,1,2,10,8,19,17,13,18,3,15,5,11,12)_{\adfTFAEI}$,

$(0,7,4,12,17,14,19,10,11,3,6,20,8,13)_{\adfTFAEI}$,

\adfLgap
$(0,18,11,14,10,8,19,5,3,7,12,4,17,9)_{\adfTFAFH}$,

$(11,5,2,7,16,13,6,10,4,15,8,1,3,9)_{\adfTFAFH}$,

\adfLgap
$(0,3,10,15,9,5,12,16,2,13,14,19,11,4)_{\adfTFAGG}$,

$(2,11,17,18,6,4,1,13,5,15,20,16,10,3)_{\adfTFAGG}$,

\adfLgap
$(0,11,7,1,3,2,18,8,5,12,10,17,15,19)_{\adfTFBBK}$,

$(5,12,18,4,17,19,16,6,15,10,1,7,2,8)_{\adfTFBBK}$,

\adfLgap
$(0,11,17,14,4,6,8,5,1,9,10,16,12,19)_{\adfTFBCJ}$,

$(17,15,6,1,13,18,2,14,16,19,20,12,9,4)_{\adfTFBCJ}$,

\adfLgap
$(0,17,13,5,7,3,10,15,16,2,1,19,12,8)_{\adfTFBDI}$,

$(5,9,18,2,17,6,3,1,7,16,8,19,14,15)_{\adfTFBDI}$,

\adfLgap
$(0,6,15,2,7,17,4,3,10,8,11,18,1,5)_{\adfTFBEH}$,

$(12,15,2,13,19,16,4,7,8,1,9,5,11,20)_{\adfTFBEH}$,

\adfLgap
$(0,9,1,12,14,17,2,13,7,19,5,15,8,3)_{\adfTFBFG}$,

$(17,11,16,18,14,5,13,3,19,4,1,6,9,7)_{\adfTFBFG}$,

\adfLgap
$(0,19,4,8,15,7,3,17,10,11,5,14,18,13)_{\adfTFCCI}$,

$(7,9,4,18,6,17,5,0,2,10,12,11,14,19)_{\adfTFCCI}$,

\adfLgap
$(0,4,18,7,19,13,11,8,2,3,9,16,17,14)_{\adfTFCDH}$,

$(15,16,6,7,11,2,0,8,18,10,5,13,17,12)_{\adfTFCDH}$,

\adfLgap
$(0,4,6,1,2,10,9,8,3,17,12,19,13,11)_{\adfTFCEG}$,

$(12,4,3,14,8,0,13,15,16,7,2,20,11,5)_{\adfTFCEG}$,

\adfLgap
$(0,3,7,5,10,19,18,12,8,16,11,1,14,17)_{\adfTFCFF}$,

$(17,8,11,4,5,15,12,3,7,18,10,13,19,0)_{\adfTFCFF}$,

\adfLgap
$(0,18,19,12,8,16,15,2,6,3,13,10,17,11)_{\adfTFDDG}$,

$(7,10,3,12,4,2,11,1,5,6,8,16,17,14)_{\adfTFDDG}$,

\adfLgap
$(0,3,1,12,11,19,7,17,8,2,16,13,9,15)_{\adfTFDEF}$,

$(7,16,8,4,9,13,0,3,14,20,6,2,5,11)_{\adfTFDEF}$

\adfLgap
$(0,1,5,15,12,7,4,17,14,2,8,10,13,3)_{\adfTFEEE}$,

$(19,18,8,14,15,11,12,3,9,10,2,16,4,20)_{\adfTFEEE}$

\adfLgap
\noindent under the action of the mapping $x \mapsto x + 3 \adfmod{21}$.

\ADFvfyParStart{21, \{\{21,7,3\}\}, -1, -1, -1} 

\noindent {\boldmath $K_{25}$}~
Let the vertex set be $Z_{25}$. The decompositions consist of the graphs

\adfLgap
$(0,19,3,12,7,21,1,2,6,8,11,20,9,13)_{\adfTFABL}$,

$(4,1,11,14,9,17,15,8,3,21,12,22,23,13)_{\adfTFABL}$,

$(10,2,14,20,21,19,22,15,11,5,18,12,23,9)_{\adfTFABL}$,

$(21,2,4,9,8,17,20,0,11,3,5,14,15,23)_{\adfTFABL}$,

\adfLgap
$(0,5,10,1,9,11,14,6,16,23,17,15,21,22)_{\adfTFACK}$,

$(8,10,15,23,5,19,6,3,13,4,11,7,2,9)_{\adfTFACK}$,

$(22,15,6,1,10,14,19,13,9,7,18,4,3,16)_{\adfTFACK}$,

$(23,2,18,16,6,12,9,22,14,4,10,7,8,17)_{\adfTFACK}$,

\adfLgap
$(0,4,23,16,11,12,7,21,8,3,20,6,10,19)_{\adfTFADJ}$,

$(6,3,12,8,2,15,20,9,7,0,10,4,22,1)_{\adfTFADJ}$,

$(0,6,1,4,23,3,17,8,19,2,11,21,9,14)_{\adfTFADJ}$,

$(9,18,22,20,8,13,14,16,17,7,10,2,24,0)_{\adfTFADJ}$,

\adfLgap
$(0,5,4,22,17,1,18,23,21,14,16,13,9,7)_{\adfTFAEI}$,

$(10,2,23,16,1,15,0,9,18,17,14,20,21,13)_{\adfTFAEI}$,

$(6,0,12,15,23,7,14,9,19,13,1,22,18,3)_{\adfTFAEI}$,

$(11,22,20,9,6,7,16,4,17,23,0,24,13,14)_{\adfTFAEI}$,

\adfLgap
$(0,1,19,13,9,20,12,16,11,15,2,22,14,4)_{\adfTFAFH}$,

$(14,1,7,0,13,18,16,12,10,21,19,20,5,2)_{\adfTFAFH}$,

$(14,13,3,21,18,12,2,6,0,9,22,23,7,4)_{\adfTFAFH}$,

$(18,0,1,17,21,2,23,3,16,9,14,10,5,8)_{\adfTFAFH}$,

\adfLgap
$(0,12,16,22,8,2,17,20,21,13,10,9,4,3)_{\adfTFAGG}$,

$(2,4,14,8,19,10,17,16,7,21,3,15,20,13)_{\adfTFAGG}$,

$(9,16,2,3,18,10,20,14,19,5,11,23,21,1)_{\adfTFAGG}$,

$(20,6,18,23,1,17,13,9,24,2,19,11,7,5)_{\adfTFAGG}$,

\adfLgap
$(0,15,21,14,16,2,7,3,5,22,20,8,17,10)_{\adfTFBBK}$,

$(16,4,13,22,18,0,12,23,20,17,21,6,7,19)_{\adfTFBBK}$,

$(21,13,14,3,4,20,10,11,16,8,0,19,17,9)_{\adfTFBBK}$,

$(18,4,6,23,12,9,21,10,14,19,8,7,17,1)_{\adfTFBBK}$,

\adfLgap
$(0,6,2,3,22,7,21,4,19,5,23,1,13,8)_{\adfTFBCJ}$,

$(9,11,4,13,21,8,14,1,5,18,2,22,23,0)_{\adfTFBCJ}$,

$(19,20,12,2,14,17,3,18,11,5,10,22,16,21)_{\adfTFBCJ}$,

$(0,10,9,4,1,15,12,2,24,8,19,21,22,18)_{\adfTFBCJ}$,

\adfLgap
$(0,19,12,9,23,6,11,17,21,1,3,22,15,14)_{\adfTFBDI}$,

$(6,7,4,13,15,5,16,9,18,8,11,10,14,2)_{\adfTFBDI}$,

$(19,10,13,4,8,15,16,3,7,23,22,14,12,1)_{\adfTFBDI}$,

$(11,17,2,15,9,20,5,18,7,6,14,13,8,0)_{\adfTFBDI}$,

\adfLgap
$(0,9,5,22,12,16,8,19,17,4,21,11,18,13)_{\adfTFBEH}$,

$(4,18,6,19,3,12,17,1,14,9,16,0,23,10)_{\adfTFBEH}$,

$(13,23,2,19,12,9,20,7,21,1,15,14,22,5)_{\adfTFBEH}$,

$(5,8,11,7,0,9,23,1,2,21,20,10,22,6)_{\adfTFBEH}$,

\adfLgap
$(0,23,9,16,4,22,15,10,1,3,20,19,13,8)_{\adfTFBFG}$,

$(20,8,10,13,6,9,11,12,1,22,3,17,4,21)_{\adfTFBFG}$,

$(10,12,14,7,8,4,13,22,2,11,15,1,23,6)_{\adfTFBFG}$,

$(10,17,22,12,1,16,21,14,4,15,18,19,24,9)_{\adfTFBFG}$,

\adfLgap
$(0,6,4,9,15,7,14,1,19,10,13,12,17,11)_{\adfTFCCI}$,

$(5,23,11,8,1,18,0,22,19,7,3,17,6,16)_{\adfTFCCI}$,

$(23,6,17,8,15,14,11,10,12,2,20,18,4,22)_{\adfTFCCI}$,

$(4,19,6,2,3,9,8,15,1,10,23,14,12,0)_{\adfTFCCI}$,

\adfLgap
$(0,17,5,15,16,7,22,8,1,21,20,9,2,13)_{\adfTFCDH}$,

$(21,20,9,14,2,4,12,15,13,1,16,18,17,23)_{\adfTFCDH}$,

$(9,21,23,18,12,3,4,19,1,2,5,17,6,10)_{\adfTFCDH}$,

$(7,9,19,3,11,14,5,0,18,8,4,6,23,10)_{\adfTFCDH}$,

\adfLgap
$(0,16,3,13,2,17,1,9,12,6,4,19,5,20)_{\adfTFCEG}$,

$(20,17,18,10,12,1,21,3,4,7,8,9,13,22)_{\adfTFCEG}$,

$(2,20,23,14,1,3,16,15,19,12,18,4,9,6)_{\adfTFCEG}$,

$(21,0,5,18,11,7,9,22,13,8,14,1,20,24)_{\adfTFCEG}$,

\adfLgap
$(0,22,13,4,5,6,19,23,12,21,8,2,1,10)_{\adfTFCFF}$,

$(2,9,3,4,18,11,1,12,16,23,5,8,19,0)_{\adfTFCFF}$,

$(1,6,7,0,9,20,10,21,23,17,22,19,2,4)_{\adfTFCFF}$,

$(9,5,19,13,6,1,23,8,3,10,14,2,0,22)_{\adfTFCFF}$,

\adfLgap
$(0,1,22,15,8,11,21,12,19,5,14,18,20,3)_{\adfTFDDG}$,

$(0,16,5,4,14,17,18,15,21,1,23,6,2,19)_{\adfTFDDG}$,

$(4,12,22,17,11,21,14,3,13,0,6,18,8,9)_{\adfTFDDG}$,

$(19,7,17,23,3,6,15,5,14,8,22,12,24,20)_{\adfTFDDG}$,

\adfLgap
$(0,1,17,22,19,9,13,5,23,14,12,18,2,11)_{\adfTFDEF}$,

$(13,22,21,23,18,14,6,19,4,15,2,5,20,8)_{\adfTFDEF}$,

$(13,7,19,3,21,10,4,9,6,23,22,20,15,1)_{\adfTFDEF}$,

$(1,5,10,9,11,22,14,2,12,21,0,4,18,6)_{\adfTFDEF}$

\adfLgap
$(0,2,17,20,19,7,4,18,21,23,11,9,13,15)_{\adfTFEEE}$,

$(0,1,16,3,20,5,13,10,12,19,6,23,18,7)_{\adfTFEEE}$,

$(20,2,14,22,6,16,13,7,23,4,21,9,19,3)_{\adfTFEEE}$,

$(0,21,5,12,22,1,14,19,18,3,9,6,7,4)_{\adfTFEEE}$

\adfLgap
\noindent under the action of the mapping $x \mapsto x + 5 \adfmod{25}$.

\ADFvfyParStart{25, \{\{25,5,5\}\}, -1, -1, -1} 

















\noindent {\boldmath $K_{36}$}~
Let the vertex set be $Z_{35} \cup \{\infty\}$. The decompositions consist of the graphs

\adfLgap
$(\infty,30,14,26,7,10,4,23,12,21,34,31,33,18)_{\adfTFABL}$,

$(7,29,19,34,1,8,5,3,25,4,15,11,31,30)_{\adfTFABL}$,

$(16,12,2,22,8,29,25,5,33,6,23,26,0,10)_{\adfTFABL}$,

$(5,10,22,32,2,30,16,0,11,19,24,23,33,28)_{\adfTFABL}$,

$(18,24,27,9,2,26,20,7,3,34,11,16,28,22)_{\adfTFABL}$,

$(19,2,1,4,11,21,8,0,9,33,17,18,\infty,22)_{\adfTFABL}$,

\adfLgap
$(\infty,14,27,5,19,1,28,33,0,3,34,18,24,22)_{\adfTFADJ}$,

$(24,12,30,10,34,13,5,28,32,3,6,4,0,2)_{\adfTFADJ}$,

$(24,11,15,31,18,27,26,30,20,13,3,12,23,6)_{\adfTFADJ}$,

$(9,26,27,15,32,1,0,21,17,18,33,24,25,20)_{\adfTFADJ}$,

$(4,18,20,12,32,11,0,7,26,29,22,1,27,16)_{\adfTFADJ}$,

$(15,6,28,16,26,12,7,23,5,\infty,3,4,14,29)_{\adfTFADJ}$,

\adfLgap
$(\infty,28,7,25,17,20,15,14,2,1,31,24,18,21)_{\adfTFAFH}$,

$(9,2,20,18,31,0,21,8,26,25,16,1,28,27)_{\adfTFAFH}$,

$(17,13,14,31,33,29,25,3,8,2,26,12,23,34)_{\adfTFAFH}$,

$(5,16,15,21,34,0,19,14,24,1,32,2,25,28)_{\adfTFAFH}$,

$(23,4,15,22,31,29,2,8,34,17,11,30,10,24)_{\adfTFAFH}$,

$(27,14,5,7,22,3,20,18,8,15,\infty,1,11,19)_{\adfTFAFH}$,

\adfLgap
$(\infty,11,28,6,20,15,32,17,23,13,4,21,9,19)_{\adfTFBBK}$,

$(18,11,21,14,4,17,31,3,22,24,0,15,13,33)_{\adfTFBBK}$,

$(19,4,0,26,13,11,3,14,29,24,30,33,2,12)_{\adfTFBBK}$,

$(34,0,1,18,33,7,10,26,15,3,17,30,20,13)_{\adfTFBBK}$,

$(24,12,10,17,27,15,14,5,1,32,16,25,31,11)_{\adfTFBBK}$,

$(1,24,7,12,13,2,3,8,0,27,20,6,32,\infty)_{\adfTFBBK}$,

\adfLgap
$(\infty,28,14,3,12,25,22,11,8,15,31,23,29,4)_{\adfTFBCJ}$,

$(24,8,28,32,20,23,12,34,16,14,15,5,6,26)_{\adfTFBCJ}$,

$(25,17,14,30,22,33,8,3,31,19,34,29,12,11)_{\adfTFBCJ}$,

$(12,17,27,16,24,21,15,9,7,26,31,5,19,10)_{\adfTFBCJ}$,

$(6,5,20,34,22,19,0,11,33,32,3,17,13,1)_{\adfTFBCJ}$,

$(18,5,9,0,3,16,2,\infty,1,11,14,33,20,7)_{\adfTFBCJ}$,

\adfLgap
$(\infty,10,1,25,5,19,22,18,11,30,2,27,9,32)_{\adfTFBDI}$,

$(33,0,28,22,14,13,19,34,31,8,6,17,3,12)_{\adfTFBDI}$,

$(27,30,25,8,7,22,10,9,19,1,34,28,24,20)_{\adfTFBDI}$,

$(33,27,11,14,25,31,16,30,24,13,28,22,9,1)_{\adfTFBDI}$,

$(12,18,9,15,16,8,26,30,28,5,24,29,6,21)_{\adfTFBDI}$,

$(29,6,1,16,17,12,22,24,\infty,8,0,3,20,31)_{\adfTFBDI}$,

\adfLgap
$(\infty,34,15,29,9,5,32,18,26,20,33,28,2,6)_{\adfTFBEH}$,

$(8,18,24,10,19,29,31,32,27,16,17,0,7,28)_{\adfTFBEH}$,

$(3,7,26,21,23,0,24,34,29,8,15,30,20,1)_{\adfTFBEH}$,

$(20,32,11,21,28,13,10,17,7,3,2,8,0,29)_{\adfTFBEH}$,

$(12,29,28,21,18,0,2,32,9,30,16,34,26,6)_{\adfTFBEH}$,

$(29,1,26,3,14,7,\infty,30,0,12,34,21,10,6)_{\adfTFBEH}$,

\adfLgap
$(\infty,22,32,11,28,6,29,19,15,25,24,7,2,33)_{\adfTFBFG}$,

$(6,9,26,20,11,2,4,0,5,12,3,8,15,17)_{\adfTFBFG}$,

$(0,34,18,30,24,31,32,12,16,22,10,21,19,23)_{\adfTFBFG}$,

$(1,27,11,9,24,29,17,30,12,18,32,5,3,6)_{\adfTFBFG}$,

$(4,8,1,3,29,16,21,33,25,13,10,30,14,7)_{\adfTFBFG}$,

$(2,26,18,6,0,8,14,28,20,33,13,\infty,19,30)_{\adfTFBFG}$,

\adfLgap
$(\infty,11,8,4,6,7,3,32,5,23,25,15,30,26)_{\adfTFCDH}$,

$(20,28,15,16,33,27,6,13,29,17,34,14,1,4)_{\adfTFCDH}$,

$(4,6,12,17,0,23,8,16,10,31,14,7,5,13)_{\adfTFCDH}$,

$(7,24,23,33,4,31,22,28,2,17,1,18,13,19)_{\adfTFCDH}$,

$(31,9,6,11,20,2,25,27,17,18,19,29,15,22)_{\adfTFCDH}$,

$(29,4,0,\infty,30,21,18,20,1,7,10,13,2,15)_{\adfTFCDH}$,

\adfLgap
$(\infty,19,21,32,34,31,0,14,17,5,10,7,12,30)_{\adfTFCFF}$,

$(15,22,17,10,0,29,28,12,8,25,6,30,21,16)_{\adfTFCFF}$,

$(25,14,11,34,9,1,23,15,8,17,16,29,2,28)_{\adfTFCFF}$,

$(12,6,2,13,25,28,15,19,21,32,1,26,20,8)_{\adfTFCFF}$,

$(27,24,28,8,1,19,12,31,13,9,4,5,6,33)_{\adfTFCFF}$,

$(10,3,8,13,19,9,16,2,\infty,28,24,12,18,6)_{\adfTFCFF}$,

\adfLgap
$(\infty,28,29,27,4,16,15,5,13,30,19,12,2,18)_{\adfTFDDG}$,

$(27,10,9,12,3,16,5,1,19,28,29,6,25,4)_{\adfTFDDG}$,

$(17,30,25,2,16,21,18,10,3,27,20,7,26,32)_{\adfTFDDG}$,

$(31,26,6,13,28,29,16,3,1,30,34,24,4,8)_{\adfTFDDG}$,

$(30,26,25,34,29,27,22,9,12,21,6,33,20,19)_{\adfTFDDG}$,

$(17,4,2,3,33,23,20,18,13,29,10,\infty,32,31)_{\adfTFDDG}$

\adfLgap
$(\infty,21,9,7,27,11,5,31,20,2,34,23,26,13)_{\adfTFDEF}$,

$(29,5,24,14,15,28,32,1,17,3,13,25,11,12)_{\adfTFDEF}$,

$(23,34,9,0,22,21,33,13,7,29,11,18,17,3)_{\adfTFDEF}$,

$(26,28,19,0,20,11,24,1,12,10,34,31,4,25)_{\adfTFDEF}$,

$(1,23,31,22,32,11,28,15,12,5,9,29,0,18)_{\adfTFDEF}$,

$(9,2,22,5,0,11,32,14,7,28,30,23,\infty,10)_{\adfTFDEF}$

\adfLgap
\noindent under the action of the mapping $x \mapsto x + 5 \adfmod{35}$, $\infty \mapsto \infty$.

\ADFvfyParStart{36, \{\{35,7,5\},\{1,1,1\}\}, 35, -1, -1} 

\noindent {\boldmath $K_{40}$}~
Let the vertex set be $Z_{39} \cup \{\infty\}$. The decompositions consist of the graphs

\adfLgap
$(\infty,3,37,29,1,25,32,5,34,13,26,23,21,19)_{\adfTFABL}$,

$(5,24,38,4,8,0,31,14,22,12,6,11,29,3)_{\adfTFABL}$,

$(8,32,36,38,22,15,18,1,27,0,29,10,16,7)_{\adfTFABL}$,

$(13,27,16,1,36,35,3,4,2,18,9,26,31,12)_{\adfTFABL}$,

\adfLgap
$(\infty,1,15,30,5,29,26,10,7,28,35,12,11,21)_{\adfTFADJ}$,

$(25,36,5,17,2,38,10,35,4,30,32,34,0,6)_{\adfTFADJ}$,

$(20,27,38,5,14,0,25,13,7,6,28,21,24,35)_{\adfTFADJ}$,

$(25,15,35,1,31,10,6,2,19,20,3,21,9,7)_{\adfTFADJ}$,

\adfLgap
$(\infty,2,34,32,1,9,35,33,14,18,37,16,5,30)_{\adfTFAFH}$,

$(20,13,29,27,4,38,0,5,36,9,10,34,31,24)_{\adfTFAFH}$,

$(25,3,19,33,31,36,18,21,30,35,10,11,32,0)_{\adfTFAFH}$,

$(10,14,37,17,27,33,4,19,2,5,31,8,20,36)_{\adfTFAFH}$,

\adfLgap
$(\infty,10,16,9,5,37,28,1,25,3,13,2,36,30)_{\adfTFBBK}$,

$(0,19,27,30,38,12,28,15,29,14,10,26,13,5)_{\adfTFBBK}$,

$(10,30,6,28,7,12,15,8,37,32,29,9,11,38)_{\adfTFBBK}$,

$(14,27,34,23,24,3,20,22,36,8,25,26,5,11)_{\adfTFBBK}$,

\adfLgap
$(\infty,13,27,19,6,38,16,28,12,7,32,1,0,20)_{\adfTFBCJ}$,

$(15,16,5,11,6,34,14,38,24,35,31,22,23,26)_{\adfTFBCJ}$,

$(26,0,31,13,9,10,25,22,4,15,2,32,20,38)_{\adfTFBCJ}$,

$(33,1,18,21,3,36,5,12,35,2,24,30,32,34)_{\adfTFBCJ}$,

\adfLgap
$(\infty,25,0,16,15,13,26,5,21,30,8,38,18,36)_{\adfTFBDI}$,

$(19,15,17,24,34,3,10,28,2,33,20,35,21,14)_{\adfTFBDI}$,

$(25,16,20,17,37,1,18,14,2,35,10,36,30,33)_{\adfTFBDI}$,

$(2,19,13,1,8,3,5,15,0,4,20,37,18,29)_{\adfTFBDI}$,

\adfLgap
$(\infty,11,4,6,2,22,38,32,14,9,17,31,7,16)_{\adfTFBEH}$,

$(22,15,35,30,31,32,26,0,12,1,28,21,25,9)_{\adfTFBEH}$,

$(1,2,37,27,3,35,18,23,34,32,33,12,26,0)_{\adfTFBEH}$,

$(26,22,16,29,5,14,24,34,0,9,6,25,7,36)_{\adfTFBEH}$,

\adfLgap
$(\infty,38,3,19,35,24,26,18,17,29,25,8,11,32)_{\adfTFBFG}$,

$(15,30,20,2,3,10,33,21,28,37,18,7,4,31)_{\adfTFBFG}$,

$(5,34,35,13,38,19,15,12,21,14,16,22,24,9)_{\adfTFBFG}$,

$(0,26,31,17,4,14,7,12,21,15,25,1,22,11)_{\adfTFBFG}$,

\adfLgap
$(\infty,33,21,8,13,17,26,32,31,27,4,1,25,3)_{\adfTFCDH}$,

$(22,19,31,8,4,38,9,5,20,18,25,27,7,32)_{\adfTFCDH}$,

$(31,23,4,35,29,24,13,38,15,12,30,11,0,17)_{\adfTFCDH}$,

$(37,34,17,9,4,12,0,24,20,23,2,3,27,33)_{\adfTFCDH}$,

\adfLgap
$(\infty,29,33,17,38,13,5,8,0,16,4,7,11,15)_{\adfTFCFF}$,

$(25,7,26,15,35,1,12,36,16,14,8,27,5,29)_{\adfTFCFF}$,

$(8,5,31,18,26,27,10,17,12,28,30,25,9,3)_{\adfTFCFF}$,

$(21,31,25,7,12,22,15,18,6,3,4,2,11,37)_{\adfTFCFF}$,

\adfLgap
$(\infty,7,30,13,22,14,16,32,4,17,28,24,3,8)_{\adfTFDDG}$,

$(13,32,34,29,22,6,7,12,27,30,5,26,17,15)_{\adfTFDDG}$,

$(16,20,6,34,37,3,30,23,22,10,2,12,11,5)_{\adfTFDDG}$,

$(5,6,25,9,0,36,8,4,17,30,10,18,33,29)_{\adfTFDDG}$

\adfLgap
$(\infty,23,25,31,1,27,18,29,16,14,34,19,6,24)_{\adfTFDEF}$,

$(11,16,21,17,13,23,9,33,36,27,31,30,37,26)_{\adfTFDEF}$,

$(8,2,14,33,16,21,7,28,0,26,9,36,5,10)_{\adfTFDEF}$,

$(32,12,23,22,10,16,0,5,7,17,20,27,33,4)_{\adfTFDEF}$

\adfLgap
\noindent under the action of the mapping $x \mapsto x + 3 \adfmod{39}$, $\infty \mapsto \infty$.

\ADFvfyParStart{40, \{\{39,13,3\},\{1,1,1\}\}, 39, -1, -1} 

\noindent {\boldmath $K_{51}$}~
Let the vertex set be $Z_{51}$. The decompositions consist of the graphs

\adfLgap
$(0,47,18,49,16,38,26,27,46,33,9,48,5,36)_{\adfTFABL}$,

$(24,15,8,9,23,44,19,17,27,10,14,46,43,2)_{\adfTFABL}$,

$(30,35,19,5,33,14,28,21,22,27,24,16,37,44)_{\adfTFABL}$,

$(0,25,45,37,21,42,44,31,40,17,23,28,11,26)_{\adfTFABL}$,

$(45,22,11,49,41,38,5,32,1,7,19,9,31,46)_{\adfTFABL}$,

\adfLgap
$(0,38,25,13,16,10,44,26,5,33,40,47,27,24)_{\adfTFADJ}$,

$(41,38,39,37,1,10,24,40,32,6,48,8,49,22)_{\adfTFADJ}$,

$(16,48,3,42,18,44,43,1,46,25,30,11,22,26)_{\adfTFADJ}$,

$(13,42,21,11,27,26,35,18,23,8,10,6,24,25)_{\adfTFADJ}$,

$(39,47,19,42,35,28,2,3,9,5,11,43,25,20)_{\adfTFADJ}$,

\adfLgap
$(0,38,4,48,31,17,34,20,1,42,22,14,33,44)_{\adfTFAFH}$,

$(5,30,33,23,9,38,29,16,28,13,31,1,36,4)_{\adfTFAFH}$,

$(9,39,2,45,43,16,25,6,30,19,27,28,15,48)_{\adfTFAFH}$,

$(47,11,37,15,17,44,16,40,12,6,45,9,14,49)_{\adfTFAFH}$,

$(5,23,35,4,10,8,11,9,26,1,2,31,34,39)_{\adfTFAFH}$,

\adfLgap
$(0,42,18,11,44,10,12,27,4,33,23,48,13,8)_{\adfTFBBK}$,

$(43,44,7,40,49,42,23,11,1,18,47,46,2,29)_{\adfTFBBK}$,

$(22,43,1,35,40,48,45,2,30,18,31,42,26,24)_{\adfTFBBK}$,

$(38,11,13,44,42,21,8,28,47,17,3,4,9,46)_{\adfTFBBK}$,

$(40,8,17,9,29,26,4,16,6,0,25,49,45,3)_{\adfTFBBK}$,

\adfLgap
$(0,42,19,26,47,22,27,17,7,34,23,16,46,38)_{\adfTFBCJ}$,

$(0,46,8,13,12,24,33,7,40,41,10,19,39,31)_{\adfTFBCJ}$,

$(20,45,38,42,31,39,21,5,1,46,34,36,43,33)_{\adfTFBCJ}$,

$(10,39,26,6,36,32,37,34,11,17,18,3,14,2)_{\adfTFBCJ}$,

$(41,3,5,7,26,43,17,0,6,22,8,35,32,23)_{\adfTFBCJ}$,

\adfLgap
$(0,44,7,3,42,21,24,29,49,10,46,16,17,9)_{\adfTFBDI}$,

$(46,36,37,30,28,18,2,21,38,24,19,35,29,40)_{\adfTFBDI}$,

$(20,19,0,42,29,32,47,36,49,22,30,5,6,37)_{\adfTFBDI}$,

$(24,5,7,30,1,35,39,28,20,16,26,33,29,34)_{\adfTFBDI}$,

$(18,32,20,9,46,14,8,34,6,31,25,28,5,41)_{\adfTFBDI}$,

\adfLgap
$(0,16,8,15,27,13,9,48,3,37,6,4,1,31)_{\adfTFBEH}$,

$(27,18,8,20,22,40,2,9,35,24,38,23,43,42)_{\adfTFBEH}$,

$(12,35,33,17,41,21,49,7,31,18,10,4,30,39)_{\adfTFBEH}$,

$(43,40,8,27,16,26,17,4,36,35,31,14,32,11)_{\adfTFBEH}$,

$(44,43,38,19,9,31,32,6,23,11,14,36,8,1)_{\adfTFBEH}$,

\adfLgap
$(0,11,39,36,16,28,14,29,18,47,7,13,26,19)_{\adfTFBFG}$,

$(9,16,13,2,15,10,45,41,17,34,5,40,4,37)_{\adfTFBFG}$,

$(3,12,9,45,15,39,38,37,14,48,10,8,29,1)_{\adfTFBFG}$,

$(12,39,49,34,29,20,18,37,13,36,41,10,14,46)_{\adfTFBFG}$,

$(40,32,48,6,38,44,3,29,13,4,14,0,20,8)_{\adfTFBFG}$,

\adfLgap
$(0,40,33,39,34,13,18,20,12,36,41,42,45,43)_{\adfTFCDH}$,

$(49,10,7,32,30,40,23,6,27,29,28,15,8,26)_{\adfTFCDH}$,

$(32,10,0,14,44,15,2,46,18,34,27,1,19,17)_{\adfTFCDH}$,

$(24,41,33,45,47,6,10,9,35,16,44,38,2,11)_{\adfTFCDH}$,

$(10,34,37,22,21,38,3,4,14,17,28,42,2,29)_{\adfTFCDH}$,

\adfLgap
$(0,19,22,32,26,21,30,2,31,29,48,38,24,45)_{\adfTFCFF}$,

$(34,40,2,7,3,19,35,33,43,32,31,1,24,17)_{\adfTFCFF}$,

$(46,16,9,20,32,8,23,26,5,12,36,48,44,27)_{\adfTFCFF}$,

$(36,35,4,0,43,37,45,40,27,5,6,24,9,10)_{\adfTFCFF}$,

$(17,16,5,47,11,28,19,26,8,30,24,27,25,1)_{\adfTFCFF}$,

\adfLgap
$(0,9,16,3,15,27,48,46,17,24,11,10,4,28)_{\adfTFDDG}$,

$(10,24,29,0,8,39,21,44,18,9,40,17,42,19)_{\adfTFDDG}$,

$(13,32,5,11,8,47,34,44,9,23,2,35,49,16)_{\adfTFDDG}$,

$(0,12,1,45,48,11,26,17,32,37,41,34,25,22)_{\adfTFDDG}$,

$(6,26,8,18,1,5,7,37,2,22,34,19,30,4)_{\adfTFDDG}$

\adfLgap
$(0,28,26,35,2,13,27,22,46,38,41,37,30,9)_{\adfTFDEF}$,

$(9,35,21,45,14,12,20,26,11,13,39,2,25,4)_{\adfTFDEF}$,

$(2,4,13,21,5,3,47,40,23,16,33,27,45,46)_{\adfTFDEF}$,

$(13,48,5,45,25,1,16,27,44,19,9,14,2,12)_{\adfTFDEF}$,

$(13,33,23,10,5,16,0,2,4,35,6,15,47,31)_{\adfTFDEF}$

\adfLgap
\noindent under the action of the mapping $x \mapsto x + 3\adfmod{51}$.

\ADFvfyParStart{51, \{\{51,17,3\}\}, -1, -1, -1} 

\noindent {\boldmath $K_{55}$}~
Let the vertex set be $Z_{55}$. The decompositions consist of the graphs

\adfLgap
$(0,29,41,26,15,13,44,50,28,34,32,4,22,45)_{\adfTFABL}$,

$(8,44,31,7,48,43,20,28,0,42,49,37,34,9)_{\adfTFABL}$,

$(21,39,18,19,15,30,11,27,49,16,40,53,43,50)_{\adfTFABL}$,

$(49,32,1,39,30,20,26,48,34,18,35,47,46,37)_{\adfTFABL}$,

$(37,53,28,15,36,4,40,19,50,22,5,51,6,2)_{\adfTFABL}$,

$(35,27,38,10,53,11,19,23,6,26,30,29,34,8)_{\adfTFABL}$,

$(17,43,51,11,10,7,13,28,46,16,31,36,47,23)_{\adfTFABL}$,

$(27,52,42,45,25,39,11,47,12,4,17,31,14,29)_{\adfTFABL}$,

$(51,53,35,3,49,46,17,10,5,7,28,16,43,54)_{\adfTFABL}$,

\adfLgap
$(0,8,16,17,35,36,19,44,4,10,15,22,33,26)_{\adfTFADJ}$,

$(52,53,10,20,22,33,31,49,4,50,14,19,42,29)_{\adfTFADJ}$,

$(24,4,11,46,23,8,45,44,53,3,6,36,17,12)_{\adfTFADJ}$,

$(53,0,26,40,12,36,2,47,19,11,39,41,31,23)_{\adfTFADJ}$,

$(28,49,45,15,16,2,48,23,43,1,47,51,30,6)_{\adfTFADJ}$,

$(13,7,47,41,45,10,53,9,35,50,34,52,38,23)_{\adfTFADJ}$,

$(0,11,13,3,4,26,52,48,19,17,6,22,25,16)_{\adfTFADJ}$,

$(5,49,25,2,26,27,42,4,8,41,35,13,19,15)_{\adfTFADJ}$,

$(13,1,20,12,42,54,2,9,50,19,52,32,44,41)_{\adfTFADJ}$,

\adfLgap
$(0,28,39,26,18,31,51,37,3,49,9,46,12,38)_{\adfTFAFH}$,

$(3,17,41,38,31,30,50,35,52,49,29,13,32,34)_{\adfTFAFH}$,

$(9,14,3,15,13,1,25,16,53,23,34,6,12,39)_{\adfTFAFH}$,

$(26,42,22,29,25,44,12,34,1,40,11,38,35,41)_{\adfTFAFH}$,

$(31,6,17,8,13,46,37,35,10,7,26,15,2,44)_{\adfTFAFH}$,

$(21,40,0,49,2,13,33,26,28,12,34,5,20,3)_{\adfTFAFH}$,

$(3,18,7,22,4,25,17,29,33,27,1,15,10,42)_{\adfTFAFH}$,

$(7,52,50,8,9,10,17,5,38,4,18,49,51,41)_{\adfTFAFH}$,

$(40,50,9,28,20,29,41,54,31,46,49,39,27,22)_{\adfTFAFH}$,

\adfLgap
$(0,16,53,3,43,44,13,37,18,51,4,10,25,11)_{\adfTFBBK}$,

$(29,38,10,37,48,19,9,11,20,50,40,33,8,28)_{\adfTFBBK}$,

$(52,21,7,40,16,10,49,27,20,24,36,29,50,18)_{\adfTFBBK}$,

$(26,43,47,22,2,7,13,0,51,42,39,34,45,25)_{\adfTFBBK}$,

$(18,25,26,24,29,14,7,30,51,31,4,43,39,52)_{\adfTFBBK}$,

$(22,7,6,20,48,41,47,38,53,15,24,33,0,37)_{\adfTFBBK}$,

$(2,50,17,19,22,51,23,6,30,1,31,20,15,12)_{\adfTFBBK}$,

$(23,9,46,7,28,26,37,19,42,53,41,51,34,31)_{\adfTFBBK}$,

$(35,36,51,49,43,9,29,0,47,19,7,48,34,4)_{\adfTFBBK}$,

\adfLgap
$(0,20,1,35,7,10,13,14,41,8,2,9,43,45)_{\adfTFBCJ}$,

$(14,45,27,46,22,32,23,33,53,24,25,30,1,38)_{\adfTFBCJ}$,

$(33,12,26,1,53,10,48,46,24,0,22,34,17,38)_{\adfTFBCJ}$,

$(22,34,53,16,43,14,4,50,12,10,26,37,18,23)_{\adfTFBCJ}$,

$(7,46,15,4,51,34,50,24,45,23,27,6,9,5)_{\adfTFBCJ}$,

$(19,52,36,42,53,13,0,51,49,9,16,29,15,26)_{\adfTFBCJ}$,

$(52,50,13,19,1,45,39,43,29,10,31,11,2,7)_{\adfTFBCJ}$,

$(28,44,0,53,13,41,26,16,25,10,7,27,29,24)_{\adfTFBCJ}$,

$(41,47,37,16,17,29,4,43,0,8,46,38,26,7)_{\adfTFBCJ}$,

\adfLgap
$(0,28,17,6,51,7,4,12,8,40,39,42,16,38)_{\adfTFBDI}$,

$(50,49,46,42,48,8,39,32,12,16,13,15,17,11)_{\adfTFBDI}$,

$(53,41,19,12,14,32,23,36,6,20,33,30,45,15)_{\adfTFBDI}$,

$(1,4,14,24,12,44,3,51,30,37,21,13,48,11)_{\adfTFBDI}$,

$(4,14,41,28,1,6,9,25,53,49,19,0,31,40)_{\adfTFBDI}$,

$(30,33,17,42,10,28,13,39,48,4,23,15,29,2)_{\adfTFBDI}$,

$(8,53,7,9,42,27,24,45,17,35,29,5,10,21)_{\adfTFBDI}$,

$(8,20,30,46,11,42,14,49,40,37,7,26,25,13)_{\adfTFBDI}$,

$(17,37,16,34,46,32,27,8,51,11,30,10,26,24)_{\adfTFBDI}$,

\adfLgap
$(0,33,11,5,37,18,38,51,48,1,17,35,36,22)_{\adfTFBEH}$,

$(34,35,48,38,14,7,29,37,9,5,53,24,45,52)_{\adfTFBEH}$,

$(24,2,10,3,35,23,21,1,28,37,29,47,26,43)_{\adfTFBEH}$,

$(13,29,12,31,27,36,44,42,37,17,47,21,45,19)_{\adfTFBEH}$,

$(16,5,11,30,45,14,26,19,1,29,20,9,8,25)_{\adfTFBEH}$,

$(3,14,9,51,36,29,30,33,21,4,23,17,41,50)_{\adfTFBEH}$,

$(6,18,31,17,27,5,2,26,39,7,3,48,4,34)_{\adfTFBEH}$,

$(48,7,50,15,23,46,20,27,14,16,6,39,4,47)_{\adfTFBEH}$,

$(45,3,0,20,36,42,18,8,35,37,10,46,47,49)_{\adfTFBEH}$,

\adfLgap
$(0,28,48,44,51,2,41,43,26,14,7,21,24,6)_{\adfTFBFG}$,

$(43,8,36,25,15,22,27,52,33,21,35,19,20,44)_{\adfTFBFG}$,

$(2,44,6,31,51,36,15,40,42,30,22,9,35,14)_{\adfTFBFG}$,

$(46,26,4,37,29,39,27,10,35,49,9,32,52,23)_{\adfTFBFG}$,

$(1,27,6,9,33,4,8,16,26,49,52,15,34,43)_{\adfTFBFG}$,

$(2,50,8,15,51,52,21,30,33,45,40,0,3,26)_{\adfTFBFG}$,

$(50,51,33,47,19,24,4,38,18,45,23,29,12,49)_{\adfTFBFG}$,

$(5,44,3,27,18,10,4,43,6,0,32,30,26,16)_{\adfTFBFG}$,

$(47,53,48,20,29,7,17,36,13,24,22,8,27,23)_{\adfTFBFG}$,

\adfLgap
$(0,6,35,40,27,42,22,1,11,18,47,10,52,50)_{\adfTFCDH}$,

$(30,32,45,35,47,17,13,8,18,29,44,49,2,26)_{\adfTFCDH}$,

$(21,40,45,39,3,29,1,46,12,23,25,16,30,9)_{\adfTFCDH}$,

$(5,17,19,26,33,53,10,37,36,44,24,11,30,38)_{\adfTFCDH}$,

$(19,40,47,24,51,23,32,15,11,46,13,20,26,7)_{\adfTFCDH}$,

$(32,16,37,19,34,33,1,46,51,22,26,13,7,8)_{\adfTFCDH}$,

$(3,13,28,45,27,49,19,53,12,24,17,33,18,0)_{\adfTFCDH}$,

$(42,46,39,29,31,43,34,30,49,53,36,14,28,9)_{\adfTFCDH}$,

$(40,24,11,8,43,9,35,10,19,32,42,4,28,26)_{\adfTFCDH}$,

\adfLgap
$(0,39,10,29,46,18,32,53,2,49,16,20,25,41)_{\adfTFCFF}$,

$(24,19,15,13,38,14,23,40,53,9,29,7,26,50)_{\adfTFCFF}$,

$(41,12,53,52,27,10,2,28,5,33,1,26,34,9)_{\adfTFCFF}$,

$(31,11,16,42,49,20,26,37,7,32,34,39,36,0)_{\adfTFCFF}$,

$(46,35,36,53,37,16,39,1,7,25,29,50,43,13)_{\adfTFCFF}$,

$(35,53,19,8,24,23,21,28,48,17,9,38,1,14)_{\adfTFCFF}$,

$(24,50,47,27,12,2,41,34,53,28,40,48,26,23)_{\adfTFCFF}$,

$(16,28,44,17,30,10,40,42,47,51,39,22,25,37)_{\adfTFCFF}$,

$(41,17,46,33,40,25,12,45,4,15,14,7,38,23)_{\adfTFCFF}$,

\adfLgap
$(0,6,18,47,52,49,25,12,26,15,44,38,43,3)_{\adfTFDDG}$,

$(35,43,34,45,2,41,36,21,49,32,39,44,40,22)_{\adfTFDDG}$,

$(7,44,41,23,53,3,13,29,21,17,50,16,30,35)_{\adfTFDDG}$,

$(34,17,32,21,8,22,38,2,42,10,11,20,51,41)_{\adfTFDDG}$,

$(33,16,2,12,23,27,20,35,29,52,34,48,47,44)_{\adfTFDDG}$,

$(21,19,2,41,11,22,5,48,33,16,51,0,30,49)_{\adfTFDDG}$,

$(52,27,24,14,1,22,25,35,32,30,2,44,13,49)_{\adfTFDDG}$,

$(3,9,25,48,46,38,40,43,50,33,44,1,24,31)_{\adfTFDDG}$,

$(48,11,0,16,14,1,39,38,20,54,34,50,30,43)_{\adfTFDDG}$

\adfLgap
$(0,36,49,38,46,35,43,19,47,34,50,48,53,23)_{\adfTFDEF}$,

$(52,36,38,3,43,50,46,32,53,4,31,7,19,15)_{\adfTFDEF}$,

$(34,26,37,47,46,49,53,17,11,12,7,38,24,25)_{\adfTFDEF}$,

$(52,18,1,53,24,10,2,39,6,35,38,11,0,5)_{\adfTFDEF}$,

$(36,50,0,9,17,6,35,42,26,38,28,37,40,13)_{\adfTFDEF}$,

$(22,21,18,47,53,33,0,14,24,37,5,48,9,44)_{\adfTFDEF}$,

$(49,51,41,7,35,6,20,5,34,3,39,38,37,1)_{\adfTFDEF}$,

$(43,40,6,0,23,21,12,41,4,9,16,14,44,31)_{\adfTFDEF}$,

$(5,27,50,39,52,17,19,14,45,30,23,7,24,47)_{\adfTFDEF}$

\adfLgap
\noindent under the action of the mapping $x \mapsto x + 5 \adfmod{55}$.

\ADFvfyParStart{55, \{\{55,11,5\}\}, -1, -1, -1} 

\noindent {\boldmath $K_{66}$}~
Let the vertex set be $Z_{65} \cup \{\infty\}$. The decompositions consist of the graphs

\adfLgap
$(\infty,54,56,37,27,44,1,53,18,46,14,15,49,59)_{\adfTFABL}$,

$(54,15,58,31,30,1,39,53,7,40,35,41,57,12)_{\adfTFABL}$,

$(31,45,0,40,21,27,36,6,23,29,3,55,64,17)_{\adfTFABL}$,

$(60,18,47,11,43,16,57,27,19,38,45,21,26,8)_{\adfTFABL}$,

$(48,64,28,26,33,47,16,12,37,15,43,60,10,34)_{\adfTFABL}$,

$(9,43,42,61,64,23,8,16,17,25,15,45,57,3)_{\adfTFABL}$,

$(53,27,22,20,38,0,14,6,46,64,10,37,51,34)_{\adfTFABL}$,

$(55,15,62,63,0,2,58,14,26,46,43,37,22,59)_{\adfTFABL}$,

$(43,48,27,54,35,52,14,34,22,41,15,18,16,60)_{\adfTFABL}$,

$(22,24,64,9,5,16,44,29,0,59,58,49,56,2)_{\adfTFABL}$,

$(8,32,61,12,56,1,5,47,21,6,48,\infty,45,29)_{\adfTFABL}$,

\adfLgap
$(\infty,29,55,25,22,23,50,51,54,26,21,17,46,27)_{\adfTFADJ}$,

$(53,47,10,63,17,19,14,58,64,18,4,28,37,30)_{\adfTFADJ}$,

$(40,43,12,24,39,50,63,28,61,10,44,5,33,53)_{\adfTFADJ}$,

$(0,36,9,8,55,21,28,17,64,45,5,63,24,40)_{\adfTFADJ}$,

$(61,8,13,56,33,54,27,53,4,17,3,5,47,19)_{\adfTFADJ}$,

$(23,18,57,35,22,2,40,55,44,36,20,53,17,56)_{\adfTFADJ}$,

$(62,46,22,1,19,52,5,36,11,12,0,8,57,33)_{\adfTFADJ}$,

$(52,4,53,51,45,2,36,12,17,37,14,59,0,29)_{\adfTFADJ}$,

$(61,50,52,9,5,43,14,15,29,6,56,11,41,33)_{\adfTFADJ}$,

$(26,9,59,46,44,35,37,45,50,17,31,3,18,21)_{\adfTFADJ}$,

$(21,11,45,6,17,64,2,34,24,33,10,54,62,\infty)_{\adfTFADJ}$,

\adfLgap
$(\infty,12,8,59,11,55,54,6,37,30,63,51,2,17)_{\adfTFAFH}$,

$(52,55,17,0,10,22,59,43,57,49,13,30,9,28)_{\adfTFAFH}$,

$(3,26,36,46,44,33,51,48,55,31,7,16,10,32)_{\adfTFAFH}$,

$(9,36,54,11,30,39,57,24,10,59,33,17,38,29)_{\adfTFAFH}$,

$(18,19,59,42,3,55,25,10,5,62,63,33,36,44)_{\adfTFAFH}$,

$(35,53,55,57,10,34,1,2,6,11,47,23,28,25)_{\adfTFAFH}$,

$(7,13,20,46,63,8,51,47,28,53,22,33,49,15)_{\adfTFAFH}$,

$(40,29,56,45,5,28,59,11,31,23,21,34,46,7)_{\adfTFAFH}$,

$(1,8,15,42,9,56,14,5,61,31,3,53,0,37)_{\adfTFAFH}$,

$(30,4,31,12,64,57,2,45,11,62,58,14,19,7)_{\adfTFAFH}$,

$(15,44,58,54,51,62,17,34,24,61,11,12,35,\infty)_{\adfTFAFH}$,

\adfLgap
$(\infty,45,43,34,60,18,40,37,19,58,59,47,5,39)_{\adfTFBBK}$,

$(42,14,8,52,53,33,46,48,2,30,55,0,36,7)_{\adfTFBBK}$,

$(35,43,52,47,20,2,41,31,54,9,64,14,16,1)_{\adfTFBBK}$,

$(18,10,3,58,39,64,61,62,13,0,26,40,27,24)_{\adfTFBBK}$,

$(64,53,31,35,51,27,48,42,11,36,47,40,6,18)_{\adfTFBBK}$,

$(11,24,38,29,3,39,0,6,49,47,31,13,40,45)_{\adfTFBBK}$,

$(32,21,17,7,62,54,37,5,58,53,52,38,28,20)_{\adfTFBBK}$,

$(15,27,35,14,50,7,45,8,5,38,12,48,31,4)_{\adfTFBBK}$,

$(22,1,31,54,44,60,51,59,48,29,53,4,45,61)_{\adfTFBBK}$,

$(47,28,52,19,1,8,36,29,20,22,2,61,41,24)_{\adfTFBBK}$,

$(30,38,6,41,26,5,9,39,11,32,\infty,31,50,4)_{\adfTFBBK}$,

\adfLgap
$(\infty,15,61,23,5,0,35,22,28,32,60,54,41,21)_{\adfTFBCJ}$,

$(14,16,4,10,42,55,29,59,63,22,33,9,58,41)_{\adfTFBCJ}$,

$(1,17,40,36,19,29,6,7,54,59,42,20,0,63)_{\adfTFBCJ}$,

$(34,6,56,12,16,13,53,33,50,2,20,48,45,24)_{\adfTFBCJ}$,

$(19,62,55,28,27,8,22,38,25,9,43,31,7,59)_{\adfTFBCJ}$,

$(47,38,7,41,44,59,22,30,1,15,20,58,61,5)_{\adfTFBCJ}$,

$(8,44,43,34,63,51,30,19,52,55,40,36,57,2)_{\adfTFBCJ}$,

$(26,46,18,12,41,42,37,8,59,60,3,47,56,45)_{\adfTFBCJ}$,

$(28,9,41,13,36,50,1,48,58,31,29,2,21,55)_{\adfTFBCJ}$,

$(40,32,52,2,58,49,56,15,55,43,41,48,25,39)_{\adfTFBCJ}$,

$(10,27,12,44,19,3,8,41,7,61,4,24,9,\infty)_{\adfTFBCJ}$,

\adfLgap
$(\infty,64,24,28,7,55,21,2,51,5,60,61,41,32)_{\adfTFBDI}$,

$(56,50,8,19,30,37,42,43,31,7,52,39,62,32)_{\adfTFBDI}$,

$(34,63,24,53,4,9,1,54,33,36,57,29,5,27)_{\adfTFBDI}$,

$(52,35,56,12,63,8,21,61,46,17,15,10,22,60)_{\adfTFBDI}$,

$(40,10,1,28,43,41,17,20,26,4,5,8,32,59)_{\adfTFBDI}$,

$(3,62,57,35,21,43,5,40,32,63,33,56,46,30)_{\adfTFBDI}$,

$(36,7,63,43,38,13,47,37,59,0,19,27,53,57)_{\adfTFBDI}$,

$(45,15,63,59,36,54,23,29,32,39,4,19,64,30)_{\adfTFBDI}$,

$(39,56,57,48,19,17,36,6,58,42,5,46,41,23)_{\adfTFBDI}$,

$(43,55,19,35,39,1,30,58,54,37,31,44,46,10)_{\adfTFBDI}$,

$(13,63,64,21,25,18,41,24,12,\infty,0,44,36,29)_{\adfTFBDI}$,

\adfLgap
$(\infty,44,21,22,12,53,47,5,49,17,54,59,20,61)_{\adfTFBEH}$,

$(7,45,3,41,22,17,40,34,43,48,6,62,18,32)_{\adfTFBEH}$,

$(35,44,31,23,32,62,54,4,8,29,3,38,22,24)_{\adfTFBEH}$,

$(28,25,17,47,48,58,46,13,6,42,49,51,33,34)_{\adfTFBEH}$,

$(57,18,35,7,54,23,31,51,21,61,62,37,5,34)_{\adfTFBEH}$,

$(2,24,30,60,42,54,43,44,9,41,56,25,31,5)_{\adfTFBEH}$,

$(53,36,22,39,35,45,25,16,42,21,31,51,59,41)_{\adfTFBEH}$,

$(16,15,12,32,52,25,28,54,37,61,18,58,41,4)_{\adfTFBEH}$,

$(18,60,38,51,35,21,30,16,45,46,0,28,31,9)_{\adfTFBEH}$,

$(38,60,45,44,19,4,28,2,50,52,39,46,49,13)_{\adfTFBEH}$,

$(16,20,4,27,53,45,18,43,\infty,34,10,50,62,19)_{\adfTFBEH}$,

\adfLgap
$(\infty,53,44,2,26,42,3,41,5,6,31,8,62,60)_{\adfTFBFG}$,

$(10,21,28,39,24,35,58,41,44,61,62,12,55,51)_{\adfTFBFG}$,

$(19,25,10,13,57,45,17,24,2,64,1,59,6,9)_{\adfTFBFG}$,

$(63,14,34,30,60,46,5,2,33,25,36,49,17,61)_{\adfTFBFG}$,

$(8,63,2,5,62,48,20,47,53,17,0,60,35,15)_{\adfTFBFG}$,

$(34,28,62,4,59,0,38,24,20,56,43,54,33,18)_{\adfTFBFG}$,

$(25,54,12,15,21,31,0,26,16,62,55,36,22,27)_{\adfTFBFG}$,

$(42,10,31,24,64,56,47,49,6,43,48,5,58,26)_{\adfTFBFG}$,

$(24,43,46,51,8,10,42,52,32,6,21,23,4,9)_{\adfTFBFG}$,

$(54,49,30,53,14,56,22,16,38,37,13,21,26,8)_{\adfTFBFG}$,

$(5,58,18,49,62,15,57,12,9,52,17,42,46,\infty)_{\adfTFBFG}$,

\adfLgap
$(\infty,59,55,35,31,6,53,62,58,27,11,61,45,2)_{\adfTFCDH}$,

$(18,20,51,23,29,34,55,5,16,19,31,50,2,21)_{\adfTFCDH}$,

$(6,62,3,26,63,38,15,12,35,28,37,17,32,27)_{\adfTFCDH}$,

$(40,63,45,51,6,17,10,59,44,4,3,2,38,0)_{\adfTFCDH}$,

$(40,27,16,24,18,33,13,1,62,64,36,45,55,41)_{\adfTFCDH}$,

$(52,2,4,13,62,25,29,61,26,43,15,0,48,41)_{\adfTFCDH}$,

$(27,34,53,4,14,48,13,8,46,3,60,57,32,56)_{\adfTFCDH}$,

$(57,19,55,51,63,45,37,30,17,54,3,47,35,29)_{\adfTFCDH}$,

$(12,16,46,36,45,1,64,53,29,52,19,10,49,11)_{\adfTFCDH}$,

$(20,13,56,49,9,29,45,34,38,28,12,24,6,8)_{\adfTFCDH}$,

$(6,5,29,30,27,26,39,58,19,62,4,23,\infty,34)_{\adfTFCDH}$,

\adfLgap
$(\infty,52,44,43,53,59,33,50,47,17,39,32,61,26)_{\adfTFCFF}$,

$(61,48,24,53,62,17,2,43,55,10,31,52,44,42)_{\adfTFCFF}$,

$(15,21,20,48,57,29,31,19,5,4,64,55,37,2)_{\adfTFCFF}$,

$(23,0,14,40,46,42,2,64,34,45,56,1,51,4)_{\adfTFCFF}$,

$(27,37,13,48,28,44,1,5,20,14,41,33,35,8)_{\adfTFCFF}$,

$(50,21,60,7,44,63,8,52,4,53,51,17,10,46)_{\adfTFCFF}$,

$(61,47,35,36,54,34,1,9,25,13,60,27,11,14)_{\adfTFCFF}$,

$(60,26,32,13,25,4,46,52,48,41,36,56,15,23)_{\adfTFCFF}$,

$(15,23,28,3,44,59,35,26,38,48,24,12,39,57)_{\adfTFCFF}$,

$(42,5,3,36,55,35,37,14,24,1,48,4,56,28)_{\adfTFCFF}$,

$(34,5,48,32,23,19,53,46,\infty,9,10,16,7,62)_{\adfTFCFF}$,

\adfLgap
$(\infty,7,30,2,51,64,47,37,42,6,43,0,50,41)_{\adfTFDDG}$,

$(31,30,0,5,12,53,59,51,57,52,10,37,44,20)_{\adfTFDDG}$,

$(33,31,45,12,18,7,52,28,10,26,9,4,2,42)_{\adfTFDDG}$,

$(15,22,26,20,4,41,53,64,37,61,23,55,2,18)_{\adfTFDDG}$,

$(56,18,60,54,8,19,39,5,52,14,53,24,25,55)_{\adfTFDDG}$,

$(42,6,64,7,8,27,63,31,53,19,10,54,4,30)_{\adfTFDDG}$,

$(39,48,21,2,23,36,27,24,6,11,49,60,20,40)_{\adfTFDDG}$,

$(47,31,59,46,24,38,56,8,61,3,48,53,39,41)_{\adfTFDDG}$,

$(24,29,37,20,38,14,18,25,60,58,1,51,5,57)_{\adfTFDDG}$,

$(43,50,57,56,36,40,42,23,58,14,63,15,29,4)_{\adfTFDDG}$,

$(18,47,48,49,14,\infty,26,55,52,46,4,16,51,50)_{\adfTFDDG}$

\adfLgap
$(\infty,40,63,31,46,45,23,17,1,16,24,37,20,58)_{\adfTFDEF}$,

$(25,7,6,46,12,29,37,15,28,17,57,27,23,34)_{\adfTFDEF}$,

$(29,8,38,57,5,34,64,36,54,58,24,22,40,11)_{\adfTFDEF}$,

$(58,52,33,9,25,11,57,3,19,46,39,59,44,34)_{\adfTFDEF}$,

$(6,54,39,2,1,48,12,16,29,49,10,58,25,55)_{\adfTFDEF}$,

$(3,6,8,51,57,18,39,5,61,16,50,27,55,11)_{\adfTFDEF}$,

$(55,4,0,41,3,64,61,34,15,60,35,42,10,25)_{\adfTFDEF}$,

$(24,41,38,15,13,28,20,21,32,12,3,17,39,25)_{\adfTFDEF}$,

$(14,46,31,2,44,17,37,61,50,7,8,28,52,4)_{\adfTFDEF}$,

$(35,10,59,0,12,28,34,8,38,32,6,36,16,30)_{\adfTFDEF}$,

$(11,37,32,47,63,3,15,44,\infty,18,8,56,58,27)_{\adfTFDEF}$

\adfLgap
\noindent under the action of the mapping $x \mapsto x + 5 \adfmod{65}$, $\infty \mapsto \infty$.

\ADFvfyParStart{66, \{\{65,13,5\},\{1,1,1\}\}, 65, -1, -1} 

\noindent {\boldmath $K_{70}$}~
Let the vertex set be $Z_{69} \cup \{\infty\}$. The decompositions consist of the graphs

\adfLgap
$(\infty,42,47,55,9,34,6,58,22,4,31,12,2,45)_{\adfTFABL}$,

$(17,41,44,10,53,60,32,54,20,45,3,35,64,8)_{\adfTFABL}$,

$(56,1,39,62,48,37,34,24,10,42,3,15,26,65)_{\adfTFABL}$,

$(28,48,32,24,0,33,53,16,42,55,64,10,9,17)_{\adfTFABL}$,

$(3,10,68,2,19,29,49,25,27,50,4,39,21,15)_{\adfTFABL}$,

$(31,1,62,32,53,55,8,21,37,15,65,50,68,49)_{\adfTFABL}$,

$(11,23,64,9,63,15,6,35,10,16,28,57,49,14)_{\adfTFABL}$,

\adfLgap
$(\infty,45,5,48,33,28,59,10,17,61,56,13,47,43)_{\adfTFADJ}$,

$(25,23,3,53,5,43,48,38,6,49,57,39,63,17)_{\adfTFADJ}$,

$(34,51,22,19,18,23,65,35,30,38,60,63,52,12)_{\adfTFADJ}$,

$(5,3,2,26,41,19,25,67,14,39,46,55,22,44)_{\adfTFADJ}$,

$(8,37,21,58,9,24,33,49,26,43,18,64,56,65)_{\adfTFADJ}$,

$(10,24,20,54,53,66,17,5,41,40,55,23,30,3)_{\adfTFADJ}$,

$(40,9,30,26,15,36,55,5,57,2,58,19,64,43)_{\adfTFADJ}$,

\adfLgap
$(\infty,40,54,21,44,49,59,2,0,11,64,3,12,60)_{\adfTFAFH}$,

$(1,68,19,4,49,61,60,67,50,47,65,44,43,29)_{\adfTFAFH}$,

$(22,64,18,30,54,38,58,33,29,36,11,53,8,23)_{\adfTFAFH}$,

$(40,65,24,55,25,39,2,15,30,33,60,32,18,52)_{\adfTFAFH}$,

$(56,39,57,26,22,58,20,67,30,37,9,4,38,29)_{\adfTFAFH}$,

$(30,52,8,42,47,14,43,48,9,19,36,56,34,26)_{\adfTFAFH}$,

$(1,24,49,56,19,21,64,47,18,5,62,36,30,43)_{\adfTFAFH}$,

\adfLgap
$(\infty,47,67,65,12,59,26,30,36,5,42,4,63,1)_{\adfTFBBK}$,

$(8,42,18,63,49,59,5,29,31,37,46,38,34,54)_{\adfTFBBK}$,

$(30,29,58,34,15,55,66,68,18,19,6,48,41,10)_{\adfTFBBK}$,

$(7,31,4,43,21,65,54,53,10,27,60,9,55,63)_{\adfTFBBK}$,

$(35,20,32,14,28,12,37,38,68,52,67,56,47,34)_{\adfTFBBK}$,

$(3,46,25,12,0,64,32,67,20,62,33,68,15,2)_{\adfTFBBK}$,

$(56,38,33,30,4,28,46,3,1,31,12,29,57,18)_{\adfTFBBK}$,

\adfLgap
$(\infty,5,28,60,26,59,12,45,36,38,9,31,43,46)_{\adfTFBCJ}$,

$(31,22,23,12,8,24,7,63,66,56,68,50,28,51)_{\adfTFBCJ}$,

$(28,45,62,2,13,49,0,5,18,50,55,46,29,19)_{\adfTFBCJ}$,

$(7,38,25,45,40,18,44,24,30,34,67,28,43,27)_{\adfTFBCJ}$,

$(42,56,23,24,3,52,46,19,43,14,59,10,35,41)_{\adfTFBCJ}$,

$(65,4,23,56,18,51,7,6,30,64,68,2,40,12)_{\adfTFBCJ}$,

$(9,65,21,32,24,63,61,8,1,38,39,0,42,35)_{\adfTFBCJ}$,

\adfLgap
$(\infty,19,0,65,39,60,52,14,1,56,32,63,55,31)_{\adfTFBDI}$,

$(47,20,55,10,32,59,42,0,2,68,53,34,19,57)_{\adfTFBDI}$,

$(11,60,59,57,52,66,64,22,21,54,25,8,10,5)_{\adfTFBDI}$,

$(62,54,15,55,57,45,29,4,22,65,48,56,27,3)_{\adfTFBDI}$,

$(43,40,36,4,48,61,60,28,25,51,23,35,54,20)_{\adfTFBDI}$,

$(9,57,53,55,45,42,29,58,52,3,62,51,38,10)_{\adfTFBDI}$,

$(48,47,41,32,28,38,37,3,19,8,7,16,52,29)_{\adfTFBDI}$,

\adfLgap
$(\infty,62,15,50,65,37,32,1,45,54,39,4,58,28)_{\adfTFBEH}$,

$(16,49,43,42,40,21,66,48,18,8,32,53,41,67)_{\adfTFBEH}$,

$(6,21,32,3,9,27,35,44,7,45,58,49,60,20)_{\adfTFBEH}$,

$(2,38,21,35,62,40,54,55,34,18,41,32,43,3)_{\adfTFBEH}$,

$(3,46,65,23,37,34,24,64,28,33,0,42,44,50)_{\adfTFBEH}$,

$(30,10,2,51,26,1,11,34,57,37,36,23,16,62)_{\adfTFBEH}$,

$(1,18,46,50,47,65,25,13,11,15,20,7,38,6)_{\adfTFBEH}$,

\adfLgap
$(\infty,67,10,54,1,63,37,35,59,57,39,34,14,19)_{\adfTFBFG}$,

$(22,14,40,52,67,57,17,25,63,43,10,34,7,30)_{\adfTFBFG}$,

$(26,22,19,9,17,55,21,13,7,59,41,65,18,3)_{\adfTFBFG}$,

$(2,22,35,50,0,12,47,51,55,57,48,4,38,26)_{\adfTFBFG}$,

$(18,12,57,56,2,15,21,44,38,8,14,28,66,33)_{\adfTFBFG}$,

$(2,35,30,61,57,58,11,8,1,47,3,62,34,21)_{\adfTFBFG}$,

$(58,15,52,36,63,49,66,8,0,68,28,53,44,18)_{\adfTFBFG}$,

\adfLgap
$(\infty,12,28,23,17,51,54,21,41,59,19,5,48,63)_{\adfTFCDH}$,

$(2,1,43,62,26,42,3,50,51,46,60,31,9,33)_{\adfTFCDH}$,

$(6,37,64,47,29,41,50,52,43,40,1,39,27,68)_{\adfTFCDH}$,

$(15,27,65,43,8,1,5,63,3,7,61,4,53,56)_{\adfTFCDH}$,

$(19,45,35,46,37,43,64,55,21,31,32,7,44,1)_{\adfTFCDH}$,

$(2,42,32,59,39,22,56,25,51,47,11,5,7,49)_{\adfTFCDH}$,

$(5,60,66,35,3,16,40,0,6,34,42,9,65,50)_{\adfTFCDH}$,

\adfLgap
$(\infty,39,63,40,53,26,48,59,7,49,21,28,62,58)_{\adfTFCFF}$,

$(35,61,63,29,7,46,4,44,28,3,17,30,42,59)_{\adfTFCFF}$,

$(37,7,49,4,29,47,11,27,63,59,9,35,30,33)_{\adfTFCFF}$,

$(15,51,4,67,19,21,60,1,44,57,43,3,9,27)_{\adfTFCFF}$,

$(14,43,37,28,64,15,24,44,50,34,3,41,62,32)_{\adfTFCFF}$,

$(58,21,66,20,36,51,47,57,65,59,2,40,50,55)_{\adfTFCFF}$,

$(48,53,27,56,43,61,5,14,28,4,25,42,44,29)_{\adfTFCFF}$,

\adfLgap
$(\infty,7,57,0,48,13,58,2,8,51,35,68,53,60)_{\adfTFDDG}$,

$(37,17,62,57,39,25,30,31,54,63,48,23,21,19)_{\adfTFDDG}$,

$(19,21,35,59,17,11,34,27,63,2,37,38,6,56)_{\adfTFDDG}$,

$(59,63,31,46,24,47,57,67,50,39,11,42,15,18)_{\adfTFDDG}$,

$(32,17,64,30,43,9,40,7,26,44,23,6,55,18)_{\adfTFDDG}$,

$(45,7,12,32,28,19,2,57,16,39,25,67,5,18)_{\adfTFDDG}$,

$(20,25,40,11,14,67,61,22,1,32,2,42,34,16)_{\adfTFDDG}$

\adfLgap
$(\infty,68,7,50,48,8,65,21,29,0,24,52,46,2)_{\adfTFDEF}$,

$(32,63,14,41,15,21,68,16,6,61,48,20,1,31)_{\adfTFDEF}$,

$(39,18,56,43,10,12,3,35,37,68,62,0,14,25)_{\adfTFDEF}$,

$(14,40,4,5,28,45,60,56,20,30,55,46,67,16)_{\adfTFDEF}$,

$(62,41,2,43,65,7,23,0,46,63,52,10,17,36)_{\adfTFDEF}$,

$(53,57,38,1,18,22,6,65,61,66,0,36,67,52)_{\adfTFDEF}$,

$(0,68,1,27,33,51,22,42,7,55,57,10,13,47)_{\adfTFDEF}$

\adfLgap
\noindent under the action of the mapping $x \mapsto x + 3 \adfmod{69}$, $\infty \mapsto \infty$.

\ADFvfyParStart{70, \{\{69,23,3\},\{1,1,1\}\}, 69, -1, -1} 

~\\
\noindent {\bf Proof of Lemma \ref{lem:theta15 multipartite}}

\noindent {\boldmath $K_{15,15}$}~
Let the vertex set be $Z_{30}$ partitioned according to residue class modulo 2.
The decompositions consist of the graphs

\adfLgap
$(0,1,3,4,5,10,17,2,11,18,7,16,27,14)_{\adfTFACK}$,

\adfLgap
$(0,1,3,4,9,12,7,14,5,10,21,8,25,16)_{\adfTFAEI}$,

\adfLgap
$(0,1,3,4,9,2,11,14,13,18,29,6,25,10)_{\adfTFAGG}$,

\adfLgap
$(0,1,3,2,5,8,7,12,21,4,17,6,25,10)_{\adfTFCCI}$,

\adfLgap
$(0,1,3,2,5,8,15,6,11,18,7,16,29,14)_{\adfTFCEG}$

\adfLgap
$(0,1,3,2,7,8,9,12,19,6,21,10,29,14)_{\adfTFEEE}$

\adfLgap
\noindent under the action of the mapping $x \mapsto x + 2 \adfmod{30}$.

\ADFvfyParStart{30, \{\{30,15,2\}\}, -1, \{\{15,\{0,1\}\}\}, -1} 

\noindent {\boldmath $K_{15,20}$}~
Let the vertex set be $\{0, 1, \dots, 34\}$ partitioned into $\{0, 1, \dots, 19\}$
and $\{20, 21, \dots, 34\}$. The decompositions consist of the graphs

\adfLgap
$(0,28,30,5,22,8,27,9,29,15,21,12,23,13)_{\adfTFACK}$,

$(11,31,28,10,27,4,29,17,21,5,33,2,20,18)_{\adfTFACK}$,

$(9,32,22,15,26,12,28,18,29,8,31,5,25,11)_{\adfTFACK}$,

$(32,10,7,20,16,30,2,21,19,27,6,31,3,33)_{\adfTFACK}$,

\adfLgap
$(0,33,23,6,21,8,20,5,29,15,30,12,31,18)_{\adfTFAEI}$,

$(22,17,7,33,1,24,8,25,0,32,13,30,19,26)_{\adfTFAEI}$,

$(28,6,3,32,15,20,0,29,13,22,5,33,14,25)_{\adfTFAEI}$,

$(33,11,15,25,13,31,6,32,10,29,4,27,2,28)_{\adfTFAEI}$,

\adfLgap
$(0,23,26,11,28,4,25,1,24,16,29,18,32,17)_{\adfTFAGG}$,

$(20,10,7,33,9,25,18,27,14,21,1,32,15,22)_{\adfTFAGG}$,

$(22,11,4,21,0,29,7,27,9,30,1,33,2,25)_{\adfTFAGG}$,

$(1,22,29,12,25,16,27,6,28,7,21,15,23,2)_{\adfTFAGG}$,

\adfLgap
$(0,26,21,4,29,12,24,3,23,10,32,8,33,9)_{\adfTFCCI}$,

$(3,23,20,2,31,1,27,9,30,6,28,16,22,15)_{\adfTFCCI}$,

$(13,25,24,14,20,10,31,9,23,19,33,7,22,4)_{\adfTFCCI}$,

$(6,28,25,7,21,2,27,15,26,17,31,4,33,1)_{\adfTFCCI}$,

\adfLgap
$(0,27,29,12,22,8,21,5,24,6,30,11,33,9)_{\adfTFCEG}$,

$(29,1,3,26,8,23,7,31,2,30,18,22,4,28)_{\adfTFCEG}$,

$(33,18,3,23,14,26,9,31,16,24,13,28,15,25)_{\adfTFCEG}$,

$(7,26,25,6,33,19,22,10,20,9,29,4,31,13)_{\adfTFCEG}$

\adfLgap
$(0,30,25,12,33,19,28,6,26,4,21,9,24,18)_{\adfTFEEE}$,

$(11,30,23,19,32,9,22,6,27,7,21,8,28,5)_{\adfTFEEE}$,

$(23,14,15,22,19,20,6,25,18,26,12,28,13,30)_{\adfTFEEE}$,

$(32,17,0,20,1,29,8,24,14,25,13,31,15,28)_{\adfTFEEE}$

\adfLgap
\noindent under the action of the mapping $x \mapsto x + 4\adfmod{20}$ for $x < 20$,
$x \mapsto (x - 20 + 3 \adfmod{15}) + 20$ for $x \ge 20$.

\ADFvfyParStart{35, \{\{20,5,4\},\{15,5,3\}\}, -1, \{\{20,\{0\}\},\{15,\{1\}\}\}, -1} 

\noindent {\boldmath $K_{15,25}$}~
Let the vertex set be $\{0, 1, \dots, 39\}$ partitioned into $\{0, 1, \dots, 24\}$
and $\{25, 26, \dots, 39\}$. The decompositions consist of the graphs

\adfLgap
$(0,34,32,12,27,17,29,5,38,18,35,22,37,9)_{\adfTFACK}$,

$(32,2,1,37,16,36,20,30,11,27,22,26,19,34)_{\adfTFACK}$,

$(27,9,7,30,12,37,11,26,24,25,3,35,16,38)_{\adfTFACK}$,

$(1,30,34,15,25,20,28,10,35,8,26,4,33,3)_{\adfTFACK}$,

$(3,27,31,8,37,4,30,18,28,1,33,24,35,20)_{\adfTFACK}$,

\adfLgap
$(0,36,37,18,38,14,25,17,27,15,33,3,31,1)_{\adfTFAEI}$,

$(23,38,36,2,30,3,34,5,37,6,35,19,31,4)_{\adfTFAEI}$,

$(24,31,32,9,27,16,33,6,28,12,35,2,26,5)_{\adfTFAEI}$,

$(32,17,20,29,0,26,6,30,9,28,2,25,3,36)_{\adfTFAEI}$,

$(1,26,28,4,33,13,39,15,36,22,31,23,29,11)_{\adfTFAEI}$,

\adfLgap
$(0,34,38,19,37,9,25,1,27,11,36,4,30,7)_{\adfTFAGG}$,

$(1,28,26,23,38,10,36,2,31,0,29,5,25,8)_{\adfTFAGG}$,

$(25,12,4,26,18,33,13,36,23,30,11,29,22,35)_{\adfTFAGG}$,

$(23,32,35,24,36,3,34,21,28,12,30,20,33,4)_{\adfTFAGG}$,

$(19,36,25,0,39,12,32,6,29,17,34,10,26,16)_{\adfTFAGG}$,

\adfLgap
$(0,35,32,14,33,17,31,4,30,22,29,19,28,9)_{\adfTFCCI}$,

$(29,17,11,34,3,32,13,30,16,33,10,28,12,25)_{\adfTFCCI}$,

$(16,33,36,4,31,7,28,3,37,5,27,21,32,20)_{\adfTFCCI}$,

$(8,29,33,5,37,21,36,2,35,10,34,15,32,6)_{\adfTFCCI}$,

$(14,37,33,19,28,23,29,18,36,3,26,22,39,11)_{\adfTFCCI}$,

\adfLgap
$(0,30,31,5,25,16,27,2,29,22,35,18,32,15)_{\adfTFCEG}$,

$(24,30,37,3,36,20,27,18,31,4,32,19,38,16)_{\adfTFCEG}$,

$(18,26,38,13,33,14,27,15,28,10,25,1,35,6)_{\adfTFCEG}$,

$(22,37,28,17,38,7,27,19,26,0,35,16,36,18)_{\adfTFCEG}$,

$(16,30,31,17,28,8,34,11,39,19,35,4,37,22)_{\adfTFCEG}$

\adfLgap
$(0,29,27,24,26,13,37,11,31,14,30,1,28,3)_{\adfTFEEE}$,

$(6,31,27,23,34,0,30,7,37,17,38,9,36,24)_{\adfTFEEE}$,

$(18,30,26,22,29,21,35,2,36,23,25,9,34,10)_{\adfTFEEE}$,

$(12,38,30,16,28,5,31,2,26,15,36,18,34,21)_{\adfTFEEE}$,

$(29,14,24,25,13,27,11,35,20,33,12,39,10,32)_{\adfTFEEE}$

\adfLgap
\noindent under the action of the mapping $x \mapsto x + 5 \adfmod{25}$ for $x < 25$,
$x \mapsto (x - 25 + 3 \adfmod{15}) + 25$ for $x \ge 25$.

\ADFvfyParStart{40, \{\{25,5,5\},\{15,5,3\}\}, -1, \{\{25,\{0\}\},\{15,\{1\}\}\}, -1} 

\noindent {\boldmath $K_{5,5,5}$}~
Let the vertex set be $Z_{15}$ partitioned according to residue class modulo 3.
The decompositions consist of the graphs

\adfLgap
$(0,1,11,2,12,10,9,14,7,5,6,4,8,3)_{\adfTFABL}$,

$(1,2,6,5,4,9,13,8,12,14,0,10,11,7)_{\adfTFABL}$,

$(2,3,10,13,11,12,4,14,6,8,0,7,9,5)_{\adfTFABL}$,

$(0,13,5,4,3,7,6,11,9,8,10,14,1,12)_{\adfTFABL}$,

$(3,13,14,11,4,2,9,1,8,7,12,5,10,6)_{\adfTFABL}$,

\adfLgap
$(0,1,14,13,11,5,10,6,8,4,2,9,7,12)_{\adfTFADJ}$,

$(1,2,8,12,10,5,13,6,14,9,4,11,3,7)_{\adfTFADJ}$,

$(2,3,0,10,8,12,14,1,6,5,7,11,9,13)_{\adfTFADJ}$,

$(3,4,1,9,5,10,14,7,6,11,0,8,13,12)_{\adfTFADJ}$,

$(0,4,13,2,6,7,8,9,10,11,12,5,3,14)_{\adfTFADJ}$

\adfLgap
$(0,1,13,11,12,2,6,8,7,5,10,9,4,14)_{\adfTFAFH}$,

$(1,2,5,12,14,9,7,11,4,6,13,3,8,10)_{\adfTFAFH}$,

$(2,3,13,14,10,6,7,9,11,0,4,8,12,1)_{\adfTFAFH}$,

$(0,5,10,12,13,8,9,2,4,3,11,7,14,6)_{\adfTFAFH}$,

$(3,5,14,0,7,12,4,10,11,6,8,1,9,13)_{\adfTFAFH}$


\ADFvfyParStart{15, \{\{15,1,15\}\}, -1, \{\{5,\{0,1,2\}\}\}, -1} 

\adfLgap
\noindent as well as the graphs

\adfLgap
$(0,3,1,2,4,11,7,8,10,5,9,14,6,13)_{\adfTFBBK}$,

\adfLgap
$(0,3,1,2,7,5,4,11,12,8,10,14,6,13)_{\adfTFBCJ}$,

\adfLgap
$(0,3,1,2,4,11,10,14,9,5,6,13,8,7)_{\adfTFBDI}$,

\adfLgap
$(0,3,1,2,4,5,10,8,12,11,7,14,9,13)_{\adfTFBEH}$,

\adfLgap
$(0,3,1,2,4,9,14,10,11,12,5,13,8,7)_{\adfTFBFG}$,

\adfLgap
$(0,3,1,2,4,6,13,5,10,8,12,14,7,11)_{\adfTFCDH}$,

\adfLgap
$(0,3,1,2,4,8,6,13,11,5,10,12,7,14)_{\adfTFCFF}$,

\adfLgap
$(0,3,1,2,10,4,6,14,5,7,11,9,8,13)_{\adfTFDDG}$

\adfLgap
$(0,3,1,2,7,8,12,14,10,5,6,4,11,13)_{\adfTFDEF}$

\adfLgap
\noindent under the action of the mapping $x \mapsto x + 3 \adfmod{15}$

\ADFvfyParStart{15, \{\{15,5,3\}\}, -1, \{\{5,\{0,1,2\}\}\}, -1} 

\noindent {\boldmath $K_{5,5,5,5}$}~
Let the vertex set be $Z_{20}$ partitioned according to residue class modulo 4.
The decompositions consist of the graphs

\adfLgap
$(0,1,2,3,4,9,6,7,5,8,14,11,16,10)_{\adfTFABL}$,

$(2,7,9,4,11,0,10,5,15,1,8,17,3,14)_{\adfTFABL}$,

\adfLgap
$(0,1,2,3,4,5,6,8,11,9,7,10,19,12)_{\adfTFADJ}$,

$(2,7,5,10,1,8,3,13,6,0,11,17,4,14)_{\adfTFADJ}$,

\adfLgap
$(0,1,2,3,4,7,10,5,6,8,14,9,11,18)_{\adfTFAFH}$,

$(2,7,8,1,4,11,16,9,0,10,19,5,3,17)_{\adfTFAFH}$,

\adfLgap
$(0,1,2,3,5,4,10,7,6,8,15,9,12,11)_{\adfTFBBK}$,

$(2,3,5,8,7,13,0,9,14,4,15,6,17,10)_{\adfTFBBK}$,

\adfLgap
$(0,1,2,3,4,5,7,6,8,9,14,11,13,10)_{\adfTFBCJ}$,

$(2,3,8,7,13,9,16,5,11,4,15,6,0,10)_{\adfTFBCJ}$,

\adfLgap
$(0,1,2,3,4,10,5,7,6,8,9,12,19,14)_{\adfTFBDI}$,

$(2,3,5,8,1,6,11,16,7,13,19,9,0,10)_{\adfTFBDI}$,

\adfLgap
$(0,1,2,3,4,5,7,6,8,13,10,11,14,19)_{\adfTFBEH}$,

$(2,3,8,9,0,7,13,15,6,1,4,14,5,12)_{\adfTFBEH}$,

\adfLgap
$(0,1,2,3,4,5,7,6,9,12,10,13,8,14)_{\adfTFBFG}$,

$(2,3,8,7,0,10,19,6,13,11,1,15,4,17)_{\adfTFBFG}$,

\adfLgap
$(0,1,2,3,5,4,7,6,8,15,10,9,12,11)_{\adfTFCDH}$,

$(2,3,5,10,8,1,12,9,18,7,13,11,16,6)_{\adfTFCDH}$,

\adfLgap
$(0,1,2,3,5,4,7,8,6,9,10,13,11,14)_{\adfTFCFF}$,

$(2,3,7,12,8,1,10,4,14,15,5,19,13,16)_{\adfTFCFF}$,

\adfLgap
$(0,1,2,3,4,5,6,8,7,9,11,14,17,10)_{\adfTFDDG}$,

$(2,3,7,1,0,8,19,9,16,11,18,13,4,14)_{\adfTFDDG}$

\adfLgap
$(0,1,2,3,4,5,6,8,11,7,9,14,17,10)_{\adfTFDEF}$,

$(2,3,7,1,8,11,9,0,10,19,5,4,18,12)_{\adfTFDEF}$

\adfLgap
\noindent under the action of the mapping $x \mapsto x + 4 \adfmod{20}$.

\ADFvfyParStart{20, \{\{20,5,4\}\}, -1, \{\{5,\{0,1,2,3\}\}\}, -1} 
\end{document}